\documentclass[letterpaper,leqno,10pt,twoside]{amsart}
\usepackage{amsmath,amsthm,amsfonts,amssymb,euscript,mathrsfs}
\usepackage{mathtools}
\usepackage{tikz}
\usetikzlibrary{arrows.meta,positioning,calc}
\usepackage[margin=0.79in]{geometry}
\usepackage[toc,page]{appendix}
\usepackage[shortlabels]{enumitem}
\usepackage{booktabs,array}
\usepackage{indentfirst}
\usepackage{float}
\usepackage{dsfont}

\numberwithin{equation}{section}

\usepackage[dvipsnames]{xcolor}

\usepackage{hyperref}
\hypersetup{colorlinks=true, pdfstartview=Fit,linkcolor=blue!70!black,citecolor=red!70!black, urlcolor=green!60!black}
\definecolor{labelkey}{rgb}{0.6,0,0}

\newtheorem{proposition}{Proposition}[section]
\newtheorem{theorem}[proposition]{Theorem}
\newtheorem{lemma}[proposition]{Lemma}
\newtheorem{corollary}[proposition]{Corollary}

\theoremstyle{definition}
\newtheorem{definition}[proposition]{Definition}
\newtheorem{assumption}[proposition]{Assumption}
\theoremstyle{remark}
\newtheorem{remark}{Remark}[section]

\setlist[enumerate]{label=\textup{(\roman*)},leftmargin=*,itemsep=4pt}
\newtagform{enumerated}{\hspace{\leftmargin}(}{)}
\AddToHook{env/enumerate/begin}{\usetagform{enumerated}}

\let\standardeqref\eqref
\renewcommand{\eqref}[1]{{\usetagform{default}\standardeqref{#1}}}
\allowdisplaybreaks[2]

\newcommand{\proofstep}[1]{\par\smallskip\noindent\underline{\textit{#1}}\enspace\ignorespaces}
\newcommand{\para}[1]{\paragraph{\underline{#1}}}

\def\ue{e}

\newcommand{\m}{\mu}
\newcommand{\mh}{\m^{\frac{1}{2}}}
\newcommand{\mih}{\m^{-\frac{1}{2}}}
\newcommand{\phase}{\mathcal O}
\newcommand{\sig}{\sigma}
\newcommand{\pk}{P}
\newcommand{\Vv}{V_v}
\newcommand{\Wv}{W_v}
\newcommand{\Wspace}{W}
\newcommand{\trace}{\tau}
\newcommand{\trp}{\trace_+}
\newcommand{\trm}{\trace_-}
\newcommand{\trpm}{\trace_\pm}
\newcommand{\Pp}{P_+}
\newcommand{\Pm}{P_-}
\newcommand{\graz}{\rho}
\newcommand{\dkin}{d_{\mathrm{kin}}}
\newcommand{\Ckin}{C_{\mathrm{kin}}}
\newcommand{\Ckinloc}{C_{\mathrm{kin},\mathrm{loc}}}
\newcommand{\Lop}{L}
\newcommand{\Kop}{K}
\newcommand{\glift}{\bar g}
\newcommand{\sourceR}{\mathcal R}
\newcommand{\bdproj}{p_{\partial}}
\newcommand{\vradius}{r_0}
\newcommand{\locA}{\widehat A}
\newcommand{\locb}{\widehat b}
\newcommand{\locB}{\widehat B}
\newcommand{\locc}{\widehat c}
\newcommand{\fluxbase}{\mathfrak j}
\newcommand{\flux}[1]{\fluxbase_{#1}}
\NewDocumentCommand{\abs}{o m}{\IfNoValueTF{#1}{\left\lvert#2\right\rvert}{{\mathopen{#1\lvert}#2\mathclose{#1\rvert}}}}
\NewDocumentCommand{\bracket}{o m}{\IfNoValueTF{#1}{\left\langle#2\right\rangle}{{\mathopen{#1\langle}#2\mathclose{#1\rangle}}}}

\newcommand{\R}{\mathbb R}
\newcommand{\dd}{\,\mathrm d}
\newcommand{\la}{\bracket{v}}
\newcommand{\Om}{\Omega}
\newcommand{\Vspace}{V}
\newcommand{\Y}{Y}
\newcommand{\Yhat}{\widehat Y}
\newcommand{\Qc}{\mathcal Q}
\newcommand{\Qplus}{\Qc^+}
\newcommand{\fluxhat}{\widehat{\fluxbase}}
\newcommand{\fluxlift}{\overline{\fluxbase}}
\newcommand{\Ga}{\Gamma}
\newcommand{\Sigmap}{\Sigma_+}
\newcommand{\Sigmam}{\Sigma_-}
\newcommand{\Af}{a_F}
\newcommand{\Bf}{B_F}
\newcommand{\Cf}{C_F}

\newcommand{\logconst}{\mathcal{D}}

\newcommand{\liftpatch}{\mathscr U}
\newcommand{\locOp}{\mathscr L}

\newcommand{\compamp}{\mathfrak r}
\newcommand{\gbarrier}{\mathscr M}
\newcommand{\spatialmult}{\vartheta}
\newcommand{\inflowweight}{w}
\newcommand{\inflowenergy}{\mathscr E}
\newcommand{\kinback}{\mathcal C}

\NewDocumentCommand{\vnorm}{o m}{\IfNoValueTF{#1}{\left\lvert#2\right\rvert}{{\mathopen{#1\lvert}#2\mathclose{#1\rvert}}}}
\NewDocumentCommand{\norm}{o m}{\IfNoValueTF{#1}{\left\lVert#2\right\rVert}{{\mathopen{#1\lVert}#2\mathclose{#1\rVert}}}}
\NewDocumentCommand{\tnorm}{o m}{\IfNoValueTF{#1}{\left\lvert\!\left\lvert\!\left\lvert#2\right\rvert\!\right\rvert\!\right\rvert}{{\mathopen{#1\lvert\!#1\lvert\!#1\lvert}#2\mathclose{#1\rvert\!#1\rvert\!#1\rvert}}}}
\NewDocumentCommand{\vgnorm}{o m m m}{\IfNoValueTF{#1}{\vnorm{#4}}{\vnorm[#1]{#4}}_{L^{#2,#3}_v}}
\NewDocumentCommand{\gnorm}{o m m m}{\IfNoValueTF{#1}{\norm{#4}}{\norm[#1]{#4}}_{L^{#2,#3}}}
\NewDocumentCommand{\tgnorm}{o m m m}{\IfNoValueTF{#1}{\tnorm{#4}}{\tnorm[#1]{#4}}_{L^{#2,#3}_{t,x,v}}}
\NewDocumentCommand{\bgnorm}{o m m m m}{\IfNoValueTF{#1}{\norm{#5}}{\norm[#1]{#5}}_{L^{#2,#3}_{\Sigma_{#4}}}}
\NewDocumentCommand{\tbgnorm}{o m m m m}{\IfNoValueTF{#1}{\tnorm{#5}}{\tnorm[#1]{#5}}_{L^{#2,#3}_{t,\Sigma_{#4}}}}

\DeclarePairedDelimiterX{\ip}[2]{\langle}{\rangle}{#1,#2}
\newcommand{\ind}{\mathds 1}
\DeclareMathOperator{\divv}{div}
\DeclareMathOperator{\tr}{tr}
\DeclareMathOperator{\supp}{supp}
\DeclareMathOperator{\dist}{dist}

\title[Landau equation in bounded domains]{Global Well-Posedness of Landau Equation in Bounded Domains}

\author{Zhimeng Ouyang}
\address{School of Mathematical Sciences, Peking University, Beijing 100871, China}
\email{ouyangzm@math.pku.edu.cn}

\author{Lei Wu}
\address{Department of Mathematics, Lehigh University, Bethlehem, PA 18015, USA}
\email{lew218@lehigh.edu}

\date{}
\subjclass[2020]{35Q20, 35K65, 35B65, 82C40}
\keywords{Landau equation, Coulomb interaction, in-flow, diffuse-reflection, kinetic trace, grazing set, strong solution, uniqueness}

\begin{document}

\begin{abstract}
We prove global well-posedness of the Landau equation near a Maxwellian in bounded $C^3$ domains with in-flow or diffuse-reflection boundary conditions for $-3\leq\gamma\leq1$, including Coulomb interactions. Small compatible Gaussian perturbations generate global strong solutions, with uniqueness and continuous dependence among weak solutions. 
The proof combines energy estimates, the comparison principle, interior and weighted boundary estimates for velocity derivatives up to second order, and a logarithmic estimate for solution difference.
\end{abstract}

\maketitle

\setcounter{tocdepth}{1}
\makeatletter
\def\l@subsection{\@tocline{2}{0pt}{2.5pc}{5pc}{}}
\makeatother
\tableofcontents


\section{Introduction}\label{sec:introduction}

The Landau equation \cite{Landau1936} describes a dilute gas or plasma whose collisions are dominated by small deflections and arises from the Boltzmann equation in the grazing collision limit \cite{Villani2002}. In a bounded domain, the boundary condition describes the exchange of particles with the exterior or their reflection at a wall. The in-flow condition specifies the incoming distribution, and diffuse-reflection describes re-emission according to the wall Maxwellian. These boundary laws determine the conservation properties and the relaxation of the gas toward equilibrium.

\subsection{Problem Formulation}\label{sec:intro-model}

Let $\Om\subset\R^3$ be a bounded connected domain with $C^3$ boundary and outward unit normal $n$.

The particle density $F(t,x,v)\geq0$, at time $t\geq0$, position $x\in\Om$ and velocity $v\in\R^3$, solves
\begin{align}
\Y F=Q(F,F),\qquad \Y:=\partial_t+v\cdot\nabla_x,\qquad F(0,x,v)=F_0(x,v). \label{eq:full}
\end{align}
Fix an interaction exponent $\gamma\in[-3,1]$. The collision kernels are
\begin{align}
a_{\gamma,ij}(v)&=\abs{v}^{\gamma+2}\bigg(\delta_{ij}-\frac{v_i v_j}{\abs{v}^2}\bigg),\qquad b_{\gamma,i}(v)=-2\abs{v}^\gamma v_i, \label{eq:kernels}\\
c_\gamma(v)&=\begin{cases}-2(\gamma+3)\abs{v}^\gamma,&-3<\gamma\leq1,\\-8\pi\delta_0,&\gamma=-3,\end{cases} \label{eq:kernel-c}
\end{align}
and the Landau collision operator is
\begin{align}
Q_\gamma(G,H)=\divv_v\Bigl(\bigl(a_\gamma*G\bigr)\nabla_vH-\bigl(b_\gamma*G\bigr)H\Bigr). \label{eq:collision}
\end{align}
Convolutions are taken in velocity at fixed $(t,x)$, and repeated indices are summed. The kernels satisfy the distributional identities
\begin{align}
b_{\gamma,i}=\partial_{v_j}a_{\gamma,ij},\qquad c_\gamma=\partial_{v_i}b_{\gamma,i}. \label{eq:kernel-identities}
\end{align}
For smooth $G,H$ with Gaussian decay, Lemma~\ref{cvel:distribution} therefore gives
\begin{align}
Q_\gamma(G,H)=\bigl(a_\gamma*G\bigr):D_v^2H-\bigl(c_\gamma*G\bigr)H. \label{eq:collision-nondivergence}
\end{align}
In the Coulomb case $\gamma=-3$, the distribution $c_{-3}=-8\pi\delta_0$ is a multiple of the Dirac mass at zero. Thus the zeroth-order term in \eqref{eq:collision-nondivergence} is
\begin{align}
-\bigl(c_{-3}*G\bigr)H=8\pi GH. \label{eq:coulomb-local-product}
\end{align}
The range $-3\leq\gamma<-2$ corresponds to very soft potentials, and $-2\leq\gamma<0$ to moderately soft potentials. The case $\gamma=0$ is the Maxwell-molecule case, and $0<\gamma\leq1$ corresponds to hard potentials. The exponent $\gamma$ is fixed throughout the paper. We suppress its subscript on $a,b,c,Q$ and on the linearized operators introduced below.

The Maxwellian is
\begin{align}
\m (v):=(2\pi)^{-\frac{3}{2}}\ue^{-\frac{\abs{v}^2}{2}}. \label{eq:I:display:002}
\end{align}
The phase boundary is divided into
\begin{align}
\Sigma_\pm:=\{(x,v)\in\partial\Om\times\R^3:\ \pm n(x)\cdot v>0\},\qquad \Sigma_0:=\{(x,v)\in\partial\Om\times\R^3:n(x)\cdot v=0\}. \label{eq:phase-boundary}
\end{align}
The sets $\Sigma_+$, $\Sigma_-$ and $\Sigma_0$ are the outgoing, incoming and grazing parts. We write $\trace F$ for the kinetic trace determined by the Green identity in Section~\ref{sec:intro-classes} and $\trpm F$ for its restrictions to $\Sigma_\pm$; Section~\ref{sec:trace} gives its construction and integrability.

We consider two boundary conditions. The in-flow condition prescribes the incoming density,
\begin{align}
\trm F=\m+\mh g\quad\text{on }(0,\infty)\times\Sigma_-, \label{eq:inflow-full}
\end{align}
where $g$ is given and $\m+\mh g\geq0$.

For diffuse-reflection, the wall is stationary and has normalized temperature one. Particles reaching the wall are re-emitted with the Maxwellian velocity distribution $\m$, with its amplitude chosen to conserve the normal mass flux \cite{Cercignani1990,Cercignani.Illner.Pulvirenti1994}. The boundary condition is
\begin{align}
\trm F(t,x,v)=\sqrt{2\pi}\,\m(v)\int_{n(x)\cdot w>0}\trp F(t,x,w)\bigl(n(x)\cdot w\bigr)\dd w. \label{eq:diffuse-full}
\end{align}
The Maxwellian normalization and the resulting flux identity are
\begin{align}
\int_{n(x)\cdot v>0}\m(v)\bigl(n(x)\cdot v\bigr)\dd v&=\frac{1}{\sqrt{2\pi}}, \label{eq:wall-normalization}\\
\int_{n(x)\cdot v<0}\trm F\,\abs[\big]{n(x)\cdot v}\dd v&=\int_{n(x)\cdot v>0}\trp F\,\bigl(n(x)\cdot v\bigr)\dd v. \label{eq:wall-flux-balance}
\end{align}
Consequently, the diffuse-reflection problem conserves total mass.

Near the Maxwellian we write
\begin{align}
F=\m +\mh f,\qquad f_0:=\mih (F_0-\m), \label{eq:normalization}
\end{align}
and define
\begin{align}
\Ga(h_1,h_2):=\mih Q(\mh h_1,\mh h_2),\qquad \Lop h:=-\Ga(\mh,h)-\Ga(h,\mh). \label{eq:linearized-operators}
\end{align}
Since $\Y\m =0$ and $Q(\m,\m)=0$, problem \eqref{eq:full} is equivalent to
\begin{align}
\Y f+\Lop f=\Ga(f,f),\qquad f(0)=f_0. \label{eq:normalized-problem}
\end{align}
The in-flow law \eqref{eq:inflow-full} becomes
\begin{align}
\trm f=g. \label{eq:inflow-normalized}
\end{align}
For $(t,x,v)\in(0,\infty)\times\Sigma_-$, the diffuse-reflection condition \eqref{eq:diffuse-full} becomes
\begin{align}
\trm f(t,x,v)=\Pm\trp f(t,x,v):=\sqrt{2\pi}\,\mh(v)\int_{n(x)\cdot w>0}\trp f(t,x,w)\mh(w)\bigl(n(x)\cdot w\bigr)\dd w. \label{eq:diffuse-normalized}
\end{align}

\subsection{Notation}\label{sec:intro-notation}

\para{General conventions} Fix $\alpha\in\left(0,\frac{1}{4}\right)$. We write $\la:=\bigl(1+\abs{v}^2\bigr)^{\frac{1}{2}}$ and $\phase:=\Om\times\R^3$. Unless otherwise indicated, constants and smallness thresholds depend only on $\Om$, $\gamma$ and $\alpha$. The letter $C$ denotes a finite positive constant that may change from line to line. Additional dependence is indicated by subscripts or stated locally. For nonnegative quantities, $A\lesssim B$ means $A\leq CB$, and $A\simeq B$ means $A\lesssim B\lesssim A$.

\para{Boundary} The notation $L^p_{\Sigma_\pm}$ means $L^p(\Sigma_\pm;\dd\gamma)$, where the boundary measure and the $L^2$ norms are
\begin{align}
\dd\gamma:=\abs[\big]{n(x)\cdot v}\dd S_x\dd v,\qquad \norm{h}_{L^2_{\Sigma_\pm}}^2:=\int_{\Sigma_\pm}\abs{h}^2\dd\gamma. \label{eq:I:display:011}
\end{align}
We write $d(x):=\dist(x,\partial\Om)$ and fix $d_0>0$ such that every $x$ in the collar $\{d(x)<d_0\}$ has a unique nearest boundary point $\bdproj(x)$. In the collar we use the local grazing weight
\begin{align}
\graz (x,v):=\bigl(d(x)\bigr)^{\frac{1}{3}}+\abs[\Big]{n\bigl(\bdproj(x)\bigr)\cdot v}. \label{eq:rho}
\end{align}
\para{Kinetic distance and H\"older spaces} For $Z=(t,x,v)$ and $Z_0=(t_0,x_0,v_0)$ in $\R\times\R^3\times\R^3$, define
\begin{align}
\dkin(Z,Z_0):=\abs{t-t_0}^{\frac{1}{2}}+\abs[\big]{x-x_0-(t-t_0)v_0}^{\frac{1}{3}}+\abs{v-v_0}. \label{eq:kinetic-distance}
\end{align}
The spatial difference measures the departure from free transport at velocity $v_0$, and the time, spatial and velocity scales are $r^2$, $r^3$ and $r$.

For $0<\beta<1$ and $\mathcal U\subset\R\times\R^3\times\R^3$, define
\begin{align}
[u]_{\Ckin^\beta(\mathcal U)}:=\sup_{\substack{Z,Z'\in\mathcal U\\Z\ne Z'}}\frac{\abs[\big]{u(Z)-u(Z')}}{\bigl(\dkin(Z,Z')\bigr)^\beta}. \label{eq:kinetic-holder-difference}
\end{align}
The space $\Ckin^\beta(\mathcal U)$ consists of functions with finite norm
\begin{align}
\norm{u}_{\Ckin^\beta(\mathcal U)}:=\sup_{Z\in\mathcal U}\abs[\big]{u(Z)}+[u]_{\Ckin^\beta(\mathcal U)}. \label{eq:kinetic-holder-norm}
\end{align}
\begin{remark}[Kinetic distance and H\"older norm]\label{rem:kinetic-holder-space}
The distance \eqref{eq:kinetic-distance} is generally nonsymmetric, but
\begin{align}
\frac{1}{2}\dkin(Z,Z_0)\leq\dkin(Z_0,Z)\leq2\dkin(Z,Z_0). \label{eq:kinetic-distance-reversal}
\end{align}
It is equivalent to the symmetric kinetic metric in \cite[(2.6)]{Ros-Oton.Weidner2025}. Thus, for $0<\beta<1$, \eqref{eq:kinetic-holder-norm} defines the same H\"older space, with equivalent norms, as \cite[Definition~2.6]{Ros-Oton.Weidner2025}. These definitions correspond to the second-order scaling $s=1$ in \cite[Definitions~2.1 and~2.3, Remark~1.11]{Imbert.Silvestre2021}.
\end{remark}
For incoming or outgoing boundary data, the norm $\norm{g}_{\Ckin^\beta\bigl((0,T)\times\Sigma_\pm\bigr)}$ is given by \eqref{eq:kinetic-holder-difference}--\eqref{eq:kinetic-holder-norm} with $\mathcal U=(0,T)\times\Sigma_\pm$. For local estimates we use
\begin{align}
u\in\Ckinloc^\beta(\mathcal U)\quad\Longleftrightarrow\quad\norm{u}_{\Ckin^\beta(K)}<\infty\quad\text{for every compact }K\subset\mathcal U. \label{eq:kinetic-holder-local}
\end{align}
In particular, $\Ckinloc^\beta\bigl([0,T]\times\overline\Om\times\R^3\bigr)$ includes the estimates on compact sets that meet the initial time or the phase boundary. The higher-order spaces used for derivative estimates extend this convention; their polynomial characterization is given in Section~\ref{lin:section-order-two}.

\para{Velocity spaces} Define
\begin{align}
\sig:=a*\m. \label{eq:sigma-definition}
\end{align}
This matrix is positive definite, with eigenvalue comparable to $\la^\gamma$ in the direction of $v$ and to $\la^{\gamma+2}$ in the orthogonal directions (Lemma~\ref{cvel:convolution}).

We use single bars for velocity norms, double bars for bulk or boundary norms, and triple bars when time is also integrated. A bulk space $L^p$ without a domain subscript means $L^p(\Om\times\R^3)$, while $L_v^p$ and $L_x^p$ indicate velocity and position separately. We retain explicit domains for boundary, mixed and local norms. For a matrix $M=(M_{ij})$, $\abs{M}=\left(\sum_{i,j}\abs{M_{ij}}^2\right)^{\frac{1}{2}}$, and the analogous sum defines the norm of a finite-dimensional tensor.

Let $\Vv$ be the completion of $C_c^\infty(\R^3)$ under the norm
\begin{align}
\vnorm{h}_{\Vv}^2:=\int_{\R^3}\left[\sig^{ij}(v)\partial_i h\partial_jh+\la^{\gamma+2}\abs{h}^2\right]\dd v. \label{eq:V}
\end{align}
The velocity space $\Wv:=L^2_v\cap \Vv$ has the sum norm
\begin{align}
\vnorm{h}_{\Wv}^2:=\vnorm{h}_{L^2_v}^2+\vnorm{h}_{\Vv}^2. \label{eq:velocity-intersection}
\end{align}
The corresponding spaces in $(x,v)$ are
\begin{align}
\Vspace:=L^2(\Om;\Vv),\qquad \Wspace:=L^2\cap\Vspace,\qquad \norm{f}_{\Wspace}^2:=\norm{f}_{L^2}^2+\norm{f}_{\Vspace}^2. \label{eq:spaces}
\end{align}
In particular, $\vnorm[\big]{f(t,x,\cdot)}_{\Vv}^2$ is the velocity energy density at $(t,x)$, and $\norm[\big]{f(t)}_{\Vspace}^2$ is its integral over $\Om$.

For $\gamma\geq-2$ we have $\Vv\subset L^2_v$. For $\gamma<-2$ the weight $\la^{\gamma+2}$ decays and this inclusion fails, so the $L^2$ norm is kept as a separate part of the energy for every $\gamma$. We use the continuous dense inclusions $\Wv\subset L^2_v\subset \Wv'$ and $\Wspace\subset L^2\subset\Wspace'$, where the dualities agree with the $L^2$ pairing. The dual space $\Vv'$ is paired with $\Vv$ through the extension of the $L^2$ pairing $\ip*{\cdot}{\cdot}$ from compactly supported functions, and the same convention is used for $\Vspace'$ and $\Vspace$.

\para{Gaussian norms} For $1\leq p\leq\infty$ and $\zeta\in\R$, the second superscript of $L^{p,\zeta}$ records the Gaussian weight. The bulk and velocity norms are
\begin{align}
\gnorm{p}{\zeta}{f}:=\norm[\Big]{\ue^{\zeta\abs{v}^2}f}_{L^p},\qquad \vgnorm{p}{\zeta}{h}:=\vnorm[\Big]{\ue^{\zeta\abs{v}^2}h}_{L_v^p}. \label{eq:Gaussian-norms}
\end{align}
Using the boundary measure in \eqref{eq:I:display:011}, we define
\begin{align}
\bgnorm{p}{\zeta}{\pm}{h}:=\norm[\Big]{\ue^{\zeta\abs{v}^2}h}_{L^p_{\Sigma_\pm}}. \label{eq:Xalpha}
\end{align}
We use the usual Bochner notation $L^q(0,T;L^{p,\zeta})$ for mixed time and phase space. When the same exponent is used in time and all remaining variables, $\tgnorm{p}{\zeta}{f}$ and $\tbgnorm{p}{\zeta}{\pm}{h}$ denote the corresponding norms with measures $\dd t\dd x\dd v$ and $\dd t\dd\gamma$ on the stated time interval.

The exponent $\alpha$ measures the Gaussian decay of the perturbation $f$. The corresponding weight for the physical perturbation is given by
\begin{align}
\gnorm{\infty}{\alpha+\frac{1}{4}}{F-\m}=(2\pi)^{-\frac{3}{4}}\gnorm{\infty}{\alpha}{f}. \label{eq:Gaussian-normalization}
\end{align}
Since $\alpha>0$, the Gaussian bound gives
\begin{align}
\gnorm{\infty}{\alpha}{f}\leq\delta\quad\Longrightarrow\quad \abs{f}\leq\delta\ue^{-\alpha\abs{v}^2},\qquad \norm{f}_{L^\infty_xL^2_v}\lesssim\delta. \label{eq:gaussian-energy}
\end{align}
Integration against the boundary measure gives
\begin{align}
\norm{g}_{L^2_{\Sigma_-}}\lesssim\bgnorm{\infty}{\alpha}{-}{g}. \label{eq:incoming-norm}
\end{align}

\subsection{Solution Classes and Kinetic Traces}\label{sec:intro-classes}

Our weak formulation follows the kinetic Sobolev and solution--trace formulations of Silvestre \cite[Section~4 and Proposition~4.3]{Silvestre2022} and Zhu \cite[Definition~2.1 and Remarks~2.2--2.3]{Zhu2024}. These formulations control $\nabla_vf$ and the transport derivative in a velocity dual space. For Landau, the trace is determined by the Green identity in Carrapatoso--Mischler's renormalized trace theorem \cite[Theorem~2.8]{Carrapatoso.Mischler2024}, based on Mischler's boundary renormalization theory \cite{Mischler2010}.

The Gaussian bound \eqref{eq:weak-gaussian-bound} is part of our solution class and will be established in Section~\ref{sec:energy}. Lemma~\ref{trace:gaussian} transfers it to the trace and gives the integrability stated below.

\para{Collision terms and weak formulation} Fix $T>0$. We assume
\begin{align}
f\in C\bigl([0,T];L^2\bigr)\cap L^2(0,T;\Vspace), \label{eq:intro-trace-space}
\end{align}
and impose the Gaussian bound
\begin{align}
f\in L^\infty(0,T;L^{\infty,\alpha}). \label{eq:weak-gaussian-bound}
\end{align}
By \eqref{eq:gaussian-energy}, the Gaussian bound also gives $f\in L^\infty(0,T;L_x^\infty L_v^2)$.

For almost every $t$, we define the nonlinear collision term by its action on $\varphi\in C_c^\infty(\Om\times\R^3)$:
\begin{align}
\ip[\big]{\Ga(f,f)}{\varphi}=-\int_{\Om\times\R^3}\Biggl\{\bigl(a*(\mh f)\bigr)\left(\nabla_vf-\frac{v}{2}f\right)-\bigl(b*(\mh f)\bigr)f\Biggr\}\cdot\left(\nabla_v\varphi+\frac{v}{2}\varphi\right)\dd x\dd v. \label{eq:intro-collision-pairing}
\end{align}
This is the integration-by-parts form of \eqref{eq:collision}; the same pairing defines both terms of $\Lop f$ in \eqref{eq:linearized-operators}.

Corollary~\ref{cvel:spatial} shows that these pairings extend continuously to $\varphi\in\Vspace$ and satisfy
\begin{align}
\tnorm{\Lop f}_{L^2(0,T;\Vspace')}+\tnorm[\big]{\Ga(f,f)}_{L^2(0,T;\Vspace')}\lesssim \Bigl(1+\tnorm{f}_{L^\infty(0,T;L^\infty_xL^2_v)}\Bigr)\tnorm{f}_{L^2(0,T;\Vspace)}<\infty. \label{eq:intro-collision-dual}
\end{align}
We interpret \eqref{eq:normalized-problem} by requiring that, for every $\varphi\in C_c^\infty\bigl((0,T)\times\Om\times\R^3\bigr)$,
\begin{align}
-\int_0^T\bigl(f,\Y\varphi\bigr)_{L^2}\dd t+\int_0^T\ip*{\Lop f}{\varphi}\dd t=\int_0^T\ip[\big]{\Ga(f,f)}{\varphi}\dd t. \label{eq:intro-interior-weak}
\end{align}
Every term is finite by \eqref{eq:intro-trace-space}--\eqref{eq:intro-collision-dual}.

\para{Kinetic trace} Suppose that $f$ satisfies \eqref{eq:intro-trace-space} and $\Y f\in L^2(0,T;\Vspace')$. Lemma~\ref{trace:CM} gives a unique kinetic trace $\trace f$, determined by the linear Green identity \eqref{trace:linear-Green}, with
\begin{align}
(n\cdot v)\trace f\in L^2_{\mathrm{loc}}\left((0,T)\times\partial\Om\times\R^3;\dd t\dd S_x\dd v\right). \label{eq:intro-trace-local}
\end{align}
Thus the local trace estimate has the weight $\abs{n\cdot v}^2$. The trace is unique up to a set of zero $\dd t\dd\gamma$ measure, and its restrictions to $\Sigma_\pm$ are denoted by $\trpm f$.

For the Landau equation, the Gaussian bound \eqref{eq:weak-gaussian-bound} and collision estimate \eqref{eq:intro-collision-dual} justify
\begin{align}
\Y f=\Ga(f,f)-\Lop f\in L^2(0,T;\Vspace'), \label{eq:intro-transport-dual}
\end{align}
so every $f$ satisfying \eqref{eq:intro-trace-space}, \eqref{eq:weak-gaussian-bound} and \eqref{eq:intro-interior-weak} has this trace. Renormalization then transfers the Gaussian bound to the trace, as proved in Lemma~\ref{trace:gaussian}:
\begin{align}
\tnorm{\trpm f}_{L^\infty(0,T;L^{\infty,\alpha}_{\Sigma_\pm})}\leq\tnorm{f}_{L^\infty(0,T;L^{\infty,\alpha})}. \label{eq:intro-trace-gaussian}
\end{align}
Integrating \eqref{eq:intro-trace-gaussian} over the finite time interval, the spatial boundary and all velocities gives
\begin{align}
\trace f\in L^2\left((0,T)\times\partial\Om\times\R^3;\dd t\dd S_x\dd v\right). \label{eq:intro-trace-unweighted}
\end{align}
The Gaussian decay also makes the diffuse-reflection integral
\begin{align}
\flux{f}(t,x):=\int_{n(x)\cdot w>0}\trp f(t,x,w)\mh(w)\bigl(n(x)\cdot w\bigr)\dd w \label{eq:diffuse-flux}
\end{align}
absolutely convergent for almost every $(t,x)\in(0,T)\times\partial\Om$. Integration with the flux measure $\dd\gamma=\abs{n\cdot v}\dd S_x\dd v$ further yields
\begin{align}
\trpm f\in L^2(0,T;L^2_{\Sigma_\pm}),\qquad \flux{f}\in L^\infty\bigl((0,T)\times\partial\Om\bigr). \label{eq:intro-trace-full}
\end{align}
In particular, $\Pm\trp f=\sqrt{2\pi}\mh\flux{f}\in L^2(0,T;L^2_{\Sigma_-})$, as required in the boundary energy identities.

Lemma~\ref{trace:CM}(iii) identifies the trace of the full density as
\begin{align}
\trpm F=\m+\mh\trpm f. \label{eq:intro-density-trace}
\end{align}
The boundary laws for $F$ and $f$ are therefore equivalent as equalities of these traces.

\begin{definition}[Weak solution]\label{def:energy}
Let $f_0\in L^2(\Om\times\R^3)$. In the in-flow case, let $g\in L^2\bigl((0,T)\times\Sigma_-;\dd t\dd\gamma\bigr)$ be given. A weak solution on $[0,T]$ is a function $f$ satisfying \eqref{eq:intro-trace-space}, \eqref{eq:weak-gaussian-bound} and \eqref{eq:intro-interior-weak}, with
\begin{align}
f(0)=f_0\qquad\text{in }L^2(\Om\times\R^3), \label{eq:intro-initial-attainment}
\end{align}
and one of the following boundary conditions:
\begin{enumerate}
\item In the in-flow case, $\trm f=g$ almost everywhere on $(0,T)\times\Sigma_-$ with respect to $\dd t\dd\gamma$, as in \eqref{eq:inflow-normalized}.
\item In the diffuse-reflection case, $\trm f=\Pm\trp f$ almost everywhere on $(0,T)\times\Sigma_-$ with respect to $\dd t\dd\gamma$, as in \eqref{eq:diffuse-normalized}.
\end{enumerate}
Here $\trpm f$ are the traces determined by \eqref{trace:linear-Green}, with the integrability stated in \eqref{eq:intro-trace-unweighted} and \eqref{eq:intro-trace-full}. Both boundary identities are equalities in $L^2\bigl((0,T)\times\Sigma_-;\dd t\dd\gamma\bigr)$. The associated density is $F=\m+\mh f$, with initial data $F_0=\m+\mh f_0$. We call the solution nonnegative when $F\geq0$ almost everywhere.
\end{definition}

The definition and Lemma~\ref{trace:CM}(iv) imply the following Green identity. For $0\leq s\leq t\leq T$ and every $\varphi\in C^\infty\bigl([0,T]\times\overline\Om\times\R^3\bigr)$ with compact support in velocity,
\begin{align}
&\bigl(f(t),\varphi(t)\bigr)_{L^2}-\bigl(f(s),\varphi(s)\bigr)_{L^2}-\int_s^t\bigl(f,\Y\varphi\bigr)_{L^2}\dd r \label{eq:green-weak}\\
&\qquad+\int_s^t\int_{\partial\Om\times\R^3}\trace f\,\varphi\,(n\cdot v)\dd S_x\dd v\dd r=\int_s^t\ip[\big]{\Ga(f,f)-\Lop f}{\varphi}\dd r. \notag
\end{align}
All terms are finite by \eqref{eq:intro-collision-dual} and \eqref{eq:intro-trace-full}.

Under \eqref{eq:intro-trace-space} and \eqref{eq:weak-gaussian-bound}, the formulations \eqref{eq:intro-interior-weak} and \eqref{eq:green-weak} are equivalent by Lemma~\ref{trace:CM}, with the same initial and boundary conditions as in Definition~\ref{def:energy}.

\begin{definition}[Strong solution]\label{def:strong}
We say that $f$ is a strong solution on $[0,T]$ if it is a weak solution in the sense of Definition~\ref{def:energy} and satisfies
\begin{align}
\int_0^T\int_{\Om\times\R^3}\left[\sum_{i,j=1}^3\abs{\partial_{v_i v_j}f}^2+\abs{\Y f}^2+\abs{\Lop f}^2+\abs[\big]{\Ga(f,f)}^2\right]\dd x\dd v\dd t<\infty. \label{eq:I:display:016}
\end{align}
Here $\partial_{v_i v_j}f$ and $\Y f$ are distributional derivatives represented by $L^2$ functions, and the collision terms are the distributions defined by \eqref{eq:intro-collision-pairing}--\eqref{eq:intro-collision-dual} and \eqref{eq:linearized-operators}. The weak identity \eqref{eq:intro-interior-weak} then gives
\begin{align}
\Y f+\Lop f=\Ga(f,f)\qquad\text{in }L^2\bigl((0,T)\times\Om\times\R^3\bigr). \label{eq:intro-strong-equation}
\end{align}
The initial and boundary conditions are those of Definition~\ref{def:energy}.
\end{definition}

Definition~\ref{def:strong} controls $\Y f$ as a whole and permits the pointwise grazing singularities quantified in Theorems~\ref{main:inflow}(ii) and \ref{main:diffuse}(ii), whose squares are integrable as required by \eqref{eq:I:display:016}.

\subsection{Main Results}\label{sec:intro-results}

\para{In-flow} The data for \eqref{eq:inflow-full} satisfy the following assumptions.

\begin{assumption}\label{in:data:regularity}
Let $f_0\in C^5(\overline\Om\times\R^3)$ and $g\in C^5\bigl([0,\infty)\times\overline\Om\times\R^3\bigr)$, with the incoming data given by the values of $g$ on $\Sigma_-$.

Assume that there is $C_0<\infty$ and, for every $T>0$, a constant $C_T<\infty$ such that
\begin{align}
\sum_{\abs{k}\leq5}\gnorm{\infty}{\alpha}{\partial_{x,v}^k f_0} &\leq C_0, \label{in:eq:data-derivatives}\\
\sum_{a+\abs{b}+\abs{c}\leq5}\tnorm{\partial_t^a\partial_x^b\partial_v^c g}_{L^\infty(0,T;L^{\infty,\alpha}_{\Sigma_-})} &\leq C_T. \label{in:eq:data-derivatives-g}
\end{align}
The constants in these derivative bounds need not be small.

Assume $\m+\mh f_0\geq0$ in $\Om\times\R^3$ and $\m+\mh g\geq0$ on $[0,\infty)\times\Sigma_-$. The initial and incoming data satisfy the compatibility conditions
\begin{align}
\trm f_0&=g(0),\label{in:eq:compat0}\\
\trm a_0&=\partial_tg(0),\qquad a_0:=-v\cdot\nabla_xf_0-\Lop f_0+\Ga(f_0,f_0).\label{in:eq:compat1}
\end{align}
Both equalities are classical identities on $\Sigma_-$ at $t=0$.
\end{assumption}

\begin{remark}\label{rem:inflow-data}
Under the remaining hypotheses of Assumption~\ref{in:data:regularity}, Theorem~\ref{main:inflow} also holds for $g$ defined only on $[0,\infty)\times\overline{\Sigma_-}$, with continuous mixed derivatives satisfying, for every $T>0$,
\begin{align}
\sum_{a+\abs{b}+\abs{c}\leq5,\,\abs{b}\leq2}\tnorm{\partial_t^a\partial_{\mathrm{tan}}^b\partial_v^c g}_{L^\infty(0,T;L^{\infty,\alpha}_{\Sigma_-})}\leq C_T. \label{in:eq:data-derivatives-g-weaker}
\end{align}
Here $\partial_{\mathrm{tan}}^b$ denotes derivatives in a fixed finite collection of $C^3$ boundary charts, with $v$ held fixed; the bound holds in each chart. These are the derivatives used in Lemmas~\ref{cb:inflow-lift} and \ref{cb:inflow-initial}. The latter uses \eqref{in:eq:compat0} for the first-order estimate at $t=0$ and \eqref{in:eq:compat1} for the second-order estimate.

For $f_0\in C_c^\infty(\Om\times\R^3)$, both compatibility conditions hold with the persistent incoming data
\begin{align}
g(t,x,v):=\varepsilon_b\chi(t)\beta(x)\mh(v), \label{in:eq:switch-on-example}
\end{align}
where $\chi\in C^\infty\bigl([0,\infty)\bigr)$ vanishes near $t=0$ and equals one for large times, $\beta\in C^2(\partial\Om)$, and $1+\varepsilon_b\chi\beta\geq0$. This example satisfies \eqref{in:eq:data-derivatives-g-weaker}; the Gaussian amplitudes in \eqref{in:eq:data} can be small independently of the derivative bounds.
\end{remark}

\begin{theorem}\label{main:inflow}
There are $\delta_*>0$, $C_*\geq1$, $\lambda_*>0$ and $\varepsilon_*>0$ such that
\begin{align}
C_*\varepsilon_*\leq\frac{\delta_*}{2}. \label{in:eq:thresholds}
\end{align}
If the data satisfy Assumption~\ref{in:data:regularity} and
\begin{align}
\gnorm{\infty}{\alpha}{f_0}+\tnorm{g}_{L^\infty(0,\infty;L^{\infty,\alpha}_{\Sigma_-})}\leq\varepsilon_*, \label{in:eq:data}
\end{align}
then there is a global weak solution $f$, with $F=\m+\mh f\geq0$, of \eqref{eq:normalized-problem} with initial data $f_0$ and in-flow condition \eqref{eq:inflow-normalized}, in the sense of Definition~\ref{def:energy}. It has the following properties.
\begin{enumerate}
\item \textbf{Global estimates.} The solution satisfies
\begingroup\usetagform{enumerated}
\begin{align}
&\tnorm{f}_{L^\infty(0,\infty;L^2)}+\tnorm{f}_{L^\infty(0,\infty;L^{\infty,\alpha})} \label{in:eq:main-uniform}\leq C_*\left(\gnorm{\infty}{\alpha}{f_0}+\tnorm{g}_{L^\infty(0,\infty;L^{\infty,\alpha}_{\Sigma_-})}\right)\leq\frac{\delta_*}{2}, \\
&\norm[\big]{f(t)}_{L^2}^2\leq C_*\ue^{-2\lambda_* t}\norm{f_0}_{L^2}^2+C_*\int_0^t\ue^{-2\lambda_*(t-s)}\norm[\big]{g(s)}_{L^2_{\Sigma_-}}^2\dd s. \label{in:eq:main-energy-forcing}
\end{align}
\endgroup
For every $T>0$,
\begingroup\usetagform{enumerated}
\begin{align}
\int_0^T\Bigl(\norm[\big]{f(t)}_{\Vspace}^2+\norm[\big]{\trp f(t)}_{L^2_{\Sigma_+}}^2\Bigr)\dd t\leq C_*\left(\norm{f_0}_{L^2}^2+\int_0^T\norm[\big]{g(t)}_{L^2_{\Sigma_-}}^2\dd t\right). \label{in:eq:main-dissipation}
\end{align}
\endgroup
The constants $C_*$, $\lambda_*$, $\delta_*$ and $\varepsilon_*$ are independent of $T$ and of the derivative bounds in Assumption~\ref{in:data:regularity}.
\item \textbf{Boundary regularity and strong solutions.} For every $T>0$,
\begingroup\usetagform{enumerated}
\begin{align}
f\in\Ckinloc^{\frac{1}{2}}\left([0,T]\times\overline\Om\times\R^3\right). \label{in:eq:main-holder}
\end{align}
\endgroup
There are $C_T>0$ and $\alpha_T\in(0,\alpha)$ such that, for $0<t<T$ and $0<d(x)<d_0$,
\begingroup\usetagform{enumerated}
\begin{align}
\ue^{\alpha_T\abs{v}^2}\abs{\nabla_vf}&\leq C_T\bigl(1+\graz^{-\frac{1}{2}}\bigr), \label{in:eq:main-gradient}\\
\ue^{\alpha_T\abs{v}^2}\bigl(\abs{D_v^2f}+\abs{\Y f}\bigr)&\leq C_T\bigl(1+\graz^{-\frac{3}{2}}\bigr). \label{in:eq:main-hessian}
\end{align}
\endgroup
For $d(x)\geq d_0$, these bounds hold with $\graz$ replaced by $1$. In particular, $f$ is a strong solution on the time interval $[0,T]$ in the sense of Definition~\ref{def:strong}, and for every $m\geq0$,
\begingroup\usetagform{enumerated}
\begin{align}
\la^mD_v^2f,\quad\la^m\Y f,\quad\la^m\Lop f,\quad\la^m\Ga(f,f)\in L^2\bigl((0,T)\times\Om\times\R^3\bigr). \label{in:eq:main-strong}
\end{align}
\endgroup
The constants $C_T$, $\alpha_T$ and the norms in \eqref{in:eq:main-strong} may depend on $T$ and the derivative bounds of the data. The norms in \eqref{in:eq:main-strong} may also depend on $m$.
\item \textbf{Uniqueness and continuous dependence.} Let $\widetilde f$ be a weak solution with the same initial and incoming data. If, for every $T>0$,
\begingroup\usetagform{enumerated}
\begin{align}
\tnorm[\big]{\widetilde f}_{L^\infty(0,T;L^{\infty,\alpha})}\leq\delta_*, \label{in:eq:comparison-class}
\end{align}
\endgroup
then $\widetilde f=f$.

Fix $T>0$. Let $f_j$ be weak solutions with incoming data $g_j$ and $\tnorm{f_j}_{L^\infty(0,T;L^{\infty,\alpha})}\leq\delta_*$. If $f_j(0)\to f_0$ in $L^2$ and $g_j\to g$ in $L^2\bigl((0,T)\times\Sigma_-;\dd t\dd\gamma\bigr)$, then
\begingroup\usetagform{enumerated}
\begin{align}
f_j&\to f\quad\text{in }C\bigl([0,T];L^2\bigr)\cap L^2(0,T;\Vspace), \label{in:eq:main-stability}\\
\trpm f_j&\to\trpm f\quad\text{in }L^2\bigl((0,T)\times\Sigma_\pm;\dd t\dd\gamma\bigr). \label{in:eq:main-stability-trace}
\end{align}
\endgroup
\end{enumerate}
\end{theorem}

Theorem~\ref{main:inflow} gives bounded global dynamics under persistent incoming forcing. When the incoming perturbation decays, the solution relaxes to the Maxwellian.

\begin{corollary}\label{in:cor:decay}
Under the hypotheses of Theorem~\ref{main:inflow}, assume in addition that
\begin{align}
\bgnorm[\big]{\infty}{\alpha}{-}{g(t)}\leq A \ue^{-\omega t},\qquad 0\leq A<\infty,\quad \omega>0. \label{in:eq:inflow-decay}
\end{align}
There is $0<\lambda_\omega\leq\omega$ such that
\begin{align}
\norm[\big]{f(t)}_{L^2}+\gnorm[\big]{\infty}{\alpha}{f(t)} \lesssim_\omega\left(\gnorm{\infty}{\alpha}{f_0}+A\right)\ue^{-\lambda_\omega t}, \label{in:eq:main-decay}
\end{align}
and, for every $t\geq0$,
\begin{align}
\int_t^\infty\Bigl(\norm[\big]{f(s)}_{\Vspace}^2+\norm[\big]{\trp f(s)}_{L^2_{\Sigma_+}}^2\Bigr)\dd s\lesssim_\omega\left(\gnorm{\infty}{\alpha}{f_0}+A\right)^2\ue^{-2\lambda_\omega t}. \label{in:eq:decay-dissipation-flux}
\end{align}
The decay rate and implicit constants are independent of the derivative bounds of the data. If only $\bgnorm[\big]{\infty}{\alpha}{-}{g(t)}\to0$ as $t\to\infty$ is assumed, then
\begin{align}
\norm[\big]{f(t)}_{L^2}+\gnorm[\big]{\infty}{\alpha}{f(t)}\to0\qquad(t\to\infty). \label{eq:I:display:029}
\end{align}
\end{corollary}

\para{Diffuse-reflection} The data for \eqref{eq:diffuse-full} satisfy the following assumptions.

\begin{assumption}\label{def:data}
Let $f_0\in C^5(\overline\Om\times\R^3)$. Assume that there is $C_0<\infty$ such that
\begin{align}
\sum_{\abs{k}+\abs{l}\leq5}\gnorm{\infty}{\alpha}{\partial_x^k\partial_v^l f_0}\leq C_0. \label{eq:data-derivatives}
\end{align}
The constant $C_0$ need not be small.

Assume $\m+\mh f_0\geq0$ and
\begin{align}
\int_{\Om\times\R^3}\mh f_0\dd x\dd v=0. \label{eq:mass0}
\end{align}
Define the initial time derivative by
\begin{align}
a_0:=-v\cdot\nabla_xf_0-\Lop f_0+\Ga(f_0,f_0). \label{eq:a0}
\end{align}
The initial data and this derivative satisfy the compatibility conditions
\begin{align}
\trm f_0=\Pm\trp f_0,\qquad \trm a_0=\Pm\trp a_0. \label{eq:compatibility}
\end{align}
\end{assumption}

The wall exchanges momentum and energy with the gas, so \eqref{eq:mass0} is the only moment condition. For example, $f_0=\varepsilon_b\beta(x)\psi(v)$ with $\beta\in C_c^\infty(\Om)$, $\int_\Om\beta\dd x=0$, $\psi\in C_c^\infty(\R^3)$ and sufficiently small $\varepsilon_b$ satisfies Assumption~\ref{def:data}: $f_0$ and $a_0$ vanish near $\partial\Om\times\R^3$, giving \eqref{eq:compatibility}.

\begin{theorem}\label{main:diffuse}
There are $\delta_*>0$, $C_*\geq1$, $\lambda_*>0$ and $\varepsilon_*>0$ such that
\begin{align}
C_*\varepsilon_*\leq\frac{\delta_*}{2}. \label{eq:diffuse-thresholds}
\end{align}
If the initial data satisfy Assumption~\ref{def:data} and
\begin{align}
\gnorm{\infty}{\alpha}{f_0}\leq\varepsilon_*, \label{eq:small-data}
\end{align}
then there is a global weak solution $f$, with $F=\m+\mh f\geq0$, of \eqref{eq:normalized-problem} with initial data $f_0$ and diffuse-reflection condition \eqref{eq:diffuse-normalized}, in the sense of Definition~\ref{def:energy}. It has the following properties.
\begin{enumerate}
\item \textbf{Global estimates.} Set
\begingroup\usetagform{enumerated}
\begin{align}
p_\gamma:=\begin{cases}1,&-2\leq\gamma\leq1,\\\frac{2}{\abs{\gamma}},&-3\leq\gamma<-2.\end{cases} \label{full:decay-exponent}
\end{align}
\endgroup
The solution satisfies
\begingroup\usetagform{enumerated}
\begin{align}
&\int_{\Om\times\R^3}\mh f(t)\dd x\dd v=0,\qquad t\geq0, \label{eq:main-mass}\\
&\tnorm{f}_{L^\infty(0,\infty;L^2)}+\tnorm{f}_{L^\infty(0,\infty;L^{\infty,\alpha})}\leq C_*\gnorm{\infty}{\alpha}{f_0}\leq\frac{\delta_*}{2}, \label{eq:main-decay}\\
&\norm[\big]{f(t)}_{L^2}^2\leq C_*\gnorm{\infty}{\alpha}{f_0}^2\ue^{-2\lambda_* t^{p_\gamma}},\qquad t\geq0, \label{eq:main-energy-decay}\\
&\int_0^\infty\Bigl(\norm{f}_{\Vspace}^2+\norm{\trp f}_{L^2_{\Sigma_+}}^2+\norm{\trm f}_{L^2_{\Sigma_-}}^2\Bigr)\dd t\leq C_*\gnorm{\infty}{\alpha}{f_0}^2. \label{eq:main-decay-dissipation}
\end{align}
\endgroup
The constants $C_*$, $\lambda_*$, $\delta_*$ and $\varepsilon_*$ are independent of $T$ and of the derivative bounds in Assumption~\ref{def:data}.
\item \textbf{Boundary regularity and strong solutions.} For every $T>0$ there is $\kappa_T\in\left(\frac{2}{3},1\right)$ such that
\begingroup\usetagform{enumerated}
\begin{align}
f\in\Ckinloc^{\frac{1}{2}}\left([0,T]\times\overline\Om\times\R^3\right),\qquad\flux{f}\in C^{\kappa_T}\left([0,T]\times\partial\Om\right). \label{eq:main-holder}
\end{align}
\endgroup
There are $C_T>0$ and $\alpha_T\in(0,\alpha)$ such that, for $0<t<T$ and $0<d(x)<d_0$,
\begingroup\usetagform{enumerated}
\begin{align}
\ue^{\alpha_T\abs{v}^2}\abs{\nabla_vf}&\leq C_T\bigl(1+\graz^{-\frac{1}{2}}\bigr), \label{eq:main-gradient}\\
\ue^{\alpha_T\abs{v}^2}\bigl(\abs{D_v^2f}+\abs{\Y f}\bigr)&\leq C_T\bigl(1+\graz^{-\frac{3}{2}}\bigr). \label{eq:main-hessian}
\end{align}
\endgroup
For $d(x)\geq d_0$, these bounds hold with $\graz$ replaced by $1$. In particular, $f$ is a strong solution on the time interval $[0,T]$ in the sense of Definition~\ref{def:strong}, and for every $m\geq0$,
\begingroup\usetagform{enumerated}
\begin{align}
\la^mD_v^2f,\quad\la^m\Y f,\quad\la^m\Lop f,\quad\la^m\Ga(f,f)\in L^2\bigl((0,T)\times\Om\times\R^3\bigr). \label{eq:main-graph}
\end{align}
\endgroup
The constants $C_T$, $\alpha_T$, $\kappa_T$ and the norms in \eqref{eq:main-graph} may depend on $T$ and the derivative bounds of the data. The norms in \eqref{eq:main-graph} may also depend on $m$.
\item \textbf{Uniqueness and continuous dependence.} Let $\widetilde f$ be a weak solution with the same initial data and diffuse-reflection condition. If, for every $T>0$,
\begingroup\usetagform{enumerated}
\begin{align}
\tnorm[\big]{\widetilde f}_{L^\infty(0,T;L^{\infty,\alpha})}\leq\delta_*, \label{eq:diffuse-comparison-class}
\end{align}
\endgroup
then $\widetilde f=f$.

Fix $T>0$. Let $f_j$ be weak solutions with diffuse-reflection and $\tnorm{f_j}_{L^\infty(0,T;L^{\infty,\alpha})}\leq\delta_*$. If $f_j(0)\to f_0$ in $L^2$, then
\begingroup\usetagform{enumerated}
\begin{align}
f_j&\to f\quad\text{in }C\bigl([0,T];L^2\bigr)\cap L^2(0,T;\Vspace), \label{eq:main-stability}\\
\trpm f_j&\to\trpm f\quad\text{in }L^2\bigl((0,T)\times\Sigma_\pm;\dd t\dd\gamma\bigr), \label{eq:main-stability-trace}\\
\flux{f_j}&\to\flux{f}\quad\text{in }L^2\bigl((0,T)\times\partial\Om\bigr). \label{eq:main-stability-flux}
\end{align}
\endgroup
\end{enumerate}
\end{theorem}

\begin{corollary}[Decay under diffuse-reflection]\label{main:diffuse-decay}
Under the hypotheses of Theorem~\ref{main:diffuse}, let $p_\gamma$ be given by \eqref{full:decay-exponent}. For every $\alpha'$ with $0<\alpha'<\alpha$, there is $\lambda_{\alpha'}>0$ such that
\begin{align}
\gnorm[\big]{\infty}{\alpha'}{f(t)}+\norm[\big]{f(t)}_{L^2}\lesssim_{\alpha'}\gnorm{\infty}{\alpha}{f_0} \ue^{-\lambda_{\alpha'}t^{p_\gamma}},\qquad t\geq0. \label{eq:main-decay-rate}
\end{align}
For $-2\leq\gamma\leq1$, with $\lambda_*$ as in Theorem~\ref{main:diffuse},
\begin{align}
\gnorm[\big]{\infty}{\alpha}{f(t)}+\norm[\big]{f(t)}_{L^2}\lesssim\gnorm{\infty}{\alpha}{f_0}\ue^{-\lambda_* t},\qquad t\geq0. \label{full:same-weight-decay}
\end{align}
All constants are independent of the derivative bounds in \eqref{eq:data-derivatives}.
\end{corollary}

The different decay rates reflect the boundary laws. The weighted energy and Gaussian estimates in Propositions~\ref{ci:energy-response} and \ref{ci:gaussian-bootstrap} give the exponential in-flow estimate of Corollary~\ref{in:cor:decay} for every $\gamma\in[-3,1]$. The diffuse rates in Corollary~\ref{main:diffuse-decay} follow from Proposition~\ref{cex:CM-input} and the Gaussian interpolation in Corollary~\ref{cex:decay-traces}.

\subsection{Background and Related Results}\label{sec:intro-background}

Landau \cite{Landau1936} introduced the equation in 1936 to describe collisions between charged particles interacting through Coulomb forces. Early mathematical work studied its derivation from the Boltzmann equation when grazing collisions dominate, including the linearized limit proved by Desvillettes \cite{Desvillettes1992}. The analysis then developed around weak existence, regularization and convergence to equilibrium. We refer to Villani \cite{Villani2002} for a broad account. We will first recall the theory without a boundary, then discuss bounded domains and the trace and regularity estimates relevant to the present problem.

\para{The Landau equation without a boundary} For the spatially homogeneous equation, Villani \cite{Villani1998} developed the notion of an $H$-solution using entropy dissipation. Desvillettes--Villani \cite{Desvillettes.Villani2000} proved global existence, uniqueness and instantaneous smoothing for hard potentials. Fournier--Gu\'erin \cite{Fournier.Guerin2009} obtained stability in the quadratic Wasserstein distance and well-posedness for soft potentials, with global results for moderately soft potentials and local results for very soft potentials under suitable moment and integrability assumptions. In the Coulomb case, Fournier \cite{Fournier2010} proved uniqueness and stability for bounded solutions. Desvillettes \cite{Desvillettes2015} controlled a weighted $H_v^1$ norm of $F^{\frac{1}{2}}$ by entropy dissipation, and Carrapatoso--Desvillettes--He \cite{Carrapatoso.Desvillettes.He2017} obtained quantitative convergence to Maxwellian equilibrium under moment assumptions.

More recently, Guillen--Silvestre \cite{Guillen.Silvestre2025} proved monotonicity of the Fisher information and global smooth solvability for suitable smooth, decaying data, including Coulomb interactions. Golding--Gualdani--Loher \cite{Golding.Gualdani.Loher2025} obtained global smooth solutions for polynomially weighted $L_v^{\frac{3}{2}}$ data satisfying additional moment and lower-bound assumptions.

For the spatially inhomogeneous equation, Villani \cite{Villani1996} proved weak sequential stability and global existence of renormalized solutions with a defect measure. Near a Maxwellian, Guo \cite{Guo2002(=)} constructed global classical solutions in a periodic box, including Coulomb interactions. Strain--Guo \cite{Strain.Guo2006,Strain.Guo2008} obtained quantitative relaxation for soft potentials through velocity-weighted energy estimates. Semigroup methods extended the perturbative theory in periodic domains to other weighted spaces \cite{Carrapatoso.Tristani.Wu2016,Carrapatoso.Tristani.Wu2017Erratum,Carrapatoso.Mischler2017}, while Kim--Guo--Hwang \cite{Kim.Guo.Hwang2020} developed an $L^2$--$L^\infty$ approach. Recent work of Chen--Nguyen--Yang \cite{Chen.Nguyen.Yang2025} treats small perturbations in the whole space with a scaling-critical norm for very soft potentials.

For large inhomogeneous data, Henderson--Snelson--Tarfulea \cite{Henderson.Snelson.Tarfulea2019,Henderson.Snelson.Tarfulea2020(=)} established local existence for soft potentials, including Coulomb interactions, first with Sobolev regularity and Gaussian decay and then for bounded data with polynomial velocity decay under additional positivity assumptions. Conditional regularity connects these local theories to continuation criteria. Cameron--Silvestre--Snelson \cite{Cameron.Silvestre.Snelson2018} obtained upper bounds and Gaussian propagation for moderately soft potentials under bounds on the hydrodynamic densities. Henderson--Snelson \cite{Henderson.Snelson2020} proved conditional $C^\infty$ regularity, with further moment assumptions for very soft potentials, and Snelson \cite{Snelson2020} established Gaussian bounds and conditional regularity for hard potentials.

\para{Landau equations in bounded domains} Duan--Liu--Sakamoto--Strain \cite{Duan.Liu.Sakamoto.Strain2020} constructed unique global mild solutions in a finite channel with in-flow or specular reflection, for the full range $-3\leq\gamma\leq1$. Their method combines velocity estimates with Fourier summability in the tangential variables and Sobolev regularity in the normal variable. Liu--Yang--Ma \cite{Liu.Yang.Ma2020} obtained global strong solutions with in-flow in a finite channel under anisotropic Sobolev assumptions on the initial and incoming data.

Guo--Hwang--Jang--Ouyang \cite{Guo.Hwang.Jang.Ouyang2020,Guo.Hwang.Jang.Ouyang2021Correction} established global stability near a Maxwellian for Coulomb Landau with specular reflection in smooth bounded domains. Their boundary regularity argument uses simultaneous changes of position and velocity and reflection across the boundary. Dong--Guo--Ouyang \cite{Dong.Guo.Ouyang2020} developed an $L^2$--$S^p$ approach for the Vlasov--Poisson--Landau system. For the relativistic Vlasov--Maxwell--Landau system with specular reflection and perfectly conducting electromagnetic boundary conditions, local well-posedness and global stability were established in \cite{Dong.Guo.Ouyang.Yastrzhembskiy2024,Dong.Guo.Yastrzhembskiy2024}.

The reflection argument is based on the specular identity
\begin{align}
\trm f(t,x,v)=\trp f\Bigl(t,x,v-2\bigl(n(x)\cdot v\bigr)n(x)\Bigr). \label{eq:intro-specular-identity}
\end{align}
After localization and flattening, \eqref{eq:intro-specular-identity} matches the incoming trace with the reflected outgoing trace. Reflection then extends the localized equation to the whole space without a boundary-supported term. Whole-space kinetic $S^p$ estimates and Sobolev embedding give $\nabla_vf\in L^\infty$. For in-flow and diffuse-reflection, the two traces generally fail to match under this reflection. The extension has a jump across the flattened boundary and its equation contains a boundary-supported term, so the same reduction to the whole-space $S^p$ theory does not apply. The estimates must retain the incoming boundary condition and account for the behavior near grazing velocities.

For in-flow, Deng--Duan \cite{DengDuan2021} proved exponential decay for linearized soft-potential Landau in bounded domains. A weight depending on both position and velocity converts transport into additional damping. Bernou--Carrapatoso--Mischler--Tristani \cite{Bernou.Carrapatoso.Mischler.Tristani2023} established hypocoercivity for linearized Landau and other collision models with Maxwell reflection. Carrapatoso--Mischler \cite{Carrapatoso.Mischler2024} constructed global weak near-equilibrium solutions for nonlinear Landau with general Maxwell accommodation and every $\gamma\in[-3,1]$. Their weighted $L^\infty$ estimates include Gaussian weights and quantitative decay. This result supplies the diffuse-reflection weak solution used in Section~\ref{cex:section}. The regularity and nonlinear uniqueness conclusions are developed in Sections~\ref{cb:section} and \ref{sec:uniqueness}.

Kim--Seo \cite[Theorems~2.1--2.5]{Kim.Seo2026} recently established stationary states and global nonlinear well-posedness for Coulomb Landau with prescribed in-flow and non-isothermal diffuse-reflection. Their results include uniqueness in a weak solution class and exponential relaxation under suitable convergence assumptions on the boundary data. The constant-wall diffuse problem is included, so the Coulomb well-posedness conclusions overlap with ours.

\para{Kinetic traces} The weak construction and the nonlinear difference estimate both require boundary traces and a Green identity. Transport and Vlasov trace theories were developed in \cite{Bardos1970,Ukai1986,Boyer2005,Mischler2000}. Mischler \cite{Mischler2010} established renormalized traces and Maxwell boundary conditions for kinetic equations with velocity diffusion. Variational kinetic Sobolev spaces and weak well-posedness were studied by Albritton--Armstrong--Mourrat--Novack \cite{Albritton.Armstrong.Mourrat.Novack2024}.

For kinetic Fokker--Planck equations, Silvestre \cite[Proposition~4.3]{Silvestre2022} obtains a local trace estimate with weight $\min\{\abs{n\cdot v},\abs{n\cdot v}^2\}$. Zhu \cite[Definition~2.1 and Remarks~2.2--2.3]{Zhu2024} formulates weak solutions together with traces determined by the Green identity. Recent results of Niebel--Valentini \cite{Niebel.Valentini2026} establish sharp trace estimates and counterexamples that depend on the boundary regularity and the velocity model.

Definition~\ref{def:energy} combines this weak solution and trace framework with the Gaussian $L^\infty$ bound \eqref{eq:weak-gaussian-bound}. The Landau equation gives $\Y f\in L^2(0,T;\Vspace')$, and the trace theorem applies in Lemma~\ref{trace:CM}. Renormalization transfers the Gaussian bound to $\trpm f$ in Lemma~\ref{trace:gaussian}. Its velocity integrability then gives both the unweighted boundary $L^2$ estimate \eqref{eq:intro-trace-unweighted} and the full flux estimate \eqref{eq:intro-trace-full}.

\para{Boundary regularity} Once the weak solution and its trace are defined, the next issue is regularity up to the boundary. Interior H\"older estimates with rough coefficients were established in \cite{Golse.Imbert.Mouhot.Vasseur2019,Guerand.Mouhot2022}. Silvestre \cite{Silvestre2022} proved H\"older estimates up to the boundary with prescribed incoming data. Zhu \cite{Zhu2024} developed weak existence, uniqueness and boundary regularity for linear kinetic Fokker--Planck equations with several boundary laws, and subsequently studied sharp a priori in-flow boundary estimates \cite[Theorems~1.3--1.4]{Zhu2025}, including their application to nonlinear Landau in \cite[Section~1.2.7]{Zhu2025}. Optimal specular regularity was studied by Ros-Oton--Weidner \cite{Ros-Oton.Weidner2025}. Nongrazing in-flow boundary expansions and Schauder estimates are developed in \cite[Section~4]{Ros-Oton.Weidner2025} and \cite[Theorem~1.3]{Zhu2025}; Section~\ref{lin:section-order-two} will adapt the Ros-Oton--Weidner argument to the coefficient assumptions and velocity dependence needed here. Kim--Weidner \cite{Kim.Weidner2026} obtained sharp $\Ckin^{\frac{1}{2}}$ regularity at grazing for in-flow and diffuse-reflection, followed by an analysis of mixed Maxwell reflection \cite{Kim.Weidner2026(=)}.

The sharp H\"older exponent and the local first-order argument in Proposition~\ref{cb:inflow-gradient} are consequences of this linear boundary theory. The passage from H\"older regularity of $f$ to regularity of the Landau coefficients and then to velocity derivatives also appears in the periodic theory of \cite{Kim.Guo.Hwang2020}. Here we will verify the hypotheses for the weak Landau solution, normalize its velocity-dependent diffusion, and track the Gaussian bounds through the changes of variables. For diffuse-reflection, we must also estimate the outgoing flux that determines the incoming data. Sections~\ref{lin:section} and \ref{cb:section} carry out these verifications and the corresponding second-order argument.

\para{Uniqueness} The specular-reflection theory closes nonlinear difference estimates using the $L^\infty$ control of $\nabla_vf$ obtained from kinetic regularity. In the periodic and finite-channel settings, Duan--Liu--Sakamoto--Strain \cite{Duan.Liu.Sakamoto.Strain2020} combine Landau trilinear estimates with the algebra property of Fourier $\ell^1$ summability. The finite-channel proof uses its product geometry to take a global tangential Fourier transform. For curved in-flow and diffuse-reflection boundaries, the pointwise velocity gradient can be singular at grazing, and its spatial dependence must enter the difference estimate.

Our $L^2$ difference estimate uses the logarithmic velocity integral of the grazing gradient bound, followed by an exponent optimization related to Yudovich's argument \cite{Yudovich1963}; Section~\ref{sec:intro-method} gives the calculation. Logarithmic comparison estimates also occur in kinetic theory. Loeper \cite{Loeper2006} treated Vlasov--Poisson with bounded spatial density, while Fournier \cite{Fournier2010} proved uniqueness and Wasserstein stability for homogeneous Coulomb Landau using a modulus of the form $s(1+\abs{\ln s})$.

\para{Contributions of this paper} Theorems~\ref{main:inflow} and \ref{main:diffuse} establish global strong solutions for small compatible Gaussian perturbations, and uniqueness and continuous dependence among weak solutions with small Gaussian amplitude, for $-3\leq\gamma\leq1$ in bounded $C^3$ domains. The smallness thresholds are independent of the derivative bounds of the data. For in-flow, the construction combines transport damping with Gaussian estimates and allows persistent, time-dependent incoming forcing. For diffuse-reflection, the weak solution of \cite{Carrapatoso.Mischler2024} is shown to satisfy the boundary and strong regularity asserted here.

The grazing estimates give weighted $L^2$ bounds for $D_v^2f$, $\Y f$, $\Lop f$ and $\Ga(f,f)$, as stated in \eqref{in:eq:main-strong} and \eqref{eq:main-graph}. They also give the logarithmic bound \eqref{eq:intro-log-bound} for $\vnorm{f(t,x,\cdot)}_{\Vv}^2$, which closes the $L^2$ difference estimate and proves uniqueness in the weak solution class with small Gaussian amplitude. For diffuse-reflection, the proof establishes regularity of the outgoing flux and constructs the incoming lift needed for the derivative estimates. The compatibility conditions extend these estimates through $t=0$.

Compared with the closely related work of Kim--Seo \cite{Kim.Seo2026}, we treat the full potential range and focus on sharp grazing estimates and weighted $L^2$ second-order regularity through $t=0$. We assume Gaussian velocity decay, while their data norms use polynomial weights. Our smallness conditions concern only the Gaussian amplitudes, independently of the finite derivative bounds; the data sizes in \cite[(2.15) and (13.85)]{Kim.Seo2026} also include derivative norms. For diffuse-reflection, our $C^3$ domain, $C^5$ initial data and compatibility through the first time derivative require fewer derivatives, while their in-flow theory permits lower domain and data regularity. Our in-flow result also allows persistent time-dependent forcing without an $L^2$-in-time condition relative to a stationary profile. Our proof combines analytic energy and boundary regularity estimates with a logarithmic difference estimate. Kim--Seo combine stochastic representations for linear solvability and traces with analytic kinetic integrability and finite-exponent derivative estimates.

\subsection{Strategy of the Proof}\label{sec:intro-method}

We will prove the three conclusions of the main theorems in order. In Section~\ref{sec:energy}, we will construct global weak solutions with uniform Gaussian bounds. In Sections~\ref{lin:section} and \ref{cb:section}, we will establish boundary regularity and the strong-solution property. In Section~\ref{sec:uniqueness}, we will compare the constructed solution with weak solutions to prove uniqueness and continuous dependence. We will use the collision estimates of Section~\ref{cvel:section} and the trace identities of Section~\ref{sec:trace} throughout.

\para{Weak solutions and global estimates} We will construct solutions satisfying \eqref{eq:intro-trace-space}--\eqref{eq:weak-gaussian-bound}. The collision estimate \eqref{eq:intro-collision-dual} gives $\Y f\in L^2(0,T;\Vspace')$, so the trace and energy identities of Section~\ref{sec:trace} apply to these solutions and their differences.

For in-flow, fix $\psi$ with small Gaussian norm and $\m+\mh\psi\geq0$, and solve
\begin{align}
\Y F_\psi=Q(\m+\mh\psi,F_\psi),\qquad F_\psi(0)=\m+\mh f_0,\qquad \trm F_\psi=\m+\mh g. \label{eq:intro-linear-construction}
\end{align}
Lemma~\ref{ci:frozen-approximation} constructs nonnegative regularized densities. Lemma~\ref{ci:frozen-inflow-construction} passes to the limit and gives $F_\psi\geq0$, together with finite-time $L^2$, Gaussian and trace estimates for $f_\psi:=\mih(F_\psi-\m)$. Lemma~\ref{trace:interior-holder} gives equicontinuity of $f_\psi$ on compact subsets, hence local uniform compactness. The Gaussian bound controls the time intervals near $0,T$, the spatial strip near $\partial\Om$ and the velocity tails by \eqref{ci:fixed-point-tails}, yielding compactness in $L^2((0,T)\times\phase)$. Separately, the energy estimate \eqref{ci:map-continuity-energy} for $f_{\psi_j}-f_\psi$ proves continuity of $\psi\mapsto f_\psi$. For sufficiently small $T$, the finite-time estimates preserve a closed convex set with a prescribed Gaussian bound, and Schauder's fixed-point theorem gives local weak existence in Proposition~\ref{ci:local-existence}.

To continue this local weak solution globally, we need $L^2$ and Gaussian bounds uniform in the length of the time interval. Testing with $\left(1-\kappa\frac{x\cdot v}{\la}\right)f$, for fixed small $\kappa>0$, gives the positive transport term $\frac{\kappa}{2}\int_{\phase}\frac{\abs{v}^2}{\la}f^2$. Together with collision coercivity and the finite-dimensional kernel estimate in Lemma~\ref{ci:velocity-transport-control}, this yields the forced decay bound \eqref{in:eq:main-energy-forcing}.
The spatially dependent Gaussian
\begin{align}
\gbarrier(x,v):=\exp\left(-\Big(\alpha+\frac{1}{4}\Big)\abs{v}^2+\frac{x\cdot v}{\la}\right), \label{eq:intro-tilted-Gaussian}
\end{align}
used in the proof of Lemma~\ref{ci:tilted-barrier} satisfies $\bigl(\Y-Q(F,\cdot)\bigr)\gbarrier\geq\gbarrier$ outside a fixed velocity ball. On bounded velocity sets, Lemma~\ref{ci:local-inflow-bound} applies the local in-flow $L^2$ to $L^\infty$ estimate to $(\pm f-k)_+$, where $k$ bounds the initial and incoming values. Its source contribution has a small factor $r^2$, with $r$ the local kinetic scale. For large velocities, the Gaussian bound for $\mh\Kop f$ from Section~\ref{cvel:section} makes the source small relative to $\gbarrier$. Proposition~\ref{ci:gaussian-bootstrap} combines these two estimates and absorbs their small multiples of $\tnorm{f}_{L^\infty(0,T;L^{\infty,\alpha})}$. The resulting Gaussian bound and \eqref{in:eq:main-energy-forcing} are uniform in $T$. Proposition~\ref{ci:global-energy} then repeats the local construction on successive intervals. Persistent small in-flow preserves the Gaussian bound, and decaying in-flow gives relaxation.

For diffuse-reflection, Proposition~\ref{cex:CM-input} starts from the global weak existence theorem of Carrapatoso--Mischler \cite{Carrapatoso.Mischler2024}. It verifies Definition~\ref{def:energy}, identifies the trace and proves nonnegativity using \eqref{ci:negative-part-energy} and the diffuse boundary condition. The global $\Vspace$ estimate follows from \cite[Proposition~6.2]{Carrapatoso.Mischler2024}, and Corollary~\ref{cex:decay-traces} supplies the decay and integrated trace estimates. Thus the diffuse-reflection argument starts with a global construction, while the in-flow argument combines local existence with the uniform estimates needed for continuation.

\para{Boundary regularity and strong solutions} We next prove that the weak solutions constructed above satisfy Definition~\ref{def:strong}. This requires estimates for $\nabla_vf$, $D_v^2f$ and $\Y f$ up to the boundary and through $t=0$.

Lemmas~\ref{cb:normalization}--\ref{cb:normalized-bounds} give the boundary flattening map, transformed weak equation and uniform coefficient bounds after the change of variables \eqref{lin:inverse-map}. A finite boundary covering and interior estimates reduce the argument to the local equations treated in Section~\ref{lin:section}. Polynomial losses in velocity are absorbed by decreasing the Gaussian exponent.

For in-flow, Lemma~\ref{cb:inflow-incoming-regularity} will verify the incoming H\"older bound from the data assumptions. Lemma~\ref{cb:inflow-coefficients} will first obtain a positive H\"older exponent for $f$ and then estimate the Landau coefficients by convolution. Proposition~\ref{cb:inflow-gradient} will use these coefficient bounds to obtain sharp H\"older regularity and the first-order estimate. For second derivatives, Lemma~\ref{cb:inflow-lift} will supply a lift $\glift$, and Proposition~\ref{cb:strong-upgrade} will estimate $u:=f-\glift$ with zero incoming trace.

For diffuse-reflection, the incoming value depends on the outgoing flux $\flux{f}$ defined in \eqref{eq:diffuse-flux}. Lemma~\ref{cb:rough-flux-coefficients} first uses rough outgoing estimates to obtain a positive H\"older exponent for this flux. The resulting regularity of the incoming data yields H\"older continuity of the solution and its coefficients up to the boundary. Proposition~\ref{cb:flux} then applies the second-order estimate of Proposition~\ref{lin:order-two} on outgoing cylinders with empty incoming boundary and improves the flux regularity to $C^\kappa$ with $\kappa>\frac{2}{3}$. This supplies the incoming-data bound required for the same sharp H\"older and first-order estimates. It also permits the characteristic lift of Lemma~\ref{lift:characteristic}, whose second velocity derivatives and transport derivative have the regularity required for Proposition~\ref{cb:graph}.

The compatibility conditions extend these arguments through $t=0$: Lemmas~\ref{cb:inflow-initial} and \ref{cb:initial} construct negative-time extensions whose residuals satisfy the local hypotheses. This gives the first- and second-order bounds \eqref{in:eq:main-gradient}--\eqref{in:eq:main-hessian} and \eqref{eq:main-gradient}--\eqref{eq:main-hessian}, in the weighted $L^\infty$ sense. To obtain the strong-solution property, we will integrate the squared second-order bound. In normal position and velocity coordinates,
\begin{align}
\int_0^{d_0}\int_{-1}^1\frac{\dd s\dd r}{\bigl(r^{\frac{1}{3}}+\abs{s}\bigr)^3}\lesssim \int_0^{d_0}r^{-\frac{2}{3}}\dd r<\infty. \label{eq:intro-strategy-integrability}
\end{align}
The Gaussian factor controls the remaining velocity integrals, giving weighted $L^2$ bounds for the derivatives. Lemma~\ref{cvel:weighted-collisions} then estimates $\Lop f$ and $\Ga(f,f)$ separately. Proposition~\ref{cb:integration} thus gives \eqref{in:eq:main-strong} and \eqref{eq:main-graph}, and the weak equation becomes an $L^2$ equality as required by Definition~\ref{def:strong}.

\para{Continuous dependence and uniqueness} We will first estimate the difference between the constructed solution $f$ and a weak solution $\widetilde f$ in the small Gaussian class. The first-order bounds \eqref{in:eq:main-gradient} and \eqref{eq:main-gradient} supply the regularity of $f$ used in this argument. For $h:=\widetilde f-f$, the difference equation is
\begin{align}
\Y h+\Bigl(\Lop-\Ga(\widetilde f,\cdot)\Bigr)h=\Ga(h,f). \label{eq:intro-difference}
\end{align}
The smallness condition and \eqref{cvel:principal} bound $\ip[\big]{\Lop h-\Ga(\widetilde f,h)}{h}$ below by $c\norm{h}_{\Vspace}^2-C\norm{h}_{L^2}^2$, with $c>0$. The derivative in the second argument of the remaining collision term falls on $f$.

For incoming data $g$ and $\widetilde g$, the difference satisfies $\trm h=\widetilde g-g$. The collision estimates and Proposition~\ref{trace:green} give the following inequality, where the velocity norms in the integral are evaluated at $(t,x)$:
\begin{align}
\frac{\dd}{\dd t}\norm{h}_{L^2}^2+c\norm{h}_{\Vspace}^2\leq C\norm{h}_{L^2}^2+C\int_\Om\vnorm{f}_{\Vv}^2\vnorm{h}_{L_v^2}^2\dd x+\norm[\big]{\widetilde g-g}_{L^2_{\Sigma_-}}^2. \label{eq:intro-strategy-difference-energy}
\end{align}
For diffuse-reflection, \eqref{trace:projection-identity} gives nonnegative boundary dissipation, so the last term is absent. Fix $T>0$. The first-order bound supplied by Proposition~\ref{cb:integration} and the Gaussian amplitude yield
\begin{align}
\vnorm{f(t,x,\cdot)}_{\Vv}^2\lesssim_T 1+\ln\left(\frac{2\operatorname{diam}(\Om)}{d(x)}\right),\qquad (t,x)\in[0,T]\times\Om\ \text{a.e.} \label{eq:intro-log-bound}
\end{align}
Lemma~\ref{uniq:velocity-log} obtains this estimate by integrating $\graz^{-1}$ in the normal velocity.

The logarithm has $L_x^p$ norm $\lesssim p$. H\"older's inequality and Lemma~\ref{uniq:logarithm} therefore give
\begin{align}
\int_\Om\left(1+\ln\left(\frac{2\operatorname{diam}(\Om)}{d(x)}\right)\right)\vnorm{h(t,x,\cdot)}_{L_v^2}^2\dd x\lesssim p\norm{h(t)}_{L_x^\infty L_v^2}^{\frac{2}{p}}\norm{h(t)}_{L^2}^{2-\frac{2}{p}},\qquad p\geq2. \label{eq:intro-strategy-log-product}
\end{align}
The Gaussian bounds control the mixed norm on the right. To recover a full factor $\norm{h(t)}_{L^2}^2$, choose $p=\ln\left(\frac{C_T}{\norm{h(t)}_{L^2}^2}\right)$ whenever $\norm{h(t)}_{L^2}>0$, with $C_T$ fixed and sufficiently large. Then $p\geq2$ and $\left(\frac{C_T}{\norm{h(t)}_{L^2}^2}\right)^{\frac1p}=\ue$. Substitution into \eqref{eq:intro-strategy-difference-energy} gives
\begin{align}
\frac{\dd}{\dd t}\norm{h(t)}_{L^2}^2\leq C_T\norm{h(t)}_{L^2}^2\ln\left(\frac{C_T}{\norm{h(t)}_{L^2}^2}\right)+\norm[\big]{\widetilde g(t)-g(t)}_{L^2_{\Sigma_-}}^2. \label{eq:intro-Osgood}
\end{align}
The logarithmic term is defined to be zero when $h(t)=0$. After adding the remaining incoming discrepancy, a logarithm and an integrating factor give a quantitative bound for $\tnorm[\big]{\widetilde f-f}_{C([0,T];L^2)}$ in terms of the initial and incoming differences; see Propositions~\ref{uniq:osgood} and \ref{uniq:osgood-diffuse}. Equal data imply $\widetilde f=f$. For diffuse-reflection, the same argument uses only the initial difference.

The first-order estimate is required of the reference solution $f$; the $L^2$ conditions in Definition~\ref{def:strong} alone do not supply \eqref{eq:intro-log-bound}. The comparison solution enters through its weak equation, its trace identity and its small Gaussian bound. Retaining the dissipation in \eqref{eq:intro-strategy-difference-energy} gives convergence in $L^2(0,T;\Vspace)$. Lemmas~\ref{uniq:bulk} and \ref{uniq:bulk-diffuse} also give convergence of the transport derivatives in $L^2(0,T;\Vspace')$ and of both full traces. In the diffuse case, a test with the normal velocity controls the flux. Section~\ref{sec:assembly} combines these conclusions with existence and regularity to prove both main theorems.

\subsection*{Acknowledgements}
Lei Wu was supported in part by the National Science Foundation under grant DMS-2405161.


\section{The Landau Collision Operator}\label{cvel:section}

\subsection{Collision Identities and \texorpdfstring{$\sigma$}{} Properties}\label{sec:cvel-kernels}

We use $H^{-1}(B_R):=(H_0^1(B_R))'$. Thus $T\in H^{-1}_{\mathrm{loc}}(\R^3)$ means that its action on $C_c^\infty(B_R)$ extends continuously to $H_0^1(B_R)$ for every $R>0$.

\begin{lemma}[Collision identities]\label{cvel:distribution}
The kernels in \eqref{eq:kernels}--\eqref{eq:kernel-c} satisfy $\partial_j a_{ij}=b_i$ and $\partial_i b_i=c$ in $\mathcal D'(\R^3)$. Suppose that $\ue^{\kappa\abs{v}^2}G\in L^\infty_v$ for some $\kappa>0$.
\begin{enumerate}
\item For $H\in H^1_{\mathrm{loc}}(\R^3)$, the divergence formula \eqref{eq:collision} defines $Q(G,H)\in H^{-1}_{\mathrm{loc}}(\R^3)$.
\item If $H\in H^2_{\mathrm{loc}}(\R^3)$, then $Q(G,H)\in L^2_{\mathrm{loc}}(\R^3)$ and
\begingroup\usetagform{enumerated}
\begin{align}
Q(G,H)=(a*G):D_v^2H-(c*G)H. \label{cvel:Q}
\end{align}
\endgroup
At $\gamma=-3$, the last term in \eqref{cvel:Q} is $8\pi GH$.
\end{enumerate}
\end{lemma}

\begin{proof}
The identities $\partial_j a_{ij}=b_i$ and $\partial_i b_i=c$, and the collision formulas for smooth functions, are given in \cite[(1.4)--(1.7)]{Carrapatoso.Mischler2024} for the kernels \eqref{eq:kernels}--\eqref{eq:kernel-c}. To obtain (i)--(ii) under the stated local Sobolev assumptions on $H$, we first prove that
\begin{align}
a*G\in W^{1,\infty}_{\mathrm{loc}}(\R^3),\qquad b*G,\ c*G\in L^\infty_{\mathrm{loc}}(\R^3). \label{cvel:local-coefficients}
\end{align}
For $v\neq0$, direct differentiation of \eqref{eq:kernels} gives
\begin{align}
\partial_k a_{ij}(v)&=(\gamma+2)\abs{v}^{\gamma}v_k\delta_{ij}-\gamma\abs{v}^{\gamma-2}v_kv_iv_j-\abs{v}^{\gamma}\bigl(\delta_{ik}v_j+v_i\delta_{jk}\bigr), \label{cvel:kernel-derivative}\\
\abs{a(v)}&\lesssim \abs{v}^{\gamma+2},\qquad \abs{\nabla_va(v)}+\abs{b(v)}\lesssim \abs{v}^{\gamma+1}. \label{cvel:kernel-growth}
\end{align}
Since $\gamma+2\geq-1$ and $\gamma+1\geq-2$, the kernels in \eqref{cvel:kernel-growth} are locally integrable. To identify their derivatives across $v=0$, integrate by parts outside $B_\varepsilon$. For $\varphi\in C_c^\infty(\R^3)$, the boundary term satisfies
\begin{align}
\abs{\int_{\abs{v}=\varepsilon}a_{ij}(v)\varphi(v)\frac{v_k}{\varepsilon}\dd S_v}\lesssim \vnorm{\varphi}_{L_v^\infty}\varepsilon^{\gamma+4}\to0. \label{cvel:kernel-derivative-boundary}
\end{align}
Thus \eqref{cvel:kernel-derivative} is the weak derivative of $a$, and $\partial_j a_{ij}=b_i$ in distributions. For $\divv_v b$, the corresponding term is
\begin{align}
\int_{\abs{v}=\varepsilon}b(v)\cdot\frac{v}{\abs{v}}\,\varphi(v)\dd S_v=-2\varepsilon^{\gamma+3}\int_{\mathbb S^2}\varphi(\varepsilon\omega)\dd S_\omega. \label{cvel:kernel-flux}
\end{align}
It tends to zero for $\gamma>-3$ and to $-8\pi\varphi(0)$ for $\gamma=-3$, proving $\divv_v b=c$ with the two cases in \eqref{eq:kernel-c}.

We will use the Gaussian convolution bound
\begin{align}
\int_{\R^3}\abs{v-w}^s\ue^{-\kappa\abs{w}^2}\dd w\lesssim_{s,\kappa}\la^s\qquad(s>-3,\ \kappa>0). \label{cvel:Gaussian-convolution-bound}
\end{align}
For $\abs{v}\leq2$, split at $\abs{v-w}=1$ and use local integrability and Gaussian decay. For $\abs{v}>2$ and $s\geq0$, use $\abs{v-w}^s\lesssim_s\abs{v}^s+\abs{w}^s$ and integrate against the Gaussian. For $-3<s<0$, split the integral into $\abs{v-w}\geq\frac{\abs{v}}{2}$ and its complement. The first part is $\lesssim_{s,\kappa}\abs{v}^s$. On the complement, $\abs{w}>\frac{\abs{v}}{2}$, so the integral is $\lesssim_s\abs{v}^{s+3}\ue^{-\frac{\kappa}{4}\abs{v}^2}\lesssim_{s,\kappa}\abs{v}^s$.

Taking $s=\gamma+2$ and $s=\gamma+1$ in \eqref{cvel:Gaussian-convolution-bound} proves the bounds for $a*G$, $(\nabla_va)*G$ and $b*G$. The same estimate justifies differentiation under convolution in distributions, so $\partial_k(a*G)=(\partial_ka)*G$. For $\gamma>-3$, use $c(v)=-2(\gamma+3)\abs{v}^{\gamma}$ and take $s=\gamma$. For $\gamma=-3$, use $c*G=-8\pi G$. This proves \eqref{cvel:local-coefficients}. The kernel identities also give
\begin{align}
\partial_i(a*G)^{ij}=(b*G)^j,\qquad \divv_v(b*G)=c*G\qquad\text{in }\mathcal D'(\R^3). \label{cvel:convolved-identities}
\end{align}
Suppose first that $H\in H^1_{\mathrm{loc}}(\R^3)$. The divergence expression \eqref{eq:collision} defines a distribution by
\begin{align}
\ip*{Q(G,H)}{\varphi}=-\int_{B_R}\bigl((a*G)\nabla_vH-(b*G)H\bigr)\cdot\nabla_v\varphi\dd v\qquad\bigl(\varphi\in C_c^\infty(B_R)\bigr). \label{cvel:local-collision-pairing}
\end{align}
By \eqref{cvel:local-coefficients} and the Cauchy--Schwarz inequality, the pairing extends to $H_0^1(B_R)$ with
\begin{align}
\vnorm{Q(G,H)}_{H^{-1}(B_R)}\leq\bigl(\vnorm{a*G}_{L^\infty(B_R)}+\vnorm{b*G}_{L^\infty(B_R)}\bigr)\vnorm{H}_{H^1(B_R)}. \label{cvel:local-negative-bound}
\end{align}
This proves (i).

For (ii), fix $R>0$ and choose smooth $H_j\to H$ in $H^2(B_R)$. The product rule and \eqref{cvel:convolved-identities} give
\begin{align}
\divv_v\bigl((a*G)\nabla_vH_j-(b*G)H_j\bigr)=(a*G):D_v^2H_j-(c*G)H_j. \label{cvel:local-product-expansion}
\end{align}
The left side converges to $Q(G,H)$ in $H^{-1}(B_R)$ by \eqref{cvel:local-negative-bound}; for the right side, \eqref{cvel:local-coefficients} gives
\begin{align}
&\vnorm{(a*G):D_v^2(H_j-H)-(c*G)(H_j-H)}_{L^2(B_R)} \label{cvel:local-second-order-limit}\\
&\qquad\leq\bigl(\vnorm{a*G}_{L^\infty(B_R)}+\vnorm{c*G}_{L^\infty(B_R)}\bigr)\vnorm{H_j-H}_{H^2(B_R)}\to0. \notag
\end{align}
Thus \eqref{cvel:Q} holds in $L^2(B_R)$ for every $R>0$. At $\gamma=-3$, $c*G=-8\pi G$ gives the stated formula.
\end{proof}

\begin{remark}\label{cvel:local-versus-global}
Under the assumptions of Lemma~\ref{cvel:distribution}, $a*G$ and $c*G$ are bounded for $-3\leq\gamma\leq-2$, so $H\in H^2(\R^3)$ implies $Q(G,H)\in L^2(\R^3)$. For $\gamma>-2$, $a*G$ can grow and this conclusion does not follow from $H\in H^2(\R^3)$ alone. Lemma~\ref{cvel:weighted-collisions} supplies the global bounds under the weighted assumptions used for strong solutions.
\end{remark}

The velocity norm in \eqref{eq:V} uses the following classical properties of $\sig$.

\begin{lemma}[Properties of $\sig$]\label{cvel:convolution}
The matrix $\sig$ in \eqref{eq:sigma-definition} is smooth and uniformly positive definite on each bounded velocity set. For $v\neq0$, write $\xi=\xi_\parallel+\xi_\perp$, where the two components are parallel and orthogonal to $v$. Then
\begin{align}
\xi\cdot\sig\xi\simeq\la^\gamma\abs{\xi_\parallel}^2+\la^{\gamma+2}\abs{\xi_\perp}^2. \label{cvel:eigenvalues}
\end{align}
\end{lemma}

\begin{proof}
Smoothness, positivity and \eqref{cvel:eigenvalues} follow from \cite[Lemma~2.1]{Carrapatoso.Mischler2024}, which treats $a*\m$ for $-3\leq\gamma\leq1$. This is exactly $\sig$ in \eqref{eq:sigma-definition}, with the same Maxwellian and interaction range. See also \cite[Lemma~3]{Guo2002(=)}.
\end{proof}

\begin{remark}\label{cvel:sigma-norm-comparison}
After the fixed change $v=\sqrt2\,w$ and Maxwellian normalization, the unweighted $\sig$-norm in \cite[(9)]{Guo2002(=)} becomes, up to a fixed positive factor,
\begin{align}
\vnorm{h}_{\sig,0}^2:=\int_{\R^3}\left[\sig^{ij}\partial_i h\partial_jh+\frac{1}{4}(v\cdot\sig v)\abs{h}^2\right]\dd v. \label{cvel:sigma-norm}
\end{align}
The factor $\frac14$ comes from $\nabla_v\mh=-\frac v2\mh$; the change of variables only modifies constants in the cited estimates. To compare \eqref{cvel:sigma-norm} with \eqref{eq:V}, Lemma~\ref{cvel:convolution} gives
\begin{align}
v\cdot\sig v\simeq\abs{v}^2\la^\gamma,\qquad \la^{\gamma+2}\simeq v\cdot\sig v\quad\text{on }\{\abs{v}\geq1\}. \label{cvel:zero-order-comparison}
\end{align}
For $h\in C_c^\infty(\R^3)$, the Poincar\'e inequality on $B_2$, with an $L^2(B_2\setminus B_1)$ term, gives
\begin{align}
\vnorm{h}_{L^2(B_1)}^2\lesssim \Bigl(\vnorm{\nabla_vh}_{L^2(B_2)}^2+\vnorm{h}_{L^2(B_2\setminus B_1)}^2\Bigr)\lesssim \vnorm{h}_{\sig,0}^2. \label{cvel:local-zero-order-control}
\end{align}
The second inequality uses local ellipticity and \eqref{cvel:zero-order-comparison} on $B_2\setminus B_1$. Together with the comparison outside $B_1$, it yields
\begin{align}
\vnorm{h}_{\Vv}\simeq\vnorm{h}_{\sig,0},\qquad -3\leq\gamma\leq1. \label{cvel:classical-norm}
\end{align}
This local argument also appears in \cite[Corollary~1]{Guo2002(=)}. By completion, \eqref{cvel:classical-norm} transfers Guo's estimates with $\theta=0$ to $\Vv$. For $\theta\neq0$, the additional velocity weights in \cite[(9)]{Guo2002(=)} must be retained.
\end{remark}

\subsection{Properties of \texorpdfstring{$\Lop$}{} and \texorpdfstring{$\Ga$}{}}\label{sec:cvel-linear}

We will first record the approximation used below in $\Vv$ and $\Wv$. Let $\chi\in C_c^\infty(B_2)$ satisfy $0\leq\chi\leq1$ and $\chi=1$ on $B_1$, and set $\chi_R(v):=\chi\left(\frac{v}{R}\right)$ for $R\geq1$. By \eqref{eq:V} and \eqref{cvel:eigenvalues}, every $h\in\Vv$ satisfies
\begin{align}
\vnorm[\big]{(1-\chi_R)h}_{\Vv}^2\lesssim \int_{\abs{v}>R}\left[\sig^{ij}\partial_ih\partial_jh+\la^{\gamma+2}\abs{h}^2\right]\dd v\to0. \label{cvel:W-density-cutoff}
\end{align}
If $h\in\Wv$, we also have $\vnorm{(1-\chi_R)h}_{L_v^2}\to0$. On each fixed ball the $\Vv$ norm is equivalent to the $H^1$ norm for functions supported there. Mollifying $\chi_Rh$ and then taking $R\to\infty$ therefore gives a sequence with
\begin{align}
h_n\in C_c^\infty(\R^3),\qquad h_n\to h\text{ in }\Vv,\qquad h_n\to h\text{ in }\Wv\text{ if }h\in\Wv. \label{cvel:W-density}
\end{align}
The bounds in the following lemmas define the collision pairings for nonsmooth arguments by approximation in the stated norms. We use this convention throughout the paper. Pairings are between $\Vv'$ and $\Vv$ unless an $L_v^2$ inner product is indicated.

For $h\in C_c^\infty(\R^3)$, write
\begin{align}
\Lop h=-\mih Q(\m,\mh h)-\Kop h,\qquad \Kop h:=\mih Q(\mh h,\m)=\Ga(h,\mh). \label{eq:I:display:052}
\end{align}
Both expressions belong to $L_v^2$. Their action on test functions is given by integration.

\begin{lemma}[Linear estimates]\label{cvel:linear-estimates}
For $h,k\in\Vv$, we have $\Lop h,\Kop h\in\Vv'$ and
\begin{align}
\abs[\big]{\ip*{\Lop h}{k}}+\abs[\big]{\ip*{\Kop h}{k}}\lesssim \vnorm{h}_{\Vv}\vnorm{k}_{\Vv}. \label{cvel:linear-boundedness}
\end{align}
There are $c_0,C_0,R_0>0$ such that
\begin{align}
\ip*{\Lop h}{h}\geq c_0\vnorm{h}_{\Vv}^2-C_0\vnorm{h}_{L^2(B_{R_0})}^2\qquad(h\in\Vv). \label{cvel:localized-garding}
\end{align}
\end{lemma}

\begin{proof}
For $h,k\in C_c^\infty(\R^3)$, \eqref{cvel:linear-boundedness} and \eqref{cvel:localized-garding} follow from \cite[Lemma~1 and Lemma~5, (23)]{Guo2002(=)} with velocity-weight parameter $\theta=0$. Remark~\ref{cvel:sigma-norm-comparison} matches the Maxwellian and norm. Approximation by \eqref{cvel:W-density} extends the pairings and both estimates to $\Vv$; the local $L^2(B_{R_0})$ convergence follows from local ellipticity.

The expression for $\Kop$ follows from \eqref{cvel:Q} and $D_v^2\m=(v\otimes v-I)\m$:
\begin{align}
\Kop h=\mh\Bigl[(v\otimes v-I):\bigl(a*(\mh h)\bigr)-c*(\mh h)\Bigr]. \label{cvel:K}
\end{align}
For $h\in C_c^\infty(\R^3)$, the convolutions have at most polynomial growth, so $\Kop h\in L_v^2$. The first term in \eqref{eq:I:display:052} has compact support and belongs to $L_v^2$ by Lemma~\ref{cvel:distribution}.
\end{proof}

\begin{lemma}[Null space and coercivity]\label{cvel:entropy-coercivity}
For $h,k\in\Vv$,
\begin{align}
\ip*{\Lop h}{k}=\ip*{\Lop k}{h},\qquad \ip*{\Lop h}{h}\geq0. \label{cvel:form-symmetry}
\end{align}
The null space of $\Lop:\Vv\to\Vv'$ is
\begin{align}
\mathcal N:=\operatorname{span}\{\mh,v_1\mh,v_2\mh,v_3\mh,\abs{v}^2\mh\}. \label{cvel:nullspace}
\end{align}
Choose an $L_v^2$ orthonormal basis $e_1,\ldots,e_5$ of $\mathcal N$ and define
\begin{align}
\pk h:=\sum_{j=1}^5\left(\int_{\R^3}h(v)e_j(v)\dd v\right)e_j\qquad(h\in\Vv). \label{cvel:projection-definition}
\end{align}
This agrees with the $L_v^2$ orthogonal projection on $\Wv$ and satisfies $\vnorm{\pk h}_{\Vv}\lesssim\vnorm{h}_{\Vv}$. There is $\lambda_\gamma>0$ such that
\begin{align}
\ip*{\Lop h}{h}\geq\lambda_\gamma\vnorm[\big]{(I-\pk)h}_{\Vv}^2\qquad(h\in\Vv). \label{cvel:micro}
\end{align}
In particular, $\ip*{\Lop h}{h}=0$ if and only if $h\in\mathcal N$.
\end{lemma}

\begin{proof}
For smooth functions, \cite[Lemma~4 and Lemma~5, (24)]{Guo2002(=)} give symmetry, nonnegativity, $\Lop\varphi=0$ for $\varphi\in\mathcal N$, and \eqref{cvel:micro}. We take $\theta=0$ and use Remark~\ref{cvel:sigma-norm-comparison}; the cited projection agrees with \eqref{cvel:projection-definition}. See also \cite[Theorem~1.2]{Mouhot.Strain2007} for coercivity.

For $h\in\Vv$ and $\varphi\in\mathcal N$, the Cauchy--Schwarz inequality gives
\begin{align}
\abs{\int h\varphi\dd v}\leq\vnorm[\Big]{\la^{\frac{\gamma+2}{2}}h}_{L_v^2}\vnorm[\Big]{\la^{-\frac{\gamma+2}{2}}\varphi}_{L_v^2}\lesssim_\varphi\vnorm{h}_{\Vv}. \label{cvel:projection-continuity}
\end{align}
Applying this to $e_1,\ldots,e_5$ proves continuity of $\pk$ on $\Vv$. Hence \eqref{cvel:W-density} and \eqref{cvel:linear-boundedness} extend \eqref{cvel:form-symmetry} and \eqref{cvel:micro} to $\Vv$. The latter gives $\ip*{\Lop h}{h}=0\Rightarrow h=\pk h\in\mathcal N$; the converse follows from $\Lop\varphi=0$ on $\mathcal N$.
\end{proof}

\begin{lemma}[Trilinear estimate]\label{cvel:trilinear-form}
For $g\in L_v^2$ and $h,k\in\Vv$, we have
\begin{align}
\abs[\Big]{\ip[\big]{\Ga(g,h)}{k}}\lesssim \vnorm{g}_{L_v^2}\vnorm{h}_{\Vv}\vnorm{k}_{\Vv}. \label{cvel:trilinear}
\end{align}
In particular, $\Ga(g,h)\in\Vv'$ and
\begin{align}
\vnorm[\big]{\Ga(g,h)}_{\Vv'}\lesssim \vnorm{g}_{L_v^2}\vnorm{h}_{\Vv}. \label{cvel:velocity-dual}
\end{align}
\end{lemma}

\begin{proof}
The bound \eqref{cvel:trilinear} follows from \cite[Lemma~3.4]{Duan.Liu.Sakamoto.Strain2020}. We apply its Landau estimate with velocity weight $1$ to $g,h,k\in C_c^\infty(\R^3)$, regarded as functions independent of $x$. The first argument then has the $L_v^2$ norm, and the norms of the other two arguments are equivalent to $\Vv$ by \eqref{cvel:classical-norm}. This gives \eqref{cvel:trilinear} for smooth arguments. Density in $L_v^2$ and \eqref{cvel:W-density} give the estimate for $g\in L_v^2$ and $h,k\in\Vv$. Taking the supremum over $\vnorm{k}_{\Vv}\leq1$ gives \eqref{cvel:velocity-dual}.

For $h\in L_v^2$, take $h_n\in C_c^\infty(\R^3)$ with $h_n\to h$ in $L_v^2$. For every $N\geq0$, the Cauchy--Schwarz inequality gives
\begin{align}
\vnorm[\big]{\la^N\mh(h_n-h)}_{L_v^1}\leq\vnorm[\Big]{\la^N\mh}_{L_v^2}\vnorm{h_n-h}_{L_v^2}\to0. \label{cvel:K-convolution-convergence}
\end{align}
The kernels $a$ and $c$ are locally integrable and have at most polynomial growth when $\gamma>-3$, while $c=-8\pi\delta_0$ when $\gamma=-3$. Thus the right side of \eqref{cvel:K} converges in $\mathcal D'(\R^3)$, and \eqref{cvel:velocity-dual} identifies its limit as
\begin{align}
\Kop h=\Ga(h,\mh)\qquad(h\in L_v^2). \label{cvel:K-bilinear-identity}
\end{align}
For $h\in\Wv$, equations~\eqref{cvel:W-density}, \eqref{cvel:linear-boundedness} and \eqref{cvel:velocity-dual} also give
\begin{align}
\Lop h=-\Ga(\mh,h)-\Ga(h,\mh)\qquad(h\in\Wv). \label{cvel:linear-bilinear-identity}
\end{align}
\end{proof}

\begin{corollary}[Spatial collision estimates]\label{cvel:spatial}
The following holds.
\begin{enumerate}[beginpenalty=10000]
\item For $h\in\Vspace=L^2(\Om;\Vv)$, we have $\Lop h\in\Vspace'$ and
\begingroup\usetagform{enumerated}
\begin{align}
\norm{\Lop h}_{\Vspace'}\lesssim \norm{h}_{\Vspace}. \label{cvel:spatial-L}
\end{align}
\endgroup
\item For $g\in L_x^\infty L_v^2$ and $h\in\Vspace$, we have $\Ga(g,h)\in\Vspace'$ and
\begingroup\usetagform{enumerated}
\begin{align}
\norm[\big]{\Ga(g,h)}_{\Vspace'}\lesssim \norm{g}_{L_x^\infty L_v^2}\norm{h}_{\Vspace}. \label{cvel:spatial-dual}
\end{align}
\endgroup
\item Let $h\in L^2$ and $f\in\Vspace$ satisfy
\begingroup\usetagform{enumerated}
\begin{align}
\int_\Om\vnorm[\big]{h(x,\cdot)}_{L_v^2}^2\vnorm[\big]{f(x,\cdot)}_{\Vv}^2\dd x<\infty. \label{cvel:spatial-reversed-assumption}
\end{align}
\endgroup
Then $\Ga(h,f)\in\Vspace'$ and
\begingroup\usetagform{enumerated}
\begin{align}
\norm[\big]{\Ga(h,f)}_{\Vspace'}^2\lesssim \int_\Om\vnorm[\big]{h(x,\cdot)}_{L_v^2}^2\vnorm[\big]{f(x,\cdot)}_{\Vv}^2\dd x. \label{cvel:spatial-reversed}
\end{align}
\endgroup
\end{enumerate}
\end{corollary}

\begin{proof}
Define the spatial pairings by integrating the velocity pairings over $\Om$. Equation~\eqref{cvel:linear-boundedness} and the Cauchy--Schwarz inequality in $x$ give \eqref{cvel:spatial-L}. Similarly, \eqref{cvel:trilinear} gives
\begin{align}
\abs[\Big]{\ip[\big]{\Ga(g,h)}{k}}\lesssim\left(\int_\Om\vnorm{g(x,\cdot)}_{L_v^2}^2\vnorm{h(x,\cdot)}_{\Vv}^2\dd x\right)^{\frac{1}{2}}\norm{k}_{\Vspace}. \label{cvel:spatial-proof}
\end{align}
The right side is finite under (ii), and under (iii) with $(g,h)$ replaced by $(h,f)$. Taking the supremum over $\norm{k}_{\Vspace}\leq1$ proves \eqref{cvel:spatial-dual} and \eqref{cvel:spatial-reversed}.
\end{proof}

We record the collision invariants in the normalization $\Ga$.

\begin{lemma}[Collision invariants]\label{cvel:conservation}
For $g\in L_v^2$ and $h\in\Vv$,
\begin{align}
\ip[\big]{\Ga(g,h)}{\mh}=0. \label{cvel:mass}
\end{align}
For $g\in\Wv$ and $\varphi\in\{1,v_1,v_2,v_3,\abs{v}^2\}$,
\begin{align}
\ip[\big]{\Ga(g,g)}{\varphi\mh}=0. \label{cvel:invariants}
\end{align}
The pairings are well defined by Lemma~\ref{cvel:trilinear-form}, since $\varphi\mh\in\Vv$.
\end{lemma}

\begin{proof}
For $g,h\in C_c^\infty(\R^3)$, integration of the divergence formula \eqref{eq:collision} gives \eqref{cvel:mass}. Applying \cite[Section~1.2, (1.18)]{Carrapatoso.Mischler2024} to the smooth compactly supported density $\mh g$ gives \eqref{cvel:invariants} for the five stated polynomials. Both identities extend by \eqref{cvel:trilinear} and \eqref{cvel:W-density}, since $\varphi\mh\in\Vv$; the approximation of $g$ in the second identity is in both $L_v^2$ and $\Vv$.
\end{proof}

\subsection{Coefficients of the Nonlinear Equation}\label{sec:cvel-coefficients}

For a velocity function $G$, write $a_G:=a*G$, $b_G:=b*G$ and $c_G:=c*G$. For $F=\m+\mh f$, in particular,
\begin{align}
\Af{}=a*F,\qquad b_F=b*F,\qquad c_F=c*F. \label{cvel:coefficient-convolutions}
\end{align}
The scalar equation for $f$ uses
\begin{align}
\Bf{}:=-b_F-\Af{}v,\qquad \Cf{}:=\frac{1}{4}v\cdot\Af{}v-\frac{1}{2}\tr\Af{}-c_F. \label{cvel:scalar-coefficients}
\end{align}
We will bound these coefficients and the source $\Kop f$, identify the differential expressions for $\Lop$ and $\Ga$, and prove their weighted $L^2$ bounds.

\begin{lemma}\label{cvel:coefficients}
Let $0<\delta\leq1$ and $F=\m+\mh f$, where
\begin{align}
\vnorm{f}_{L_v^\infty}\leq\delta. \label{cvel:coefficient-bounded}
\end{align}
Then the following holds.
\begin{enumerate}
\item There is $\delta_A=\delta_A(\gamma)\in(0,1]$ such that $\delta\leq\delta_A$ implies
\begingroup\usetagform{enumerated}
\begin{align}
\frac{1}{2}\sig\leq\Af{}\leq\frac{3}{2}\sig \label{cvel:ellipticity}
\end{align}
\endgroup
where $\sig$ is defined in \eqref{eq:sigma-definition}. In particular,
\begingroup\usetagform{enumerated}
\begin{align}
0\leq v\cdot\Af{}v\leq\frac{3}{2}v\cdot\sig v\lesssim \la^{\gamma+2}. \label{eq:I:display:194}
\end{align}
\endgroup
\item For every $\delta\leq1$, we have $\Af{}\in W^{1,\infty}_{\mathrm{loc}}(\R^3)$, $\nabla_v\Af{}=(\nabla_va)*F$ in $\mathcal D'(\R^3)$, and
\begingroup\usetagform{enumerated}
\begin{align}
\abs[\big]{\Af{}(v)}\lesssim \la^{\gamma+2},\qquad \abs[\big]{\nabla_v\Af{}(v)}+\abs[\big]{b_F(v)}\lesssim \la^{\gamma+1}. \label{cvel:rough-coefficients}
\end{align}
\endgroup
The remaining coefficient satisfies
\begingroup\usetagform{enumerated}
\begin{align}
\abs[\big]{c_F(v)}\lesssim \la^\gamma\qquad(-3\leq\gamma\leq1). \label{cvel:coefficient-c}
\end{align}
\endgroup
At $\gamma=-3$, we also have
\begingroup\usetagform{enumerated}
\begin{align}
c_F=-8\pi F. \label{cvel:coefficient-Coulomb}
\end{align}
\endgroup
\item For $v=r\omega$, where $r\geq2$ and $\abs{\omega}=1$, we have
\begingroup\usetagform{enumerated}
\begin{align}
\tr\Af{}(r\omega)&=2r^{\gamma+2}\int_{\R^3}F(w)\dd w+O(r^{\gamma+1}), \label{eq:I:display:195}\\
(r\omega)\cdot\Af{}(r\omega)(r\omega)&=r^{\gamma+2}\int_{\R^3}\bigl(\abs{w}^2-(\omega\cdot w)^2\bigr)F(w)\dd w+O(r^{\gamma+1}), \label{ci:radial-moment-expansion}\\
(r\omega)\cdot b_F(r\omega)&=-2r^{\gamma+2}\int_{\R^3}F(w)\dd w+O(r^{\gamma+1}). \label{ci:weighted-drift-expansion}
\end{align}
\endgroup
The error constants are uniform in $\omega$ and in $f$ satisfying \eqref{cvel:coefficient-bounded}. The moments satisfy
\begingroup\usetagform{enumerated}
\begin{align}
\abs{\int_{\R^3}F(w)\dd w-1}+\abs{\int_{\R^3}\bigl(\abs{w}^2-(\omega\cdot w)^2\bigr)F(w)\dd w-2}\lesssim \delta. \label{cvel:coefficient-moments}
\end{align}
\endgroup
\end{enumerate}
\end{lemma}

\begin{proof}
The estimate \eqref{cvel:ellipticity} follows from \cite[Lemma~2.2, (2.9)]{Carrapatoso.Mischler2024}, and the absolute bounds for $\Af{}$, $b_F$ and $c_F$ follow from (2.6)--(2.8) of the same lemma. We apply these estimates with the weight $\la^{10}$. They require boundedness in this weighted $L_v^\infty$ norm and a weight exponent greater than $8+\gamma$. The exponent condition holds because $10>8+\gamma$ throughout $-3\leq\gamma\leq1$. For the perturbation, \eqref{cvel:coefficient-bounded} gives
\begin{align}
\vnorm[\Big]{\la^{10}\mh f}_{L_v^\infty}\leq\delta\vnorm[\Big]{\la^{10}\mh}_{L_v^\infty}\lesssim\delta. \label{cvel:coefficient-weighted-bound}
\end{align}
In particular, (2.9) of that reference gives
\begin{align}
\abs[\Big]{\xi\cdot\bigl(a*(\mh f)\bigr)(v)\xi}\leq C_\gamma\delta\,\xi\cdot\sig(v)\xi\qquad(\xi\in\R^3). \label{cvel:perturbed-covariance}
\end{align}
Since $\Af{}=\sig+a*(\mh f)$, choosing $\delta_A\leq\min\left\{1,\frac{1}{2C_\gamma}\right\}$ proves \eqref{cvel:ellipticity}. Together with the radial eigenvalue bound in \eqref{cvel:eigenvalues}, this gives \eqref{eq:I:display:194}.

For the absolute bounds, $\abs{F}\leq\m+\mh$ gives $\vnorm[\Big]{\la^{10}F}_{L_v^\infty}\lesssim1$. The cited estimates (2.6)--(2.8) yield the bounds for $\Af{}$, $b_F$ and, when $\gamma>-3$, $c_F$. At $\gamma=-3$, \eqref{eq:kernel-c} gives $\abs{c_F}=8\pi\abs{F}\lesssim\la^{-3}$ and \eqref{cvel:coefficient-Coulomb}.

The proof of Lemma~\ref{cvel:distribution} gives $\nabla_v\Af{}=(\nabla_va)*F$. Using \eqref{cvel:kernel-growth} and \eqref{cvel:Gaussian-convolution-bound} with $s=\gamma+1$ gives
\begin{align}
\abs[\big]{\nabla_v\Af{}(v)}\lesssim\int_{\R^3}\abs{v-w}^{\gamma+1}\ue^{-\frac{1}{4}\abs{w}^2}\dd w\lesssim\la^{\gamma+1}. \label{cvel:coefficient-derivative-bound}
\end{align}
This completes \eqref{cvel:rough-coefficients}.

To prove part~(iii), write the convolutions explicitly, with $w$ as the integration variable:
\begin{align}
\tr\Af{}(v)&=2\int_{\R^3}\abs{v-w}^{\gamma+2}F(w)\dd w, \label{cvel:trace-convolution}\\
v\cdot\Af{}(v)v&=\int_{\R^3}\abs{v-w}^{\gamma}\abs{v\times w}^2F(w)\dd w, \label{cvel:radial-convolution}\\
v\cdot b_F(v)&=-2\int_{\R^3}\abs{v-w}^{\gamma}\bigl(v\cdot(v-w)\bigr)F(w)\dd w. \label{cvel:b-moment-convolution}
\end{align}
Let $v=r\omega$. On $\left\{\abs{w}\leq\frac{r}{2}\right\}$, the mean value theorem gives
\begin{align}
\abs{\abs{r\omega-w}^{\gamma+2}-r^{\gamma+2}}&\lesssim r^{\gamma+1}\abs{w}, \label{cvel:trace-expansion-error}\\
\abs{r^2\abs{r\omega-w}^{\gamma}-r^{\gamma+2}}\abs{\omega\times w}^2&\lesssim r^{\gamma+1}\abs{w}^3, \label{cvel:radial-expansion-error}\\
\abs{\abs{r\omega-w}^{\gamma}\bigl(r^2-r\omega\cdot w\bigr)-r^{\gamma+2}}&\lesssim r^{\gamma+1}\abs{w}. \label{cvel:b-expansion-error}
\end{align}
Using $\abs{F(w)}\lesssim\ue^{-\frac{1}{4}\abs{w}^2}$, the integrals of these errors are $\lesssim r^{\gamma+1}$.

On $\left\{\abs{w}>\frac{r}{2}\right\}$, use $\abs{v\times w}=\abs{v\times(w-v)}\leq r\abs{v-w}$ and $\abs{v\cdot(v-w)}\leq r\abs{v-w}$. Thus the absolute values of the three omitted convolution integrals are bounded by
\begin{align}
\int_{\abs{w}>\frac{r}{2}}\bigl((1+r^2)\abs{r\omega-w}^{\gamma+2}+r\abs{r\omega-w}^{\gamma+1}\bigr)\abs{F(w)}\dd w\lesssim r^5\ue^{-\frac{r^2}{32}}. \label{cvel:moment-complement}
\end{align}
Factor $\ue^{-\frac{r^2}{32}}$ from the Gaussian bound and apply \eqref{cvel:Gaussian-convolution-bound} to the remaining Gaussian with $s=\gamma+2,\gamma+1$. The omitted portions of the leading moments are bounded by a polynomial in $r$ times $\ue^{-\frac{r^2}{32}}$. Combining \eqref{cvel:trace-expansion-error}--\eqref{cvel:moment-complement} proves \eqref{eq:I:display:195}--\eqref{ci:weighted-drift-expansion}. Finally, the corresponding moments of $\m$ are $1$ and $2$, and
\begin{align}
\int_{\R^3}(1+\abs{w}^2)\abs{F(w)-\m(w)}\dd w\leq\delta\int_{\R^3}(1+\abs{w}^2)\mh(w)\dd w\lesssim \delta. \label{cvel:moment-perturbation}
\end{align}
This proves \eqref{cvel:coefficient-moments}.
\end{proof}

\begin{lemma}[Bounds for the coefficients and source terms]\label{cvel:scalar-bounds}
Let $0<\alpha<\frac{1}{4}$ and $0<\delta\leq1$. Suppose that
\begin{align}
\vnorm[\Big]{\ue^{\alpha\abs{v}^2}f}_{L_v^\infty}\leq\delta. \label{cvel:coefficient-Gaussian}
\end{align}
Let $F=\m+\mh f$, and define $\Bf{},\Cf{}$ by \eqref{cvel:scalar-coefficients}. Then
\begin{align}
\vnorm[\Big]{\la^{-4}\Bf{}}_{L_v^\infty}+\vnorm[\Big]{\la^{-5}\Cf{}}_{L_v^\infty}\lesssim1, \label{cvel:scalar-coefficient-growth}\\
\vnorm[\Big]{\la^{-5}\ue^{\alpha\abs{v}^2}\bigl(\Cf{}f+\Kop f\bigr)}_{L_v^\infty}\lesssim \delta. \label{cvel:scalar-source}
\end{align}
\end{lemma}

\begin{proof}
By \eqref{cvel:rough-coefficients}--\eqref{cvel:coefficient-c} and $\gamma\leq1$,
\begin{align}
\abs{b_F}+\abs{\Af{}v}\lesssim\la^4,\qquad \abs{v\cdot\Af{}v}+\abs{\tr\Af{}}+\abs{c_F}\lesssim\la^5. \label{cvel:scalar-growth-proof}
\end{align}
This proves \eqref{cvel:scalar-coefficient-growth}. To bound $\Kop f$, we use \cite[Lemma~2.2, (2.6)--(2.8)]{Carrapatoso.Mischler2024} in the formula \eqref{cvel:K}. Its weighted $L_v^\infty$ hypothesis for $\mh f$ is \eqref{cvel:coefficient-weighted-bound}, and the exponent condition was checked in the proof of Lemma~\ref{cvel:coefficients}. Also, \eqref{cvel:coefficient-Gaussian} implies $f\in L_v^2$, so \eqref{cvel:K-bilinear-identity} identifies this formula with $\Kop f$. Therefore
\begin{align}
\abs[\big]{\Kop f(v)}\lesssim\delta\la^5\ue^{-\frac{1}{4}\abs{v}^2}\lesssim\delta\la^5\ue^{-\alpha\abs{v}^2}. \label{cvel:scalar-K-bound}
\end{align}
At $\gamma=-3$, the term involving $c$ in \eqref{cvel:K} is $8\pi\m f$ and satisfies the same bound. Combining \eqref{cvel:scalar-coefficient-growth}, \eqref{cvel:coefficient-Gaussian} and \eqref{cvel:scalar-K-bound} proves \eqref{cvel:scalar-source}. The corresponding Gaussian estimate for $G\mapsto Q(G,\m)$ is also given in \cite[Lemma~2.3]{Carrapatoso.Mischler2024}.
\end{proof}

Under \eqref{cvel:coefficient-bounded}, identities \eqref{cvel:convolved-identities} and bound \eqref{cvel:rough-coefficients} give $\divv_v\Af{}=b_F$ and $\Af{}\in W^{1,\infty}_{\mathrm{loc}}$. For $h\in H^1_{\mathrm{loc}}(\R^3)$, we therefore use the distributional convention
\begin{align}
\Af{}:D_v^2h:=\divv_v(\Af{}\nabla_vh)-b_F\cdot\nabla_vh\quad\text{in }H^{-1}_{\mathrm{loc}}(\R^3). \label{cvel:weak-contraction}
\end{align}
If $D_v^2h\in L^2_{\mathrm{loc}}$, this agrees almost everywhere with the usual matrix contraction.

\begin{lemma}[Collision identities in distributions]\label{cvel:scalar}
Let $0<\alpha<\frac{1}{4}$, $0<\delta\leq1$, and let $f\in\Wv$ satisfy \eqref{cvel:coefficient-Gaussian}. Set $F=\m+\mh f$, with $\Af{},\Bf{},\Cf{}$ defined by \eqref{cvel:coefficient-convolutions} and \eqref{cvel:scalar-coefficients}. Then
\begin{align}
-\Lop f+\Ga(f,f)=\divv_v(\Af{}\nabla_vf)+\Bf{}\cdot\nabla_vf+\Cf{}f+\Kop f\quad\text{in }\mathcal D'(\R^3). \label{cvel:scalar-div}
\end{align}
Every term on the right belongs to $H^{-1}_{\mathrm{loc}}(\R^3)$. Equivalently, with convention \eqref{cvel:weak-contraction},
\begin{align}
-\Lop f+\Ga(f,f)=\Af{}:D_v^2f-(\Af{}v)\cdot\nabla_vf+\Cf{}f+\Kop f\quad\text{in }\mathcal D'(\R^3). \label{cvel:scalar-nondiv}
\end{align}
If also $D_v^2f\in L^2_{\mathrm{loc}}$, both identities hold in $L^2_{\mathrm{loc}}$ and hence almost everywhere.
\end{lemma}

\begin{proof}
The definition \eqref{eq:V} and the local ellipticity of $\sig$ imply $f\in H^1_{\mathrm{loc}}(\R^3)$. For every $R>0$, Lemma~\ref{cvel:scalar-bounds} and \eqref{cvel:rough-coefficients} give
\begin{align}
\Af{}\nabla_vf,\quad\Bf{}\cdot\nabla_vf,\quad\Cf{}f,\quad\Kop f\in L^2(B_R). \label{cvel:scalar-local-products}
\end{align}
Consequently, the right side of \eqref{cvel:scalar-div} belongs to $H^{-1}(B_R)$.

By \eqref{cvel:W-density}, choose $f_j\in C_c^\infty(\R^3)$ such that $f_j\to f$ in $\Wv$. Keep $F$ fixed. Direct differentiation gives
\begin{align}
\mih\partial_{ij}(\mh f_j)=\partial_{ij}f_j-\frac{1}{2}\bigl(v_i\partial_jf_j+v_j\partial_if_j\bigr)+\left(\frac{1}{4}v_iv_j-\frac{1}{2}\delta_{ij}\right)f_j. \label{eq:I:display:071}
\end{align}
Using \eqref{cvel:Q}, the symmetry of $\Af{}$ and $\partial_i\Af{}^{ij}=b_F^j$ from \eqref{cvel:convolved-identities}, we obtain
\begin{align}
\mih Q(F,\mh f_j)=\divv_v(\Af{}\nabla_vf_j)+\Bf{}\cdot\nabla_vf_j+\Cf{}f_j. \label{cvel:scalar-approximation}
\end{align}
Lemma~\ref{cvel:distribution} and multiplication by the smooth function $\mih$ on $B_R$ give
\begin{align}
\vnorm[\Big]{\mih Q(F,\mh(f_j-f))}_{H^{-1}(B_R)}\lesssim_R\vnorm{f_j-f}_{H^1(B_R)}\to0. \label{cvel:scalar-local-convergence}
\end{align}
The bound follows from \eqref{cvel:rough-coefficients}--\eqref{cvel:coefficient-c}, and local ellipticity gives $H^1(B_R)$ convergence from $\Wv$ convergence.

To identify this limit with $\Ga(\mh,f)+\Ga(f,f)$ in $\Vv'$, we also approximate the first argument. The local integrability of $a,b$ and the Gaussian factor $\mh$ give
\begin{align}
\vnorm[\Big]{a*(\mh(f_k-f))}_{L^2(B_R)}+\vnorm[\Big]{b*(\mh(f_k-f))}_{L^2(B_R)}\lesssim_R\vnorm{f_k-f}_{L_v^2}\to0. \label{cvel:scalar-first-convergence}
\end{align}
For the part $\abs{v-w}\leq1$, this is Young's inequality with $a\ind_{B_1},b\ind_{B_1}\in L^1$. On $\abs{v-w}>1$ with $v\in B_R$, the kernels have polynomial growth, and Cauchy--Schwarz in $w$ gives the same bound because $\mh(w)$ decays exponentially. Thus, for fixed $j$ and $\varphi\in C_c^\infty(\R^3)$, the divergence formula \eqref{cvel:local-collision-pairing} gives
\begin{align}
\ip*{\Ga(f_k,f_j)}{\varphi}\to\ip*{\mih Q(\mh f,\mh f_j)}{\varphi}. \label{cvel:scalar-first-pairing}
\end{align}
On the other hand, \eqref{cvel:trilinear} gives $\Ga(f_k,f_j)\to\Ga(f,f_j)$ in $\Vv'$. Applying the same argument to the first argument $\mh$ shows that
\begin{align}
\mih Q(F,\mh f_j)=\Ga(\mh,f_j)+\Ga(f,f_j)\quad\text{in }\mathcal D'(\R^3). \label{cvel:scalar-form-identification}
\end{align}
Finally, \eqref{cvel:trilinear} gives
\begin{align}
\vnorm[\Big]{\Ga(\mh,f_j-f)+\Ga(f,f_j-f)}_{\Vv'}\lesssim \bigl(\vnorm[\Big]{\mh}_{L_v^2}+\vnorm{f}_{L_v^2}\bigr)\vnorm{f_j-f}_{\Vv}\to0. \label{cvel:scalar-form-convergence}
\end{align}
Passing to the limit in \eqref{cvel:scalar-approximation} by \eqref{cvel:scalar-local-convergence} and \eqref{cvel:scalar-form-convergence}, and using \eqref{cvel:linear-bilinear-identity}, proves \eqref{cvel:scalar-div}. In particular, bilinearity and $Q(\m,\m)=0$ give
\begin{align}
-\Lop f+\Ga(f,f)=\mih Q(F,\mh f)+\Kop f=\mih Q(F,F)\quad\text{in }\mathcal D'(\R^3). \label{cvel:scalar-bilinearity}
\end{align}
Convention \eqref{cvel:weak-contraction} gives
\begin{align}
\divv_v(\Af{}\nabla_vf)=\Af{}:D_v^2f+b_F\cdot\nabla_vf. \label{cvel:scalar-product-rule}
\end{align}
Substitution in \eqref{cvel:scalar-div} proves \eqref{cvel:scalar-nondiv}. If $D_v^2f\in L^2_{\mathrm{loc}}$, the product rule identifies \eqref{cvel:scalar-product-rule} in $L^2(B_R)$, and \eqref{cvel:scalar-local-products} proves the stated $L^2$ identities.
\end{proof}

\begin{lemma}[Weighted $L^2$ collision estimates]\label{cvel:weighted-collisions}
Let $0\leq t_0<T<\infty$, $0<\alpha<\frac{1}{4}$ and $0<\delta\leq1$. Suppose that
\begin{align}
\tnorm[\Big]{\ue^{\alpha\abs{v}^2}f}_{L^\infty((t_0,T)\times\phase)}\leq\delta. \label{cvel:weighted-Gaussian}
\end{align}
If, for every $m\geq0$,
\begin{align}
\la^mf,\quad\la^m\nabla_vf,\quad\la^mD_v^2f\in L^2\bigl((t_0,T)\times\Om\times\R^3\bigr), \label{cvel:scalar-graph-assumption}
\end{align}
then, for every $m\geq0$,
\begin{align}
\la^m\Lop f,\quad\la^m\Ga(f,f)\in L^2\bigl((t_0,T)\times\Om\times\R^3\bigr). \label{cvel:scalar-graph-conclusion}
\end{align}
\end{lemma}

\begin{proof}
The assumption \eqref{cvel:scalar-graph-assumption} implies $f(t,x,\cdot)\in\Wv\cap H^2_{\mathrm{loc}}(\R^3)$ for almost every $(t,x)$. Apply Lemma~\ref{cvel:scalar} at those points. The approximation in its proof also identifies each operator separately, and \eqref{cvel:Q}, \eqref{eq:I:display:071} give
\begin{align}
\Ga(f,f)&=\bigl(a*(\mh f)\bigr):D_v^2f-\Bigl[\bigl(a*(\mh f)\bigr)v\Bigr]\cdot\nabla_vf \label{cvel:graph-expansion}\\
&\quad+\left[\frac{1}{4}v\cdot\bigl(a*(\mh f)\bigr)v-\frac{1}{2}\tr\bigl(a*(\mh f)\bigr)-c*(\mh f)\right]f.\notag
\end{align}
The corresponding linear formula is
\begin{align}
\Lop f=-\sig:D_v^2f+(\sig v)\cdot\nabla_vf-\left(\frac{1}{4}v\cdot\sig v-\frac{1}{2}\tr\sig-c*\m\right)f-\Kop f. \label{cvel:L-expansion}
\end{align}
The coefficient bounds and \eqref{cvel:scalar-K-bound} imply, almost everywhere in $(t_0,T)\times\phase$,
\begin{align}
\abs{\Lop f}+\abs{\Ga(f,f)}\lesssim\la^5\bigl(\abs{D_v^2f}+\abs{\nabla_vf}+\abs{f}\bigr)+\delta\la^5\ue^{-\frac{1}{4}\abs{v}^2}. \label{cvel:scalar-graph-pointwise}
\end{align}
At $\gamma=-3$, use $c*(\mh f)=-8\pi\mh f$ in \eqref{cvel:graph-expansion}. Multiplying \eqref{cvel:scalar-graph-pointwise} by $\la^m$ proves \eqref{cvel:scalar-graph-conclusion}: the derivative terms belong to $L^2$ by \eqref{cvel:scalar-graph-assumption} with $m+5$, and the Gaussian term has every polynomial moment since $T<\infty$ and $\Om$ is bounded.
\end{proof}


\section{Kinetic Traces and Energy Identities}\label{sec:trace}

We will apply the kinetic trace theorem of \cite[Theorem~2.8]{Carrapatoso.Mischler2024}, based on \cite[Section~4]{Mischler2010}, and use Gaussian bounds to obtain the boundary integrability and $L^2$ continuity needed for the energy identities.

We consider
\begin{align}
\Y f+\Lop f=S,\qquad S\in L^2(0,T;\Vspace'), \label{trace:linearized}
\end{align}
with $\Lop$ from \eqref{eq:linearized-operators}. Corollary~\ref{cvel:spatial} gives
\begin{align}
\norm{\Lop f}_{\Vspace'}\lesssim\norm{f}_{\Vspace},\qquad \norm[\big]{\Ga(f,f)}_{\Vspace'}\lesssim\norm{f}_{L^\infty_xL^2_v}\norm{f}_{\Vspace}. \label{trace:dual-bounds}
\end{align}
Thus the weak Landau equation \eqref{eq:intro-interior-weak}, together with \eqref{eq:intro-trace-space}--\eqref{eq:weak-gaussian-bound}, is a case of \eqref{trace:linearized} with $S=\Ga(f,f)$. The same estimates apply to the difference equation and to the linear equations in Section~\ref{sec:energy}.

\subsection{The Landau Trace and Its Gaussian Bound}\label{sec:trace-definition}

\begin{lemma}[Kinetic trace for the linearized Landau equation]\label{trace:CM}
Let $T>0$ and $S\in L^2(0,T;\Vspace')$. Suppose that
\begin{align}
f\in L^\infty(0,T;L^2)\cap L^2(0,T;\Vspace). \label{trace:bulk-hypotheses}
\end{align}
Assume that $f$ solves \eqref{trace:linearized}, meaning that
\begin{align}
-\int_0^T\bigl(f,\Y\varphi\bigr)_{L^2}\dd t+\int_0^T\ip*{\Lop f}{\varphi}\dd t=\int_0^T\ip*{S}{\varphi}\dd t \label{trace:linearized-weak}
\end{align}
for every $\varphi\in C_c^\infty\bigl((0,T)\times\Om\times\R^3\bigr)$. The following conclusions hold.
\begin{enumerate}
\item The function $f$ has a representative and a unique kinetic trace $\trace f$ such that
\begin{align}
f\in C\bigl([0,T];L^2_{\mathrm{loc}}(\Om\times\R^3)\bigr),\qquad (n\cdot v)\trace f\in L^2_{\mathrm{loc}}\bigl((0,T)\times\partial\Om\times\R^3\bigr). \label{trace:local-trace}
\end{align}
The trace is determined, up to a set of zero $\dd t\dd\gamma$ measure, by
\begin{align}
-\int_0^T\bigl(f,\Y\varphi\bigr)_{L^2}\dd t+\int_0^T\int_{\partial\Om\times\R^3}\trace f\,\varphi(n\cdot v)\dd S_x\dd v\dd t=\int_0^T\ip*{S-\Lop f}{\varphi}\dd t \label{trace:compact-Green}
\end{align}
for every smooth $\varphi$ on $(0,T)\times\overline\Om\times\R^3$ with compact support in time and velocity. The boundary norm in \eqref{trace:local-trace} uses the measure $\dd t\dd S_x\dd v$.
\item Let $f$ and $\widetilde f$ satisfy the hypotheses with sources $S$ and $\widetilde S$, respectively. Then $h:=\widetilde f-f$ satisfies \eqref{trace:linearized-weak} with source $\widetilde S-S$, and its trace is
\begin{align}
\trpm h=\trpm\widetilde f-\trpm f. \label{trace:difference-trace}
\end{align}
\item For the corresponding density $F:=\m+\mh f$, its kinetic trace is given by \eqref{eq:intro-density-trace}. The boundary conditions \eqref{eq:inflow-full} and \eqref{eq:diffuse-full} are equivalent to \eqref{eq:inflow-normalized} and \eqref{eq:diffuse-normalized}, respectively, whenever the boundary integrals are well-defined.
\item Suppose additionally that
\begin{align}
f\in C\bigl([0,T];L^2\bigr). \label{trace:energy-hypotheses}
\end{align}
For every $0\leq s\leq t\leq T$ and every smooth $\varphi$ on $[0,T]\times\overline\Om\times\R^3$ with compact velocity support, the Green identity is
\begin{align}
\Bigl[\bigl(f(r),\varphi(r)\bigr)_{L^2}\Bigr]_{r=s}^{r=t}-\int_s^t\bigl(f,\Y\varphi\bigr)_{L^2}\dd r+\int_s^t\int_{\partial\Om\times\R^3}\trace f\,\varphi(n\cdot v)\dd S_x\dd v\dd r =\int_s^t\ip*{S-\Lop f}{\varphi}\dd r. \label{trace:linear-Green}
\end{align}
\end{enumerate}
\end{lemma}

\begin{proof}
We apply \cite[Theorem~2.8(1)--(2)]{Carrapatoso.Mischler2024}. The equation in that theorem has a velocity-diffusion term on the left. We therefore rewrite \eqref{trace:linearized-weak}, in distributions, as
\begin{align}
\Y f-\divv_v(\sig\nabla_vf)=S-\Lop f-\divv_v(\sig\nabla_vf). \label{trace:fixed-equation}
\end{align}
For this application, we verify, for every $R<\infty$,
\begin{align}
\sig(v)=\sig(v)^{\mathsf T},\qquad C_R^{-1}I\leq\sig(v)\leq C_RI\quad(v\in B_R), \label{trace:local-matrix-assumption}
\end{align}
and
\begin{align}
f\in L^2\bigl((0,T)\times\Om;H^1(B_R)\bigr),\qquad S-\Lop f-\divv_v(\sig\nabla_vf)\in L^2\bigl((0,T)\times\Om;H^{-1}(B_R)\bigr). \label{trace:local-function-source}
\end{align}
Here $H^{-1}(B_R)$ is the dual of $H^1_0(B_R)$.

The matrix $\sig$ in \eqref{eq:sigma-definition} satisfies \eqref{trace:local-matrix-assumption} by Lemma~\ref{cvel:convolution}. Its lower eigenvalue bound and the zeroth-order term of \eqref{eq:V} give
\begin{align}
\tnorm{f}_{L^2((0,T)\times\Om;H^1(B_R))}\lesssim_R\tnorm{f}_{L^2(0,T;\Vspace)}. \label{trace:local-H1-bound}
\end{align}
For $\varphi\in\Vspace$, define the divergence in \eqref{trace:fixed-equation} by
\begin{align}
\ip[\big]{\divv_v(\sig\nabla_vf)}{\varphi}=-\int_{\phase}\sig\nabla_vf\cdot\nabla_v\varphi\dd x\dd v. \label{trace:diffusion-pairing}
\end{align}
Cauchy--Schwarz, \eqref{eq:V} and \eqref{trace:dual-bounds} give
\begin{align}
\tnorm{S-\Lop f-\divv_v(\sig\nabla_vf)}_{L^2(0,T;\Vspace')}\lesssim\tnorm{S}_{L^2(0,T;\Vspace')}+\tnorm{f}_{L^2(0,T;\Vspace)}. \label{trace:fixed-source}
\end{align}
For $\varphi\in L^2(\Om;H^1_0(B_R))$, its extension by zero for $v\notin B_R$, denoted by $\ind_{B_R}\varphi$, satisfies
\begin{align}
\norm{\ind_{B_R}\varphi}_{\Vspace}\lesssim_R\norm{\varphi}_{L^2(\Om;H^1(B_R))}. \label{trace:local-test-bound}
\end{align}
Testing the source in \eqref{trace:fixed-source} against these extensions yields
\begin{align}
\tnorm{S-\Lop f-\divv_v(\sig\nabla_vf)}_{L^2((0,T)\times\Om;H^{-1}(B_R))}\lesssim_R\left(\tnorm{S}_{L^2(0,T;\Vspace')}+\tnorm{f}_{L^2(0,T;\Vspace)}\right). \label{trace:local-source-bound}
\end{align}
This proves \eqref{trace:local-function-source}. Thus \cite[Theorem~2.8(1)--(2)]{Carrapatoso.Mischler2024} gives \eqref{trace:local-trace}. We record the identity from part~(1) in the notation of \eqref{trace:fixed-equation}. For $\beta\in C^2(\R)$ with $\beta',\beta''\in L^\infty(\R)$, its equation~(2.25), after one integration by parts in velocity, reads
\begin{align}
&\left[\int_{\phase}\beta\bigl(f(r)\bigr)\varphi(r)\dd x\dd v\right]_{r=s}^{r=t}+\int_s^t\int_{\partial\Om\times\R^3}\beta(\trace f)\varphi(n\cdot v)\dd S_x\dd v\dd r \label{trace:beta-Green}\\
&\quad=\int_s^t\int_{\phase}\beta(f)\Y\varphi\dd x\dd v\dd r+\int_s^t\ip[\Big]{S-\Lop f-\divv_v(\sig\nabla_vf)}{\beta'(f)\varphi}\dd r\notag\\
&\qquad-\int_s^t\int_{\phase}\sig\nabla_vf\cdot\nabla_v\bigl(\beta'(f)\varphi\bigr)\dd x\dd v\dd r.\notag
\end{align}
Here $0\leq s\leq t\leq T$, and $\varphi$ is smooth on $[0,T]\times\overline\Om\times\R^3$ with compact velocity support. The time values and the boundary trace in \eqref{trace:beta-Green} are those supplied by the cited theorem. The Sobolev chain rule gives
\begin{align}
\nabla_v\bigl(\beta'(f)\varphi\bigr)=\beta''(f)\varphi\nabla_vf+\beta'(f)\nabla_v\varphi. \label{trace:beta-test-gradient}
\end{align}
The bounds on $\beta',\beta''$ and \eqref{trace:local-H1-bound} imply $\beta'(f)\varphi\in L^2(0,T;\Vspace)$, so the source pairing is well-defined. By \eqref{trace:diffusion-pairing},
\begin{align}
\ip[\Big]{S-\Lop f-\divv_v(\sig\nabla_vf)}{\beta'(f)\varphi}=\ip*{S-\Lop f}{\beta'(f)\varphi} +\int_{\phase}\sig\nabla_vf\cdot\nabla_v\bigl(\beta'(f)\varphi\bigr)\dd x\dd v. \label{trace:beta-source-cancellation}
\end{align}
This cancels the last line of \eqref{trace:beta-Green}. Taking $s=0$, $t=T$, $\beta(\xi)=\xi$ and $\varphi$ compactly supported in time gives \eqref{trace:compact-Green}.

Subtracting \eqref{trace:compact-Green} for two traces shows that their difference times $n\cdot v$ vanishes as a boundary distribution. This proves uniqueness $\dd t\dd\gamma$-almost everywhere. The same subtraction for $(\widetilde f,\widetilde S)$ and $(f,S)$ identifies the trace of $\widetilde f-f$ and proves (ii).

For (iii), use $\mh\varphi$ in \eqref{trace:compact-Green}. Since $\Y(\mh\varphi)=\mh\Y\varphi$, this identifies the trace of $\mh f$. Adding the ordinary trace of $\m$ gives
\begin{align}
\trace(\mh f)=\mh\trace f,\qquad \trace F=\m+\mh\trace f. \label{trace:density-multiplication}
\end{align}
Substitution into \eqref{eq:inflow-full} and \eqref{eq:diffuse-full}, using \eqref{eq:wall-normalization}, proves the stated equivalence of the boundary conditions.

For (iv), take $\beta(\xi)=\xi$ in \eqref{trace:beta-Green} and use \eqref{trace:beta-source-cancellation}. Uniqueness identifies the boundary trace with $\trace f$, and \eqref{trace:energy-hypotheses} identifies the time values with $f(s),f(t)$, giving \eqref{trace:linear-Green}.
\end{proof}

The Gaussian bound gives the full boundary $L^2(\dd t\dd\gamma)$ estimate and makes the diffuse-reflection integral \eqref{eq:diffuse-flux} well-defined. It also supplies the time continuity required in Lemma~\ref{trace:CM}(iv).

\begin{lemma}[Gaussian bounds and continuity]\label{trace:gaussian}
Let $f,S$ satisfy the hypotheses of Lemma~\ref{trace:CM}. Suppose additionally that, for $M\geq0$ and $\alpha>0$,
\begin{align}
\tnorm{f}_{L^\infty(0,T;L^{\infty,\alpha})}\leq M. \label{trace:Gaussian-assumption}
\end{align}
Let $\trace f$ be the trace in Lemma~\ref{trace:CM}. Then the following hold.
\begin{enumerate}
\item The kinetic traces satisfy
\begin{align}
\tnorm{\trpm f}_{L^\infty(0,T;L^{\infty,\alpha}_{\Sigma_\pm})}\leq M. \label{trace:Gaussian-boundary}
\end{align}
In particular,
\begin{align}
\tnorm{\trpm f}_{L^2\bigl((0,T)\times\Sigma_\pm;\dd t\dd\gamma\bigr)}^2\leq\frac{\pi M^2T\abs{\partial\Om}}{8\alpha^2}. \label{trace:Gaussian-full-traces}
\end{align}
\item The scalar integral $\flux{f}$ in \eqref{eq:diffuse-flux} is well-defined for almost every $(t,x)\in(0,T)\times\partial\Om$, and
\begin{align}
\tnorm{\flux{f}}_{L^\infty\bigl((0,T)\times\partial\Om\bigr)}\leq\frac{\pi M}{2(2\pi)^{\frac{3}{4}}\left(\alpha+\frac{1}{4}\right)^2}. \label{trace:scalar-bound}
\end{align}
\item The representative in Lemma~\ref{trace:CM}(i) satisfies
\begin{align}
f\in C\bigl([0,T];L^2\bigr),\qquad \gnorm[\big]{\infty}{\alpha}{f(t)}\leq M\quad(0\leq t\leq T). \label{trace:Gaussian-time-continuity}
\end{align}
Thus Lemma~\ref{trace:CM}(iv) applies.
\item If $\lambda\geq0$ and $\tnorm[\Big]{\ue^{\lambda t}f}_{L^\infty(0,T;L^{\infty,\alpha})}\leq M$, then \eqref{trace:Gaussian-boundary} and \eqref{trace:scalar-bound} hold with $\trpm f$ and $\flux{f}$ replaced by $\ue^{\lambda t}\trpm f$ and $\ue^{\lambda t}\flux{f}$, respectively.
\end{enumerate}
\end{lemma}

\begin{proof}
\proofstep{The equation after multiplication by the Gaussian weight} For each $R<\infty$, the product rule and \eqref{trace:local-H1-bound} give
\begin{align}
\tnorm[\Big]{\ue^{\alpha\abs{v}^2}f}_{L^2((0,T)\times\Om;H^1(B_R))}\lesssim_R\tnorm{f}_{L^2((0,T)\times\Om;H^1(B_R))}\lesssim_R\tnorm{f}_{L^2(0,T;\Vspace)}. \label{trace:weighted-H1}
\end{align}

Testing \eqref{trace:linearized-weak} against $\ue^{\alpha\abs{v}^2}\varphi$, for $\varphi\in C_c^\infty((0,T)\times\Om\times\R^3)$, and using $\Y\ue^{\alpha\abs{v}^2}=0$ gives, in distributions, the equation required by \cite[Theorem~2.8]{Carrapatoso.Mischler2024}:
\begin{align}
\Y\bigl(\ue^{\alpha\abs{v}^2}f\bigr)-\divv_v\Bigl(\sig\nabla_v\bigl(\ue^{\alpha\abs{v}^2}f\bigr)\Bigr) =\ue^{\alpha\abs{v}^2}(S-\Lop f)-\divv_v\Bigl(\sig\nabla_v\bigl(\ue^{\alpha\abs{v}^2}f\bigr)\Bigr). \label{trace:weighted-equation}
\end{align}

We check its source in $H^{-1}(B_R)$ directly. For $\varphi\in L^2(\Om;H^1_0(B_R))$, the zero extension of the weighted test satisfies
\begin{align}
\norm[\Big]{\ind_{B_R}\ue^{\alpha\abs{v}^2}\varphi}_{\Vspace}\lesssim_R\norm{\varphi}_{L^2(\Om;H^1(B_R))}. \label{trace:weighted-test-bound}
\end{align}
The source pairing is defined as in \eqref{trace:diffusion-pairing}. By \eqref{trace:weighted-test-bound}, \eqref{trace:dual-bounds} and Cauchy--Schwarz for the diffusion term, \eqref{trace:weighted-H1} gives
\begin{align}
&\tnorm[\Big]{\ue^{\alpha\abs{v}^2}(S-\Lop f)-\divv_v\Bigl(\sig\nabla_v\bigl(\ue^{\alpha\abs{v}^2}f\bigr)\Bigr)}_{L^2((0,T)\times\Om;H^{-1}(B_R))} \label{trace:weighted-source-terms}\\
&\qquad\lesssim_R\left(\tnorm{S}_{L^2(0,T;\Vspace')}+\tnorm{f}_{L^2(0,T;\Vspace)}\right).\notag
\end{align}
The matrix bound \eqref{trace:local-matrix-assumption}, the velocity regularity \eqref{trace:weighted-H1} and the source bound \eqref{trace:weighted-source-terms} verify the hypotheses of \cite[Theorem~2.8(1)--(2)]{Carrapatoso.Mischler2024}. Therefore $\ue^{\alpha\abs{v}^2}f$ has a kinetic trace and satisfies the identity \eqref{trace:beta-Green} with the function and source from \eqref{trace:weighted-equation}.

\proofstep{The Gaussian bound on the trace} Let $\beta\in C^2(\R)$ with $\beta',\beta''\in L^\infty(\R)$, and let $\varphi$ be smooth up to $\partial\Om$ with compact support in time and velocity. Apply \eqref{trace:beta-Green} to \eqref{trace:weighted-equation} on $(0,T)$. The time terms vanish. The diffusion terms cancel as in \eqref{trace:beta-source-cancellation}, giving
\begin{align}
&\int_0^T\int_{\partial\Om\times\R^3}\beta\Bigl(\trace\bigl(\ue^{\alpha\abs{v}^2}f\bigr)\Bigr)\varphi(n\cdot v)\dd S_x\dd v\dd t \label{trace:weighted-beta-Green}\\
&\quad=\int_0^T\int_{\phase}\beta\bigl(\ue^{\alpha\abs{v}^2}f\bigr)\Y\varphi\dd x\dd v\dd t+\int_0^T\ip[\Big]{S-\Lop f}{\ue^{\alpha\abs{v}^2}\beta'\bigl(\ue^{\alpha\abs{v}^2}f\bigr)\varphi}\dd t.\notag
\end{align}
The test in the last pairing belongs to $L^2(0,T;\Vspace)$ by the Sobolev chain rule and \eqref{trace:weighted-H1}, since its velocity support is compact. First take $\beta(\xi)=\xi$ in \eqref{trace:weighted-beta-Green} and compare with \eqref{trace:compact-Green} tested against $\ue^{\alpha\abs{v}^2}\varphi$. Uniqueness of the trace gives
\begin{align}
\trace\bigl(\ue^{\alpha\abs{v}^2}f\bigr)=\ue^{\alpha\abs{v}^2}\trace f\quad\dd t\dd\gamma\text{-almost everywhere}. \label{trace:weighted-trace}
\end{align}

Choose a nonnegative $\beta\in C^2(\R)$ with $\beta,\beta',\beta''\in L^\infty(\R)$, with $\beta=\beta'=\beta''=0$ on $[-M,M]$ and $\beta>0$ outside this interval. For a smooth $\varphi$ compactly supported in time and velocity, choose $R$ containing its velocity support. By \eqref{trace:Gaussian-assumption} and the Sobolev chain rule,
\begin{align}
\beta'\bigl(\ue^{\alpha\abs{v}^2}f\bigr)\varphi=0\quad\text{in }L^2\bigl((0,T)\times\Om;H^1_0(B_R)\bigr), \label{trace:vanishing-source-test}
\end{align}
since $\beta'$ and $\beta''$ vanish on $[-M,M]$. The first term on the right of \eqref{trace:weighted-beta-Green} vanishes because $\beta\bigl(\ue^{\alpha\abs{v}^2}f\bigr)=0$. By \eqref{trace:vanishing-source-test} and \eqref{trace:weighted-test-bound}, the test in its duality pairing is zero in $L^2(0,T;\Vspace)$, so the second term also vanishes. Using \eqref{trace:weighted-trace}, we obtain
\begin{align}
\int_0^T\int_{\partial\Om\times\R^3}\beta\bigl(\ue^{\alpha\abs{v}^2}\trace f\bigr)\varphi(n\cdot v)\dd S_x\dd v\dd t=0. \label{trace:weighted-renormalization}
\end{align}
Arbitrary boundary tests give $\beta\bigl(\ue^{\alpha\abs{v}^2}\trace f\bigr)=0$ away from $n\cdot v=0$. The choice of $\beta$ therefore implies \eqref{trace:Gaussian-boundary}.

\proofstep{The boundary integrals} For each $x\in\partial\Om$, rotate velocity so that $n(x)\cdot v=v_3$. For every $\kappa>0$,
\begin{align}
\int_{v_3>0}\ue^{-\kappa\abs{v}^2}v_3\dd v=\left(\int_{\R^2}\ue^{-\kappa\abs{v'}^2}\dd v'\right)\left(\int_0^\infty s\ue^{-\kappa s^2}\dd s\right)=\frac{\pi}{2\kappa^2}. \label{trace:Gaussian-flux-integral}
\end{align}
Taking $\kappa=2\alpha$ and integrating over $(0,T)\times\partial\Om$ proves \eqref{trace:Gaussian-full-traces}. Taking $\kappa=\alpha+\frac14$ and using $\mh=(2\pi)^{-\frac34}\ue^{-\frac{\abs{v}^2}{4}}$ proves absolute convergence of \eqref{eq:diffuse-flux} and the bound \eqref{trace:scalar-bound}.

\proofstep{Continuity in $L^2$} For $\varrho>0$, write $\Om_\varrho:=\{x\in\Om:d(x)>\varrho\}$. Fix $t\in[0,T]$ and choose times $t_j\to t$ from the full-measure set on which $\gnorm{\infty}{\alpha}{f(t_j)}\leq M$. Writing $\xi_+:=\max\{\xi,0\}$, the Gaussian bound at $t_j$ and \eqref{trace:local-trace} give, for each $m\geq1$,
\begin{align}
\norm[\Big]{\bigl(\abs{f(t)}-M\ue^{-\alpha\abs{v}^2}\bigr)_+}_{L^2(\Om_{\frac1m}\times B_m)}\leq\norm[\big]{f(t)-f(t_j)}_{L^2(\Om_{\frac1m}\times B_m)}\to0\quad(j\to\infty). \label{trace:Gaussian-every-time}
\end{align}
Exhausting $\Om\times\R^3$ by $\Om_{\frac1m}\times B_m$ proves $\gnorm[\big]{\infty}{\alpha}{f(t)}\leq M$ for every $t\in[0,T]$.

Now let $t_j\to t$ be any sequence in $[0,T]$. Splitting the integral into $\Om_{\frac1m}\times B_m$, the region $\abs{v}>m$, and the boundary strip gives
\begin{align}
\norm[\big]{f(t_j)-f(t)}_{L^2}^2&\leq\norm[\big]{f(t_j)-f(t)}_{L^2(\Om_{\frac1m}\times B_m)}^2+4M^2\abs{\Om}\int_{\abs{v}>m}\ue^{-2\alpha\abs{v}^2}\dd v \label{trace:global-time-continuity}\\
&\quad+4M^2\abs[\big]{\Om\setminus\Om_{\frac1m}}\int_{\R^3}\ue^{-2\alpha\abs{v}^2}\dd v.\notag
\end{align}
For fixed $m$, the first term tends to zero as $j\to\infty$ by \eqref{trace:local-trace}. The other two terms tend to zero as $m\to\infty$. This proves \eqref{trace:Gaussian-time-continuity}.

\proofstep{The time weight} For (iv), use $\ue^{\lambda t}\varphi$ in \eqref{trace:linearized-weak}. The identity $\Y(\ue^{\lambda t}\varphi)=\ue^{\lambda t}(\Y\varphi+\lambda\varphi)$ gives, in distributions,
\begin{align}
\Y\bigl(\ue^{\lambda t}f\bigr)+\Lop\bigl(\ue^{\lambda t}f\bigr)=\ue^{\lambda t}(S+\lambda f). \label{trace:time-weighted-equation}
\end{align}
The term $\lambda f$ belongs to $L^2(0,T;\Vspace')$. To see this, the Cauchy--Schwarz inequality and \eqref{eq:V} give, for $\varphi\in\Vspace$,
\begin{align}
\abs{\int_{\phase}f\varphi\dd x\dd v}\leq\norm[\Big]{\la^{-\frac{\gamma+2}{2}}f}_{L^2}\norm[\Big]{\la^{\frac{\gamma+2}{2}}\varphi}_{L^2}\leq\norm[\Big]{\la^{-\frac{\gamma+2}{2}}f}_{L^2}\norm{\varphi}_{\Vspace}. \label{trace:time-weighted-source}
\end{align}
The Gaussian bound makes the first factor bounded in time. Since $T<\infty$, the right-hand side of \eqref{trace:time-weighted-equation} belongs to $L^2(0,T;\Vspace')$. Comparing \eqref{trace:compact-Green} for $\ue^{\lambda t}f$ with the same identity for $f$ tested against $\ue^{\lambda t}\varphi$ gives $\trace(\ue^{\lambda t}f)=\ue^{\lambda t}\trace f$. Parts (i)--(ii), applied to $\ue^{\lambda t}f$, prove (iv).
\end{proof}

\subsection{Energy Identities}\label{sec:trace-green}

To express the boundary term in the $L^2$ identity, define the operator on outgoing velocities by
\begin{align}
\Pp\trp f(t,x,v):=\sqrt{2\pi}\,\mh(v)\flux{f}(t,x),\qquad (x,v)\in\Sigma_+. \label{eq:projection}
\end{align}
The diffuse-reflection condition gives the following boundary identities.

\begin{lemma}\label{trace:projection-lemma}
Fix a boundary point $x\in\partial\Om$, and use the velocity norms
\begin{align}
\vnorm{f}_{L^2_{v,\pm}}^2:=\int_{\{\pm n(x)\cdot v>0\}}\abs[\big]{f(v)}^2\abs[\big]{n(x)\cdot v}\dd v. \label{trace:velocity-flux-norms}
\end{align}
Suppose that $\trp f(t,x,\cdot)\in L^2_{v,+}$ and $f$ satisfies the diffuse-reflection condition \eqref{eq:diffuse-normalized}. Then
\begin{align}
\vnorm{\Pp\trp f}_{L^2_{v,+}}^2&=\vnorm{\trm f}_{L^2_{v,-}}^2=\sqrt{2\pi}\,\abs{\flux{f}}^2, \label{trace:projection-identity-a}\\
\vnorm{\trp f}_{L^2_{v,+}}^2-\vnorm{\trm f}_{L^2_{v,-}}^2&=\vnorm[\big]{(I-\Pp)\trp f}_{L^2_{v,+}}^2. \label{trace:projection-identity}
\end{align}
\end{lemma}

\begin{proof}
The Cauchy--Schwarz inequality and \eqref{eq:wall-normalization} give
\begin{align}
\abs{\flux{f}}^2\leq\vnorm{\trp f}_{L^2_{v,+}}^2\int_{n(x)\cdot v>0}\m(v)\bigl(n(x)\cdot v\bigr)\dd v=\frac{1}{\sqrt{2\pi}}\vnorm{\trp f}_{L^2_{v,+}}^2. \label{trace:projection-flux-bound}
\end{align}
Thus the integral \eqref{eq:diffuse-flux} is finite and $\Pp\trp f\in L^2_{v,+}$. The definition \eqref{eq:projection} then yields
\begin{align}
\int_{n(x)\cdot v>0}\bigl(\trp f-\Pp\trp f\bigr)\mh(v)\bigl(n(x)\cdot v\bigr)\dd v=0. \label{trace:projection-orthogonality}
\end{align}
Since $\Pp\trp f$ is proportional to $\mh$, orthogonality gives
\begin{align}
\vnorm{\trp f}_{L^2_{v,+}}^2=\vnorm{\Pp\trp f}_{L^2_{v,+}}^2+\vnorm[\big]{(I-\Pp)\trp f}_{L^2_{v,+}}^2. \label{trace:projection-Pythagoras}
\end{align}
Squaring $\trm f=\sqrt{2\pi}\mh\flux{f}$ and \eqref{eq:projection}, then using \eqref{eq:wall-normalization}, gives \eqref{trace:projection-identity-a}. Substitution into \eqref{trace:projection-Pythagoras} proves \eqref{trace:projection-identity}.
\end{proof}

The identities may be integrated in $(t,x)$. By linearity of the boundary condition, they also apply to differences of solutions.

We next apply the trace theorem to \eqref{trace:linearized} with the additional time continuity \eqref{trace:energy-hypotheses}. An incoming trace in $L^2(\dd t\dd\gamma)$ then gives the outgoing trace in the same space and the following energy identity. Lemma~\ref{trace:gaussian} supplies both properties when a Gaussian bound is available.

\begin{proposition}\label{trace:green}
Let $f\in C([0,T];L^2)\cap L^2(0,T;\Vspace)$ and $S\in L^2(0,T;\Vspace')$. Suppose that $f$ solves \eqref{trace:linearized} in the sense of \eqref{trace:linearized-weak}, and let $\trace f$ be its trace from Lemma~\ref{trace:CM}. The following identities hold for every $0\leq s\leq t\leq T$.
\begin{enumerate}
\item \textup{In-flow boundary.} Suppose that $\trm f=g$ with $g\in L^2\bigl((0,T)\times\Sigmam;\dd t\dd\gamma\bigr)$. Then $\trp f\in L^2\bigl((0,T)\times\Sigmap;\dd t\dd\gamma\bigr)$ and
\begin{align}
\frac12\norm[\big]{f(t)}_{L^2}^2+\frac12\int_s^t\norm[\big]{\trp f(r)}_{L^2_{\Sigma_+}}^2\dd r=\frac12\norm[\big]{f(s)}_{L^2}^2+\frac12\int_s^t\norm[\big]{g(r)}_{L^2_{\Sigma_-}}^2\dd r +\int_s^t\ip*{\Y f}{f}\dd r. \label{trace:full-Green}
\end{align}
\item \textup{Diffuse-reflection boundary.} Suppose that the integral \eqref{eq:diffuse-flux} is well-defined and
\begin{align}
\trm f=\Pm\trp f\in L^2\bigl((0,T)\times\Sigmam;\dd t\dd\gamma\bigr). \label{trace:diffuse-energy-assumption}
\end{align}
Then $\trp f\in L^2\bigl((0,T)\times\Sigmap;\dd t\dd\gamma\bigr)$ and
\begin{align}
\frac12\norm[\big]{f(t)}_{L^2}^2+\frac12\int_s^t\norm[\big]{(I-\Pp)\trp f(r)}_{L^2_{\Sigma_+}}^2\dd r=\frac12\norm[\big]{f(s)}_{L^2}^2+\int_s^t\ip*{\Y f}{f}\dd r. \label{trace:energy-Green}
\end{align}
\end{enumerate}
In both identities, $\Y f=S-\Lop f\in L^2(0,T;\Vspace')$, and the last integral is a duality pairing. Under \eqref{eq:weak-gaussian-bound}, Lemma~\ref{trace:gaussian} gives the incoming $L^2$ bound and the existence of the integral \eqref{eq:diffuse-flux}.
\end{proposition}

\begin{proof}
We apply \cite[Theorem~2.8(3), (2.25)]{Carrapatoso.Mischler2024} to \eqref{trace:fixed-equation} to obtain the local outgoing $L^2$ bound and permit $\beta(\xi)=\frac12\xi^2$ in \eqref{trace:beta-Green}. Lemma~\ref{trace:CM} verifies its coefficient and source assumptions. The remaining hypotheses are
\begin{align}
f(0)\in L^2_{\mathrm{loc}}(\overline\Om\times\R^3),\qquad \trm f\in L^2_{\mathrm{loc}}\bigl((0,T)\times\Sigma_-;\dd t\dd\gamma\bigr), \label{trace:local-incoming-hypotheses}
\end{align}
which follow from the assumed continuity in $L^2$ and the incoming-trace bound.

For a smooth $\varphi$ on $[0,T]\times\overline\Om\times\R^3$ with compact velocity support, \eqref{trace:local-H1-bound} gives $f\varphi\in L^2(0,T;\Vspace)$. Its weak derivative is $\varphi\nabla_vf+f\nabla_v\varphi$, so \eqref{trace:diffusion-pairing} gives
\begin{align}
\ip[\big]{S-\Lop f-\divv_v(\sig\nabla_vf)}{f\varphi}=\ip*{S-\Lop f}{f\varphi}+\int_{\phase}\varphi\,\sig\nabla_vf\cdot\nabla_vf\dd x\dd v +\int_{\phase}f\,\sig\nabla_vf\cdot\nabla_v\varphi\dd x\dd v. \label{trace:cutoff-cancellation}
\end{align}
Substituting this into \eqref{trace:beta-Green} with $\beta(\xi)=\frac12\xi^2$ cancels the diffusion terms and gives
\begin{align}
&\left[\frac12\int_{\phase}\bigl(f(r)\bigr)^2\varphi(r)\dd x\dd v\right]_{r=s}^{r=t}+\frac12\int_s^t\int_{\partial\Om\times\R^3}(\trace f)^2\varphi(n\cdot v)\dd S_x\dd v\dd r \label{trace:local-quadratic-Green}\\
&\qquad=\frac12\int_s^t\int_{\phase}f^2\Y\varphi\dd x\dd v\dd r+\int_s^t\ip*{S-\Lop f}{f\varphi}\dd r.\notag
\end{align}
The time terms in \eqref{trace:beta-Green} are defined by the assumed continuity of $f$ in $L^2$, so \eqref{trace:local-quadratic-Green} holds for every $0\leq s\leq t\leq T$. Equation \eqref{trace:local-quadratic-Green} is linear in $\varphi$ and therefore also holds for signed test functions.

Choose a smooth radial function $\chi$ that equals one on $\{\abs{v}\leq1\}$, vanishes on $\{\abs{v}\geq2\}$, and is nonincreasing in $\abs{v}$, with $0\leq\chi\leq1$. Set $\chi_R(v):=\chi\bigl(\frac{v}{R}\bigr)$ for $R\geq1$. Then $\chi_R\uparrow1$ and $\abs{\nabla_v\chi_R}\lesssim\frac1R$. Taking $\varphi=\chi_R^2$ in \eqref{trace:local-quadratic-Green} gives
\begin{align}
&\frac12\norm[\big]{\chi_Rf(t)}_{L^2}^2+\frac12\int_s^t\int_{\Sigmap}\chi_R^2\abs{\trp f}^2\dd\gamma\dd r \label{trace:cutoff-Green}\\
&\qquad=\frac12\norm[\big]{\chi_Rf(s)}_{L^2}^2+\frac12\int_s^t\int_{\Sigmam}\chi_R^2\abs{\trm f}^2\dd\gamma\dd r+\int_s^t\ip*{\Y f}{\chi_R^2f}\dd r,\notag
\end{align}
where $\Y\chi_R^2=0$ and $\Y f=S-\Lop f$.

To pass to the limit in the pairing, apply the cutoff estimate \eqref{cvel:W-density-cutoff} with $\chi^2$ and integrate in $(t,x)$:
\begin{align}
\tnorm{(1-\chi_R^2)f}_{L^2(0,T;\Vspace)}^2\lesssim\int_0^T\int_{\Om}\int_{\abs{v}\geq R}\left(\sig\nabla_vf\cdot\nabla_vf+\la^{\gamma+2}\abs{f}^2\right)\dd v\dd x\dd r\to0. \label{trace:cutoff-V-convergence}
\end{align}
The integral tends to zero because $f\in L^2(0,T;\Vspace)$. Consequently,
\begin{align}
\abs[\Big]{\int_s^t\ip*{\Y f}{\chi_R^2f-f}\dd r}\leq\tnorm{\Y f}_{L^2(0,T;\Vspace')}\tnorm{(1-\chi_R^2)f}_{L^2(0,T;\Vspace)}\to0. \label{trace:cutoff-dual-limit}
\end{align}
In \eqref{trace:cutoff-Green}, $\chi_Rf(s)\to f(s)$ and $\chi_Rf(t)\to f(t)$ in $L^2$. The incoming integral converges by the assumed $L^2$ bound. Monotone convergence applies to the outgoing term. The limits of the other terms are finite, so the outgoing trace belongs to $L^2(\dd t\dd\gamma)$ and setting $\trm f=g$ proves (i).

For (ii), apply (i) with $g=\trm f$, using \eqref{trace:diffuse-energy-assumption}. Both traces are now in $L^2(\dd t\dd\gamma)$, so integrating \eqref{trace:projection-identity} over $(s,t)\times\partial\Om$ and substituting into \eqref{trace:full-Green} gives \eqref{trace:energy-Green}. The final trace bound follows from \eqref{trace:Gaussian-full-traces} and \eqref{trace:scalar-bound}.
\end{proof}


\section{Global Existence of Weak Solutions}\label{sec:energy}

We will first construct local weak solutions under in-flow and then extend them using estimates uniform in time. We will treat diffuse-reflection separately in Section~\ref{cex:section}, where we will verify that the global solution supplied by \cite[Theorem~1.1]{Carrapatoso.Mischler2024} satisfies Definition~\ref{def:energy} and the estimates needed below.

Throughout the in-flow construction, fix $0<\alpha<\frac{1}{4}$ and $\beta:=\frac{1}{2}\left(\alpha+\frac{1}{4}\right)$. For the linear problem, assume
\begin{align}
\tnorm{\psi}_{L^\infty(0,T;L^{\infty,\alpha})}\leq\delta,\qquad G:=\m+\mh\psi. \label{ci:map-coefficient}
\end{align}
For $\delta$ sufficiently small, Lemma~\ref{cvel:coefficients} gives $\frac{1}{2}\sig\leq a*G\leq\frac{3}{2}\sig$ and the coefficient bounds used below. We will denote the solution of \eqref{ci:frozen-problem} by $F_\psi$ and set $f_\psi:=\mih(F_\psi-\m)$. By $Q(\m,\m)=0$ and \eqref{eq:linearized-operators}, its perturbation equation is
\begin{align}
\Y f_\psi+\Lop f_\psi=\Ga(\psi,f_\psi)-\Ga(f_\psi,\mh)+\Ga(\psi,\mh),\qquad f_\psi(0)=f_0,\qquad \trm f_\psi=g. \label{ci:map-problem}
\end{align}
At a fixed point $f_\psi=\psi$, the last two terms cancel, giving $\Y f+\Lop f=\Ga(f,f)$. Throughout this section, $F=\m+\mh f$ denotes the full density.

\para{Local existence}
Section~\ref{sec:energy-linear} will construct $F_\psi\geq0$ by regularization, with Gaussian and finite-time energy bounds for $f_\psi$. In Proposition~\ref{ci:local-existence}, the H\"older estimate of Lemma~\ref{trace:interior-holder} will give local equicontinuity; Arzel\`a--Ascoli and the Gaussian tails will then give compactness in spacetime $L^2$. An energy estimate for $f_{\psi_j}-f_\psi$ will give continuity of the map. Schauder's theorem will yield a local solution on an interval uniform for data in a fixed small Gaussian ball.

\para{Uniform estimates and continuation}
Sections~\ref{sec:energy-decay}--\ref{sec:energy-Gaussian} will prove the uniform energy and Gaussian estimates \eqref{ci:forced-energy-decay}--\eqref{ci:integrated-energy} and \eqref{ci:uniform-Gaussian}. Proposition~\ref{ci:global-energy} will then restart the local construction on successive intervals. Strong regularity and uniqueness for the data of Theorem~\ref{main:inflow} will follow in Sections~\ref{cb:section} and \ref{sec:uniqueness}.

\subsection{Linear In-flow Equation and Finite-Time Estimates}\label{sec:energy-linear}

For $G$ in \eqref{ci:map-coefficient}, we will solve
\begin{align}
\Y F=Q(G,F),\qquad F(0)=\m+\mh f_0,\qquad \trm F=\m+\mh g. \label{ci:frozen-problem}
\end{align}
Lemma~\ref{ci:frozen-inflow-construction} constructs the weak solution $f=f_\psi$ used in Section~\ref{sec:energy-construction}.

We will construct $F$ by first solving, for $0<\varepsilon\leq1$,
\begin{align}
\Y F_\varepsilon&=Q(G,F_\varepsilon)+\varepsilon\left(\Delta_vF_\varepsilon+\divv_v(vF_\varepsilon)\right), \qquad F_\varepsilon(0)=\m+\mh f_0,\qquad \trm F_\varepsilon=\m+\mh g. \label{ci:regularized-problem}
\end{align}
The additional operator gives $a*G+\varepsilon I\geq\varepsilon I$ and preserves the Maxwellian:
\begin{align}
\Delta_v\m+\divv_v(v\m)=0. \label{eq:I:display:216}
\end{align}
For $f_\varepsilon:=\mih(F_\varepsilon-\m)$, the corresponding perturbation equation is
\begin{align}
\Y f_\varepsilon+\Lop f_\varepsilon&=\Ga(\psi,f_\varepsilon)-\Ga(f_\varepsilon,\mh)+S+\varepsilon\left(\Delta_v f_\varepsilon+\left(\frac{3}{2}-\frac{\abs{v}^2}{4}\right)f_\varepsilon\right). \label{ci:regularized-normalized}
\end{align}
For \eqref{ci:regularized-problem}, the source is
\begin{align}
S:=\Ga(\psi,\mh). \label{ci:regularized-source}
\end{align}
For $f\in C_c^\infty(\phase)$, integration by parts gives
\begin{align}
-\varepsilon\ip[\bigg]{\Delta_vf+\left(\frac{3}{2}-\frac{\abs{v}^2}{4}\right)f}{f}=\varepsilon\norm{\nabla_vf}_{L^2}^2+\frac{\varepsilon}{4}\norm{\abs{v}f}_{L^2}^2-\frac{3\varepsilon}{2}\norm{f}_{L^2}^2. \label{ci:regularizer-energy}
\end{align}
Thus the additional terms supply the $\varepsilon$-weighted bounds for $\nabla_vf_\varepsilon$ and $\abs{v}f_\varepsilon$ in \eqref{ci:regularized-energy}.

We will obtain estimates uniform in $\varepsilon$ and pass to $\varepsilon=0$ by weak convergence with $\psi$ fixed. Equation~\eqref{ci:regularized-normalized} is understood through
\begin{align}
&-\int f_\varepsilon\Y\varphi+\int_0^T\ip[\Big]{\Lop f_\varepsilon-\Ga(\psi,f_\varepsilon)+\Ga(f_\varepsilon,\mh)}{\varphi}\dd t \label{ci:linear-weak-form}\\
&\qquad=\int S\varphi-\varepsilon\int\nabla_v f_\varepsilon\cdot\nabla_v\varphi+\varepsilon\int\left(\frac{3}{2}-\frac{\abs{v}^2}{4}\right)f_\varepsilon\varphi.\notag
\end{align}
Here the unmarked integrals are over $(0,T)\times\phase$, and $\varphi\in C_c^\infty\bigl((0,T)\times\phase\bigr)$. At $\varepsilon=0$ with the source \eqref{ci:regularized-source}, this is the weak formulation of \eqref{ci:map-problem}, with the subscript $\psi$ suppressed. Initial values and incoming traces are imposed separately in the spaces specified below.

For the limiting equation, we seek $f$ in the class
\begin{align}
f\in C\bigl([0,T];L^2\bigr)\cap L^2(0,T;\Vspace),\qquad \Y f\in L^2(0,T;\Vspace'). \label{energy:positive-regularity}
\end{align}

The linear construction uses the following estimates for a radial Gaussian and for $\Ga(f,\mh)$. As in Section~\ref{cvel:section}, $a_G=a*G$, $b_G=b*G$ and $c_G=c*G$, with convolution in $v$.

\begin{lemma}\label{ci:radial-estimate}
For $\delta>0$ sufficiently small in \eqref{ci:map-coefficient}, there are $C,c_*,R>0$ such that, for $\alpha+\frac{1}{4}\leq\zeta\leq\beta+\frac{1}{4}$ and $W_\zeta(v):=\ue^{-\zeta\abs{v}^2}$,
\begin{align}
Q(G,W_\zeta)\leq CW_\zeta, \label{ci:radial-upper}
\end{align}
and
\begin{align}
\frac{Q(G,W_\zeta)}{W_\zeta}\leq-c_*\la^{\gamma+2}\qquad\text{for }\abs{v}\geq R. \label{ci:radial-barrier}
\end{align}
The constants are uniform in $(t,x)$ and $\zeta$ and do not depend on derivatives of $\psi$.
\end{lemma}

\begin{proof}
Since $\nabla_vW_\zeta=-2\zeta vW_\zeta$ and $D_v^2W_\zeta=(4\zeta^2v\otimes v-2\zeta I)W_\zeta$, equation~\eqref{cvel:Q} gives
\begin{align}
\frac{Q(G,W_\zeta)}{W_\zeta}=4\zeta^2v\cdot a_Gv-2\zeta\tr a_G-c_G. \label{ci:radial-quotient}
\end{align}
Apply Lemma~\ref{cvel:coefficients} with $F=G$ and $f=\psi$. The hypothesis \eqref{cvel:coefficient-bounded} follows from \eqref{ci:map-coefficient}. Equations~\eqref{eq:I:display:195}, \eqref{ci:radial-moment-expansion} and \eqref{cvel:coefficient-moments}, together with $\abs{c_G}\lesssim\la^\gamma$ from \eqref{cvel:coefficient-c}, yield
\begin{align}
\frac{Q(G,W_\zeta)}{W_\zeta}=\bigl(4\zeta(2\zeta-1)+O(\delta)\bigr)\abs{v}^{\gamma+2}+O\bigl(\abs{v}^{\gamma+1}\bigr)\qquad\text{as }\abs{v}\to\infty. \label{eq:I:display:196}
\end{align}
All error constants are uniform in $(t,x)$ and $\zeta$. Since $0<\alpha<\beta<\frac{1}{4}$, the coefficient $4\zeta(2\zeta-1)$ has a strictly negative upper bound on $\left[\alpha+\frac{1}{4},\beta+\frac{1}{4}\right]$. Choose $\delta$ sufficiently small and then $R\geq2$ sufficiently large. This proves \eqref{ci:radial-barrier}. On $B_R$, the bounds \eqref{cvel:rough-coefficients}--\eqref{cvel:coefficient-c} and \eqref{ci:radial-quotient} bound $Q(G,W_\zeta)/W_\zeta$ uniformly. Together with its negative upper bound outside $B_R$, this gives \eqref{ci:radial-upper} on $\R^3$.
\end{proof}

\begin{lemma}\label{ci:Maxwellian-argument}
For $f\in L^2$ and $f\in L^{\infty,\alpha}$, respectively, the collision term with second argument $\mh$ satisfies
\begin{align}
\gnorm[\big]{2}{\beta}{\Ga(f,\mh)}\lesssim\norm{f}_{L^2},\qquad \gnorm[\big]{\infty}{\beta}{\Ga(f,\mh)}\lesssim\gnorm{\infty}{\alpha}{f}. \label{ci:Maxwellian-estimates}
\end{align}
\end{lemma}

\begin{proof}
Equations~\eqref{cvel:K} and \eqref{cvel:K-bilinear-identity} give, in distributions, for every $f\in L^2$,
\begin{align}
\Ga(f,\mh)=\mh(v)\bigl(v\otimes v-I\bigr):\bigl(a*(\mh f)\bigr)(v)-\mh(v)\bigl(c*(\mh f)\bigr)(v). \label{eq:I:display:204}
\end{align}
Since $2\gamma+4\geq-2>-3$, the Cauchy--Schwarz inequality and \eqref{cvel:Gaussian-convolution-bound}, applied with $s=2\gamma+4$, give
\begin{align}
\abs[\Big]{\bigl(a*(\mh f)\bigr)(x,v)}\lesssim\left(\int_{\R^3}\abs{v-w}^{2\gamma+4}\m(w)\dd w\right)^{\frac{1}{2}}\vnorm[\big]{f(x,\cdot)}_{L^2_v} \lesssim\la^{\gamma+2}\vnorm[\big]{f(x,\cdot)}_{L^2_v}. \label{eq:I:display:205}
\end{align}
Since $\beta<\frac{1}{4}$, the $L^{2,\beta}$ norm of the first term in \eqref{eq:I:display:204} is bounded by $C\norm{f}_{L^2}$. For $\gamma>-3$, split $c$ by the indicators $\ind_{\{\abs{v}<1\}}$ and $\ind_{\{\abs{v}\geq1\}}$. The first part belongs to $L^1$, and Young's inequality gives
\begin{align}
\vnorm[\Big]{\ue^{\beta\abs{v}^2}\mh\bigl((c\ind_{\{\abs{v}<1\}})*(\mh f)\bigr)}_{L^2_v}\lesssim\vnorm{c\ind_{\{\abs{v}<1\}}}_{L^1_v}\vnorm{f}_{L^2_v}. \label{eq:I:display:206}
\end{align}
For the far convolution, the Cauchy--Schwarz inequality gives the bound $C\la^{\max(\gamma,0)}\vnorm[\big]{f(x,\cdot)}_{L^2_v}$. The factor $\ue^{\beta\abs{v}^2}\mh$ makes its $L^2_v$ norm finite. At $\gamma=-3$, use $c=-8\pi\delta_0$ to write the last term in \eqref{eq:I:display:204} as $8\pi\m f$, which is estimated directly. Integrating in $x$ proves the $L^{2,\beta}$ estimate.

For $f\in L^{\infty,\alpha}$, apply \eqref{cvel:scalar-K-bound} by homogeneity at almost every $x$. Since $\beta<\frac{1}{4}$, the function $\la^5\ue^{-\left(\frac{1}{4}-\beta\right)\abs{v}^2}$ is bounded. Together with \eqref{cvel:K-bilinear-identity}, this proves the $L^{\infty,\beta}$ estimate in \eqref{ci:Maxwellian-estimates}.
\end{proof}

In the remainder of this subsection, take $\delta$ in \eqref{ci:map-coefficient} small enough for Lemmas~\ref{cvel:coefficients} and \ref{ci:radial-estimate}. Constants may also depend on the fixed upper bound for $\delta$, with dependence on $T$ indicated by a subscript. They are independent of derivatives of $\psi$.

\begin{lemma}[The regularized in-flow problem]\label{ci:frozen-approximation}
Let $f_0\in L^2$ and $g\in L^2\bigl((0,T)\times\Sigma_-;\dd t\dd\gamma\bigr)$. For each $0<\varepsilon\leq1$, problem \eqref{ci:regularized-problem} has a unique weak solution $F_\varepsilon=\m+\mh f_\varepsilon$ in the class
\begin{align}
f_\varepsilon\in C\bigl([0,T];L^2\bigr)\cap L^2(0,T;\Vspace),\qquad \nabla_vf_\varepsilon,\abs{v}f_\varepsilon\in L^2\bigl((0,T)\times\phase\bigr). \label{ci:regularized-class}
\end{align}
The equation holds in the sense of \eqref{ci:linear-weak-form}, with the source \eqref{ci:regularized-source}. The initial condition holds in $L^2$, and the incoming condition holds $\dd t\dd\gamma$-almost everywhere, with
\begin{align}
f_\varepsilon(0)=f_0,\qquad \trm f_\varepsilon=g,\qquad \trpm f_\varepsilon\in L^2\bigl(0,T;L^2_{\Sigma_\pm}\bigr). \label{ci:regularized-data-traces}
\end{align}
The solution satisfies
\begin{align}
&\tnorm{f_\varepsilon}_{L^\infty(0,T;L^2)}^2+\int_0^T\left(\norm{f_\varepsilon}_{\Vspace}^2+\norm{\trp f_\varepsilon}_{L^2_{\Sigma_+}}^2+\varepsilon\norm{\nabla_vf_\varepsilon}_{L^2}^2+\varepsilon\norm{\abs{v}f_\varepsilon}_{L^2}^2\right)\dd t \label{ci:regularized-energy}\\
&\qquad\lesssim_T\left(\norm{f_0}_{L^2}^2+\int_0^T\left(\norm{g}_{L^2_{\Sigma_-}}^2+\norm{\psi}_{L^2}^2\right)\dd t\right).\notag
\end{align}
The estimate is uniform in $\varepsilon$. Nonnegative initial and incoming densities give $F_\varepsilon\geq0$.
\end{lemma}

\begin{proof}
We apply \cite[proof of Theorem~2.11, Step~1]{Carrapatoso.Mischler2024} with weight $\omega=\mih$ to construct the weak solution with \eqref{ci:regularized-class}--\eqref{ci:regularized-data-traces}. We first verify the coefficient and data conditions of that construction, then prove \eqref{ci:regularized-energy} uniformly in $\varepsilon$.

\proofstep{Verification of the cited hypotheses} In the nondivergence notation of \cite[(2.29)]{Carrapatoso.Mischler2024}, equation \eqref{ci:regularized-problem} has diffusion matrix $a_G+\varepsilon I$, first-order coefficient $\varepsilon v$ and zeroth-order coefficient $3\varepsilon-c_G$, since
\begin{align}
\Y F_\varepsilon=(a_G+\varepsilon I):D_v^2F_\varepsilon+\varepsilon v\cdot\nabla_vF_\varepsilon+(3\varepsilon-c_G)F_\varepsilon. \label{ci:regularized-density-expanded}
\end{align}
Lemma~\ref{cvel:coefficients} gives, for every $R>0$,
\begin{align}
\tnorm{\abs{a_G+\varepsilon I}+\varepsilon\abs{v}+\abs{3\varepsilon-c_G}}_{L^\infty((0,T)\times\Om\times B_R)}\lesssim_R 1, \label{ci:regularized-coefficients-local}
\end{align}
and, for every $\xi\in\R^3$,
\begin{align}
\xi\cdot(a_G+\varepsilon I)\xi\geq\frac{1}{2}\xi\cdot\sig\xi+\varepsilon\abs{\xi}^2\geq\varepsilon\abs{\xi}^2. \label{ci:regularized-ellipticity}
\end{align}
The divergence formulation has first-order coefficient $\varepsilon v-b_G$, which is also locally bounded by Lemma~\ref{cvel:coefficients}. Thus \cite[(2.21) and (2.30)]{Carrapatoso.Mischler2024} hold for each fixed $\varepsilon>0$, with global ellipticity constant $\varepsilon$.

At $\omega=\mih$, the collision contribution to the coefficient defined in \cite[(2.18)]{Carrapatoso.Mischler2024} is
\begin{align}
\varpi_G:=\frac{1}{4}v\cdot a_Gv+\frac{1}{2}v\cdot b_G-\frac{1}{2}c_G. \label{eq:I:display:218}
\end{align}
The regularizing terms contribute $\varepsilon\left(\frac{3}{2}-\frac{\abs{v}^2}{4}\right)$ by \eqref{ci:regularizer-energy}.
Applying \eqref{ci:radial-moment-expansion}, \eqref{ci:weighted-drift-expansion} and \eqref{cvel:coefficient-moments} to $F=G$, and using \eqref{cvel:coefficient-c}, gives
\begin{align}
\varpi_G=\left(-\frac{1}{2}+O(\delta)\right)\abs{v}^{\gamma+2}+O\bigl(\abs{v}^{\gamma+1}\bigr)\qquad\text{as }\abs{v}\to\infty. \label{ci:weighted-potential-expansion}
\end{align}
The error constants are uniform in $(t,x)$ by \eqref{ci:map-coefficient}. Combining this expansion with \eqref{ci:regularized-coefficients-local}, and decreasing $\delta$, gives $c,C>0$ independent of $\varepsilon$ such that
\begin{align}
-C\left(\la^{\gamma+2}+\varepsilon\la^2\right)&\leq\varpi_G+\varepsilon\left(\frac{3}{2}-\frac{\abs{v}^2}{4}\right)\leq C-c\left(\la^{\gamma+2}+\varepsilon\la^2\right). \label{eq:I:display:219}
\end{align}
This is precisely \cite[(2.31)]{Carrapatoso.Mischler2024}, with its function $\varsigma(v)=c\bigl(\la^{\gamma+2}+\varepsilon\la^2\bigr)$, $\lambda_0=C$ and $\lambda_1=C/c$.

For \cite[(2.38)]{Carrapatoso.Mischler2024}, it suffices in our notation to verify
\begin{align}
\tnorm[\Big]{\mh\left(\abs{a_G+\varepsilon I}+\abs{\nabla_va_G}+\abs{D_v^2a_G}+\varepsilon\abs{v}+3\varepsilon+\abs{3\varepsilon-c_G}\right)}_{L^2((0,T)\times\phase)}<\infty. \label{ci:regularized-weighted-coefficients}
\end{align}
The other requirement there is $\la\omega^{-1}=\la\mh\in L_v^2$, which follows directly from the Gaussian decay of $\mh$. For $\gamma>-3$, direct differentiation of $a$ gives locally integrable kernels of sizes $\abs{v}^{\gamma+2}$, $\abs{v}^{\gamma+1}$ and $\abs{v}^{\gamma}$ for derivatives of orders zero, one and two. Convolution with the Gaussian bound for $G$ gives polynomial growth. Multiplication by $\mh$ therefore yields \eqref{ci:regularized-weighted-coefficients}.

At $\gamma=-3$, the same argument applies to $a_G$ and $\nabla_va_G$, while $c_G=-8\pi G$. For the second derivatives, $a_{ij}=\partial_i\partial_j\abs{v}$ and $\Delta_v^2\abs{v}=-8\pi\delta_0$ in distributions. Thus $\partial_k\partial_l(a_{ij}*G)$ is a Fourier multiplier applied to $G$, whose symbol is $-8\pi\frac{\xi_i\xi_j\xi_k\xi_l}{\abs{\xi}^4}$ for $\xi\ne0$. Plancherel's theorem gives
\begin{align}
\vnorm{D_v^2a_G}_{L_v^2}\lesssim\vnorm{G}_{L_v^2}\lesssim1\qquad\text{for almost every }(t,x). \label{ci:Coulomb-second-coefficients}
\end{align}
Since $\mh$ is bounded, this also proves \eqref{ci:regularized-weighted-coefficients} in the Coulomb case. The drift and its divergence are $\varepsilon v$ and $3\varepsilon$, so their contributions are integrable with the same Gaussian factor.

Condition~\cite[(2.36)]{Carrapatoso.Mischler2024}, used there for the $L^1$ estimate, also holds: the nonlocal operator in that condition is zero for \eqref{ci:regularized-density-expanded}, and the required coefficient upper bound is zero because
\begin{align}
\sum_{i,j}\partial_i\partial_j\bigl((a_G)_{ij}+\varepsilon\delta_{ij}\bigr)-\divv_v(\varepsilon v)+(3\varepsilon-c_G)=c_G-3\varepsilon+3\varepsilon-c_G=0. \label{ci:regularized-mass-coefficient}
\end{align}
Finally, the data satisfy
\begin{align}
\norm[\Big]{\mih(\m+\mh f_0)}_{L^2}<\infty,\qquad\tnorm[\Big]{\mih(\m+\mh g)}_{L^2((0,T)\times\Sigma_-;\dd t\dd\gamma)}<\infty. \label{ci:regularized-weighted-data}
\end{align}
These follow from $\mih(\m+\mh f_0)=\mh+f_0$ and $\mih(\m+\mh g)=\mh+g$. The cited in-flow construction now gives $\mih F_\varepsilon\in C([0,T];L^2)$ and $\mih\trpm F_\varepsilon\in L^2((0,T)\times\Sigma_\pm;\dd t\dd\gamma)$. The energy estimate at the end of that step, with the derivative norm defined in \cite[(2.35)]{Carrapatoso.Mischler2024}, together with \eqref{ci:regularized-ellipticity} and the specified $\varsigma$, gives \eqref{ci:regularized-class} after subtracting $\mh$. We prove $F_\varepsilon\geq0$ for nonnegative initial and incoming densities below.

\proofstep{The estimate for the perturbation} Fix $\varepsilon>0$ and set $f_R:=\chi_Rf_\varepsilon$, with the velocity cutoff from \eqref{trace:cutoff-Green}. The cutoff places the additional $\varepsilon$-terms in $L^2(0,T;\Vspace')$, allowing us to apply Proposition~\ref{trace:green}. Testing \eqref{ci:linear-weak-form} with $\chi_R\varphi$ gives
\begin{align}
\Y f_R=\chi_R\left[S-\Lop f_\varepsilon+\Ga(\psi,f_\varepsilon)-\Ga(f_\varepsilon,\mh)+\varepsilon\left(\Delta_v+\frac{3}{2}-\frac{\abs{v}^2}{4}\right)f_\varepsilon\right]. \label{ci:localized-regularized-equation}
\end{align}
The collision terms belong to $L^2(0,T;\Vspace')$ by Corollary~\ref{cvel:spatial}, \eqref{ci:map-coefficient} and \eqref{ci:regularized-class}. Local ellipticity \eqref{trace:local-H1-bound} gives
\begin{align}
\norm{\chi_R\Delta_vf_\varepsilon}_{\Vspace'}+\norm[\bigg]{\chi_R\bigg(\frac{3}{2}-\frac{\abs{v}^2}{4}\bigg)f_\varepsilon}_{\Vspace'}\lesssim_R\left(\norm{\nabla_vf_\varepsilon}_{L^2}+\norm{f_\varepsilon}_{L^2}\right), \label{ci:localized-diffusion-bound}
\end{align}
where the diffusion term acts on $\varphi\in\Vspace$ as $-\int\nabla_vf_\varepsilon\cdot\nabla_v(\chi_R\varphi)$. Hence $f_R\in C([0,T];L^2)\cap L^2(0,T;\Vspace)$ and $\Y f_R\in L^2(0,T;\Vspace')$. Its initial value is $\chi_Rf_0$, and testing the Green identity with $\chi_R\varphi$ gives $\trpm f_R=\chi_R\trpm f_\varepsilon$ by Lemma~\ref{trace:CM}. Proposition~\ref{trace:green}(i) applies to $f_R$, with source $\Lop f_R+\Y f_R$.

Substitute \eqref{ci:localized-regularized-equation} into \eqref{trace:full-Green}. Integration by parts in the diffusion term gives
\begin{align}
\varepsilon\ip[\big]{\chi_R\Delta_vf_\varepsilon}{f_R}=-\varepsilon\norm{\chi_R\nabla_vf_\varepsilon}_{L^2}^2-2\varepsilon\int_{\phase}\chi_Rf_\varepsilon\nabla_vf_\varepsilon\cdot\nabla_v\chi_R\dd x\dd v. \label{ci:localized-diffusion-identity}
\end{align}
The error vanishes as $R\to\infty$, since
\begin{align}
2\varepsilon\abs[\bigg]{\int_0^t\int_{\phase}\chi_Rf_\varepsilon\nabla_vf_\varepsilon\cdot\nabla_v\chi_R\dd x\dd v\dd r}\lesssim\frac{\varepsilon}{R}\tnorm{f_\varepsilon}_{L^2((0,T)\times\phase)}\tnorm{\nabla_vf_\varepsilon}_{L^2((0,T)\times\phase)}\to0. \label{ci:localized-diffusion-limit}
\end{align}
The convergence $\chi_R^2f_\varepsilon\to f_\varepsilon$ in $L^2(0,T;\Vspace)$ follows from \eqref{trace:cutoff-V-convergence} and passes all collision pairings to the limit. The bounds \eqref{ci:regularized-class}--\eqref{ci:regularized-data-traces} give convergence of the remaining terms. Thus, for every $t\in[0,T]$,
\begin{align}
&\frac{1}{2}\norm{f_\varepsilon(t)}_{L^2}^2+\int_0^t\left(\ip[\Big]{\Lop f_\varepsilon-\Ga(\psi,f_\varepsilon)}{f_\varepsilon}+\frac{1}{2}\norm{\trp f_\varepsilon}_{L^2_{\Sigma_+}}^2+\varepsilon\norm{\nabla_vf_\varepsilon}_{L^2}^2+\frac{\varepsilon}{4}\norm{\abs{v}f_\varepsilon}_{L^2}^2\right)\dd r \label{ci:regularized-quadratic}\\
&\qquad=\frac{1}{2}\norm{f_0}_{L^2}^2+\int_0^t\left(\frac{1}{2}\norm{g}_{L^2_{\Sigma_-}}^2+\frac{3\varepsilon}{2}\norm{f_\varepsilon}_{L^2}^2+\ip[\Big]{\Ga(\psi,\mh)-\Ga(f_\varepsilon,\mh)}{f_\varepsilon}\right)\dd r.\notag
\end{align}

By \eqref{ci:map-coefficient} and \eqref{eq:gaussian-energy}, $\tnorm{\psi}_{L^\infty(0,T;L_x^\infty L_v^2)}\lesssim\delta$. Apply \eqref{cvel:localized-garding}, integrated in $x$, and \eqref{cvel:spatial-dual} to the collision terms in \eqref{ci:regularized-quadratic}. Choose $\delta$ so that the resulting $C\delta\norm{f_\varepsilon}_{\Vspace}^2$ is absorbed by $\frac{c_0}{2}\norm{f_\varepsilon}_{\Vspace}^2$, where $c_0$ is from \eqref{cvel:localized-garding}. Lemma~\ref{ci:Maxwellian-argument} and the Cauchy--Schwarz inequality bound the last pairing by $C\bigl(\norm{\psi}_{L^2}^2+\norm{f_\varepsilon}_{L^2}^2\bigr)$. We obtain, in distributions in $t$,
\begin{align}
\partial_t\norm{f_\varepsilon}_{L^2}^2+c\norm{f_\varepsilon}_{\Vspace}^2+\norm{\trp f_\varepsilon}_{L^2_{\Sigma_+}}^2+\frac{\varepsilon}{4}\left(\norm{\nabla_v f_\varepsilon}_{L^2}^2+\norm{\abs{v}f_\varepsilon}_{L^2}^2\right)\leq C\left(\norm{f_\varepsilon}_{L^2}^2+\norm{\psi}_{L^2}^2\right)+\norm{g}_{L^2_{\Sigma_-}}^2. \label{ci:regularized-differential}
\end{align}
Integration and Gr\"onwall's inequality prove \eqref{ci:regularized-energy}.

\proofstep{Nonnegative densities} Set $q:=\mih F_\varepsilon=\mh+f_\varepsilon$. The preceding construction gives this $q$ the regularity \eqref{ci:regularized-class}, with $q(0)=\mh+f_0$ and $\trpm q=\mh+\trpm f_\varepsilon\in L^2(0,T;L^2_{\Sigma_\pm})$ by \eqref{ci:regularized-data-traces}. We derive the negative-part estimate with the initial and incoming terms retained, so that it also applies to the Gaussian supersolutions below and the diffuse solution in Section~\ref{cex:section}. The calculation applies for $0\leq\varepsilon\leq1$ whenever
\begin{align}
q\in C([0,T];L^2)\cap L^2(0,T;\Vspace),\qquad\trm q\in L^2(0,T;L^2_{\Sigma_-}), \label{ci:negative-part-class}
\end{align}
and
\begin{align}
\Y q-\mih Q(G,\mh q)-\varepsilon\left(\Delta_v+\frac{3}{2}-\frac{\abs{v}^2}{4}\right)q\geq0, \label{ci:normalized-supersolution}
\end{align}
where the left side belongs to $L^2((0,T)\times\Om\times B_M)$ for every $M>0$. When $\varepsilon>0$, we also require $\nabla_vq,\abs{v}q\in L^2((0,T)\times\phase)$, as in \eqref{ci:regularized-class}. At $\varepsilon=0$ all terms containing $\varepsilon$ are omitted. For $q=\mih F_\varepsilon$, equation \eqref{ci:regularized-problem} gives equality in \eqref{ci:normalized-supersolution}.

Define $q_-:=\max\{-q,0\}$. We use the renormalized identity \eqref{trace:beta-Green}. Under the initial and incoming $L^2$ assumptions \eqref{trace:local-incoming-hypotheses}, \cite[Theorem~2.8(3)]{Carrapatoso.Mischler2024} allows $\beta\in W^{2,\infty}_{\mathrm{loc}}(\R)$ with bounded $\beta''$, without requiring bounded $\beta'$. Thus it permits $\beta(s)=\frac{1}{2}(s_-)^2$. Here \eqref{ci:negative-part-class} supplies those data assumptions and local velocity $H^1$ regularity. The diffusion $a_G+\varepsilon I$ is locally uniformly elliptic, and it and the drift $-(b_G+a_Gv)$ are locally bounded by Lemma~\ref{cvel:coefficients}. The zeroth-order term and the residual in \eqref{ci:normalized-supersolution} belong locally to $L^2$, so the identity applies to this equation.

The cutoff $\chi_R^2$ from \eqref{trace:cutoff-Green} supplies the compact velocity support required in \eqref{trace:beta-Green}. Applying that identity amounts to testing \eqref{ci:normalized-supersolution} with $-\chi_R^2q_-$. The nonnegative residual contributes a nonpositive term. Integration by parts in the drift gives
\begin{align}
\frac{1}{4}v\cdot a_Gv-\frac{1}{2}\tr a_G-c_G+\frac{1}{2}\divv_v(b_G+a_Gv)=\frac{1}{4}v\cdot a_Gv+\frac{1}{2}b_G\cdot v-\frac{1}{2}c_G=\varpi_G. \label{ci:negative-part-coefficient}
\end{align}
In the next three displays, integrals are over $(0,t)\times\phase$, with any additional velocity restriction indicated. The terms differentiating $\chi_R$ are
\begin{align}
-2\int\chi_Rq_-(a_G+\varepsilon I)\nabla_vq_-\cdot\nabla_v\chi_R+\int\chi_Rq_-^2(b_G+a_Gv)\cdot\nabla_v\chi_R. \label{ci:negative-part-cutoff}
\end{align}
For the first integral, $a_G\leq C\sig$ and Cauchy--Schwarz bound its absolute value by
\begin{align}
\frac{C}{R}\left(\int_{\abs{v}>R}\left(\sig\nabla_vq_-\cdot\nabla_vq_-+\varepsilon\abs{\nabla_vq_-}^2\right)\right)^{\frac{1}{2}}\left(\int_{\abs{v}>R}\left(\la^{\gamma+2}+\varepsilon\right)q_-^2\right)^{\frac{1}{2}}\to0. \label{ci:negative-part-cutoff-diffusion}
\end{align}
For the second, $\abs{b_G+a_Gv}\lesssim\la^{\gamma+3}$ and $\abs{\nabla_v\chi_R}\lesssim\frac1R$ on $R<\abs{v}<2R$ give
\begin{align}
\abs{\int\chi_Rq_-^2(b_G+a_Gv)\cdot\nabla_v\chi_R}\lesssim\int_{\abs{v}>R}\la^{\gamma+2}q_-^2\to0. \label{ci:negative-part-cutoff-drift}
\end{align}
These limits follow from \eqref{ci:negative-part-class}, the additional integrability when $\varepsilon>0$, and the Sobolev chain rule for $q_-$. The signed zeroth-order term is integrable by \eqref{eq:I:display:219}, whose bound also holds at $\varepsilon=0$. The initial and incoming integrals converge by their $L^2$ bounds, and the outgoing integral converges monotonically. Letting $R\to\infty$ and using the upper bound in \eqref{eq:I:display:219} gives, for every $t\in[0,T]$,
\begin{align}
&\frac{1}{2}\norm{q_-(t)}_{L^2}^2+\frac{1}{2}\int_0^t\norm{(\trp q)_-}_{L^2_{\Sigma_+}}^2\dd s+\int_0^t\int_{\phase}(a_G+\varepsilon I)\nabla_vq_-\cdot\nabla_vq_-\dd x\dd v\dd s \label{ci:negative-part-energy}\\
&\qquad\leq\frac{1}{2}\norm{q_-(0)}_{L^2}^2+\frac{1}{2}\int_0^t\norm{(\trm q)_-}_{L^2_{\Sigma_-}}^2\dd s+C\int_0^t\norm{q_-}_{L^2}^2\dd s.\notag
\end{align}
For $q=\mih F_\varepsilon$, nonnegative density data give $q_-(0)=0$ and $(\trm q)_-=0$. Gr\"onwall's inequality in \eqref{ci:negative-part-energy} yields $q_-=0$, hence $F_\varepsilon\geq0$.
\end{proof}

\begin{lemma}[Linear in-flow with Gaussian data]\label{ci:frozen-inflow-construction}
Assume
\begin{align}
f_0\in L^{\infty,\alpha},\qquad g\in L^\infty\bigl(0,T;L^{\infty,\alpha}_{\Sigma_-}\bigr). \label{ci:linear-Gaussian-data}
\end{align}
Problem \eqref{ci:frozen-problem} has a unique weak solution $F=\m+\mh f$ in the class
\begin{align}
f\in C\bigl([0,T];L^2\bigr)\cap L^2(0,T;\Vspace)\cap L^\infty\bigl(0,T;L^{\infty,\alpha}\bigr). \label{ci:linear-Gaussian-class}
\end{align}
The equation holds in the sense of \eqref{ci:linear-weak-form} at $\varepsilon=0$, with the source \eqref{ci:regularized-source}. The initial condition holds in $L^2$, and the incoming condition holds $\dd t\dd\gamma$-almost everywhere, with
\begin{align}
f(0)=f_0,\qquad \trm f=g,\qquad \trpm f\in L^2\bigl(0,T;L^2_{\Sigma_\pm}\bigr). \label{ci:linear-data-traces}
\end{align}
The solution satisfies the energy identity \eqref{trace:full-Green} and, for every $0\leq t\leq T$,
\begin{align}
&\sup_{0\leq s\leq t}\norm{f(s)}_{L^2}^2+\int_0^t\left(\norm{f}_{\Vspace}^2+\norm{\trp f}_{L^2_{\Sigma_+}}^2\right)\dd s \label{ci:linear-energy} \leq C\ue^{Ct}\left(\norm{f_0}_{L^2}^2+\int_0^t\left(\norm{g}_{L^2_{\Sigma_-}}^2+\norm{\psi}_{L^2}^2\right)\dd s\right).
\end{align}
There is also the Gaussian bound
\begin{align}
\gnorm[\big]{\infty}{\alpha}{f(t)}\leq\ue^{Ct}\left(\gnorm{\infty}{\alpha}{f_0}+\tnorm{g}_{L^\infty(0,T;L^{\infty,\alpha}_{\Sigma_-})}\right)+\delta\left(\ue^{Ct}-1\right). \label{rev14:frozen-Gaussian-bound}
\end{align}
The constants in \eqref{ci:linear-energy} and \eqref{rev14:frozen-Gaussian-bound} are independent of $T$, the data and derivatives of $\psi$. Nonnegative initial and incoming densities give $F\geq0$.

\end{lemma}

\begin{proof}
Lemma~\ref{ci:frozen-approximation} supplies the regularized weak solutions and the estimate \eqref{ci:regularized-energy}. We will prove a Gaussian bound uniform in $\varepsilon$, pass to the limit as $\varepsilon\to0$, and identify the initial value and boundary traces of the limit.

\proofstep{The Gaussian estimate uniform in $\varepsilon$} The additional diffusion satisfies
\begin{align}
\frac{\left(\Delta_v+\divv_v(v\,\cdot)\right)\ue^{-\zeta\abs{v}^2}}{\ue^{-\zeta\abs{v}^2}}=(1-2\zeta)\left(3-2\zeta\abs{v}^2\right)\leq3(1-2\zeta),\qquad \alpha+\frac{1}{4}\leq\zeta\leq\beta+\frac{1}{4}. \label{ci:OU-comparison}
\end{align}
Estimate~\eqref{ci:radial-upper}, at $\zeta=\alpha+\frac{1}{4}$, and \eqref{ci:Maxwellian-estimates}, applied to $\psi$, give
\begin{align}
Q\left(G,\mh\ue^{-\alpha\abs{v}^2}\right)\leq C\mh\ue^{-\alpha\abs{v}^2},\qquad \abs[\big]{Q(G-\m,\m)}\leq C\delta\mh\ue^{-\alpha\abs{v}^2}. \label{ci:affine-Gaussian-source}
\end{align}
Take $C$ large enough for these inequalities and \eqref{ci:OU-comparison}. The scalar function
\begin{align}
\compamp(t):=\ue^{Ct}\left(\gnorm{\infty}{\alpha}{f_0}+\tnorm{g}_{L^\infty(0,T;L^{\infty,\alpha}_{\Sigma_-})}\right)+\delta\left(\ue^{Ct}-1\right) \label{ci:comparison-amplitude}
\end{align}
satisfies $\compamp'=C\compamp+C\delta$ and dominates both data norms. Set $u_\pm:=\compamp\ue^{-\alpha\abs{v}^2}\mp f_\varepsilon$. Equations \eqref{eq:I:display:216}, \eqref{ci:OU-comparison} and \eqref{ci:affine-Gaussian-source} give
\begin{align}
\Y(\mh u_\pm)-Q(G,\mh u_\pm)-\varepsilon\left(\Delta_v+\divv_v(v\,\cdot)\right)(\mh u_\pm)\geq0,\qquad u_\pm(0)\geq0,\qquad\trm u_\pm\geq0. \label{ci:Gaussian-supersolutions}
\end{align}
The functions $u_\pm$ have the regularity \eqref{ci:regularized-class}, since $\compamp\ue^{-\alpha\abs{v}^2}$ and its velocity derivatives have Gaussian decay. To apply \eqref{ci:negative-part-energy} with $q=u_\pm$, use $\Y\mh=0$ and
\begin{align}
\mih\left(\Delta_v+\divv_v(v\,\cdot)\right)(\mh u_\pm)=\left(\Delta_v+\frac{3}{2}-\frac{\abs{v}^2}{4}\right)u_\pm. \label{ci:Gaussian-normalization}
\end{align}
Dividing \eqref{ci:Gaussian-supersolutions} by $\mh$ therefore gives
\begin{align}
\Y u_\pm-\mih Q(G,\mh u_\pm)-\varepsilon\left(\Delta_v+\frac{3}{2}-\frac{\abs{v}^2}{4}\right)u_\pm\geq\bigl(\compamp'-C\compamp-C\delta\bigr)\ue^{-\alpha\abs{v}^2}=0. \label{ci:Gaussian-normalized-supersolutions}
\end{align}
The residual is locally $L^2$: the equation for $f_\varepsilon$ leaves only the Gaussian terms in these bounds and $\Ga(\psi,\mh)$. The initial and incoming negative parts vanish by \eqref{ci:Gaussian-supersolutions}. Thus \eqref{ci:negative-part-energy}, with $q=u_\pm$, gives
\begin{align}
\norm{(u_\pm)_-(t)}_{L^2}^2\lesssim\int_0^t\norm{(u_\pm)_-}_{L^2}^2\dd s. \label{ci:Gaussian-negative-part}
\end{align}
Gr\"onwall's inequality gives $u_\pm\geq0$, and therefore
\begin{align}
\gnorm[\big]{\infty}{\alpha}{f_\varepsilon(t)}\leq\compamp(t)\qquad(0\leq t\leq T). \label{ci:regularized-Gaussian-bound}
\end{align}
The bound first holds for almost every time and then at every time by the $L^2$ continuity in \eqref{ci:regularized-class}. Its constants are independent of $\varepsilon$.

\proofstep{Weak convergence and the equation} Fix $\psi$. The estimate \eqref{ci:regularized-energy} gives a sequence $\varepsilon_j\to0$ and limits $f,f_+$ such that
\begin{align}
f_{\varepsilon_j}\rightharpoonup f\text{ in }L^2(0,T;\Vspace),\qquad f_{\varepsilon_j}\overset{*}{\rightharpoonup}f\text{ in }L^\infty(0,T;L^2),\qquad \trp f_{\varepsilon_j}\rightharpoonup f_+\text{ in }L^2(0,T;L^2_{\Sigma_+}). \label{ci:linear-subsequence}
\end{align}
The bound \eqref{ci:regularized-Gaussian-bound} passes to the weak limit because $\{h:\abs{h(t,x,v)}\leq\compamp(t)\ue^{-\alpha\abs{v}^2}\}$ is a closed convex subset of spacetime $L^2$. In particular, $f$ satisfies \eqref{rev14:frozen-Gaussian-bound} almost everywhere in time.

Corollary~\ref{cvel:spatial} gives weak convergence of $\Lop f_{\varepsilon_j}$ and $\Ga(\psi,f_{\varepsilon_j})$ in $L^2(0,T;\Vspace')$. By \eqref{ci:Maxwellian-estimates} and \eqref{eq:V},
\begin{align}
\norm[\big]{\Ga(h,\mh)}_{\Vspace'}\lesssim\gnorm[\big]{2}{\beta}{\Ga(h,\mh)}\lesssim\norm{h}_{L^2}\qquad(h\in L^2). \label{ci:remainder-dual}
\end{align}
For the first inequality, pair against $\varphi\in\Vspace$ and use the boundedness of $\la^{-\frac{\gamma+2}{2}}\ue^{-\beta\abs{v}^2}$. Hence $\Ga(f_{\varepsilon_j},\mh)\rightharpoonup\Ga(f,\mh)$ in $L^2(0,T;\Vspace')$.

For $\varphi\in C_c^\infty((0,T)\times\phase)$, weak convergence in spacetime $L^2$ gives
\begin{align}
\ip[\big]{\Y(f_{\varepsilon_j}-f)}{\varphi}=-\int(f_{\varepsilon_j}-f)\Y\varphi\to0. \label{ci:weak-transport-limit}
\end{align}
The additional diffusion and zeroth-order term vanish since
\begin{align}
\abs{\varepsilon_j\int\nabla_vf_{\varepsilon_j}\cdot\nabla_v\varphi}\leq\varepsilon_j^{\frac{1}{2}}\tnorm{\varepsilon_j^{\frac{1}{2}}\nabla_vf_{\varepsilon_j}}_{L^2}\tnorm{\nabla_v\varphi}_{L^2}\to0, \label{ci:vanishing-diffusion}\\
\abs{\varepsilon_j\int\bigg(\frac{3}{2}-\frac{\abs{v}^2}{4}\bigg)f_{\varepsilon_j}\varphi}\lesssim_\varphi\varepsilon_j\tnorm{f_{\varepsilon_j}}_{L^2}\to0. \label{ci:vanishing-regularized-drift}
\end{align}
The integrals and $L^2$ norms here are over $(0,T)\times\phase$. Passing to the limit in \eqref{ci:linear-weak-form} yields
\begin{align}
\Y f+\Lop f=\Ga(\psi,f)-\Ga(f,\mh)+\Ga(\psi,\mh). \label{ci:frozen-normalized-source}
\end{align}
Corollary~\ref{cvel:spatial}, \eqref{ci:remainder-dual} and \eqref{ci:map-coefficient} put the right side in $L^2(0,T;\Vspace')$.

\proofstep{Initial value, boundary trace and continuity in $L^2$} The Gaussian bound already obtained and Lemma~\ref{trace:gaussian}(iii), applied to \eqref{ci:frozen-normalized-source}, give $f\in C\bigl([0,T];L^2\bigr)$.

For any $t\in[0,T]$, choose $t_j\to t$ from the full-measure set where \eqref{rev14:frozen-Gaussian-bound} holds. With $\compamp$ from \eqref{ci:comparison-amplitude},
\begin{align}
\norm[\Big]{\bigl(\abs[\big]{f(t)}-\compamp(t)\ue^{-\alpha\abs{v}^2}\bigr)_+}_{L^2}\leq\norm[\big]{f(t)-f(t_j)}_{L^2}+\abs[\big]{\compamp(t_j)-\compamp(t)}\norm[\Big]{\ue^{-\alpha\abs{v}^2}}_{L^2}\to0. \label{ci:Gaussian-limit-every-time}
\end{align}
Thus \eqref{rev14:frozen-Gaussian-bound} holds for every $t\in[0,T]$.

Using the outgoing limit $f_+$ in \eqref{ci:linear-subsequence}, pass to the limit in the linear Green identities for $f_{\varepsilon_j}$, using tests smooth up to $t=0$ and $\partial\Om$, with compact velocity support and zero value at $t=T$. The incoming and initial terms are fixed. Comparison with \eqref{trace:compact-Green}, first for tests compactly supported in time, identifies the outgoing limit with $\trp f$ and gives $\trm f=g$. Comparison with \eqref{trace:linear-Green} from Lemma~\ref{trace:CM}(iv), using tests supported away from the spatial boundary, then gives $f(0)=f_0$ in $L^2$. Proposition~\ref{trace:green} now yields \eqref{trace:full-Green}. Apply \eqref{cvel:localized-garding}, \eqref{cvel:spatial-dual} and \eqref{ci:remainder-dual} in this identity, as in the derivation of \eqref{ci:regularized-differential}, with $\varepsilon=0$. Gr\"onwall's inequality gives \eqref{ci:linear-energy}.

\proofstep{Uniqueness and positivity} Let $h:=\widetilde f-f$ for two weak solutions in the stated class with the same data. Then
\begin{align}
\Y h+\Lop h=\Ga(\psi,h)-\Ga(h,\mh),\qquad h(0)=0,\qquad\trm h=0. \label{ci:linear-difference}
\end{align}
Apply \eqref{cvel:localized-garding}, integrated in $x$, \eqref{cvel:spatial-dual} and Lemma~\ref{ci:Maxwellian-argument} in the energy identity \eqref{trace:full-Green}. The bound \eqref{ci:map-coefficient} absorbs the $\Ga(\psi,h)$ term for sufficiently small $\delta$, giving
\begin{align}
\norm[\big]{h(t)}_{L^2}^2+c\int_0^t\norm{h}_{\Vspace}^2\dd r+\int_0^t\norm{\trp h}_{L^2_{\Sigma_+}}^2\dd r
\lesssim\int_0^t\norm{h}_{L^2}^2\dd r. \label{ci:linear-uniqueness}
\end{align}
Gr\"onwall's inequality gives $h=0$, and the limit is independent of the subsequence. Nonnegative initial and incoming densities give $F_{\varepsilon_j}\geq0$ by Lemma~\ref{ci:frozen-approximation}. Since $F_{\varepsilon_j}\rightharpoonup F$ in spacetime $L^2$ and its nonnegative cone is weakly closed, $F\geq0$.
\end{proof}

\subsection{Local Nonlinear Existence under In-flow}\label{sec:energy-construction}

For each $\psi$ satisfying \eqref{ci:map-coefficient}, let $f_\psi$ be the unique weak solution supplied by Lemma~\ref{ci:frozen-inflow-construction}. This defines $\psi\mapsto f_\psi$ for fixed initial and incoming data. We will obtain a nonlinear weak solution as a fixed point of this map.

The following estimate gives interior compactness for the family $f_\psi$ when the Gaussian bound on $\psi$ is uniform. Its constants require no derivatives of $\psi$.

\begin{lemma}[Interior kinetic H\"older estimate]\label{trace:interior-holder}
Let $f$ be the solution supplied by Lemma~\ref{ci:frozen-inflow-construction}. Let $\mathcal U$ and $\mathcal U'$ be bounded open subsets of $(0,T)\times\phase$ with $\overline{\mathcal U'}\subset\mathcal U$. There is $\theta\in(0,1)$ such that
\begin{align}
\tnorm{f}_{\Ckin^\theta(\overline{\mathcal U'})}\lesssim_{\mathcal U,\mathcal U'}\left(\tnorm{f}_{L^\infty(\mathcal U)}+\delta\right). \label{trace:interior-bound}
\end{align}
Here $\Ckin^\theta$ is defined by \eqref{eq:kinetic-holder-difference}--\eqref{eq:kinetic-holder-norm}. The exponent may depend on $\mathcal U,\mathcal U'$, and the estimate is uniform for $\psi$ satisfying \eqref{ci:map-coefficient}.
\end{lemma}

\begin{proof}
We apply \cite[Lemma~2.24]{Ros-Oton.Weidner2025}. Lemma~\ref{ci:frozen-inflow-construction} and local ellipticity give $f\in L^\infty(\mathcal U)$ and $\nabla_vf\in L^2(\mathcal U)$. By \eqref{cvel:linear-bilinear-identity}, equation \eqref{ci:frozen-normalized-source} is $\Y f=\mih Q(G,\mh f)+\Ga(\psi,\mh)$. Using the divergence definition \eqref{eq:collision} and \eqref{cvel:convolved-identities}, valid here by the local coefficient bounds and $\nabla_vf\in L^2(\mathcal U)$, gives in distributions
\begin{align}
\Y f-\divv_v(a_G\nabla_vf)=-(b_G+a_Gv)\cdot\nabla_vf+\left(\frac{1}{4}v\cdot a_Gv-\frac{1}{2}\tr a_G-c_G\right)f+\Ga(\psi,\mh). \label{ci:interior-Landau-divergence}
\end{align}
On $\mathcal U$, Lemma~\ref{cvel:coefficients} and \eqref{ci:Maxwellian-estimates} give constants $\lambda,\Lambda,C>0$ such that
\begin{align}
\lambda\abs{\xi}^2\leq\xi\cdot a_G\xi\leq\Lambda\abs{\xi}^2,\qquad \abs{b_G+a_Gv}+\abs{\frac{1}{4}v\cdot a_Gv-\frac{1}{2}\tr a_G-c_G}\leq C,\qquad \abs{\Ga(\psi,\mh)}\leq C\delta. \label{ci:interior-ellipticity}
\end{align}
These bounds hold almost everywhere and uniformly in $\psi$. At $\gamma=-3$, the bound for $c_G$ uses $c_G=-8\pi G$. Thus \eqref{ci:interior-Landau-divergence} satisfies the ellipticity, coefficient and source assumptions of the cited lemma.

Cover $\overline{\mathcal U'}$ by finitely many smaller kinetic cylinders whose corresponding larger cylinders lie in $\mathcal U$. The cited estimate holds on each smaller cylinder. There is $\rho>0$ such that two points of $\overline{\mathcal U'}$ at kinetic distance less than $\rho$ lie in one of these cylinders. For the remaining pairs,
\begin{align}
\frac{\abs{f(Z)-f(Z')}}{\bigl(\dkin(Z,Z')\bigr)^\theta}\leq2\rho^{-\theta}\tnorm{f}_{L^\infty(\mathcal U)}\qquad\text{if }\dkin(Z,Z')\geq\rho. \label{ci:interior-separated-points}
\end{align}
Together with the estimates on the cylinders, this proves \eqref{trace:interior-bound}.
\end{proof}

\begin{proposition}[Local existence]\label{ci:local-existence}
There are $\delta_0,C_0>0$ with the following property. Let $T_0>0$, and let the data $f_0$ and $g$ satisfy $\m+\mh f_0\geq0$, $\m+\mh g\geq0$, and
\begin{align}
\delta:=\gnorm{\infty}{\alpha}{f_0}+\tnorm{g}_{L^\infty(0,T_0;L^{\infty,\alpha}_{\Sigma_-})}\leq\delta_0,\qquad T_*:=\min\bigg\{T_0,\frac{\ln\left(\frac{4}{3}\right)}{C_0}\bigg\}. \label{ci:local-data}
\end{align}
There is a weak solution $f$ on $[0,T_*]$ in the sense of Definition~\ref{def:energy}, with $f(0)=f_0$, $\trm f=g$ and $\m+\mh f\geq0$, satisfying
\begin{align}
\tnorm{f}_{L^\infty(0,T_*;L^{\infty,\alpha})}\leq2\delta. \label{ci:local-Gaussian}
\end{align}
With the time norms in the following estimate taken over $(0,T_*)$, it also satisfies
\begin{align}
\tnorm{f}_{L^\infty_tL^2}+\tnorm{f}_{L^2_t\Vspace}+\tnorm{\trp f}_{L^2_tL^2_{\Sigma_+}}\lesssim\left(\norm{f_0}_{L^2}+\tnorm{g}_{L^2_tL^2_{\Sigma_-}}\right). \label{ci:local-energy}
\end{align}
The constants are independent of $T_0$ and of the data satisfying $\delta\leq\delta_0$.
\end{proposition}

\begin{proof}
Following the compact fixed-point method of \cite[Section~8.1]{Carrapatoso.Mischler2024}, we will apply Schauder's theorem to $\psi\mapsto f_\psi$ in $L^2\bigl((0,T_*)\times\phase\bigr)$, verifying invariance, compactness and continuity in the following steps.

\proofstep{The Gaussian bound and the lifespan}
Choose $\delta_0$ small enough for Lemma~\ref{ci:frozen-inflow-construction} to apply to coefficients of Gaussian norm at most $2\delta_0$, and fix $C_0$ uniformly in \eqref{rev14:frozen-Gaussian-bound}. If $\delta=0$, take $f=0$. We henceforth assume $\delta>0$.

Consider
\begin{align}
\mathscr B:=\left\{\psi:\tnorm{\psi}_{L^\infty(0,T_*;L^{\infty,\alpha})}\leq2\delta,\quad \m+\mh\psi\geq0\text{ almost everywhere}\right\}. \label{ci:fixed-point-ball}
\end{align}
This set is convex, contains $0$, and is bounded in spacetime $L^2$. It is closed there because its pointwise bounds pass to an almost everywhere convergent subsequence. For $\psi\in\mathscr B$, Lemma~\ref{ci:frozen-inflow-construction} supplies $f_\psi$ with initial value $f_0$ and incoming trace $g$. Equation~\eqref{rev14:frozen-Gaussian-bound}, with coefficient bound $2\delta$, gives, for $0\leq t\leq T_*$,
\begin{align}
\gnorm[\big]{\infty}{\alpha}{f_\psi(t)}\leq\delta\ue^{C_0t}+2\delta\left(\ue^{C_0t}-1\right)=\delta\left(3\ue^{C_0t}-2\right)\leq2\delta. \label{ci:local-map-bound}
\end{align}
Lemma~\ref{ci:frozen-inflow-construction} also gives $\m+\mh f_\psi\geq0$. Thus $\psi\mapsto f_\psi$ maps $\mathscr B$ into itself.

\proofstep{The uniform energy estimate}
All time norms in this step are over $(0,T_*)$. Equation~\eqref{ci:linear-energy} gives
\begin{align}
\tnorm{f_\psi}_{L^\infty_tL^2}^2+\tnorm{f_\psi}_{L^2_t\Vspace}^2+\tnorm{\trp f_\psi}_{L^2_tL^2_{\Sigma_+}}^2\lesssim\left(\norm{f_0}_{L^2}^2+\tnorm{g}_{L^2_tL^2_{\Sigma_-}}^2+\tnorm{\psi}_{L^2_tL^2}^2\right)\lesssim\delta^2. \label{ci:fixed-point-energy}
\end{align}
The last inequality follows from the Gaussian data bounds, \eqref{ci:fixed-point-ball} and $T_*\leq C_0^{-1}\ln\left(\frac{4}{3}\right)$. The constant is independent of $T_0$.

\proofstep{Compactness}
Let $K$ be a compact subset of $(0,T_*)\times\Om\times\R^3$. Lemma~\ref{trace:interior-holder} and \eqref{ci:local-map-bound} give an exponent $\theta\in(0,1)$ and a bound
\begin{align}
\tnorm{f_\psi}_{\Ckin^\theta(K)}\lesssim_K\delta,\qquad \psi\in\mathscr B. \label{ci:fixed-point-Holder}
\end{align}
The H\"older estimate gives uniform boundedness and equicontinuity on $K$. Arzel\`a--Ascoli and a diagonal extraction therefore give, from any sequence in $\{f_\psi:\psi\in\mathscr B\}$, a subsequence converging uniformly on every compact subset of $(0,T_*)\times\phase$.

For $\Om_r:=\{x\in\Om:d(x)>r\}$, $R>1$ and $0<s<\frac{T_*}{2}$, \eqref{ci:local-map-bound} also gives
\begin{align}
\int_{((0,T_*)\times\phase)\setminus([s,T_*-s]\times\Om_r\times B_R)}\abs{f_\psi}^2\dd t\dd x\dd v \lesssim\delta^2\left(s+T_*\abs{\Om\setminus\Om_r}+T_*\int_{\abs{v}>R}\ue^{-2\alpha\abs{v}^2}\dd v\right). \label{ci:fixed-point-tails}
\end{align}
Choose $s,r$ small and $R$ large to make this bound uniformly small. The local uniform convergence then makes the subsequence Cauchy in spacetime $L^2$, proving relative compactness of the image.

\proofstep{Continuity of the map}
Suppose that $\psi_j\to\psi$ in $L^2((0,T_*)\times\phase)$, with $\psi_j,\psi\in\mathscr B$, and set $h_j:=f_{\psi_j}-f_\psi$. Both solutions are supplied by Lemma~\ref{ci:frozen-inflow-construction}. Their difference therefore satisfies
\begin{align}
\Y h_j+\Lop h_j=\Ga(\psi_j,h_j)-\Ga(h_j,\mh)+\Ga(\psi_j-\psi,f_\psi)+\Ga(\psi_j-\psi,\mh),\qquad h_j(0)=0,\qquad\trm h_j=0. \label{ci:map-difference-equation}
\end{align}
The Gaussian bounds give $\vnorm{\psi_j-\psi}_{L_v^2}^2\lesssim\delta^2$ almost everywhere in $(t,x)$. Estimate~\eqref{cvel:spatial-reversed} gives
\begin{align}
\tnorm[\big]{\Ga(\psi_j-\psi,f_\psi)}_{L^2_t\Vspace'}^2\lesssim\int_0^{T_*}\int_\Om\vnorm{\psi_j-\psi}_{L_v^2}^2\vnorm{f_\psi}_{\Vv}^2\dd x\dd t. \label{ci:map-forcing-bound}
\end{align}
For $M>0$, split this integral where $\vnorm{f_\psi}_{\Vv}^2\leq M$ and where it exceeds $M$. Its right side is bounded by
\begin{align}
CM\tnorm{\psi_j-\psi}_{L^2_{t,x,v}}^2+C\delta^2\int_0^{T_*}\int_{\{x:\vnorm{f_\psi(t,x,\cdot)}_{\Vv}^2>M\}}\vnorm{f_\psi}_{\Vv}^2\dd x\dd t. \label{ci:map-forcing-truncation}
\end{align}
The first term tends to zero for fixed $M$. Since $f_\psi\in L^2_t\Vspace$, the second tends to zero as $M\to\infty$, uniformly in $j$. Consequently, \eqref{ci:map-forcing-bound} tends to zero. Also, \eqref{ci:remainder-dual} gives $\Ga(\psi_j-\psi,\mh)\to0$ in $L^2_t\Vspace'$.

Apply \eqref{trace:full-Green} to $h_j$. The bound \eqref{cvel:localized-garding}, integrated in $x$, controls $\Lop h_j$. The estimate \eqref{cvel:spatial-dual} and the small Gaussian bound on $\psi_j$ absorb $\Ga(\psi_j,h_j)$ into the velocity dissipation. Estimate~\eqref{ci:Maxwellian-estimates} bounds the pairing with $\Ga(h_j,\mh)$ by $C\norm{h_j}_{L^2}^2$. Young's inequality and Gr\"onwall's inequality yield
\begin{align}
\tnorm{h_j}_{L^\infty_tL^2}^2+\tnorm{h_j}_{L^2_t\Vspace}^2+\tnorm{\trp h_j}_{L^2_tL^2_{\Sigma_+}}^2\lesssim\left(\tnorm[\big]{\Ga(\psi_j-\psi,f_\psi)}_{L^2_t\Vspace'}^2+\tnorm{\psi_j-\psi}_{L^2_{t,x,v}}^2\right)\to0. \label{ci:map-continuity-energy}
\end{align}
All time norms in this step are over $(0,T_*)$. This proves continuity of $\psi\mapsto f_\psi$ in spacetime $L^2$.

\proofstep{The fixed point and its energy bound}
Schauder's fixed-point theorem gives $f\in\mathscr B$ with $f_f=f$. Substitution into \eqref{ci:frozen-normalized-source} cancels the two terms with second argument $\mh$ and gives $\Y f+\Lop f=\Ga(f,f)$. Its traces and time continuity follow from the linear construction, and \eqref{ci:local-map-bound} gives \eqref{ci:local-Gaussian}.

Apply \eqref{ci:linear-energy} on $[0,t]$ with $\psi=f$. Since $T_*\leq C_0^{-1}\ln\left(\frac{4}{3}\right)$, its constant is uniform for $t\leq T_*$, and
\begin{align}
\sup_{0\leq s\leq t}\norm{f(s)}_{L^2}^2+\int_0^t\left(\norm{f}_{\Vspace}^2+\norm{\trp f}_{L^2_{\Sigma_+}}^2\right)\dd s\lesssim\left(\norm{f_0}_{L^2}^2+\int_0^t\norm{g}_{L^2_{\Sigma_-}}^2\dd s+\int_0^t\norm{f}_{L^2}^2\dd s\right). \label{ci:local-energy-proof}
\end{align}
Gr\"onwall's inequality bounds the $L^\infty_tL^2$ norm by the data. Substitution into \eqref{ci:local-energy-proof} bounds the other two terms, proving \eqref{ci:local-energy} and the energy requirements of Definition~\ref{def:energy}.
\end{proof}

\subsection{Uniform \texorpdfstring{$L^2$}{L2} Estimates under In-flow}\label{sec:energy-decay}

We will prove estimates uniform in time for the nonlinear solution constructed in Proposition~\ref{ci:local-existence} by weighted energy estimates. The collision estimate \eqref{cvel:micro} controls $(I-\pk)f$. We will choose a bounded positive weight whose transport derivative supplies an additional term that, together with this collision estimate, controls $\pk f$. The use of a space--velocity weight follows \cite{DengDuan2021}. The broader principle of choosing test functions so that transport controls $\pk f$ also appears in \cite[Lemmas~3.3 and~6.1]{Esposito.Guo.Kim.Marra2013}. Here an explicit weight and the finite-dimensional estimate below suffice.

Set
\begin{align}
\spatialmult(x,v):=\frac{x\cdot v}{\la},\qquad R_\Om:=\sup_{x\in\Om}\abs{x}. \label{ci:transport-weight}
\end{align}
Since $\Y\spatialmult=\frac{\abs{v}^2}{\la}$, testing with $(1-\kappa\spatialmult)f$ produces $\frac{\kappa}{2}\int_{\phase}\frac{\abs{v}^2}{\la}\abs{f}^2$. The following inequality combines this term with the collision control of $(I-\pk)f$.

\begin{lemma}\label{ci:velocity-transport-control}
For $f\in L^2\cap\Vspace$ such that $\la^{-\frac{1}{2}}\abs{v}f\in L^2$, we have
\begin{align}
\norm{f}_{\Vspace}^2+\norm{f}_{L^2}^2\lesssim\biggl(\norm[\big]{(I-\pk)f}_{\Vspace}^2+\norm[\Big]{\la^{-\frac{1}{2}}\abs{v}f}_{L^2}^2\biggr). \label{ci:velocity-transport-norm}
\end{align}
\end{lemma}

\begin{proof}
The annular $L^2$ norm is a norm on $\mathcal N$ in \eqref{cvel:nullspace}: a polynomial multiplying $\mh$ cannot vanish on $\{1<\abs{v}<2\}$ unless it is zero. Finite-dimensional norm equivalence, integrated in $x$, gives
\begin{align}
\norm{\pk f}_{\Vspace}^2+\norm{\pk f}_{L^2}^2\lesssim\int_\Om\int_{1<\abs{v}<2}\abs{\pk f}^2\dd v\dd x. \label{ci:annulus-nullspace}
\end{align}
On $\{1<\abs{v}<2\}$, both $\la^{\gamma+2}$ and $\frac{\abs{v}^2}{\la}$ have positive lower bounds. Hence the identity $\pk f=f-(I-\pk)f$ implies
\begin{align}
\int_\Om\int_{1<\abs{v}<2}\abs{\pk f}^2\dd v\dd x\lesssim\norm[\big]{(I-\pk)f}_{\Vspace}^2+\norm[\Big]{\la^{-\frac{1}{2}}\abs{v}f}_{L^2}^2. \label{ci:projection-annulus-bound}
\end{align}
Together with $f=(I-\pk)f+\pk f$, these estimates prove the $\Vspace$ bound.

For the $L^2$ bound, split at $\abs{v}=1$, using $\la^{\gamma+2}\gtrsim1$ inside $B_1$ and $\frac{\abs{v}^2}{\la}\geq\frac{1}{\sqrt2}$ outside:
\begin{align}
\norm{f}_{L^2}^2\lesssim\norm[\big]{(I-\pk)f}_{\Vspace}^2+\norm{\pk f}_{L^2}^2+\norm[\Big]{\la^{-\frac{1}{2}}\abs{v}f}_{L^2}^2. \label{ci:velocity-L2-splitting}
\end{align}
Substituting \eqref{ci:annulus-nullspace}--\eqref{ci:projection-annulus-bound} proves the $L^2$ bound.
\end{proof}

\begin{proposition}\label{ci:energy-response}
There exist $\delta,\lambda_*>0$ with the following property. Let $f$ be a weak solution of the in-flow problem in Definition~\ref{def:energy} on $[0,T]$, with initial value $f_0\in L^2$ and incoming value $g\in L^2\bigl((0,T)\times\Sigma_-;\dd t\dd\gamma\bigr)$. Suppose that
\begin{align}
\tnorm{f}_{L^\infty(0,T;L^{\infty,\alpha})}\leq\delta. \label{ci:energy-coefficient-bound}
\end{align}
Then, for every $t\in[0,T]$,
\begin{align}
\norm[\big]{f(t)}_{L^2}^2\lesssim\ue^{-2\lambda_*t}\norm{f_0}_{L^2}^2+\int_0^t\ue^{-2\lambda_*(t-r)}\norm[\big]{g(r)}_{L^2_{\Sigma_-}}^2\dd r, \label{ci:forced-energy-decay}
\end{align}
\begin{align}
\int_0^t\left(\norm[\big]{f(r)}_{\Vspace}^2+\norm[\big]{\trp f(r)}_{L^2_{\Sigma_+}}^2\right)\dd r\lesssim\norm{f_0}_{L^2}^2+\int_0^t\norm[\big]{g(r)}_{L^2_{\Sigma_-}}^2\dd r. \label{ci:integrated-energy}
\end{align}
The constants are independent of $T$ and of derivatives of the data.
\end{proposition}

\begin{proof}
\proofstep{Construction of the weight}
The definition \eqref{ci:transport-weight} gives
\begin{align}
\abs{\spatialmult}\leq R_\Om,\qquad \abs{\nabla_v\spatialmult}\leq\frac{2R_\Om}{\la},\qquad \Y\spatialmult=\frac{\abs{v}^2}{\la}. \label{ci:transport-multiplier-bounds}
\end{align}
Combining the derivative bound with $\sig\leq C\la^{\gamma+2}I$ from \eqref{cvel:eigenvalues}, we obtain $\sig^{ij}\partial_i\spatialmult\partial_j\spatialmult\lesssim\la^\gamma$. Expanding $\nabla_v(\spatialmult f)=\spatialmult\nabla_vf+f\nabla_v\spatialmult$ in \eqref{eq:V} therefore gives
\begin{align}
\norm{\spatialmult f}_{\Vspace}^2\leq2R_\Om^2\norm{f}_{\Vspace}^2+2\int_{\phase}\sig^{ij}\partial_i\spatialmult\partial_j\spatialmult\abs{f}^2\dd x\dd v\lesssim\norm{f}_{\Vspace}^2. \label{ci:transport-multiplier-V}
\end{align}
Here $\la^\gamma\leq\la^{\gamma+2}$ bounds the second integral by the zero-order term of the $\Vspace$ norm. Choose $0<\kappa\leq\frac{1}{2(1+R_\Om)}$, to be decreased below, and define
\begin{align}
\inflowweight_\kappa:=1-\kappa\spatialmult,\qquad \inflowenergy_\kappa(t):=\frac{1}{2}\int_{\phase}\inflowweight_\kappa\abs[\big]{f(t)}^2\dd x\dd v. \label{ci:weighted-inflow-energy}
\end{align}
In particular,
\begin{align}
\frac{1}{2}\leq\inflowweight_\kappa\leq\frac{3}{2},\qquad \frac{1}{4}\norm[\big]{f(t)}_{L^2}^2\leq\inflowenergy_\kappa(t)\leq\frac{3}{4}\norm[\big]{f(t)}_{L^2}^2. \label{ci:weighted-energy-equivalence}
\end{align}

\proofstep{The weighted energy identity}
Apply \eqref{trace:local-quadratic-Green} with $S=\Ga(f,f)$ and $\varphi=\inflowweight_\kappa\chi_R^2$, where $\chi_R$ is the velocity cutoff from \eqref{trace:cutoff-Green}. Since $\chi_R$ depends only on $v$,
\begin{align}
\Y\chi_R=0,\qquad \Y\bigl(\inflowweight_\kappa\chi_R^2\bigr)=-\kappa\frac{\abs{v}^2}{\la}\chi_R^2. \label{ci:weighted-transport-derivative}
\end{align}
We pass to the limit $R\to\infty$ in each term of the Green identity.

For the initial and final terms, $f(0),f(t)\in L^2$ and the boundedness of $\inflowweight_\kappa$ give convergence by dominated convergence for every $t\in[0,T]$.

For the boundary terms, Proposition~\ref{trace:green} gives $\trpm f\in L^2(\dd t\dd\gamma)$. Since $\inflowweight_\kappa$ is bounded and $\trm f=g$, dominated convergence gives the incoming and outgoing integrals with weight $\inflowweight_\kappa$.

For the collision pairing, \eqref{eq:intro-collision-dual} gives $\Ga(f,f)-\Lop f\in L^2(0,T;\Vspace')$. Multiplication by $\inflowweight_\kappa$ is bounded on $\Vspace$ by \eqref{ci:transport-multiplier-V}, so \eqref{trace:cutoff-V-convergence} yields
\begin{align}
\inflowweight_\kappa\chi_R^2f\to\inflowweight_\kappa f\quad\text{in }L^2(0,T;\Vspace). \label{ci:weighted-cutoff-convergence}
\end{align}
Consequently,
\begin{align}
\abs[\Big]{\int_0^t\ip[\big]{\Ga(f,f)-\Lop f}{\inflowweight_\kappa(\chi_R^2-1)f}\dd r}&\lesssim(1+\delta)\tnorm{f}_{L^2(0,T;\Vspace)}\tnorm{\inflowweight_\kappa(\chi_R^2-1)f}_{L^2(0,T;\Vspace)}\to0. \label{ci:weighted-cutoff-pairing}
\end{align}
For the transport integral, the Gaussian assumption \eqref{ci:energy-coefficient-bound} gives
\begin{align}
\tnorm[\Big]{\la^{-\frac{1}{2}}\abs{v}f}_{L^\infty(0,T;L^2)}^2\leq\abs{\Om}\delta^2\int_{\R^3}\frac{\abs{v}^2}{\la}\ue^{-2\alpha\abs{v}^2}\dd v\lesssim\delta^2. \label{ci:weighted-transport-integrability}
\end{align}
Thus $\frac{\abs{v}^2}{\la}\abs{f}^2\in L^1((0,T)\times\phase)$. Moving the transport integral to the left side of the Green identity and using dominated convergence, we obtain
\begin{align}
-\frac12\int_0^t\int_{\phase}\abs{f}^2\Y\bigl(\inflowweight_\kappa\chi_R^2\bigr)\dd x\dd v\dd r&=\frac{\kappa}{2}\int_0^t\int_{\phase}\frac{\abs{v}^2}{\la}\chi_R^2\abs{f}^2\dd x\dd v\dd r \label{ci:weighted-transport-limit}\\
&\to\frac{\kappa}{2}\int_0^t\norm[\Big]{\la^{-\frac{1}{2}}\abs{v}f(r)}_{L^2}^2\dd r.\notag
\end{align}
This is the positive transport contribution in the weighted energy identity. All time integrands in the limiting identity belong to $L^1(0,T)$, so $\inflowenergy_\kappa$ is absolutely continuous. Differentiating gives, for almost every $t\in(0,T)$,
\begin{align}
\frac{\dd}{\dd t}\inflowenergy_\kappa(t)+\frac{1}{2}\int_{\Sigma_+}\inflowweight_\kappa\abs{\trp f}^2\dd\gamma+\frac{\kappa}{2}\norm[\Big]{\la^{-\frac{1}{2}}\abs{v}f}_{L^2}^2+\ip*{\Lop f}{\inflowweight_\kappa f}&=\frac{1}{2}\int_{\Sigma_-}\inflowweight_\kappa\abs{g}^2\dd\gamma+\ip[\big]{\Ga(f,f)}{\inflowweight_\kappa f}. \label{ci:weighted-inflow-identity}
\end{align}

\proofstep{The collision terms and the uniform estimates}
By \eqref{cvel:nullspace}, $\Lop f=\Lop(I-\pk)f$. Hence the additional pairing in the weighted identity contains the factor $(I-\pk)f$. Equations \eqref{cvel:linear-boundedness}, \eqref{cvel:micro} and \eqref{ci:transport-multiplier-V} give
\begin{align}
\ip*{\Lop f}{\inflowweight_\kappa f}&=\ip*{\Lop f}{f}-\kappa\ip[\big]{\Lop(I-\pk)f}{\spatialmult f} \label{ci:weighted-collision-bound}\\
&\geq\lambda_\gamma\norm[\big]{(I-\pk)f}_{\Vspace}^2-C\kappa\norm[\big]{(I-\pk)f}_{\Vspace}\norm{f}_{\Vspace}\geq\frac{3\lambda_\gamma}{4}\norm[\big]{(I-\pk)f}_{\Vspace}^2-C\kappa^2\norm{f}_{\Vspace}^2.\notag
\end{align}
Since $\norm{\inflowweight_\kappa f}_{\Vspace}\lesssim\norm{f}_{\Vspace}$ uniformly in the chosen range of $\kappa$, equations \eqref{cvel:spatial-dual}, \eqref{eq:gaussian-energy} and \eqref{ci:energy-coefficient-bound} imply
\begin{align}
\abs[\Big]{\ip[\big]{\Ga(f,f)}{\inflowweight_\kappa f}}\lesssim\delta\norm{f}_{\Vspace}^2. \label{ci:weighted-nonlinear-bound}
\end{align}
Inserting these estimates into \eqref{ci:weighted-inflow-identity} yields
\begin{align}
\frac{\dd}{\dd t}\inflowenergy_\kappa(t)+\frac{1}{4}\norm{\trp f}_{L^2_{\Sigma_+}}^2+\frac{3\lambda_\gamma}{4}\norm[\big]{(I-\pk)f}_{\Vspace}^2+\frac{\kappa}{2}\norm[\Big]{\la^{-\frac{1}{2}}\abs{v}f}_{L^2}^2&\leq\frac{3}{4}\norm{g}_{L^2_{\Sigma_-}}^2+C(\kappa^2+\delta)\norm{f}_{\Vspace}^2. \label{ci:weighted-precoercive}
\end{align}
Choose $\kappa>0$ sufficiently small in the range fixed above, and then $\delta>0$ sufficiently small, so that
\begin{align}
C(\kappa^2+\delta)\leq\frac14\min\{\lambda_\gamma,\kappa\}. \label{ci:weighted-smallness}
\end{align}
By \eqref{ci:velocity-transport-norm}, this absorbs the error in \eqref{ci:weighted-precoercive} and gives
\begin{align}
\frac{\dd}{\dd t}\inflowenergy_\kappa(t)+\frac14\norm{\trp f}_{L^2_{\Sigma_+}}^2+\frac{\lambda_\gamma}{2}\norm[\big]{(I-\pk)f}_{\Vspace}^2+\frac{\kappa}{4}\norm[\Big]{\la^{-\frac{1}{2}}\abs{v}f}_{L^2}^2\leq\frac34\norm{g}_{L^2_{\Sigma_-}}^2. \label{ci:weighted-absorbed-energy}
\end{align}
By \eqref{ci:velocity-transport-norm}, half of the last two terms controls $\norm{f}_{\Vspace}^2$, and the other half controls $\norm{f}_{L^2}^2$. Using \eqref{ci:weighted-energy-equivalence}, we obtain $c,\lambda_*>0$ such that
\begin{align}
\frac{\dd}{\dd t}\inflowenergy_\kappa(t)+2\lambda_*\inflowenergy_\kappa(t)+c\norm{f}_{\Vspace}^2+\frac{1}{4}\norm{\trp f}_{L^2_{\Sigma_+}}^2\leq\frac{3}{4}\norm{g}_{L^2_{\Sigma_-}}^2. \label{ci:weighted-energy-decay}
\end{align}
Multiply \eqref{ci:weighted-energy-decay} by $\ue^{2\lambda_*t}$, integrate from $0$ to $t$, and use \eqref{ci:weighted-energy-equivalence}. This proves \eqref{ci:forced-energy-decay}. Integrating \eqref{ci:weighted-energy-decay} without the exponential factor proves \eqref{ci:integrated-energy}.
\end{proof}

The same estimates apply after translating time. On an interval starting at $s$, one uses $f(s)$ as the initial value and $g(s+\cdot)$ as the incoming data, with the same constants.

\subsection{Uniform Gaussian Estimates under In-flow}\label{sec:energy-Gaussian}

We will obtain the uniform Gaussian bound by combining two estimates. On bounded velocity sets, local $L^2$ to $L^\infty$ estimates control $f$ using the energy bound from Section~\ref{sec:energy-decay}. At large velocities, a Gaussian supersolution controls $f$ using the data and its bound on bounded velocity sets. Lemmas~\ref{ci:local-inflow-bound} and \ref{ci:tilted-barrier} establish these respective estimates. Proposition~\ref{ci:gaussian-bootstrap} combines them and absorbs their small multiples of the full Gaussian norm, giving bounds uniform in time and decay when the incoming data decay.

Throughout this subsection, $f$ is a weak solution in the sense of Definition~\ref{def:energy}, with initial value $f_0$ and incoming value $g$, and
\begin{align}
\tnorm{f}_{L^\infty(0,T;L^{\infty,\alpha})}\leq\delta,\qquad 0<\delta\leq\delta_A. \label{ci:Gaussian-assumption}
\end{align}
Here $\delta_A$ is the threshold in Lemma~\ref{cvel:coefficients}(i), so \eqref{cvel:ellipticity} gives $\frac12\sig\leq\Af{}\leq\frac32\sig$. We also choose $\delta$ small enough for Proposition~\ref{ci:energy-response}. Recall that \eqref{eq:intro-trace-space} gives $f\in C([0,T];L^2)\cap L^2(0,T;\Vspace)$, while \eqref{ci:forced-energy-decay}--\eqref{ci:integrated-energy} give the quantitative $L^2$ and dissipation bounds. The former estimate controls the $L^\infty_tL^2$ terms below and yields $L^2$ decay when the incoming data decay. We will decrease $\delta$ further when needed.

We will use the divergence equation \eqref{cvel:scalar-div}, which Lemma~\ref{cvel:scalar} justifies in distributions for the present weak solution. Its coefficients are defined in \eqref{cvel:coefficient-convolutions} and \eqref{cvel:scalar-coefficients}; Lemmas~\ref{cvel:coefficients} and \ref{cvel:scalar-bounds} give their bounds. The trace is the one in Definition~\ref{def:energy}, by Section~\ref{sec:trace-definition}.

\begin{lemma}[Estimate on bounded velocity sets]\label{ci:local-inflow-bound}
For $R\geq1$, $0<r\leq\frac{1}{2}$ and $0\leq\lambda\leq\frac{1}{2}$, we have
\begin{align}
\tnorm[\Big]{\ue^{\lambda t}f}_{L^\infty((0,T)\times\Om\times B_R)}\leq&\norm{f_0}_{L^\infty}+\tnorm[\Big]{\ue^{\lambda t}g}_{L^\infty((0,T)\times\Sigma_-)} \label{ci:bounded-velocity-estimate}\\
&\quad+C_Rr^{-6}\tnorm[\Big]{\ue^{\lambda t}f}_{L^\infty(0,T;L^2)}+C_Rr^2\tnorm[\Big]{\ue^{\lambda t}f}_{L^\infty(0,T;L^{\infty,\alpha})}.\notag
\end{align}
The constant is independent of $r,T,\lambda$ and of derivatives of the data. If $T>4r^2$, then
\begin{align}
\tnorm{f}_{L^\infty((4r^2,T)\times\Om\times B_R)}\leq\tnorm{g}_{L^\infty((0,T)\times\Sigma_-)}+C_Rr^{-6}\tnorm{f}_{L^\infty(0,T;L^2)} +C_Rr^2\tnorm{f}_{L^\infty(0,T;L^{\infty,\alpha})}. \label{ci:bounded-velocity-late}
\end{align}
\end{lemma}

\begin{proof}
We prove \eqref{ci:bounded-velocity-estimate} by the argument of \cite[Theorem~1.1]{Silvestre2022}, whose main ingredient is the local $L^2$ to $L^\infty$ estimate in \cite[Proposition~6.1]{Silvestre2022}, applied to the truncations \eqref{ci:bounded-velocity-truncations}. We verify the coefficient and source bounds \eqref{ci:bounded-velocity-coefficients}--\eqref{ci:bounded-velocity-source}, the $L^2$ regularity \eqref{ci:bounded-velocity-truncation-data}, and the weak subsolution inequality \eqref{ci:bounded-velocity-subsolution}. We then extend the truncations by zero across the initial time and spatial boundary and rescale Proposition~6.1 to obtain \eqref{ci:bounded-velocity-scaled}.

\proofstep{Coefficient and source bounds}
Multiplying \eqref{cvel:scalar-div} by $\ue^{\lambda t}$ gives
\begin{align}
\Y(\ue^{\lambda t}f)-\divv_v\bigl(\Af{}\nabla_v(\ue^{\lambda t}f)\bigr)-\Bf{}\cdot\nabla_v(\ue^{\lambda t}f)=\ue^{\lambda t}\bigl((\Cf{}+\lambda)f+\Kop f\bigr). \label{ci:bounded-velocity-weighted-equation}
\end{align}
Equations \eqref{cvel:ellipticity}, \eqref{cvel:rough-coefficients} and \eqref{cvel:scalar-coefficient-growth} give, on $\abs{v}<R+1$,
\begin{align}
C_R^{-1}I\leq\Af{}\leq C_RI,\qquad \tnorm{\abs{\Bf{}}+\abs{\Cf{}}}_{L^\infty((0,T)\times\Om\times B_{R+1})}\lesssim_R1. \label{ci:bounded-velocity-coefficients}
\end{align}
Linearity of $\Kop$ and \eqref{cvel:scalar-K-bound} bound the source in \eqref{ci:bounded-velocity-weighted-equation}:
\begin{align}
\tnorm[\Big]{\ue^{\lambda t}\bigl((\Cf{}+\lambda)f+\Kop f\bigr)}_{L^\infty((0,T)\times\Om\times B_{R+1})}\lesssim_R\tnorm[\Big]{\ue^{\lambda t}f}_{L^\infty(0,T;L^{\infty,\alpha})}. \label{ci:bounded-velocity-source}
\end{align}
These bounds are uniform for $0\leq\lambda\leq\frac12$ and $T>0$.

\proofstep{Positive-part subsolutions}
Set
\begin{align}
k:=\norm{f_0}_{L^\infty}+\tnorm[\Big]{\ue^{\lambda t}g}_{L^\infty((0,T)\times\Sigma_-)},\qquad u_\pm:=\bigl(\pm\ue^{\lambda t}f-k\bigr)_+. \label{ci:bounded-velocity-truncations}
\end{align}
The Lipschitz continuity of the positive-part map, \eqref{eq:intro-trace-space} and \eqref{trace:local-H1-bound} give
\begin{align}
u_\pm\in C([0,T];L^2),\quad \nabla_vu_\pm\in L^2((0,T)\times\Om\times B_{R+1}),\quad u_\pm(0)=0,\quad\trm u_\pm=0. \label{ci:bounded-velocity-truncation-data}
\end{align}
We will prove the following inequality in the weak sense of \eqref{ci:bounded-velocity-subsolution}, for the nonnegative test functions specified in \eqref{ci:bounded-velocity-tests}:
\begin{align}
\Y u_\pm-\divv_v(\Af{}\nabla_vu_\pm)-\Bf{}\cdot\nabla_vu_\pm\leq\ue^{\lambda t}\abs{(\Cf{}+\lambda)f+\Kop f}. \label{ci:bounded-velocity-positive-equation}
\end{align}

We use the renormalized identity \eqref{trace:beta-Green} from Section~\ref{sec:trace-definition}. To verify its hypotheses for $\ue^{\lambda t}f$, observe that
\begin{align}
\Y(\ue^{\lambda t}f)+\Lop(\ue^{\lambda t}f)=\ue^{\lambda t}\Ga(f,f)+\lambda\ue^{\lambda t}f. \label{ci:bounded-velocity-linearized-equation}
\end{align}
The weighted function belongs to $C([0,T];L^2)\cap L^2(0,T;\Vspace)$ by \eqref{eq:intro-trace-space}, since $T<\infty$. The first source term belongs to $L^2(0,T;\Vspace')$ by \eqref{trace:dual-bounds} and the Gaussian bound. For the second, the zeroth-order term in \eqref{eq:V} gives
\begin{align}
\norm{f}_{\Vspace'}\leq\norm[\Big]{\la^{-\frac{\gamma+2}{2}}f}_{L^2}\lesssim\gnorm{\infty}{\alpha}{f}, \label{ci:bounded-velocity-Gaussian-dual}
\end{align}
so it also belongs to $L^2(0,T;\Vspace')$. Thus \eqref{trace:beta-Green} applies with
\begin{align}
\beta\in C^2(\R),\qquad \beta',\beta''\in L^\infty(\R). \label{ci:bounded-velocity-renormalization-class}
\end{align}
The trace is $\ue^{\lambda t}\trace f$ by uniqueness in \eqref{trace:linear-Green}. Substituting \eqref{ci:bounded-velocity-weighted-equation} into the source pairing in \eqref{trace:beta-Green} gives the identity with diffusion $\Af{}$ and first-order coefficient $\Bf{}$. Its integrals are finite by \eqref{trace:local-H1-bound} and \eqref{ci:bounded-velocity-coefficients}--\eqref{ci:bounded-velocity-source}.

The admissible regularity of $\beta$ is essential here. In the proof of Lemma~\ref{ci:frozen-approximation}, \cite[Theorem~2.8(3)]{Carrapatoso.Mischler2024} allows $\beta(\xi)=\frac12(\xi^-)^2$: its weak second derivative is bounded, and that part of the theorem removes the boundedness requirement on $\beta'$. Here the desired function is $\beta(\xi)=\xi_+$, whose distributional second derivative is the Dirac mass at zero. It does not belong to $W^{2,\infty}_{\mathrm{loc}}$, so even the relaxed class does not permit its direct use. We therefore approximate it by smooth convex functions satisfying \eqref{ci:bounded-velocity-renormalization-class}. Choose a smooth nondecreasing function $\eta$, equal to zero on $(-\infty,0]$ and one on $[1,\infty)$, and define
\begin{align}
\beta_\varepsilon(\xi):=\int_0^\xi\eta\left(\frac{s}{\varepsilon}\right)\dd s,\qquad 0\leq\beta_\varepsilon'\leq1,\quad\beta_\varepsilon''\geq0,\quad0\leq\xi_+-\beta_\varepsilon(\xi)\leq\varepsilon. \label{ci:bounded-velocity-positive-approximation}
\end{align}
For each $\varepsilon>0$, both derivatives are bounded and $\beta_\varepsilon=0$ on $(-\infty,0]$. Apply \eqref{trace:beta-Green} with $\xi\mapsto\beta_\varepsilon(\pm\xi-k)$. For $u_{\pm,\varepsilon}:=\beta_\varepsilon(\pm\ue^{\lambda t}f-k)$, the resulting equation in distributions is
\begin{align}\label{ci:bounded-velocity-convexity}
&\Y u_{\pm,\varepsilon}-\divv_v(\Af{}\nabla_vu_{\pm,\varepsilon})-\Bf{}\cdot\nabla_vu_{\pm,\varepsilon}\\
&=\pm\beta_\varepsilon'(\pm\ue^{\lambda t}f-k)\ue^{\lambda t}\bigl((\Cf{}+\lambda)f+\Kop f\bigr)-\beta_\varepsilon''(\pm\ue^{\lambda t}f-k)\Af{}\nabla_v(\ue^{\lambda t}f)\cdot\nabla_v(\ue^{\lambda t}f)\notag\\
&\leq\ue^{\lambda t}\abs{(\Cf{}+\lambda)f+\Kop f}.\notag
\end{align}
Convexity and positivity of $\Af{}$ make the second term nonpositive. In \eqref{trace:beta-Green}, take
\begin{align}
0\leq\varphi\in C_c^\infty\bigl((-\infty,T)\times\R^3\times B_{R+1}\bigr). \label{ci:bounded-velocity-tests}
\end{align}
We restrict this test function to $[0,T]\times\overline\Om\times B_{R+1}$. It may be nonzero at $t=0$ and on the incoming, outgoing and grazing spatial boundary. The initial and incoming terms vanish by \eqref{ci:bounded-velocity-truncation-data}, while the outgoing term satisfies
\begin{align}
\int_0^T\int_{\Sigma_+}\beta_\varepsilon\bigl(\pm\ue^{\lambda t}\trp f-k\bigr)\varphi\dd\gamma\dd t\geq0. \label{ci:bounded-velocity-outgoing-sign}
\end{align}

Discard this outgoing term and the nonpositive diffusion term in \eqref{ci:bounded-velocity-convexity}. As $\varepsilon\to0$, the Sobolev chain rule and dominated convergence give $u_{\pm,\varepsilon}\to u_\pm$ and $\nabla_vu_{\pm,\varepsilon}\to\nabla_vu_\pm$ locally in $L^2$. Bounding the source by the right side of \eqref{ci:bounded-velocity-positive-equation} yields, for every nonnegative $\varphi$ in \eqref{ci:bounded-velocity-tests},
\begin{align}
-\int u_\pm\Y\varphi+\int\Af{}\nabla_vu_\pm\cdot\nabla_v\varphi-\int\Bf{}\cdot\nabla_vu_\pm\,\varphi\leq\int\ue^{\lambda t}\abs{(\Cf{}+\lambda)f+\Kop f}\,\varphi. \label{ci:bounded-velocity-subsolution}
\end{align}
The integrals are over $(0,T)\times\Om\times B_{R+1}$. Nonnegative smooth approximation in $C^1$, with support in a fixed compact subset, extends the inequality to nonnegative $C_c^1$ tests; all integrals converge by \eqref{ci:bounded-velocity-truncation-data} and the local coefficient and source bounds.

\proofstep{Extension across the initial time and spatial boundary}
Proposition~6.1 of \cite{Silvestre2022} requires a subsolution defined on a full kinetic cylinder. To apply it near the spatial boundary or initial time, we must extend $u_\pm$ beyond $(0,T)\times\Om$. The cylinders we use are
\begin{align}
\kinback_\rho(t_0,x_0,v_0):=\left\{(t,x,v):t_0-\rho^2<t<t_0,\ \abs{x-x_0-(t-t_0)v_0}<\rho^3,\ \abs{v-v_0}<\rho\right\}. \label{ci:backward-cylinders}
\end{align}
For $t<0$ or $x\notin\Om$, extend $u_\pm$ by zero, the diffusion $\Af{}$ and first-order coefficient $\Bf{}$ in \eqref{ci:bounded-velocity-weighted-equation} by $I$ and $0$, respectively, and the right side of \eqref{ci:bounded-velocity-positive-equation} by zero. We keep the same notation for these extensions. The coefficient and source bounds \eqref{ci:bounded-velocity-coefficients}--\eqref{ci:bounded-velocity-source} remain valid on $(-\infty,T)\times\R^3\times B_{R+1}$.

The zero initial value in \eqref{ci:bounded-velocity-truncation-data} removes the time-boundary term at $t=0$. For the spatial extension, we have
\begin{align}
\nabla_v(\ind_\Om u_\pm)=\ind_\Om\nabla_vu_\pm. \label{ci:bounded-velocity-extension-gradient}
\end{align}
Thus the extended $u_\pm$ and $\nabla_vu_\pm$ belong to $L^2$ on every cylinder contained in $(-\infty,T)\times\R^3\times B_{R+1}$. A nonnegative $C_c^1$ test supported in such a full cylinder restricts to an admissible test in \eqref{ci:bounded-velocity-subsolution}. The integrals outside $(0,T)\times\Om$ vanish, so that inequality also holds for the extensions. This verifies the weak subsolution inequality across the initial time and the incoming, outgoing and grazing spatial boundary.

In the notation of \cite[Definition~5.3 and Proposition~6.1]{Silvestre2022}, the diffusion matrix, first-order coefficient and source on the original domain are
\begin{align}
(a,b,G)=\left(\Af{},-\Bf{},\ue^{\lambda t}\abs{(\Cf{}+\lambda)f+\Kop f}\right). \label{ci:bounded-velocity-coefficient-identification}
\end{align}
Their extensions are $a=I$, $b=0$ and $G=0$.

\proofstep{Scaling and conclusion}
Fix $0<t_0<T$, $x_0\in\R^3$, $\abs{v_0}\leq R$ and $0<r\leq\frac12$. The condition $2r\leq1$ places the velocities of $\kinback_{2r}(t_0,x_0,v_0)$ in $B_{R+1}$, where \eqref{ci:bounded-velocity-coefficients}--\eqref{ci:bounded-velocity-source} hold. We apply \cite[Proposition~6.1]{Silvestre2022} to this cylinder after the change of variables
\begin{align}
t=t_0+(2r)^2s,\qquad x=x_0+(2r)^3y+(2r)^2sv_0,\qquad v=v_0+2rw. \label{ci:bounded-velocity-scaling-map}
\end{align}
Writing $Q_\rho:=\kinback_\rho(0,0,0)$, this map sends $Q_1$ to $\kinback_{2r}(t_0,x_0,v_0)$ and $Q_{\frac12}$ to $\kinback_r(t_0,x_0,v_0)$. Let hatted functions denote pullbacks under \eqref{ci:bounded-velocity-scaling-map}, and write $\Yhat=\partial_s+w\cdot\nabla_y$. The extended weak inequality becomes
\begin{align}
\Yhat\widehat u_\pm-\divv_w(\widehat a\nabla_w\widehat u_\pm)+2r\widehat b\cdot\nabla_w\widehat u_\pm\leq(2r)^2\widehat G\qquad\text{in }Q_1. \label{ci:bounded-velocity-scaled-equation}
\end{align}
The diffusion has the same ellipticity bounds, and $\norm[\big]{2r\widehat b}_{L^\infty(Q_1)}\lesssim_R1$, uniformly in $r$. The Jacobian of \eqref{ci:bounded-velocity-scaling-map} is $(2r)^{14}$, so the $L^2$ norm gains the factor $(2r)^{-7}$. Applying \cite[Proposition~6.1]{Silvestre2022} to \eqref{ci:bounded-velocity-scaled-equation}, returning to the original variables and using \eqref{ci:bounded-velocity-source}, we obtain
\begin{align}
\tnorm{u_\pm}_{L^\infty(\kinback_r\cap\{x\in\Om\})}\lesssim_R r^{-7}\tnorm{u_\pm}_{L^2(\kinback_{2r}\cap\{x\in\Om\})}+r^2\tnorm[\Big]{\ue^{\lambda t}f}_{L^\infty(0,T;L^{\infty,\alpha})}. \label{ci:bounded-velocity-scaled}
\end{align}
Since $0\leq u_\pm\leq\ue^{\lambda t}\abs{f}$ for $t>0$, $u_\pm=0$ for $t<0$, and the time interval of $\kinback_{2r}$ has length $4r^2$,
\begin{align}
\tnorm{u_\pm}_{L^2(\kinback_{2r}\cap\{x\in\Om\})}\leq2r\tnorm[\Big]{\ue^{\lambda t}f}_{L^\infty(0,T;L^2)}. \label{ci:bounded-velocity-time-factor}
\end{align}
Substitution in \eqref{ci:bounded-velocity-scaled} gives the factor $r^{-6}$. Covering $(0,T)\times\Om\times B_R$ by the smaller cylinders and using $\abs{\ue^{\lambda t}f}\leq k+\max\{u_+,u_-\}$ proves \eqref{ci:bounded-velocity-estimate}. The time extension includes cylinders meeting $t=0$ with the same constants.

For \eqref{ci:bounded-velocity-late}, set $\lambda=0$ and $k=\tnorm{g}_{L^\infty((0,T)\times\Sigma_-)}$. The preceding argument applies to cylinders with $t_0>4r^2$, since their larger cylinders lie at positive times. Equations \eqref{ci:bounded-velocity-scaled}--\eqref{ci:bounded-velocity-time-factor} then give \eqref{ci:bounded-velocity-late} without an initial-data term.
\end{proof}

\begin{lemma}[Gaussian estimate for large velocities]\label{ci:tilted-barrier}
For $\delta$ sufficiently small, there is $R_0\geq1$ such that, for every $R\geq R_0$, there is $\varepsilon_R>0$ for which, whenever $0\leq\lambda\leq\frac12$ and $0\leq\tau<T$, we have
\begin{align}
\tnorm[\Big]{\ind_{\{\abs{v}\geq R+1\}}\ue^{\lambda t}f}_{L^\infty(\tau,T;L^{\infty,\alpha})}\leq &C\ue^{-\frac{\tau}{2}}\gnorm{\infty}{\alpha}{f_0}+C\tnorm[\Big]{\ue^{\lambda t}g}_{L^\infty(0,T;L^{\infty,\alpha}_{\Sigma_-})} \label{ci:exterior-Gaussian-estimate}\\
&+C\ue^{\alpha(R+1)^2}\tnorm[\Big]{\ue^{\lambda t}f}_{L^\infty((0,T)\times\Om\times B_{R+1})}+\varepsilon_R\tnorm[\Big]{\ue^{\lambda t}f}_{L^\infty(0,T;L^{\infty,\alpha})}.\notag
\end{align}
The constant $C$ is independent of $R,T,\lambda,\tau$ and of derivatives of the data. The numbers $\varepsilon_R$ can be chosen so that $\varepsilon_R\to0$ as $R\to\infty$.

\end{lemma}

\begin{proof}
\proofstep{Construction of the Gaussian supersolution}
With $\spatialmult$ and $R_\Om$ from \eqref{ci:transport-weight}, define
\begin{align}
\gbarrier(x,v):=\ue^{-\left(\alpha+\frac{1}{4}\right)\abs{v}^2+\spatialmult(x,v)}. \label{eq:I:display:193}
\end{align}
The bound $\abs{\spatialmult}\leq R_\Om$ gives
\begin{align}
\gbarrier\simeq\ue^{-\left(\alpha+\frac14\right)\abs{v}^2}. \label{rev14:barrier-Gaussian-comparison}
\end{align}
Lemma~\ref{ci:radial-estimate}, with $G=F$ and $\zeta=\alpha+\frac14$, gives \eqref{ci:radial-barrier} for the radial factor $W_{\alpha+\frac14}$. Its smallness requirement is met by decreasing the standing $\delta$. Since $\gbarrier=\ue^{\spatialmult}W_{\alpha+\frac14}$, we must also estimate the terms in which velocity derivatives act on $\ue^{\spatialmult}$. The product rule gives
\begin{align}
\frac{Q(F,\gbarrier)}{\gbarrier}=\frac{Q(F,W_{\alpha+\frac{1}{4}})}{W_{\alpha+\frac{1}{4}}}-(4\alpha+1)(\Af{}v)\cdot\nabla_v\spatialmult+\Af{}:\left(D_v^2\spatialmult+\nabla_v\spatialmult\otimes\nabla_v\spatialmult\right). \label{eq:I:display:198}
\end{align}
By \eqref{eq:I:display:194}, \eqref{cvel:rough-coefficients} and the derivatives of \eqref{ci:transport-weight},
\begin{align}
\abs{(\Af{}v)\cdot\nabla_v\spatialmult}\leq(v\cdot\Af{}v)^{\frac12}(\nabla_v\spatialmult\cdot\Af{}\nabla_v\spatialmult)^{\frac12}\lesssim\la^{\gamma+1}. \label{ci:exterior-weight-errors}
\end{align}
\begin{align}
\abs{\Af{}:\bigl(D_v^2\spatialmult+\nabla_v\spatialmult\otimes\nabla_v\spatialmult\bigr)}\lesssim\la^\gamma. \label{ci:exterior-weight-second-derivatives}
\end{align}
Here $\abs{\nabla_v\spatialmult}\lesssim\la^{-1}$ and $\abs{D_v^2\spatialmult}\lesssim\la^{-2}$. Thus \eqref{ci:radial-barrier} implies
\begin{align}
\frac{Q(F,\gbarrier)}{\gbarrier}\leq-\frac{c_*}{2}\la^{\gamma+2}\qquad(\abs{v}\geq R_0),\qquad \frac{\Y\gbarrier}{\gbarrier}=\frac{\abs{v}^2}{\la}. \label{ci:exterior-barrier-calculation}
\end{align}
Increasing $R_0$ therefore gives, throughout $-3\leq\gamma\leq1$,
\begin{align}
\Y\gbarrier-Q(F,\gbarrier)\geq\gbarrier\qquad(\abs{v}\geq R_0). \label{ci:barrier-inequality}
\end{align}

\proofstep{Source estimate and comparison functions}
The collision identity \eqref{cvel:scalar-bilinearity} gives
\begin{align}
\Y\bigl(\ue^{\lambda t}\mh f\bigr)=Q\bigl(F,\ue^{\lambda t}\mh f\bigr)+\lambda\ue^{\lambda t}\mh f+\ue^{\lambda t}\mh\Kop f. \label{ci:exterior-density-equation}
\end{align}
By \eqref{cvel:scalar-K-bound} and \eqref{rev14:barrier-Gaussian-comparison},
\begin{align}
\frac{\ue^{\lambda t}\mh\abs{\Kop f}}{\gbarrier}\lesssim\la^5\ue^{-\left(\frac{1}{4}-\alpha\right)\abs{v}^2}\tnorm[\Big]{\ue^{\lambda t}f}_{L^\infty(0,T;L^{\infty,\alpha})}. \label{ci:exterior-source-smallness}
\end{align}
This includes $\gamma=-3$, since the part of $\Kop f$ containing $c$ equals $8\pi\m f$. The coefficient on the right tends to zero uniformly as $\abs{v}\to\infty$ because $\alpha<\frac{1}{4}$.

Choose a constant $C>0$, independent of $R,T,\lambda$, sufficiently large, and define
\begin{align}
\compamp(t)&:=C\ue^{-\frac{t}{2}}\gnorm{\infty}{\alpha}{f_0}+C\tnorm[\Big]{\ue^{\lambda t}g}_{L^\infty(0,T;L^{\infty,\alpha}_{\Sigma_-})}+C\ue^{\alpha(R+1)^2}\tnorm[\Big]{\ue^{\lambda t}f}_{L^\infty((0,T)\times\Om\times B_{R+1})} \label{ci:exterior-comparison-constant}\\
&\quad+C\vnorm[\Big]{\ind_{\{\abs{v}\geq R\}}\la^5\ue^{-\left(\frac{1}{4}-\alpha\right)\abs{v}^2}}_{L_v^\infty}\tnorm[\Big]{\ue^{\lambda t}f}_{L^\infty(0,T;L^{\infty,\alpha})}.\notag
\end{align}
Since $0\leq\lambda\leq\frac{1}{2}$, the contribution of its first term to $\compamp'+(1-\lambda)\compamp$ is nonnegative. The last term and \eqref{ci:exterior-source-smallness} therefore give
\begin{align}
\compamp'(t)+(1-\lambda)\compamp(t)\geq\frac{\ue^{\lambda t}\mh\abs{\Kop f}}{\gbarrier}\qquad(\abs{v}>R). \label{ci:exterior-amplitude-inequality}
\end{align}
Set $h_\pm:=\compamp(t)\gbarrier\pm\ue^{\lambda t}\mh f$. By \eqref{ci:exterior-density-equation}, \eqref{ci:barrier-inequality} and \eqref{ci:exterior-amplitude-inequality}, respectively, we have, for $\abs{v}>R$,
\begin{align}
\Y h_\pm-Q(F,h_\pm)-\lambda h_\pm&=\compamp'\gbarrier+\compamp\bigl(\Y\gbarrier-Q(F,\gbarrier)-\lambda\gbarrier\bigr)\pm\ue^{\lambda t}\mh\Kop f \label{ci:exterior-supersolutions}\\
&\geq\bigl(\compamp'+(1-\lambda)\compamp\bigr)\gbarrier-\ue^{\lambda t}\mh\abs{\Kop f}\geq0.\notag
\end{align}

By \eqref{rev14:barrier-Gaussian-comparison}, $\mh/\gbarrier\lesssim\ue^{\alpha\abs{v}^2}$. Therefore the first and second terms of \eqref{ci:exterior-comparison-constant} dominate the initial and incoming values, since
\begin{align}
\frac{\mh\abs{f_0}}{\gbarrier}\lesssim\gnorm{\infty}{\alpha}{f_0},\qquad \frac{\ue^{\lambda t}\mh\abs{g}}{\gbarrier}\lesssim\tnorm[\Big]{\ue^{\lambda t}g}_{L^\infty(0,T;L^{\infty,\alpha}_{\Sigma_-})}. \label{ci:exterior-data-domination}
\end{align}
On $R<\abs{v}<R+1$, the third term dominates the solution because
\begin{align}
\frac{\ue^{\lambda t}\mh\abs{f}}{\gbarrier}\lesssim\ue^{\alpha(R+1)^2}\tnorm[\Big]{\ue^{\lambda t}f}_{L^\infty((0,T)\times\Om\times B_{R+1})}. \label{ci:exterior-annulus-domination}
\end{align}
With $C$ in \eqref{ci:exterior-comparison-constant} chosen to satisfy these bounds, we obtain
\begin{align}
h_\pm(0)\geq0,\qquad\trm h_\pm\geq0,\qquad h_\pm\geq0\quad\text{on }(0,T)\times\Om\times\{R<\abs{v}<R+1\}. \label{ci:exterior-comparison-data}
\end{align}

\proofstep{The squared negative part}
We prove $h_\pm\geq0$ for $\abs{v}>R$. Define
\begin{align}
q_\pm:=\ue^{-\lambda t}\mih h_\pm=\ue^{-\lambda t}\compamp(t)\mih\gbarrier\pm f. \label{ci:exterior-normalized-functions}
\end{align}
The function $\mih\gbarrier$ and its velocity derivatives have Gaussian decay by \eqref{eq:I:display:193} and the derivative bounds following \eqref{ci:exterior-weight-second-derivatives}. Since $\compamp$ is smooth in time, \eqref{eq:intro-trace-space} gives $q_\pm\in C([0,T];L^2)\cap L^2(0,T;\Vspace)$. Their incoming traces are $\ue^{-\lambda t}\compamp\mih\gbarrier\pm g$, which belong to $L^2(0,T;L^2_{\Sigma_-})$ by the Gaussian bound for the data.

Multiplying the first equality in \eqref{ci:exterior-supersolutions} by $\ue^{-\lambda t}\mih$ gives
\begin{align}
\Y q_\pm-\mih Q(F,\mh q_\pm)&=\ue^{-\lambda t}\mih\bigl[\compamp'\gbarrier+\compamp\bigl(\Y\gbarrier-Q(F,\gbarrier)-\lambda\gbarrier\bigr)\bigr]\pm\Kop f\geq0\quad(\abs{v}>R). \label{ci:exterior-normalized-supersolutions}
\end{align}
Write $q$ for either $q_+$ or $q_-$, set $q^-:=\max\{-q,0\}$, and denote the corresponding right side of \eqref{ci:exterior-normalized-supersolutions} by $\mathcal S$. By \eqref{eq:I:display:198}, the local coefficient bounds and \eqref{cvel:scalar-K-bound}, we have $\mathcal S\in L^2((0,T)\times\Om\times B_M)$ for every $M>0$. The divergence definition \eqref{eq:collision} and \eqref{cvel:convolved-identities}, with the local velocity $H^1$ regularity of $q$, give in distributions
\begin{align}
\Y q-\divv_v(\Af{}\nabla_vq)+(b_F+\Af{}v)\cdot\nabla_vq=\left(\frac14v\cdot\Af{}v-\frac12\tr\Af{}-c_F\right)q+\mathcal S. \label{ci:exterior-negative-equation}
\end{align}
We next verify that \eqref{trace:beta-Green} from Section~\ref{sec:trace-definition} applies to \eqref{ci:exterior-negative-equation} with $\beta(s)=\frac12(s^-)^2$. The local velocity $H^1$ regularity follows from \eqref{trace:local-H1-bound}. Lemma~\ref{cvel:coefficients} gives locally uniformly elliptic $\Af{}$ and locally bounded $b_F+\Af{}v$ and zeroth-order coefficient. Hence the right side of \eqref{ci:exterior-negative-equation} is locally $L^2$. The initial and incoming values satisfy \eqref{trace:local-incoming-hypotheses} by the preceding bounds.

Under these hypotheses, \cite[Theorem~2.8(3)]{Carrapatoso.Mischler2024} allows $\beta\in W^{2,\infty}_{\mathrm{loc}}(\R)$ with bounded $\beta''$. Our choice satisfies $\beta'(s)=-s^-$ and $\beta''(s)=\ind_{\{s<0\}}$, so it is admissible. Uniqueness of the kinetic trace identifies its boundary values with those of $q$ above, and $C_tL^2$ continuity identifies its time values. We may therefore use \eqref{trace:beta-Green} for $q$ with this squared negative part and every smooth $\varphi$ with compact velocity support.

\proofstep{The cutoff identity}
Keep $R$ fixed. Choose $\eta_R\in C^\infty(\R^3)$, with $0\leq\eta_R\leq1$, equal to zero on $B_R$ and one outside $B_{R+1}$. Let $\chi_L$ be the radial cutoff from \eqref{trace:cutoff-Green}, with $L>R+1$, and take $\varphi=\eta_R^2\chi_L^2$. The factor $\chi_L^2$ gives compact velocity support, while $\varphi=0$ on $B_R$ and \eqref{ci:exterior-normalized-supersolutions} give $\mathcal S\varphi\geq0$.

The positive factor in \eqref{ci:exterior-normalized-functions} transfers \eqref{ci:exterior-comparison-data} to $q$, so
\begin{align}
q^-(0)=0,\qquad(\trm q)^-=0,\qquad q^-=\nabla_vq^-=0\quad\text{on }(0,T)\times\Om\times\{R<\abs{v}<R+1\}. \label{ci:exterior-negative-data}
\end{align}
The gradient vanishes by the Sobolev chain rule.

Apply \eqref{trace:beta-Green} with $\beta(s)=\frac12(s^-)^2$ and the chosen $\varphi$. We have $\Y\varphi=0$, and the initial and incoming terms vanish by \eqref{ci:exterior-negative-data}. Every term differentiating $\eta_R$ contains $q^-$ or $\nabla_vq^-$ on this annulus and is therefore zero. Integrating the term $(b_F+\Af{}v)\cdot\nabla_vq$ by parts gives the coefficient identity \eqref{ci:negative-part-coefficient} with $G=F$. Thus, with $\varpi_F$ defined by \eqref{eq:I:display:218} and unmarked integrals over $(0,t)\times\phase$, we obtain
\begin{align}
&\frac12\norm{\eta_R\chi_Lq^-(t)}_{L^2}^2+\frac12\int_0^t\norm{\eta_R\chi_L(\trp q)^-}_{L^2_{\Sigma_+}}^2\dd s+\int\eta_R^2\chi_L^2\Af{}\nabla_vq^-\cdot\nabla_vq^- \label{ci:exterior-negative-cutoff-identity}\\
&\quad=\int\eta_R^2\chi_L^2\varpi_F\big(q^-\big)^2-\int\eta_R^2\chi_L^2\mathcal S q^-\notag\\
&\qquad-2\int\eta_R^2\chi_Lq^-\Af{}\nabla_vq^-\cdot\nabla_v\chi_L+\int\eta_R^2\chi_L\big(q^-\big)^2(b_F+\Af{}v)\cdot\nabla_v\chi_L.\notag
\end{align}

\proofstep{The limit $L\to\infty$}
We first estimate the two terms differentiating $\chi_L$, keeping $R$ fixed. By \eqref{cvel:ellipticity}, $\Af{}\leq C\sig$. Together with \eqref{cvel:eigenvalues} and $\abs{\nabla_v\chi_L}\lesssim L^{-1}$, this gives
\begin{align}
\abs{2\int\eta_R^2\chi_Lq^-\Af{}\nabla_vq^-\cdot\nabla_v\chi_L}\lesssim\frac1L\left(\int_{\abs{v}>L}\sig\nabla_vq^-\cdot\nabla_vq^-\right)^{\frac12}\left(\int_{\abs{v}>L}\la^{\gamma+2}\big(q^-\big)^2\right)^{\frac12}\to0. \label{ci:exterior-negative-cutoff-diffusion}
\end{align}
By \eqref{cvel:rough-coefficients}, $\abs{b_F+\Af{}v}\lesssim\la^{\gamma+3}$. Since $\nabla_v\chi_L$ is supported on $L<\abs{v}<2L$, we also have
\begin{align}
\abs{\int\eta_R^2\chi_L\big(q^-\big)^2(b_F+\Af{}v)\cdot\nabla_v\chi_L}\lesssim\int_{\abs{v}>L}\la^{\gamma+2}\big(q^-\big)^2\to0. \label{ci:exterior-negative-cutoff-first-order}
\end{align}
Both limits follow from $q\in L^2(0,T;\Vspace)$ and the Sobolev chain rule.

In \eqref{ci:exterior-negative-cutoff-identity}, the source term is nonpositive because $\mathcal S\geq0$ on the support of $\eta_R$. We discard it and the nonnegative outgoing and diffusion terms. The standing smallness of $F-\m$ allows \eqref{eq:I:display:219} with $G=F$ and $\varepsilon=0$, so $\varpi_F\leq C$. We then let $L\to\infty$. Dominated convergence in the time-slice and spacetime $L^2$ integrals, together with \eqref{ci:exterior-negative-cutoff-diffusion} and \eqref{ci:exterior-negative-cutoff-first-order}, gives
\begin{align}
\norm{\eta_Rq^-(t)}_{L^2}^2\lesssim\int_0^t\norm{\eta_Rq^-(s)}_{L^2}^2\dd s. \label{ci:exterior-negative-Gronwall}
\end{align}
Gr\"onwall's inequality gives $\eta_Rq^-(t)=0$ for every $t\in[0,T]$. The annular vanishing in \eqref{ci:exterior-negative-data} also holds at every time by $q\in C([0,T];L^2)$. Since $\eta_R=1$ outside $B_{R+1}$, we obtain $(q_\pm)^-=0$, hence $(h_\pm)^-=0$, for $\abs{v}>R$. Consequently,
\begin{align}
\ue^{\lambda t}\mh\abs{f}\leq\compamp(t)\gbarrier\qquad(\abs{v}>R). \label{ci:exterior-Gaussian-domination}
\end{align}
Divide by $\mh$, use \eqref{rev14:barrier-Gaussian-comparison} and \eqref{ci:exterior-comparison-constant}, and take the essential supremum over $\tau<t<T$. Since $\ue^{-t/2}\leq\ue^{-\tau/2}$, this proves \eqref{ci:exterior-Gaussian-estimate}, with
\begin{align}
\varepsilon_R:=C\vnorm[\Big]{\ind_{\{\abs{v}\geq R\}}\la^5\ue^{-\left(\frac14-\alpha\right)\abs{v}^2}}_{L_v^\infty}\to0. \label{ci:exterior-small-constant}
\end{align}
\end{proof}

\begin{proposition}\label{ci:gaussian-bootstrap}
There is $\delta>0$ with the following property. Let $f$ be a weak solution of the in-flow problem on $[0,T]$, with initial value $f_0$, incoming value $g$, and
\begin{align}
\tnorm{f}_{L^\infty(0,T;L^{\infty,\alpha})}\leq\delta. \label{eq:I:display:207}
\end{align}
Then the following holds.
\begin{enumerate}
\item Uniformly in $T$,
\begingroup\usetagform{enumerated}
\begin{align}
\tnorm{f}_{L^\infty(0,T;L^{\infty,\alpha})}\lesssim\gnorm{\infty}{\alpha}{f_0}+\tnorm{g}_{L^\infty(0,T;L^{\infty,\alpha}_{\Sigma_-})}. \label{ci:uniform-Gaussian}
\end{align}
\endgroup
\item If $\bgnorm{\infty}{\alpha}{-}{g(t)}\leq\eta_\omega\ue^{-\omega t}$ with $\omega>0$ and $\eta_\omega\geq0$, there are $C_\omega,\lambda>0$, with $\lambda$ depending on $\omega$, such that
\begingroup\usetagform{enumerated}
\begin{align}
\norm[\big]{f(t)}_{L^2}+\gnorm[\big]{\infty}{\alpha}{f(t)}\leq C_\omega\left(\gnorm{\infty}{\alpha}{f_0}+\eta_\omega\right)\ue^{-\lambda t}\qquad(0\leq t\leq T). \label{ci:decaying-Gaussian}
\end{align}
\endgroup

\item Suppose that $f$ is defined on $[0,\infty)$, that \eqref{eq:I:display:207} holds for every $T$, and that $\bgnorm{\infty}{\alpha}{-}{g(t)}\to0$ as $t\to\infty$. Then
\begingroup\usetagform{enumerated}
\begin{align}
\norm[\big]{f(t)}_{L^2}+\gnorm[\big]{\infty}{\alpha}{f(t)}\to0\qquad(t\to\infty). \label{ci:vanishing-inflow-Gaussian}
\end{align}
\endgroup
\end{enumerate}
\end{proposition}

\begin{proof}
Choose $\delta$ small enough for Proposition~\ref{ci:energy-response} and Lemma~\ref{ci:tilted-barrier}. We first prove an estimate valid for $0\leq\lambda\leq\frac{1}{2}$. The constants in this argument are independent of $T$ and $\lambda$.

\proofstep{Absorption of the Gaussian norm} We combine \eqref{ci:exterior-Gaussian-estimate} with $\tau=0$ and \eqref{ci:bounded-velocity-estimate} with $R+1$ in place of $R$. Choose $R\geq R_0$ so that $\varepsilon_R\leq\frac{1}{4}$. Since the Gaussian weight on $B_{R+1}$ is bounded by $\ue^{\alpha(R+1)^2}$, these two estimates give
\begin{align}
\tnorm[\Big]{\ue^{\lambda t}f}_{L^\infty(0,T;L^{\infty,\alpha})}\leq &C_R\left(\gnorm{\infty}{\alpha}{f_0}+\tnorm[\Big]{\ue^{\lambda t}g}_{L^\infty(0,T;L^{\infty,\alpha}_{\Sigma_-})}+r^{-6}\tnorm[\Big]{\ue^{\lambda t}f}_{L^\infty(0,T;L^2)}\right) \label{ci:Gaussian-before-absorption}\\
&\qquad+\left(\frac{1}{4}+C_Rr^2\right)\tnorm[\Big]{\ue^{\lambda t}f}_{L^\infty(0,T;L^{\infty,\alpha})}.\notag
\end{align}
Fix $0<r\leq\frac{1}{2}$ so that $C_Rr^2\leq\frac{1}{4}$. The Gaussian norm is finite by \eqref{eq:I:display:207}, since $T<\infty$, so its last occurrence can be absorbed. We have proved
\begin{align}
&\tnorm[\Big]{\ue^{\lambda t}f}_{L^\infty(0,T;L^{\infty,\alpha})}\lesssim\gnorm{\infty}{\alpha}{f_0}+\tnorm[\Big]{\ue^{\lambda t}g}_{L^\infty(0,T;L^{\infty,\alpha}_{\Sigma_-})}+\tnorm[\Big]{\ue^{\lambda t}f}_{L^\infty(0,T;L^2)}. \label{ci:Gaussian-from-energy}
\end{align}
With $R$ and $r$ fixed, the constant is independent of $T$ and $\lambda$.

For $\lambda=0$, equations \eqref{ci:forced-energy-decay}, \eqref{eq:gaussian-energy} and \eqref{eq:incoming-norm}, together with the integrability of the exponential kernel, give
\begin{align}
\tnorm{f}_{L^\infty(0,T;L^2)}\lesssim\gnorm{\infty}{\alpha}{f_0}+\tnorm{g}_{L^\infty(0,T;L^{\infty,\alpha}_{\Sigma_-})}. \label{ci:Gaussian-energy-input}
\end{align}
Substitution in \eqref{ci:Gaussian-from-energy} proves \eqref{ci:uniform-Gaussian}.

\proofstep{Exponentially decaying in-flow} Let $\lambda_*>0$ be the decay rate in \eqref{ci:forced-energy-decay}, and fix $\lambda:=\frac12\min\{\frac12,\lambda_*,\omega\}$.
After multiplying \eqref{ci:forced-energy-decay} by $\ue^{2\lambda t}$, the initial term is bounded by $\norm{f_0}_{L^2}^2$ because $\lambda<\lambda_*$. For the boundary term, \eqref{eq:incoming-norm} and $\bgnorm{\infty}{\alpha}{-}{g(s)}\leq\eta_\omega\ue^{-\omega s}$ give
\begin{align}
\ue^{2\lambda t}\int_0^t\ue^{-2\lambda_*(t-s)}\norm[\big]{g(s)}_{L^2_{\Sigma_-}}^2\dd s&\lesssim\eta_\omega^2\int_0^t\ue^{-2(\lambda_*-\lambda)(t-s)}\ue^{-2(\omega-\lambda)s}\dd s \label{ci:Gaussian-weighted-boundary}\\
&\leq\frac{\eta_\omega^2}{2(\lambda_*-\lambda)}.\notag
\end{align}
Taking square roots and using \eqref{eq:gaussian-energy} for the initial value yields
\begin{align}
\tnorm[\Big]{\ue^{\lambda t}f}_{L^\infty(0,T;L^2)}\lesssim_\omega\gnorm{\infty}{\alpha}{f_0}+\eta_\omega. \label{ci:Gaussian-weighted-energy-input}
\end{align}
Apply \eqref{ci:Gaussian-from-energy} at this $\lambda$ to obtain \eqref{ci:decaying-Gaussian}. Lemma~\ref{trace:gaussian} and $C_tL^2$ continuity give the bound at every $t$.

\proofstep{In-flow tending to zero} First, \eqref{ci:forced-energy-decay} implies $\norm{f(t)}_{L^2}\to0$. To check this directly, split its boundary integral at a fixed $s>0$. The integral over $[0,s]$ tends to zero as $t\to\infty$, and the remaining integral is at most $C\tnorm{g}_{L^\infty(s,\infty;L^{\infty,\alpha}_{\Sigma_-})}^2$. Let $t\to\infty$ and then $s\to\infty$.

The decreasing time intervals and \eqref{eq:I:display:207} give the finite limit
\begin{align}
m_\infty:=\lim_{s\to\infty}\tnorm{f}_{L^\infty(s,\infty;L^{\infty,\alpha})}. \label{ci:Gaussian-limit-amplitude}
\end{align}
Apply \eqref{ci:bounded-velocity-late} after the time translation $t\mapsto s+t$, with $R+1$ in place of $R$, and then let the terminal time tend to infinity. Letting $s\to\infty$ and using the vanishing incoming and $L^2$ norms gives
\begin{align}
\lim_{s\to\infty}\tnorm{f}_{L^\infty((s,\infty)\times\Om\times B_{R+1})}\leq C_Rr^2m_\infty. \label{ci:Gaussian-limit-local}
\end{align}

Next apply \eqref{ci:exterior-Gaussian-estimate}, with $\lambda=0$, to the same translated solution and let the terminal time tend to infinity. Its left side is the exterior Gaussian norm on $(s+\tau,\infty)$. For fixed $s$, its initial contribution $\ue^{-\tau/2}\gnorm{\infty}{\alpha}{f(s)}$ tends to zero as $\tau\to\infty$. The remaining terms are bounded by norms on $(s,\infty)$. Letting $\tau\to\infty$, then $s\to\infty$, and including the region $B_{R+1}$ gives
\begin{align}
m_\infty\leq C\ue^{\alpha(R+1)^2}\lim_{s\to\infty}\tnorm{f}_{L^\infty((s,\infty)\times\Om\times B_{R+1})}+\varepsilon_Rm_\infty. \label{ci:Gaussian-limit-exterior}
\end{align}
Choose $R$ large enough that $\varepsilon_R\leq\frac{1}{4}$, and then choose $r$ small enough that the product of the constants multiplying $r^2m_\infty$ in \eqref{ci:Gaussian-limit-local}--\eqref{ci:Gaussian-limit-exterior} is at most $\frac{1}{4}$. Thus $m_\infty\leq\frac{1}{2}m_\infty$, so $m_\infty=0$. Lemma~\ref{trace:gaussian} transfers each time-tail bound to every time in that interval by $C_tL^2$ continuity. This proves \eqref{ci:vanishing-inflow-Gaussian}.
\end{proof}

\subsection{Global Continuation under In-flow}\label{sec:energy-continuation}

We will continue the local solution by combining the uniform existence interval from Proposition~\ref{ci:local-existence} with the bounds in Propositions~\ref{ci:energy-response} and \ref{ci:gaussian-bootstrap}. These bounds keep the data small enough to restart the local construction.

\begin{proposition}[Global continuation]\label{ci:global-energy}
Let $0<\alpha<\frac{1}{4}$, and suppose that $\m+\mh f_0\geq0$ and $\m+\mh g\geq0$. There is $\varepsilon_*>0$ such that
\begin{align}
\gnorm{\infty}{\alpha}{f_0}+\tnorm{g}_{L^\infty(0,\infty;L^{\infty,\alpha}_{\Sigma_-})}\leq\varepsilon_* \label{eq:I:display:214}
\end{align}
implies the existence of a global weak solution $f$ in the sense of Definition~\ref{def:energy}, with initial value $f_0$, in-flow value $\trm f=g$, and $\m+\mh f\geq0$. It satisfies \eqref{ci:uniform-Gaussian}, \eqref{ci:forced-energy-decay}, and \eqref{ci:integrated-energy}. If the incoming perturbation satisfies the exponential bound in Proposition~\ref{ci:gaussian-bootstrap}, then \eqref{ci:decaying-Gaussian} also holds. The constants in these estimates are independent of the final time.
\end{proposition}

\begin{proof}
Choose $0<\delta\leq\delta_0$ below the thresholds of Propositions~\ref{ci:energy-response} and \ref{ci:gaussian-bootstrap}, where $\delta_0,C_0$ are from Proposition~\ref{ci:local-existence}. Let $C\geq2$ bound the sum of the right side of \eqref{ci:uniform-Gaussian} and the incoming Gaussian norm by $C\varepsilon_*$. Fix the local lifespan and the data threshold by
\begin{align}
h:=\frac{\ln\left(\frac43\right)}{C_0},\qquad C\varepsilon_*<\frac{\delta}{4}. \label{ci:global-smallness}
\end{align}
Proposition~\ref{ci:local-existence} gives a solution on $[0,h]$ with nonnegative density and Gaussian norm at most $\delta$. Suppose that it has been constructed on $[0,t_0]$, where $t_0=nh$, with this bound. Apply \eqref{ci:uniform-Gaussian} on the entire interval $[0,t_0]$. Lemma~\ref{trace:gaussian} and $L^2$ continuity give the same bound at $t_0$, so
\begin{align}
\gnorm[\big]{\infty}{\alpha}{f(t_0)}+\tnorm{g}_{L^\infty(t_0,t_0+h;L^{\infty,\alpha}_{\Sigma_-})}\leq C\varepsilon_*<\frac{\delta}{4}. \label{ci:restart-data}
\end{align}
The nonnegativity of $\m+\mh f(t_0)$ also follows by $L^2$ continuity. The time-translated local existence theorem supplies a solution on $[t_0,t_0+h]$ with initial value $f(t_0)$ and Gaussian norm at most $\delta/2$. Their $L^2$ values agree at $t_0$, so the weak formulations concatenate to give a solution on $[0,t_0+h]$ with nonnegative density and Gaussian norm at most $\delta$.

Induction constructs a single global weak solution. Propositions~\ref{ci:energy-response} and \ref{ci:gaussian-bootstrap}, applied on each entire constructed interval, give all the stated estimates with constants independent of its length.
\end{proof}

\subsection{Global Weak Existence under Diffuse-Reflection}\label{cex:section}

We will apply the global weak existence theorem of Carrapatoso--Mischler \cite[Theorem~1.1]{Carrapatoso.Mischler2024} and express its Gaussian and weighted $L^2$ estimates in the normalization $F=\m+\mh f$. Section~\ref{sec:trace} identifies the normalized traces, and \eqref{ci:negative-part-energy} gives $F\geq0$ using the diffuse boundary contraction. Corollary~\ref{cex:decay-traces} then gives decay and the global time integrals.

\begin{proposition}[Global weak existence under diffuse-reflection]\label{cex:CM-input}
There is $\varepsilon_*>0$ such that, if
\begin{align}
\m+\mh f_0\geq0,\qquad \int_{\phase}\mh f_0\dd x\dd v=0,\qquad \varepsilon:=\gnorm{\infty}{\alpha}{f_0}\leq\varepsilon_*, \label{cex:initial-data}
\end{align}
then the diffuse-reflection problem has a global weak solution $f$ in the sense of Definition~\ref{def:energy}. It satisfies $F=\m+\mh f\geq0$ and $\int_{\phase}\mh f(t)\dd x\dd v=0$ for every $t\geq0$.

The solution satisfies
\begin{align}
\tnorm{f}_{L^\infty(0,\infty;L^{\infty,\alpha})}\lesssim\varepsilon. \label{cex:CM-bounds}
\end{align}
\end{proposition}

\begin{proof}
\proofstep{The cited solution and Gaussian bound} We use the Gaussian case of \cite[Theorem~1.1]{Carrapatoso.Mischler2024}. For diffuse-reflection, under our assumptions on $\Om$ and $\gamma$, its quantitative form gives a fixed $c>0$ such that, for every sufficiently small $\eta>0$,
\begin{align}
\int_{\phase}(F_0-\m)\dd x\dd v=0,\qquad \gnorm{\infty}{\alpha+\frac14}{F_0-\m}\leq c\eta \label{cex:CM-hypotheses}
\end{align}
imply the existence of a global weak solution of \eqref{eq:full} and \eqref{eq:diffuse-full}, with initial data $F_0$, conservation of mass, and
\begin{align}
\tnorm{F-\m}_{L^\infty(0,\infty;L^{\infty,\alpha+\frac14})}\leq\eta. \label{cex:CM-threshold-bound}
\end{align}
The Gaussian exponent must lie in $(0,\frac12)$, which holds because $0<\alpha<\frac14$. Our $C^3$ boundary supplies the required $W^{2,\infty}(\R^3)$ extension of the collar distance, as described in Section~\ref{lin:boundary-change}.

We apply this statement with $F_0=\m+\mh f_0$. The mass condition is \eqref{cex:initial-data}, and \eqref{eq:Gaussian-normalization} gives
\begin{align}
\gnorm{\infty}{\alpha+\frac14}{F_0-\m}=(2\pi)^{-\frac34}\varepsilon\leq c\eta \label{cex:input-weight}
\end{align}
for sufficiently small $\varepsilon_*$. This verifies \eqref{cex:CM-hypotheses} and gives the cited solution $F=\m+\mh f$.

\proofstep{Estimates and trace statements supplied by the cited results} We collect the conclusions for the fixed point constructed in \cite[Section~8.1]{Carrapatoso.Mischler2024}, writing them throughout for $f=\mih(F-\m)$. The estimates require
\begin{align}
\tnorm[\Big]{\la^{10}\mh f}_{L^\infty((0,\infty)\times\phase)}\leq\eta_1,\qquad f_0\in L^2,\qquad\int_{\phase}\mh f_0\dd x\dd v=0.\label{cex:CM-required}
\end{align}
Here $\eta_1>0$ is chosen below the thresholds in \cite[Theorem~3.4, Proposition~6.2 and (7.1)]{Carrapatoso.Mischler2024} for the two Gaussian weights used below. The exponent $10$ satisfies its requirement $10>8+\gamma$.

First, \cite[(7.2)]{Carrapatoso.Mischler2024}, applied at this fixed point, gives the proportional Gaussian estimate
\begin{align}
\tnorm{f}_{L^\infty(0,\infty;L^{\infty,\alpha})}\lesssim\gnorm{\infty}{\alpha}{f_0}. \label{cex:CM-proportional-bound}
\end{align}
Next, \cite[Theorem~3.4 and Proposition~6.2, (6.7)]{Carrapatoso.Mischler2024} give, for every $T<\infty$,
\begin{align}
f\in C\bigl([0,T];L^2\bigr)\cap L^2(0,T;\Vspace),\qquad f(0)=f_0\ \text{in }L^2, \label{cex:normalized-energy}\\
\tnorm{f}_{L^\infty(0,T;L^2)}^2+\int_0^T\norm[\big]{f(t)}_{\Vspace}^2\dd t\lesssim\norm{f_0}_{L^2}^2, \label{cex:global-energy}
\end{align}
with a constant independent of $T$. These are the cited weighted estimates with weight $\mih$ on the density perturbation $\mh f$: the weighted $L^2$ norm is exactly $\norm f_{L^2}$, and the velocity dissipation norm is comparable to $\norm f_{\Vspace}$ by \eqref{eq:V} and \eqref{cvel:eigenvalues}. Theorem~3.4 also supplies the diffuse boundary condition for its kinetic trace.

To identify this trace in our formulation, we use Lemma~\ref{trace:CM}(i), (iii) and Lemma~\ref{trace:gaussian}(i). Besides the bulk regularity \eqref{cex:normalized-energy} and the Gaussian bound, they require the weak equation \eqref{eq:intro-interior-weak} with $S=\Ga(f,f)\in L^2(0,T;\Vspace')$. Their conclusions are uniqueness of the trace determined by \eqref{trace:compact-Green}, the density identity \eqref{eq:intro-density-trace}, and
\begin{align}
\tnorm{\trpm f}_{L^\infty(0,T;L^{\infty,\alpha}_{\Sigma_\pm})}\leq\tnorm{f}_{L^\infty(0,T;L^{\infty,\alpha})},\qquad\trpm f\in L^2\bigl((0,T)\times\Sigma_\pm;\dd t\dd\gamma\bigr). \label{cex:finite-trace}
\end{align}
The last integrability makes the diffuse boundary integral well-defined.

For nonnegativity, we will use \eqref{ci:negative-part-energy} with $G=F$ and $\varepsilon=0$. Its hypotheses are the regularity \eqref{ci:negative-part-class}, the small Gaussian coefficient bound \eqref{ci:map-coefficient} with $\psi=f$, and \eqref{ci:normalized-supersolution} with a locally $L^2$ residual. The initial and incoming terms are retained in that estimate.

\proofstep{Verification of the Gaussian and $L^2$ hypotheses} Choose the fixed point in the ball defined in \cite[Section~8.1]{Carrapatoso.Mischler2024}. Its defining bound is the first condition in \eqref{cex:CM-required}; the construction applies after decreasing $\varepsilon_*$ in \eqref{cex:initial-data}. The mass condition in \eqref{cex:CM-required} is already part of \eqref{cex:initial-data}. By \eqref{eq:gaussian-energy},
\begin{align}
\norm{f_0}_{L^2}^2\lesssim\gnorm{\infty}{\alpha}{f_0}^2=\varepsilon^2. \label{eq:I:display:235}
\end{align}
Thus all hypotheses in \eqref{cex:CM-required} hold. Estimate \eqref{cex:CM-proportional-bound} proves \eqref{cex:CM-bounds}, and \eqref{cex:normalized-energy}--\eqref{cex:global-energy} give the stated time continuity and velocity regularity, with the right-hand side of \eqref{cex:global-energy} bounded by $C\varepsilon^2$. We choose $\varepsilon_*$ small enough that \eqref{cex:CM-bounds} also gives \eqref{ci:map-coefficient} with $\psi=f$.

\proofstep{Verification of the trace hypotheses and boundary condition} Multiplying the cited equation by $\mih$ gives \eqref{eq:intro-interior-weak}. The collision estimate \eqref{trace:dual-bounds}, the Gaussian bound \eqref{cex:CM-bounds} and \eqref{cex:global-energy} imply
\begin{align}
\tnorm{\Ga(f,f)}_{L^2(0,T;\Vspace')}\lesssim\tnorm{f}_{L^\infty(0,T;L_x^\infty L_v^2)}\tnorm{f}_{L^2(0,T;\Vspace)}\lesssim\varepsilon\tnorm{f}_{L^2(0,T;\Vspace)}<\infty.\label{cex:trace-source}
\end{align}
Together with \eqref{cex:normalized-energy}, \eqref{eq:intro-interior-weak} and \eqref{cex:CM-bounds}, this verifies the trace hypotheses collected above and gives \eqref{cex:finite-trace}, with Gaussian norm bounded by $C\varepsilon$.

Testing the cited weak formulation for $F-\m$ with $\mih\varphi$ gives \eqref{trace:compact-Green} for its normalized boundary trace. Uniqueness in Lemma~\ref{trace:CM}(i) therefore identifies it with our trace. Consequently \eqref{eq:intro-density-trace} and the wall normalization \eqref{eq:wall-normalization} convert its diffuse boundary condition into
\begin{align}
\trm f=\Pm\trp f\quad\text{in }L^2\bigl((0,T)\times\Sigma_-;\dd t\dd\gamma\bigr).\label{cex:normalized-boundary}
\end{align}
The initial value in \eqref{cex:normalized-energy}, the weak equation and these trace conclusions verify Definition~\ref{def:energy}.

\proofstep{Verification of the negative-part estimate} Set $q:=\mih F=\mh+f$. Since $\mh\in\Vspace$ and its boundary traces belong to $L^2_{\Sigma_\pm}$, \eqref{cex:normalized-energy} and \eqref{cex:finite-trace} give \eqref{ci:negative-part-class}. The density equation and \eqref{cex:initial-data} give
\begin{align}
\Y q-\mih Q(F,\mh q)=0,\qquad q(0)=\mh+f_0\geq0.\label{cex:negative-part-equation}
\end{align}
Thus \eqref{ci:normalized-supersolution} has zero residual, which is locally $L^2$. The coefficient bound required by \eqref{ci:negative-part-energy} was verified in the preceding Gaussian estimate. All hypotheses of that estimate now hold, and $q_-(0)=0$.

By \eqref{eq:wall-normalization}, $\Pm\mh=\mh$. Hence \eqref{cex:normalized-boundary} gives $\trm q=\Pm\trp q$. Positivity of $\Pm$ and \eqref{trace:projection-flux-bound} imply
\begin{align}
\norm[\Big]{(\trm q)_-}_{L^2_{\Sigma_-}}^2\leq\norm[\Big]{\Pm((\trp q)_-)}_{L^2_{\Sigma_-}}^2\leq\norm[\Big]{(\trp q)_-}_{L^2_{\Sigma_+}}^2. \label{cex:negative-boundary-contraction}
\end{align}
Apply \eqref{ci:negative-part-energy}. Its initial term vanishes, and \eqref{cex:negative-boundary-contraction} absorbs the incoming term into the outgoing term. Dropping the nonnegative diffusion integral yields
\begin{align}
\norm{q_-(t)}_{L^2}^2\lesssim\int_0^t\norm{q_-(s)}_{L^2}^2\dd s. \label{cex:negative-part-Gronwall}
\end{align}
Gr\"onwall's inequality gives $q_-(t)=0$, hence $F(t)\geq0$ for every $t\geq0$.
\end{proof}

\begin{corollary}[Decay and global integral estimates]\label{cex:decay-traces}
Let $f$ be the solution of Proposition~\ref{cex:CM-input}, and let $p_\gamma$ be as in \eqref{full:decay-exponent}. For every $0<\alpha'<\alpha$, there is $\lambda_{\alpha'}>0$ such that, uniformly in time,
\begin{align}
\gnorm[\big]{\infty}{\alpha'}{f(t)}+\norm[\big]{f(t)}_{L^2}&\lesssim_{\alpha'}\varepsilon\ue^{-\lambda_{\alpha'}t^{p_\gamma}}, \label{cex:weaker-Gaussian}\\
\bgnorm[\big]{\infty}{\alpha'}{\pm}{\trpm f(t)}&\lesssim_{\alpha'}\varepsilon\ue^{-\lambda_{\alpha'}t^{p_\gamma}}\qquad\text{for almost every }t>0. \label{cex:trace-Gaussian-bound}
\end{align}
Moreover,
\begin{align}
\int_0^\infty\left(\norm[\big]{f(t)}_{L^2}^2+\norm[\big]{f(t)}_{\Vspace}^2+\sum_{\pm}\norm[\big]{\trpm f(t)}_{L^2_{\Sigma_\pm}}^2+\norm[\big]{\flux{f}(t)}_{L^2(\partial\Om)}^2\right)\dd t\lesssim_{\alpha'}\varepsilon^2. \label{cex:global-integrals}
\end{align}
If $-2\leq\gamma\leq1$, there is $\lambda_D>0$ such that, uniformly in time,
\begin{align}
\gnorm[\big]{\infty}{\alpha}{f(t)}+\norm[\big]{f(t)}_{L^2}\lesssim\varepsilon\ue^{-\lambda_Dt}. \label{cex:same-weight-decay}
\end{align}
\end{corollary}

\begin{proof}
\proofstep{Decay of the density perturbation} Proposition~\ref{cex:CM-input} verifies \eqref{cex:CM-hypotheses}. The decay estimate \cite[(1.25)]{Carrapatoso.Mischler2024} controls $F-\m$ with the same Gaussian weight $\ue^{\left(\alpha+\frac14\right)\abs{v}^2}$ when $\gamma\geq-2$, and with the polynomial weight $\la^{10}$ when $\gamma<-2$. The exponent $10$ satisfies the condition $10>8+\gamma$ in that reference, and the time exponent is $p_\gamma$ in \eqref{full:decay-exponent}. Both velocity weights are at least one. Thus \eqref{cex:input-weight} gives $\lambda_D>0$ such that
\begin{align}
\norm[\big]{\mh f(t)}_{L^\infty}\lesssim\varepsilon\ue^{-\lambda_Dt^{p_\gamma}}. \label{cex:density-decay}
\end{align}
\proofstep{Gaussian and $L^2$ decay} For $0<\alpha'<\alpha$, we have
\begin{align}
\ue^{\alpha'\abs{v}^2}\abs{f}=\Bigl(\ue^{\alpha\abs{v}^2}\abs{f}\Bigr)^{\frac{4\alpha'+1}{4\alpha+1}}\Bigl(\ue^{-\frac{\abs{v}^2}{4}}\abs{f}\Bigr)^{\frac{4(\alpha-\alpha')}{4\alpha+1}}. \label{eq:I:display:241}
\end{align}
Since $\ue^{-\frac{\abs{v}^2}{4}}\abs{f}=(2\pi)^{\frac34}\abs{\mh f}$, estimates \eqref{cex:CM-bounds} and \eqref{cex:density-decay} give the Gaussian term of \eqref{cex:weaker-Gaussian}, with $\lambda_{\alpha'}:=\frac{4(\alpha-\alpha')}{4\alpha+1}\lambda_D>0$. The $L^2$ term follows from \eqref{eq:gaussian-energy}.

\proofstep{Trace decay and global integrals} On every interval $(a,b)$ with rational endpoints $0<a<b$, estimate \eqref{cex:weaker-Gaussian} bounds the bulk Gaussian norm by a constant times $\varepsilon\ue^{-\lambda_{\alpha'}a^{p_\gamma}}$. Lemma~\ref{trace:gaussian}(i) gives the same bound for the trace on $(a,b)$. Letting $a\uparrow t$ through rational endpoints gives \eqref{cex:trace-Gaussian-bound} for almost every $t>0$.

Formula \eqref{trace:Gaussian-flux-integral} and the Cauchy--Schwarz estimate for \eqref{eq:diffuse-flux} therefore give
\begin{align}
\sum_{\pm}\norm[\big]{\trpm f(t)}_{L^2_{\Sigma_\pm}}^2+\norm[\big]{\flux{f}(t)}_{L^2(\partial\Om)}^2\lesssim_{\alpha'}\varepsilon^2\ue^{-2\lambda_{\alpha'}t^{p_\gamma}}. \label{cex:flux-decay}
\end{align}
Since $p_\gamma>0$, the function $\ue^{-2\lambda_{\alpha'}t^{p_\gamma}}$ belongs to $L^1(0,\infty)$. Squaring \eqref{cex:weaker-Gaussian} and integrating it together with \eqref{cex:flux-decay} gives the $L^2$ and boundary terms of \eqref{cex:global-integrals}. Letting $T\to\infty$ in \eqref{cex:global-energy} gives its $\Vspace$ term.

\proofstep{Decay with the initial Gaussian weight} If $\gamma\geq-2$, estimate \cite[(1.25)]{Carrapatoso.Mischler2024} retains the initial Gaussian weight and $p_\gamma=1$. The normalization \eqref{eq:Gaussian-normalization} and \eqref{eq:gaussian-energy} therefore give \eqref{cex:same-weight-decay}.
\end{proof}


\section{Local Regularity Estimates for Linear Kinetic Equations}\label{lin:section}

We will derive local regularity estimates for the weak solutions constructed in Section~\ref{sec:energy}. We will first flatten the boundary and normalize the velocity diffusion, then establish interior estimates and boundary H\"older and first-order derivative estimates for a more general linear equation. For second-order estimates at incoming points, we will subtract a boundary lift. In Section~\ref{cb:section}, we will verify the hypotheses for in-flow and diffuse-reflection and return to the original variables.

\subsection{Landau Equation and Boundary Reduction}\label{cb:section-normalization}

We will first write the Landau equation in divergence and nondivergence form using Section~\ref{sec:cvel-coefficients}, then derive the transformed weak equation and coefficient bounds under the flattening map.

\subsubsection{The Landau equation and its coefficients}\label{lin:Landau-rewriting}

Let $I\subset\R$ be an open interval and let $J\subset\subset I$ be a compact subinterval. Suppose that, on every such interval,
\begin{align}
f\in C(J;L^2)\cap L^2(J;\Vspace),\qquad \Y f\in L^2(J;\Vspace'). \label{lin:bulk-hypotheses}
\end{align}
Assume also that
\begin{align}
F:=\m+\mh f,\qquad \tnorm[\Big]{\ue^{\alpha\abs{v}^2}f}_{L^\infty(I\times\Om\times\R^3)}\leq\delta,\qquad 0<\alpha<\frac{1}{4},\quad 0<\delta\leq\delta_A, \label{cb:setting}
\end{align}
where $\delta_A$ is the constant in Lemma~\ref{cvel:coefficients}. We allow the additional source $\sourceR$ in
\begin{align}
\Y f+\Lop f=\Ga(f,f)+\sourceR. \label{cb:equation-residual}
\end{align}
Equation~\eqref{cb:equation-residual} means the weak identity \eqref{eq:intro-interior-weak} with $\int\sourceR\varphi\dd t\dd x\dd v$ added to the right-hand side. We assume that, for every $J\subset\subset I$, there are $M_J>0$ and $\alpha_J>0$ such that
\begin{align}
\tnorm[\Big]{\ue^{\alpha_J\abs{v}^2}\sourceR}_{L^\infty(J\times\Om\times\R^3)}\leq M_J. \label{lin:source-Gaussian}
\end{align}
For the in-flow first-order estimate, \eqref{lin:source-Gaussian} suffices. The second-order and improved diffuse-flux estimates will additionally assume
\begin{align}
\ue^{\alpha_J\abs{v_0}^2}[\sourceR]_{\Ckin^\beta\left(J\times\overline\Om\times B_{\bracket{v_0}^{-1}}(v_0)\right)}\leq M_J\qquad(v_0\in\R^3), \label{lin:source-holder}
\end{align}
for some $\beta\in(0,1)$ and $M_J,\alpha_J>0$ independent of $v_0$.

The source $\sourceR$ allows the negative-time extensions used to obtain estimates at $t=0$; Lemmas~\ref{cb:inflow-initial} and \ref{cb:initial} will verify \eqref{lin:source-Gaussian}--\eqref{lin:source-holder} for these extensions. Boundary lifts contribute the separate terms in \eqref{lin:Landau-source}. The kinetic trace is that of Lemma~\ref{trace:CM}; Section~\ref{cb:section} will impose the boundary law.

The coefficients $\Af{},b_F,c_F,\Bf{},\Cf{}$ are defined in \eqref{cvel:coefficient-convolutions}--\eqref{cvel:scalar-coefficients}, and $\Kop$ in \eqref{eq:I:display:052}. Assumptions \eqref{lin:bulk-hypotheses} and \eqref{cb:setting} verify the hypotheses of Lemma~\ref{cvel:scalar} for almost every $(t,x)$. Its identity \eqref{cvel:scalar-div} rewrites \eqref{cb:equation-residual} as
\begin{align}
\Y f-\divv_v(\Af{}\nabla_vf)=\Bf{}\cdot\nabla_vf+\Cf{}f+\Kop f+\sourceR. \label{cb:coefficients}
\end{align}
The equivalent identity \eqref{cvel:scalar-nondiv} gives
\begin{align}
\Y f-\Af{}:D_v^2f=-\bigl(\Af{}v\bigr)\cdot\nabla_vf+\Cf{}f+\Kop f+\sourceR. \label{cb:scalar}
\end{align}
Both equations hold in distributions, with the contraction in \eqref{cb:scalar} understood by \eqref{cvel:weak-contraction}. We use the divergence form for energy and boundary H\"older estimates, and the nondivergence form for Schauder and second-order estimates. The bounds \eqref{cvel:rough-coefficients}, \eqref{cvel:scalar-source} and \eqref{lin:source-Gaussian} make all their terms well defined under \eqref{lin:bulk-hypotheses}.

By Lemma~\ref{cvel:coefficients}, $\frac{1}{2}\sig\leq\Af{}\leq\frac{3}{2}\sig$, with $\sig$ defined in \eqref{eq:sigma-definition}. The following change of variables makes the diffusion uniformly elliptic near each fixed velocity $v_0$.

\subsubsection{Flattening map}\label{lin:boundary-change}

For a fixed point $Z_0=(t_0,x_0,v_0)$ and a variable $\zeta=(\theta,\xi,\eta)$, define the kinetic translation and dilation by
\begin{align}
T_{Z_0}(\zeta):=Z_0\circ\zeta&:=(t_0+\theta,x_0+\xi+\theta v_0,v_0+\eta), \label{eq:I:display:108}\\
S_r(\zeta)&:=(r^2\theta,r^3\xi,r\eta),\qquad r>0. \label{eq:intro-kinetic-scaling}
\end{align}
We use the kinetic cylinders of \cite[Section~1.3.2]{Zhu2024}, with two-sided time intervals as in \cite[Section~2.1]{Ros-Oton.Weidner2025}. For $\Qc_1(0):=(-1,1)\times B_1(0)\times B_1(0)$, set
\begin{align}
\Qc_r(Z_0)&:=T_{Z_0}\Bigl(S_r\bigl(\Qc_1(0)\bigr)\Bigr) =\Bigl\{(t,x,v):\abs{t-t_0}<r^2,\ \abs[\big]{x-x_0-(t-t_0)v_0}<r^3,\ \abs{v-v_0}<r\Bigr\}. \label{lin:cylinders}
\end{align}
Definition~\eqref{eq:kinetic-distance} gives
\begin{align}
\{Z:\dkin(Z,Z_0)<r\}&\subset \Qc_r(Z_0)\subset\{Z:\dkin(Z,Z_0)<3r\}. \label{lin:cylinder-distance}
\end{align}
For the variables $\zeta=(\theta,\xi,\eta)$ and a center $\zeta_0$, we use \eqref{lin:cylinders} with $(t,x,v)$ replaced by $(\theta,\xi,\eta)$. Writing $e_3:=(0,0,1)$, we set
\begin{align}
\Qplus_r(\zeta_0)&:=\Qc_r(\zeta_0)\cap\{\xi_3>0\}. \label{lin:half-cylinders}
\end{align}
The incoming, outgoing and grazing sets, with time unrestricted, are
\begin{align}
\gamma_-&:=\{\xi_3=0,\ \eta_3>0\},\qquad \gamma_+:=\{\xi_3=0,\ \eta_3<0\},\qquad \gamma_0:=\{\xi_3=\eta_3=0\}. \label{lin:boundary-sets}
\end{align}

The joint position-velocity change in \cite[Sections~7.1.1 and 7.1.3]{Guo.Hwang.Jang.Ouyang2020} motivates the following boundary-preserving map. The compact $C^3$ boundary admits finitely many $C^3$ graph neighborhoods and a tubular collar with $C^3$ signed distance. Multiplying the signed distance by a smooth cutoff supported in the collar and equal to one near $\partial\Om$ gives a $W^{2,\infty}(\R^3)$ function agreeing with $d$ on a smaller interior collar.

Fix $J\subset\subset I$ and $Z_0=(t_0,x_0,v_0)$ with $t_0\in J$, $x_0\in\overline\Om$ and $d(x_0)<\frac{d_0}{2}$. After a translation and rotation near $\bdproj(x_0)$, with the same rotation applied to $v$, we use a graph over the first two spatial variables and the signed extension of $d$ to define
\begin{align}
\phi(x)&:=\bigl(x_1,x_2,d(x)\bigr),\qquad y:=\phi(x),\qquad w:=D\phi(x)v. \label{lin:intermediate-map}
\end{align}
On the interior collar, $\nabla d(x)=-n(\bdproj(x))$, so
\begin{align}
y_3&=d(x),\qquad w_3=-n\bigl(\bdproj(x)\bigr)\cdot v. \label{lin:intermediate-normal}
\end{align}
Set $y_{\partial}:=\phi(x_0)-d(x_0)e_3$. We choose an upper triangular matrix $S_0$ with positive diagonal so that
\begin{align}
S_0S_0^T&=D\phi(x_0)\sig(v_0)D\phi(x_0)^T,\qquad e_3^TS_0=s_ne_3^T,\qquad s_n>0. \label{eq:I:display:262}
\end{align}
Such a factor is obtained by applying Cholesky factorization in reversed coordinate order to the positive-definite matrix on the right of \eqref{eq:I:display:262}; its upper triangular form gives the stated third row. The eigenvalue bounds \eqref{cvel:eigenvalues} and the uniform bounds for $D\phi$ and its inverse give
\begin{align}
\norm{S_0}&\lesssim\bracket{v_0}^{\frac{\gamma+2}{2}},\qquad \norm{S_0^{-1}}\lesssim\bracket{v_0}^{-\frac{\gamma}{2}},\qquad \bracket{v_0}^{\frac{\gamma}{2}}\lesssim s_n\lesssim\bracket{v_0}^{\frac{\gamma+2}{2}}. \label{cb:matrix-bounds}
\end{align}
For parameters $0<c_J\leq1$ and $N>0$, set $\vradius:=c_J\bracket{v_0}^{-N}$ and define
\begin{align}
t&=t_0+\vradius^2\theta,\qquad y=y_{\partial}+\vradius^3S_0\xi,\qquad w=\vradius S_0\eta. \label{lin:intermediate-scaling}
\end{align}
Combining \eqref{lin:intermediate-map} and \eqref{lin:intermediate-scaling}, the map from the final variables to the original variables is
\begin{align}
\Psi_0(\theta,\xi,\eta)&:=(t,x,v),\qquad t=t_0+\vradius^2\theta,\quad x=\phi^{-1}(y_{\partial}+\vradius^3S_0\xi),\quad v=\bigl(D\phi(x)\bigr)^{-1}\vradius S_0\eta. \label{lin:inverse-map}
\end{align}
Equivalently, the forward map is
\begin{align}
\Psi_0^{-1}(t,x,v)&=\left(\frac{t-t_0}{\vradius^2},\frac{1}{\vradius^3}S_0^{-1}\bigl(\phi(x)-y_{\partial}\bigr),\frac{1}{\vradius}S_0^{-1}D\phi(x)v\right). \label{cb:normalizing-map}
\end{align}
The image of the center is
\begin{align}
\widehat Z_0&:=(0,\xi_0,\eta_0)=\Psi_0^{-1}(Z_0)=\left(0,\frac{d(x_0)}{\vradius^3}S_0^{-1}e_3,\frac{1}{\vradius}S_0^{-1}D\phi(x_0)v_0\right). \label{lin:normalized-center}
\end{align}

\begin{figure}[htbp]
\centering
\begin{tikzpicture}[>=Latex]
\begin{scope}
\fill[blue!6] (-2,1.8)--(-2,0.45) .. controls (-0.8,-0.25) and (0.7,-0.15) .. (1.7,0.55)--(1.7,1.8)--cycle;
\draw[thick,blue!65!black] (-2,0.45) .. controls (-0.8,-0.25) and (0.7,-0.15) .. (1.7,0.55);
\node at (-0.3,1.45) {$\Om$};
\node at (-0.3,-0.55) {$\partial\Om$};
\fill (-0.3,0.65) circle (1.6pt);
\node[left] at (-0.3,0.65) {$x_0$};
\draw[->,thick] (-0.3,0.65)--(0.15,1.1) node[right] {$v$};
\node at (-0.2,-1.05) {$(t,x,v)$};
\end{scope}
\draw[->,very thick] (2,0.75)--(4.2,0.75) node[midway,above] {$\Psi_0^{-1}$};
\begin{scope}[shift={(6.4,0)}]
\fill[blue!6] (-1.7,0)--(1.7,0)--(1.7,1.8)--(-1.7,1.8)--cycle;
\draw[thick,blue!65!black] (-1.7,0)--(1.7,0);
\node at (0,1.45) {$\xi_3>0$};
\node at (0,-0.55) {$\xi_3=0$};
\fill (-0.4,0.65) circle (1.6pt);
\node[left] at (-0.4,0.65) {$\xi_0$};
\draw[->,thick] (-0.4,0.65)--(0.15,1.05) node[right] {$\eta$};
\node at (0,-1.05) {$(\theta,\xi,\eta)$};
\end{scope}
\end{tikzpicture}
\caption{The local flattening and scaling map $\Psi_0^{-1}$ in \eqref{cb:normalizing-map} sends the boundary to $\{\xi_3=0\}$.}\label{fig:boundary-change}
\end{figure}
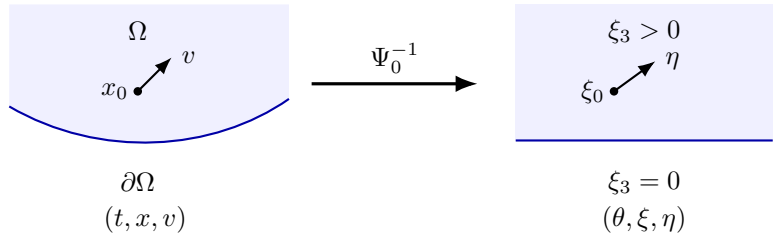

\begin{lemma}\label{cb:normalization}
The map $\phi$ in \eqref{lin:intermediate-map} is a local $C^3$ diffeomorphism with uniform bounds through order three. There exist $c_J,N>0$, independent of $Z_0$, such that $\Psi_0$ maps $\Qplus_2(\widehat Z_0)$ diffeomorphically onto an open subset of $I\times\Om\times\R^3$. It extends to $\Qc_2(\widehat Z_0)\cap\{\xi_3\geq0\}$, whose image is relatively open in $I\times\overline\Om\times\R^3$ and contains $Z_0$. At corresponding points,
\begin{align}
\xi_3&=\frac{d(x)}{\vradius^3s_n},\qquad \eta_3=-\frac{n\bigl(\bdproj(x)\bigr)\cdot v}{\vradius s_n}. \label{eq:I:display:279}
\end{align}
Within this image, the physical boundary corresponds to $\Qc_2(\widehat Z_0)\cap\{\xi_3=0\}$, and its incoming, outgoing and grazing parts correspond to $\Qc_2(\widehat Z_0)\cap\gamma_-$, $\Qc_2(\widehat Z_0)\cap\gamma_+$ and $\Qc_2(\widehat Z_0)\cap\gamma_0$, respectively. In particular, $(\xi_0)_3>0$ when $x_0\in\Om$.
\end{lemma}

\begin{proof}
After the rotation used in \eqref{lin:intermediate-map}, $\partial_{x_3}d$ is bounded away from zero on a smaller neighborhood. The inverse function theorem gives the regularity of $\phi$ and $\phi^{-1}$, and a finite boundary covering gives uniform bounds through order three.

For $\zeta\in\Qplus_2(\widehat Z_0)$, the cylinder definition gives $\abs{\theta}<4$, $\abs{\xi-\xi_0-\theta\eta_0}<8$ and $\abs{\eta-\eta_0}<2$. Substituting \eqref{lin:normalized-center} into \eqref{lin:inverse-map} and using \eqref{cb:matrix-bounds}, we obtain
\begin{align}
\abs{x-x_0}&\lesssim\vradius^2\bracket{v_0}+\vradius^3\bracket{v_0}^{\frac{\gamma+2}{2}}, \label{lin:position-variation}\\
\abs{v-v_0}&\lesssim\vradius\bracket{v_0}^{\frac{\gamma+2}{2}}+\bracket{v_0}\left(\vradius^2\bracket{v_0}+\vradius^3\bracket{v_0}^{\frac{\gamma+2}{2}}\right). \label{eq:I:display:265}
\end{align}
Choosing $N$ large and $c_J$ small keeps $x$ in the coordinate neighborhood and $t$ in $I$, and gives $\abs{v-v_0}\leq\bracket{v_0}^{-1}$. The local inverse in \eqref{cb:normalizing-map} then proves the stated neighborhood inclusion.

Since $e_3^TS_0^{-1}=s_n^{-1}e_3^T$ and $y_{\partial}\cdot e_3=0$, \eqref{cb:normalizing-map} and \eqref{lin:intermediate-normal} give \eqref{eq:I:display:279}. The signs in that formula prove all boundary correspondences.
\end{proof}

\para{Kinetic distance under $\Psi_0$}
For $\zeta,\zeta'\in\Qc_2(\widehat Z_0)\cap\{\xi_3\geq0\}$ and $Z=\Psi_0(\zeta)$, $Z'=\Psi_0(\zeta')$, we have, for a fixed $N_1$ independent of $N$,
\begin{align}
\dkin(Z,Z')&\lesssim\vradius\bracket{v_0}^{N_1}\dkin(\zeta,\zeta'), \label{lin:kinetic-change-bound}\\
\dkin(\zeta,\zeta')&\lesssim\vradius^{-1}\bracket{v_0}^{N_1}\dkin(Z,Z'). \label{lin:inverse-kinetic-change-bound}
\end{align}
We use the Taylor argument in the proof of \cite[Lemma~2.17(iii)]{Ros-Oton.Weidner2025}, keeping the factors from $S_0$ and $\vradius$ explicit. For the first inequality, Taylor's formula for $\phi^{-1}$ cancels the characteristic increment $(t-t')v'$ and gives
\begin{align}
\abs{x-x'-(t-t')v'}&\lesssim\bracket{v_0}^{N_1}\left(\vradius^3\bigl(\dkin(\zeta,\zeta')\bigr)^3+\vradius^4\bigl(\dkin(\zeta,\zeta')\bigr)^4\right), \label{lin:kinetic-position-increment}\\
\abs{v-v'}&\lesssim\bracket{v_0}^{N_1}\left(\vradius\dkin(\zeta,\zeta')+\vradius^2\bigl(\dkin(\zeta,\zeta')\bigr)^2\right). \label{lin:kinetic-velocity-increment}
\end{align}
Together with $\abs{t-t'}^{\frac12}=\vradius\abs{\theta-\theta'}^{\frac12}$, these bounds prove \eqref{lin:kinetic-change-bound}, since $\vradius\leq1$ and the kinetic diameter of $\Qc_2(\widehat Z_0)$ is bounded uniformly in its center.

For the inverse estimate, \eqref{cb:normalizing-map} and Taylor's formula give the explicit cancellation
\begin{align}
\xi-\xi'-(\theta-\theta')\eta'&=\vradius^{-3}S_0^{-1}\left[D\phi(x')\bigl(x-x'-(t-t')v'\bigr)+O(\abs{x-x'}^2)\right]. \label{lin:inverse-position-cancellation}
\end{align}
Here $\abs{x-x'}\lesssim\bigl(\dkin(Z,Z')\bigr)^3+\bracket{v_0}\bigl(\dkin(Z,Z')\bigr)^2$, while $\eta-\eta'=\vradius^{-1}S_0^{-1}\bigl(D\phi(x)v-D\phi(x')v'\bigr)$. Taking the cube root in \eqref{lin:inverse-position-cancellation} and using \eqref{cb:matrix-bounds}, the time increment, and the polynomial diameter bound from \eqref{lin:kinetic-change-bound} proves \eqref{lin:inverse-kinetic-change-bound}. All velocity powers come from \eqref{cb:matrix-bounds} and the fixed geometric bounds, so $N_1$ is independent of $N$.

For the Landau equation \eqref{cb:equation-residual}, define the transformed function, transport operator, coefficients and source by
\begin{align}
\widehat f&:=f\circ\Psi_0,\qquad \Yhat:=\partial_\theta+\eta\cdot\nabla_\xi, \label{cb:normalized-solution}\\
\locA&:=\bigl[S_0^{-1}D\phi(x)\Af{}D\phi(x)^TS_0^{-T}\bigr]\circ\Psi_0,\qquad \locc:=\vradius^2\Cf{}\circ\Psi_0, \label{lin:normalized-f-coefficients}\\
\locb&:=-\vradius\Bigl[S_0^{-1}\bigl(D\phi(x)(\Af{}v)+D^2\phi(x)[v,v]\bigr)\Bigr]\circ\Psi_0, \label{lin:Landau-coefficients}\\
\locB&:=\vradius\Bigl[S_0^{-1}\bigl(D\phi(x)\Bf{}-D^2\phi(x)[v,v]\bigr)\Bigr]\circ\Psi_0, \label{lin:normalized-drift}\\
\widehat S_f&:=\vradius^2\bigl(\Kop f+\sourceR\bigr)\circ\Psi_0. \label{lin:normalized-f-source}
\end{align}

\begin{lemma}\label{cb:normalized-equation}
Assume \eqref{lin:bulk-hypotheses}, \eqref{cb:setting}, \eqref{cb:equation-residual} and \eqref{lin:source-Gaussian}. On $\Qplus_2(\widehat Z_0)$, we have
\begin{align}
\Yhat\widehat f-\divv_\eta(\locA\nabla_\eta\widehat f)&=\locB\cdot\nabla_\eta\widehat f+\locc\widehat f+\widehat S_f. \label{lin:normalized-f-equation}
\end{align}
Precisely, $\widehat f,\nabla_\eta\widehat f\in L^2_{\mathrm{loc}}$, and every $\varphi\in C_c^\infty(\Qplus_2(\widehat Z_0))$ satisfies
\begin{align}
-\int\widehat f\Yhat\varphi+\int\locA\nabla_\eta\widehat f\cdot\nabla_\eta\varphi&=\int\bigl(\locB\cdot\nabla_\eta\widehat f+\locc\widehat f+\widehat S_f\bigr)\varphi. \label{lin:normalized-weak}
\end{align}
All integrals are with respect to $\dd\theta\dd\xi\dd\eta$. The kinetic trace is the transformed trace from Section~\ref{sec:trace}, and the corresponding boundary identity follows by changing variables in \eqref{trace:linear-Green}. In particular, if $\trm f=g$, then $\trm\widehat f=g\circ\Psi_0$ on the incoming part of the transformed region. The equivalent nondivergence form is
\begin{align}
\Yhat\widehat f-\locA:D_\eta^2\widehat f&=\locb\cdot\nabla_\eta\widehat f+\locc\widehat f+\widehat S_f, \label{lin:normalized-f-nondivergence}
\end{align}
where $\locA:D_\eta^2\widehat f:=\divv_\eta\bigl(\locA\nabla_\eta\widehat f\bigr)-\bigl(\divv_\eta\locA\bigr)\cdot\nabla_\eta\widehat f$ until second derivatives have been established.
\end{lemma}

\begin{proof}
\proofstep{The weak equation in $(t,y,w)$} Define
\begin{align}
f^\phi(t,y,w)&:=f\Bigl(t,\phi^{-1}(y),\bigl[D\phi(\phi^{-1}(y))\bigr]^{-1}w\Bigr). \label{cb:flattened-solution}
\end{align}
In the following formulas, $(y,w)=\bigl(\phi(x),D\phi(x)v\bigr)$. The chain rule gives
\begin{align}
\Y f&=(\partial_t+w\cdot\nabla_y)f^\phi+D^2\phi(x)[v,v]\cdot\nabla_wf^\phi, \label{cb:transformed-transport}\\
\divv_v(\Af{}\nabla_vf)&=\divv_w\left(D\phi(x)\Af{}D\phi(x)^T\nabla_wf^\phi\right). \label{lin:transformed-diffusion}
\end{align}
Thus the first-order coefficient on the right of the equation for $f^\phi$ is $D\phi(x)\Bf{}-D^2\phi(x)[v,v]$, and the other coefficients and scalar source are composed with the inverse map.

To verify these identities in the weak formulation, use the volume change $\dd x\dd v=j(y)\dd y\dd w$, where $j(y):=\abs[\big]{\det D\phi(\phi^{-1}(y))}^{-2}$. The transport field satisfies
\begin{align}
\divv_w\bigl(D^2\phi(x)[v,v]\bigr)&=-w\cdot\nabla_y\ln j(y). \label{lin:Jacobian-cancellation}
\end{align}
To verify this identity, differentiating in $w$ gives $2v\cdot\nabla_x\ln\abs{\det D\phi(x)}$, which is the negative of $w\cdot\nabla_y\ln j(y)$. Test the weak equation \eqref{cb:coefficients} with the composition of $\varphi/j$ with the forward change of variables. The derivative of $j$ in the transport term cancels the velocity divergence in \eqref{lin:Jacobian-cancellation}. This yields the unweighted weak equation with precisely the coefficients above. The local $L^2$ bounds on $f$ and $\nabla_vf$ justify the change of test functions by smooth approximation.

\proofstep{Scaling and trace} The function $\widehat f=f\circ\Psi_0$ equals $f^\phi$ evaluated at \eqref{lin:intermediate-scaling}. The two transport terms acquire the same factor,
\begin{align}
\Yhat\widehat f(\theta,\xi,\eta)&=\vradius^2\bigl((\partial_t+w\cdot\nabla_y)f^\phi\bigr)(t,y,w). \label{lin:normalized-transport}
\end{align}
The diffusion term has the same factor. Substituting into the intermediate equation gives \eqref{lin:normalized-f-equation} with the coefficients in \eqref{lin:normalized-f-coefficients}--\eqref{lin:normalized-f-source}. Moreover,
\begin{align}
\divv_\eta\locA&=\vradius\bigl[S_0^{-1}D\phi(x)b_F\bigr]\circ\Psi_0. \label{lin:normalized-diffusion-divergence}
\end{align}
Together with \eqref{cvel:scalar-coefficients}, this gives $\locB=\locb-\divv_\eta\locA$ and hence \eqref{lin:normalized-f-nondivergence}. The intermediate weak formulation and the constant linear scaling prove \eqref{lin:normalized-weak}.

Apply the same changes of variables and test functions to \eqref{trace:linear-Green}. The boundary factors are transformed by \eqref{eq:I:display:279}, so the resulting identity is the Green identity for \eqref{lin:normalized-f-equation}. Its trace is the composition of the original trace with $\Psi_0$; in particular, $\trm\widehat f=g\circ\Psi_0$ whenever $\trm f=g$.
\end{proof}

\begin{lemma}\label{cb:normalized-bounds}
Under the assumptions of Lemma~\ref{cb:normalized-equation}, there exist $c_J,N>0$, independent of $Z_0$, for which the preceding lemmas hold and the transformed coefficients satisfy
\begin{align}
C^{-1}I&\leq\locA\leq CI,\qquad \tnorm[\big]{\nabla_\eta\locA}_{L^\infty}+\tnorm[\big]{\locb}_{L^\infty}+\tnorm[\big]{\locB}_{L^\infty}+\tnorm[\big]{\locc}_{L^\infty}\leq C. \label{lin:normalized-coefficient-bounds}
\end{align}
Moreover, the Gaussian bounds in \eqref{cb:setting} and \eqref{lin:source-Gaussian} imply, for some $N_*,\alpha_*>0$,
\begin{align}
\tnorm[\big]{\widehat f}_{L^\infty}+\tnorm[\big]{\widehat S_f}_{L^\infty}&\lesssim_J\bracket{v_0}^{N_*}\ue^{-\alpha_*\abs{v_0}^2}. \label{cb:Gaussian-bound}
\end{align}
All transformed norms in this lemma are over $\Qplus_2(\widehat Z_0)$.

For the H\"older conclusions, fix $0<\beta<1$ and $m\geq0$. Suppose that physical neighborhoods $\mathcal U_{Z_0}$, containing the images in Lemma~\ref{cb:normalization} and fixed independently of the final choice of $\vradius$, satisfy
\begin{align}
\tnorm{\Af{}}_{\Ckin^\beta(\mathcal U_{Z_0})}+\tnorm{\Bf{}}_{\Ckin^\beta(\mathcal U_{Z_0})}+\tnorm{\Cf{}}_{\Ckin^\beta(\mathcal U_{Z_0})}&\lesssim_J\bracket{v_0}^m. \label{lin:physical-coefficient-holder}
\end{align}
Then $c_J,N$ can be chosen after these hypotheses are fixed so that all preceding conclusions hold and
\begin{align}
[\locA]_{\Ckin^\beta}+[\locb]_{\Ckin^\beta}+[\locB]_{\Ckin^\beta}+[\locc]_{\Ckin^\beta}&\lesssim_J\vradius^\beta\bracket{v_0}^M\lesssim_J1. \label{lin:normalized-coefficient-holder}
\end{align}
If, in addition, for some $a>0$,
\begin{align}
[\Kop f+\sourceR]_{\Ckin^\beta(\mathcal U_{Z_0})}&\lesssim_J\bracket{v_0}^m\ue^{-a\abs{v_0}^2}, \label{lin:physical-source-holder}
\end{align}
then
\begin{align}
[\widehat S_f]_{\Ckin^\beta}&\lesssim_J\vradius^{2+\beta}\bracket{v_0}^M\ue^{-a\abs{v_0}^2}. \label{lin:normalized-source-holder}
\end{align}
Here $M$ depends on $\beta,m$ and the fixed coefficient and geometric bounds, but is independent of the chosen $N$. The constants may depend on $J$ and the stated source and H\"older bounds on a fixed compact enlargement of $J$ in $I$.
\end{lemma}

\begin{proof}
\proofstep{Ellipticity} The range estimates \eqref{lin:position-variation}--\eqref{eq:I:display:265} and the bounds $\norm{D_v\sig(v)}\lesssim\bracket{v_0}^{\gamma+1}$ and \eqref{cb:matrix-bounds} give
\begin{align}
\norm[\Big]{S_0^{-1}\bigl[D\phi(x)\sig(v)D\phi(x)^T-D\phi(x_0)\sig(v_0)D\phi(x_0)^T\bigr]S_0^{-T}}&\lesssim\bigl(\bracket{v_0}\abs{v-v_0}+\bracket{v_0}^2\abs{x-x_0}\bigr). \label{eq:I:display:266}
\end{align}
If the norm on the left of \eqref{eq:I:display:266} is at most $\frac{1}{2}$, then \eqref{eq:I:display:262} and the bounds $\frac{1}{2}\sig\leq\Af{}\leq\frac{3}{2}\sig$ from Lemma~\ref{cvel:coefficients} give the ellipticity in \eqref{lin:normalized-coefficient-bounds}.

\proofstep{The lower-order coefficients and source} Before scaling, $\abs{\Bf{}}+\abs{D^2\phi(x)[v,v]}\lesssim\bracket{v_0}^4$ and $\abs{\Cf{}}\lesssim\bracket{v_0}^5$. Therefore
\begin{align}
\abs[\big]{\locb}+\abs[\big]{\locB}&\lesssim\vradius\bracket{v_0}^{\frac{11}{2}},\qquad \abs[\big]{\locc}\lesssim\vradius^2\bracket{v_0}^5,\qquad \abs[\big]{\nabla_\eta\locA}\lesssim\vradius\bracket{v_0}^{\frac{5}{2}}. \label{lin:scaled-lower-coefficients}
\end{align}
For the last bound, use the two factors $S_0^{-1}$, the factor $\vradius S_0$ from the velocity derivative, and \eqref{cvel:rough-coefficients}.

By \eqref{eq:I:display:265}, $\abs{v-v_0}\leq\bracket{v_0}^{-1}$ on the mapped cylinder, and hence $\abs[\big]{\abs{v}^2-\abs{v_0}^2}\lesssim1$. Substituting this bound into \eqref{cb:setting}, \eqref{cvel:scalar-K-bound} and \eqref{lin:source-Gaussian}, and using the factor $\vradius^2\leq1$ in \eqref{lin:normalized-f-source}, proves \eqref{cb:Gaussian-bound}.

\proofstep{Kinetic H\"older seminorms}
Under \eqref{lin:physical-coefficient-holder}, composition with $\Psi_0$ gains $\vradius^\beta\bracket{v_0}^{\beta N_1}$. The factors $D\phi(x)$, $D^2\phi(x)[v,v]$ and $v$ obey the same increment bound with fixed polynomial factors in $\bracket{v_0}$. Applying the product inequality for $\Ckin^\beta$ to \eqref{lin:normalized-f-coefficients}--\eqref{lin:normalized-drift} gives
\begin{align}
[\locA]_{\Ckin^\beta}+\vradius^{-1}\bigl([\locb]_{\Ckin^\beta}+[\locB]_{\Ckin^\beta}\bigr)+\vradius^{-2}[\locc]_{\Ckin^\beta}&\lesssim_J\vradius^\beta\bracket{v_0}^M. \label{lin:scaled-coefficient-holder}
\end{align}
Since $\vradius\leq1$, \eqref{lin:scaled-coefficient-holder} gives the first inequality in \eqref{lin:normalized-coefficient-holder}. Under the separate source assumption \eqref{lin:physical-source-holder}, the same composition bound and the factor $\vradius^2$ in \eqref{lin:normalized-f-source} give \eqref{lin:normalized-source-holder}. The exponent $M$ depends only on the fixed powers in \eqref{lin:physical-coefficient-holder}, \eqref{cb:matrix-bounds} and \eqref{lin:kinetic-change-bound}, enlarged to include $m+\beta N_1$.

Choose $N$ to dominate the fixed powers in \eqref{lin:position-variation}, \eqref{eq:I:display:265}, \eqref{eq:I:display:266} and \eqref{lin:scaled-lower-coefficients}, and, when the H\"older hypotheses are imposed, require $\beta N\geq M$. Then $\vradius^\beta\bracket{v_0}^M=c_J^\beta\bracket{v_0}^{M-\beta N}\leq1$. A sufficiently small $c_J$, uniform in $Z_0$, then ensures the domain and time restrictions and the ellipticity condition in \eqref{eq:I:display:266}.
\end{proof}

For later use, the velocity derivative conversion follows directly from \eqref{lin:inverse-map}. At corresponding points $(t,x,v)=\Psi_0(\theta,\xi,\eta)$,
\begin{align}
\nabla_vf(t,x,v)&=\frac{1}{\vradius}D\phi(x)^TS_0^{-T}\nabla_\eta\widehat f(\theta,\xi,\eta), \label{lin:physical-gradient-conversion}\\
D_v^2f(t,x,v)&=\frac{1}{\vradius^2}D\phi(x)^TS_0^{-T}\bigl(D_\eta^2\widehat f(\theta,\xi,\eta)\bigr)S_0^{-1}D\phi(x). \label{lin:physical-hessian-conversion}
\end{align}
The Hessian identity applies once these derivatives have been established. By \eqref{cb:matrix-bounds}, the conversion factors grow at most polynomially in $\bracket{v_0}$. The transport identity is
\begin{align}
(\Y f)\circ\Psi_0&=\frac{1}{\vradius^2}\Yhat\widehat f+\frac{1}{\vradius}\bigl[S_0^{-1}D^2\phi(x)[v,v]\bigr]\circ\Psi_0\cdot\nabla_\eta\widehat f. \label{lin:physical-transport-conversion}
\end{align}

\subsection{Interior Regularity}\label{lin:section-inputs}

We will now establish interior regularity for a more general linear kinetic equation. We will use a common weak equation and state the H\"older and Schauder estimates under separate coefficient assumptions. Here $r$ denotes a cylinder radius in the transformed variables, whereas $\vradius$ is the fixed scale in the map $\Psi_0$.

For noninteger $s\in(1,3)$, we use the kinetic H\"older space $\Ckin^s$ of \cite[Definition~2.6]{Ros-Oton.Weidner2025}, with the equivalent distance \eqref{eq:kinetic-distance}. This convention extends \eqref{eq:kinetic-holder-difference}: for $0<s<1$ the two definitions agree, and the equivalent distances give equivalent norms. We will recall the polynomial characterization in Section~\ref{lin:section-order-two}, where it enters the boundary expansion argument.

Interior H\"older estimates for kinetic equations with bounded measurable coefficients were proved in \cite{Golse.Imbert.Mouhot.Vasseur2019}. Local Schauder estimates in kinetic H\"older spaces are available, for example, in \cite[Theorem~3.9]{Imbert.Mouhot2021}. We will use the formulations in \cite[Lemmas~2.23--2.24]{Ros-Oton.Weidner2025}, which allow weak solutions and the two-sided cylinders used here.

The local equation is
\begin{align}
\Yhat\widehat f-\divv_\eta(\locA\nabla_\eta\widehat f)=\locB\cdot\nabla_\eta\widehat f+\locc\widehat f+\widehat S. \label{lin:Landau-equation}
\end{align}
It is understood in distributions. For a fixed center $\widehat Z_0=(\theta_0,\xi_0,\eta_0)$, we omit its argument from $\Qc_s$ and $\Qplus_s$.

\begin{proposition}[Interior estimates]\label{lin:interior-inputs}
Let $0<r\leq1$ and $0<\lambda\leq\Lambda<\infty$. Assume that
\begin{align}
\lambda I\leq\locA\leq\Lambda I\quad\text{in }\Qc_{2r}. \label{lin:local-ellipticity}
\end{align}
Let $\widehat f$ solve \eqref{lin:Landau-equation} in the weak sense of distributions on $\Qc_{2r}$, with
\begin{align}
\widehat f\in L^\infty(\Qc_{2r}),\qquad \nabla_\eta\widehat f\in L^2_{\mathrm{loc}}(\Qc_{2r}). \label{lin:common-local-hypotheses}
\end{align}
The following two estimates hold under their respective assumptions.
\begin{enumerate}
\item\label{lin:rough-full} Assume
\begin{align}
\locB,\locc,\widehat S\in L^\infty(\Qc_{2r}). \label{lin:interior-rough-hypotheses}
\end{align}
There is $\beta_0\in(0,1)$ such that, for $0<\beta<\beta_0$, we have $\widehat f\in\Ckin^\beta(\Qc_r)$ and
\begin{align}
r^\beta[\widehat f]_{\Ckin^\beta(\Qc_r)}\lesssim\tnorm[\big]{\widehat f}_{L^\infty(\Qc_{2r})}+r^2\tnorm[\big]{\widehat S}_{L^\infty(\Qc_{2r})}. \label{lin:rough-interior}
\end{align}
The exponent $\beta_0$ depends only on $\lambda,\Lambda$, $r\tnorm[\big]{\locB}_{L^\infty(\Qc_{2r})}$ and $r^2\tnorm[\big]{\locc}_{L^\infty(\Qc_{2r})}$; the implicit constant may also depend on $\beta$.

\item\label{lin:interior} Let $0<\beta<1$, set $\locb:=\locB+\divv_\eta\locA$, and assume
\begin{align}
\nabla_\eta\locA\in L^\infty(\Qc_{2r}),\qquad\locA,\locb,\locc,\widehat S\in\Ckin^\beta(\Qc_{2r}). \label{lin:interior-holder-hypotheses}
\end{align}
The equivalent nondivergence equation is
\begin{align}
\Yhat\widehat f-\locA:D_\eta^2\widehat f=\locb\cdot\nabla_\eta\widehat f+\locc\widehat f+\widehat S, \label{lin:Landau-nondivergence}
\end{align}
where $\locA:D_\eta^2\widehat f:=\divv_\eta(\locA\nabla_\eta\widehat f)-(\divv_\eta\locA)\cdot\nabla_\eta\widehat f$ until second derivatives are known. Then $\widehat f\in\Ckin^{2+\beta}(\Qc_r)$ and
\begin{align}
&r\tnorm[\big]{\nabla_\eta\widehat f}_{L^\infty(\Qc_r)}+r^2\left(\tnorm[\Big]{D_\eta^2\widehat f}_{L^\infty(\Qc_r)}+\tnorm[\big]{\Yhat\widehat f}_{L^\infty(\Qc_r)}\right)+r^{2+\beta}[\widehat f]_{\Ckin^{2+\beta}(\Qc_r)} \label{lin:interior-bound}\\
&\qquad\lesssim\tnorm[\big]{\widehat f}_{L^\infty(\Qc_{2r})}+r^2\tnorm[\big]{\widehat S}_{L^\infty(\Qc_{2r})}+r^{2+\beta}[\widehat S]_{\Ckin^\beta(\Qc_{2r})}.\notag
\end{align}
The implicit constant depends on $\lambda,\Lambda,\beta$ and the following scaled coefficient bounds on $\Qc_{2r}$:
\begin{align}
r^\beta[\locA]_{\Ckin^\beta}+r\tnorm[\big]{\nabla_\eta\locA}_{L^\infty}+r\tnorm[\big]{\locb}_{L^\infty}+r^{1+\beta}[\locb]_{\Ckin^\beta}+r^2\tnorm[\big]{\locc}_{L^\infty}+r^{2+\beta}[\locc]_{\Ckin^\beta}. \label{lin:interior-coefficient-dependence}
\end{align}
\end{enumerate}
\end{proposition}

\begin{proof}
\proofstep{Proof of \ref{lin:rough-full}.}
\cite[Lemma~2.24]{Ros-Oton.Weidner2025} gives a positive kinetic H\"older exponent for a weak solution of \eqref{lin:Landau-equation} with uniformly elliptic diffusion, bounded coefficients and source, $\widehat f\in L^\infty$ and $\nabla_\eta\widehat f\in L^2_{\mathrm{loc}}$. These requirements are \eqref{lin:local-ellipticity}, \eqref{lin:common-local-hypotheses} and \eqref{lin:interior-rough-hypotheses}. Under the translation and dilation $T_{\widehat Z_0}\circ S_r$, the coefficient and source tuple becomes $(\locA,r\locB,r^2\locc,r^2\widehat S)\circ T_{\widehat Z_0}\circ S_r$. The unit-scale estimate therefore gives \eqref{lin:rough-interior}, including the factor $r^2$ on the source and the stated coefficient dependence. The cited lemma uses the time intervals in \eqref{lin:cylinders}.

\proofstep{Proof of \ref{lin:interior}.}
Assumption~\eqref{lin:interior-holder-hypotheses} makes the contraction in \eqref{lin:Landau-nondivergence} well defined in distributions. Together with \eqref{lin:local-ellipticity}, it verifies the coefficient assumptions of \cite[Lemma~2.23]{Ros-Oton.Weidner2025} with $k=2$. The paragraph immediately after that lemma explicitly permits weak solutions. Thus it applies directly to $\widehat f$ and gives $\widehat f\in\Ckin^{2+\beta}(\Qc_r)$. The rescaled source has H\"older norm $r^2\tnorm[\big]{\widehat S}_{L^\infty(\Qc_{2r})}+r^{2+\beta}[\widehat S]_{\Ckin^\beta(\Qc_{2r})}$, so the cited estimate gives the seminorm term of \eqref{lin:interior-bound}.

The derivative and interpolation estimates \cite[Lemmas~2.12 and 2.10]{Ros-Oton.Weidner2025} bound $\nabla_\eta\widehat f,D_\eta^2\widehat f,\Yhat\widehat f$ by the $\Ckin^{2+\beta}$ norm on a slightly larger interior cylinder. Apply the preceding estimate on that cylinder and rescale: the three derivatives acquire the factors $r,r^2,r^2$, respectively. This gives \eqref{lin:interior-bound}, with the coefficient bounds \eqref{lin:interior-coefficient-dependence}.
\end{proof}

\subsection{Boundary H\"older and First-Order Derivative Estimates}\label{lin:section-boundary-inputs}

We will next establish boundary H\"older and first-order estimates for weak solutions of \eqref{lin:Landau-equation} on $\Qplus_{2r}$. We will first specify the weak formulation at $\xi_3=0$, then give three estimates with separate regularity assumptions on the coefficients and incoming data. Using the sets in \eqref{lin:boundary-sets}, we write
\begin{align}
\partial_-\Qplus_{2r}:=\Qc_{2r}\cap\gamma_-=\Qc_{2r}\cap\{\xi_3=0,\ \eta_3>0\}. \label{lin:local-incoming-boundary}
\end{align}
Boundary norms below are taken on this set and are zero when it is empty.

The weak formulation identifies the boundary trace as well as the interior equation. For every $\varphi\in C_c^\infty(\Qc_{2r})$, restricted to $\{\xi_3\geq0\}$, we require
\begin{align}
-\int\widehat f\Yhat\varphi+\int\locA\nabla_\eta\widehat f\cdot\nabla_\eta\varphi-\int_{\xi_3=0}\eta_3\trace\widehat f\,\varphi&=\int\bigl(\locB\cdot\nabla_\eta\widehat f+\locc\widehat f+\widehat S\bigr)\varphi. \label{lin:boundary-weak-setting}
\end{align}
The bulk integrals are over $\Qplus_{2r}$ with measure $\dd\theta\dd\xi\dd\eta$, and the boundary measure is $\dd\theta\dd\xi'\dd\eta$. The test function may have either sign and need not vanish on $\xi_3=0$. Formula~\eqref{lin:boundary-weak-setting} is the local version of \eqref{trace:compact-Green}, with outward transport flux $-\eta_3$.

\begin{proposition}[Boundary estimates]\label{lin:boundary-inputs}
Let $0<r\leq1$, $0<\lambda\leq\Lambda<\infty$, and $(\xi_0)_3\geq0$. Assume \eqref{lin:local-ellipticity} and \eqref{lin:interior-rough-hypotheses} on $\Qplus_{2r}$. Let $\widehat f$ be a weak solution of \eqref{lin:Landau-equation} in the sense of \eqref{lin:boundary-weak-setting}, with
\begin{align}
\widehat f\in L^\infty(\Qplus_{2r}),\qquad \chi\nabla_\eta\widehat f\in L^2(\Qplus_{2r})\quad\text{for every }\chi\in C_c^\infty(\Qc_{2r}). \label{lin:boundary-local-hypotheses}
\end{align}
The following conclusions hold.
\begin{enumerate}
\item\label{lin:rough} There is $\beta_0\in(0,1)$ such that, for $0<\beta<\beta_0$, the assumption
\begin{align}
\trm\widehat f\in\Ckin^\beta(\partial_-\Qplus_{2r}) \label{lin:boundary-rough-hypotheses}
\end{align}
implies $\widehat f\in\Ckin^\beta(\Qplus_r)$ and
\begin{align}
r^\beta[\widehat f]_{\Ckin^\beta(\Qplus_r)}&\lesssim\tnorm[\big]{\widehat f}_{L^\infty(\Qplus_{2r})}+r^\beta[\trm\widehat f]_{\Ckin^\beta(\partial_-\Qplus_{2r})}+r^2\tnorm[\big]{\widehat S}_{L^\infty(\Qplus_{2r})}. \label{lin:rough-bound}
\end{align}
The exponent $\beta_0$ depends only on $\lambda,\Lambda$, $r\tnorm[\big]{\locB}_{L^\infty(\Qplus_{2r})}$ and $r^2\tnorm[\big]{\locc}_{L^\infty(\Qplus_{2r})}$. The implicit constant may also depend on $\beta$ and is independent of the center velocity.

\item\label{lin:sharp} Let $0<\beta<\frac{1}{2}$ and assume
\begin{align}
\locA\in\Ckin^\beta(\Qplus_{2r}),\qquad\trm\widehat f\in\Ckin^{\frac{1}{2}+\beta}(\partial_-\Qplus_{2r}). \label{lin:boundary-sharp-hypotheses}
\end{align}
Then $\widehat f\in\Ckin^{\frac{1}{2}}(\Qplus_r)$ and
\begin{align}
r^{\frac{1}{2}}[\widehat f]_{\Ckin^{\frac{1}{2}}(\Qplus_r)}&\lesssim\left(1+r\tnorm[\big]{\locB}_{L^\infty(\Qplus_{2r})}\right)^N\Bigl[\tnorm[\big]{\widehat f}_{L^\infty(\Qplus_{2r})}+r^2\tnorm[\big]{\widehat S}_{L^\infty(\Qplus_{2r})} \label{lin:sharp-bound}\\
&\qquad+r^{\frac{1}{2}+\beta}[\trm\widehat f]_{\Ckin^{\frac{1}{2}+\beta}(\partial_-\Qplus_{2r})}\Bigr].\notag
\end{align}
The exponent $N>0$ and the implicit constant depend only on $\beta,\lambda,\Lambda$, $r^\beta[\locA]_{\Ckin^\beta(\Qplus_{2r})}$ and $r^2\tnorm[\big]{\locc}_{L^\infty(\Qplus_{2r})}$. The cylinder may meet the grazing set $\{\xi_3=\eta_3=0\}$.

\item\label{lin:gradient} Let $0<\beta<1$ and assume
\begin{align}
\Qc_{4r}\cap\{\xi_3=\eta_3=0\}=\varnothing,\qquad\locA\in\Ckin^\beta(\Qplus_{2r}),\qquad\trm\widehat f\in\Ckin^{1+\beta}(\partial_-\Qplus_{2r}). \label{lin:boundary-gradient-hypotheses}
\end{align}
Then $\widehat f\in\Ckin^{1+\beta}(\Qplus_r)$, $\nabla_\eta\widehat f\in\Ckin^\beta(\Qplus_r)$, and
\begin{align}
&r\tnorm[\big]{\nabla_\eta\widehat f}_{L^\infty(\Qplus_r)}+r^{1+\beta}[\nabla_\eta\widehat f]_{\Ckin^\beta(\Qplus_r)}+r^{1+\beta}[\widehat f]_{\Ckin^{1+\beta}(\Qplus_r)} \label{lin:gradient-bound}\\
&\qquad\lesssim\left(1+r\tnorm[\big]{\locB}_{L^\infty(\Qplus_{2r})}\right)^{N}\Bigl[\tnorm[\big]{\widehat f}_{L^\infty(\Qplus_{2r})}+r^2\tnorm[\big]{\widehat S}_{L^\infty(\Qplus_{2r})}+r^{1+\beta}[\trm\widehat f]_{\Ckin^{1+\beta}(\partial_-\Qplus_{2r})}\Bigr].\notag
\end{align}
The exponent $N>0$ and the implicit constant have the same dependence as in item~\ref{lin:sharp}, with the present value of $\beta$.

\item\label{lin:boundary-trace-identification} Under the assumptions of any one of items~\ref{lin:rough}--\ref{lin:gradient}, the continuous boundary values of $\widehat f$ on the plane $\xi_3=0$ inside $\Qc_r$ agree with the kinetic trace in \eqref{lin:boundary-weak-setting} almost everywhere for the measure $\abs{\eta_3}\,\dd\theta\,\dd\xi'\,\dd\eta$.
\end{enumerate}
\end{proposition}

\begin{proof}
We first verify the common weak-solution requirements. In \eqref{lin:boundary-weak-setting}, restrict to test functions satisfying $\supp\varphi\cap\{\xi_3=0\}\subset\{\eta_3\geq0\}$. Only the incoming trace remains in the boundary integral, giving the formulation in \cite[Definition~2.1]{Kim.Weidner2026} and \cite[Definition~2.21]{Ros-Oton.Weidner2025}. Condition~\eqref{lin:boundary-local-hypotheses} gives the required $L^2$ velocity gradient on every smaller cylinder, including its intersection with $\xi_3=0$. We apply the cited local estimates on such cylinders and cover $\Qplus_r$; the outer radius can remain a fixed multiple of $r$ strictly below $2r$.

We treat $\locc\widehat f+\widehat S$ as the scalar source. The common boundedness assumptions give
\begin{align}
r^2\tnorm[\big]{\locc\widehat f+\widehat S}_{L^\infty(\Qplus_{2r})}\leq r^2\tnorm[\big]{\locc}_{L^\infty(\Qplus_{2r})}\tnorm[\big]{\widehat f}_{L^\infty(\Qplus_{2r})}+r^2\tnorm[\big]{\widehat S}_{L^\infty(\Qplus_{2r})}. \label{lin:boundary-source-bound}
\end{align}
Under kinetic dilation, the coefficient and source tuple is $(\locA,r\locB,r^2(\locc\widehat f+\widehat S))$ composed with the dilation. The scaled diffusion seminorm is $r^\beta[\locA]_{\Ckin^\beta}$, and an incoming seminorm of order $s$ acquires the factor $r^s$.

\proofstep{Proof of \ref{lin:rough}.}
The divergence-form estimate \cite[Theorem~1.2]{Silvestre2022} gives a positive H\"older exponent under uniform ellipticity, bounded $\locB$ and scalar source, and a H\"older incoming trace. Apply the divergence-form argument in the proof of \cite[Lemma~2.25]{Ros-Oton.Weidner2025} to \eqref{lin:Landau-equation} on $\Qplus_{2r}$, with scalar source $\locc\widehat f+\widehat S$. Its hypotheses are \eqref{lin:local-ellipticity}, \eqref{lin:interior-rough-hypotheses}, \eqref{lin:boundary-source-bound} and \eqref{lin:boundary-rough-hypotheses}, together with the weak-solution conditions verified above. The conclusion is $\widehat f\in\Ckin^\beta(\Qplus_r)$ with \eqref{lin:rough-bound}. Convexity of the half-space gives uniformity in the center velocity. The kinetic Lipschitz assumption on the diffusion in the statement of the cited lemma is used to pass from its nondivergence equation to divergence form; here the weak equation already has the form \eqref{lin:Landau-equation}.

\proofstep{Proof of \ref{lin:sharp}.}
For a weak solution of \eqref{lin:Landau-equation}, \cite[Theorem~6.2]{Kim.Weidner2026} gives the $\Ckin^{\frac{1}{2}}$ estimate with bounded $\locB$ and scalar source, uniformly elliptic $\locA\in\Ckin^\beta$, and incoming trace in $\Ckin^{\frac{1}{2}+\beta}$. The additional requirements are exactly \eqref{lin:boundary-sharp-hypotheses}. Applying the theorem with equation radius $2r$, using \eqref{lin:boundary-source-bound} and the scaling above, gives \eqref{lin:sharp-bound}. No exclusion of the grazing set is required.

\proofstep{Proof of \ref{lin:gradient}.}
The estimate \cite[Lemma~7.3, (7.9)]{Kim.Weidner2026} gives the three quantities on the left of \eqref{lin:gradient-bound}. Its additional assumptions are $\locA\in\Ckin^\beta$, incoming trace in $\Ckin^{1+\beta}$, and exclusion of the grazing set from the cylinder of twice the equation radius. With equation radius $2r$, these are \eqref{lin:boundary-gradient-hypotheses}. Its additional velocity-divergence source is zero here. The common weak-solution conditions and \eqref{lin:boundary-source-bound} therefore give \eqref{lin:gradient-bound}, with the stated scaled coefficient dependence.

\proofstep{Proof of \ref{lin:boundary-trace-identification}.}
Write $\widehat f_b$ for the continuous boundary value supplied by the estimates. Let $\psi(\theta,\xi',\eta)$ be a smooth boundary test function compactly supported in $\Qc_r\cap\{\xi_3=0,\ \abs{\eta_3}>0\}$, and choose $\rho\in C_c^\infty(\R)$ with $\rho(0)=1$. In \eqref{lin:boundary-weak-setting}, take $\varphi_\delta=\psi\rho(\xi_3/\delta)$. All bulk terms tend to zero except the normal transport term, for which continuity gives
\begin{align}
-\int_{\Qplus_{2r}}\widehat f\,\psi\,\eta_3\frac{1}{\delta}\rho'\left(\frac{\xi_3}{\delta}\right)\dd\theta\dd\xi\dd\eta\to\int_{\xi_3=0}\eta_3\widehat f_b\psi\,\dd\theta\dd\xi'\dd\eta. \label{lin:continuous-trace-limit}
\end{align}
The terms containing $\nabla_\eta\widehat f$ vanish by Cauchy--Schwarz and \eqref{lin:boundary-local-hypotheses}. Hence
\begin{align}
\int_{\xi_3=0}\eta_3\left(\widehat f_b-\trace\widehat f\right)\psi\,\dd\theta\dd\xi'\dd\eta=0. \label{lin:continuous-trace-identification}
\end{align}
As in the uniqueness argument in Lemma~\ref{trace:CM}, exhaustion of the incoming and outgoing sets identifies the two boundary values almost everywhere for the measure $\abs{\eta_3}\,\dd\theta\,\dd\xi'\,\dd\eta$.
\end{proof}

\subsection{A Second-Order Derivative Estimate at Nongrazing Boundary Points}\label{lin:section-order-two}

We will obtain the second-order estimate by adapting the nongrazing boundary expansions and Schauder estimates of \cite[Lemmas~4.1, 4.3--4.4]{Ros-Oton.Weidner2025} to polynomial degree two. Related boundary Schauder estimates are given in \cite[Theorem~1.3(i) and Proposition~4.3]{Zhu2025}. We will verify the adaptation using $\nabla_\eta\locA\in L^\infty$ in \eqref{lin:coefficients} and retain the polynomial dependence on the normal center velocity. Boundary polynomial approximations will first give small remainders on neighboring interior cylinders; the interior estimate will then yield the polynomial approximation at every center. We will impose zero incoming trace; Section~\ref{cb:section} will construct the lifts needed for nonzero incoming data.

Let $\mathcal U\subset\R\times\R^3\times\R^3$ be open in the transformed variables $\widehat Z=(\theta,\xi,\eta)$. For $0<\beta<1$, the polynomial spaces used for orders $1+\beta$ and $2+\beta$ are
\begin{align}
\mathcal P_1:=\operatorname{span}\{1,\eta_1,\eta_2,\eta_3\},\qquad \mathcal P_2:=\operatorname{span}\{1,\theta,\eta_i,\eta_i\eta_j:1\leq i\leq j\leq3\}. \label{lin:polynomial-space}
\end{align}
For $k=1,2$, a bounded function $\widehat f$ belongs to $\Ckin^{k+\beta}(\mathcal U)$ if and only if there is $C\geq0$ such that every $\widehat Z_0\in\mathcal U$ admits a polynomial $P_{\widehat Z_0}\in\mathcal P_k$ satisfying
\begin{align}
\abs[\big]{\widehat f(\widehat Z)-P_{\widehat Z_0}(\widehat Z)}&\leq C\bigl(\dkin(\widehat Z,\widehat Z_0)\bigr)^{k+\beta}\qquad(\widehat Z\in\mathcal U). \label{eq:kinetic-holder-polynomial}
\end{align}
The same $C$ works for all centers. Its least value is $[\widehat f]_{\Ckin^{k+\beta}(\mathcal U)}$; adding $\tnorm[\big]{\widehat f}_{L^\infty(\mathcal U)}$ gives the norm. This is the characterization in \cite[Definition~2.6]{Ros-Oton.Weidner2025}, with the equivalent kinetic distance \eqref{eq:kinetic-distance}. For boundary data, we use the same condition with both points restricted to the boundary set.

For the fixed transformed coefficients, write
\begin{align}
\locOp:=\Yhat-\locA:D_\eta^2-\locb\cdot\nabla_\eta-\locc. \label{lin:local-operator}
\end{align}
We use the weak formulation \eqref{lin:boundary-weak-setting}, with $\locB=\locb-\divv_\eta\locA$, for the equation $\locOp\widehat u=\widehat S$. Thus the nondivergence expression has the distributional meaning specified after \eqref{lin:Landau-nondivergence}. As before, we omit a fixed cylinder center.

\begin{proposition}[Second-order boundary estimate]\label{lin:order-two}
Let $0<\beta<1$, $0<r\leq1$, $0<\lambda\leq\Lambda<\infty$, $M\geq0$, and $\widehat Z_0=(\theta_0,\xi_0,\eta_0)$ with $(\xi_0)_3\geq0$. Assume
\begin{align}
\Qc_{64r}(\widehat Z_0)\cap\{\xi_3=\eta_3=0\}=\varnothing. \label{lin:second-order-nongrazing}
\end{align}
The coefficients and source satisfy
\begin{align}
\lambda I\leq\locA\leq\Lambda I,\qquad \locA,\locb,\locc,\widehat S\in\Ckin^\beta(\Qplus_{16r}),\qquad\nabla_\eta\locA\in L^\infty(\Qplus_{16r}), \label{lin:coefficients}
\end{align}
and, with all norms in the following display taken on $\Qplus_{16r}$,
\begin{align}
r^\beta[\locA]_{\Ckin^\beta}+r\tnorm[\big]{\nabla_\eta\locA}_{L^\infty}+r\tnorm[\big]{\locb}_{L^\infty}+r^{1+\beta}[\locb]_{\Ckin^\beta}+r^2\tnorm[\big]{\locc}_{L^\infty}+r^{2+\beta}[\locc]_{\Ckin^\beta}\leq M. \label{lin:coefficient-norm}
\end{align}
Let $\widehat u$ be a weak solution, in the sense of \eqref{lin:boundary-weak-setting}, of
\begin{align}
\locOp\widehat u=\widehat S\quad\text{in }\Qplus_{16r}. \label{lin:equation}
\end{align}
Assume \eqref{lin:boundary-local-hypotheses}, with $\widehat f$ replaced by $\widehat u$ and $2r$ replaced by $16r$ in both that condition and the weak formulation, and suppose
\begin{align}
\trm\widehat u=0\quad\text{on }\partial_-\Qplus_{16r}. \label{lin:zero-incoming-trace}
\end{align}
The incoming boundary may be empty. Then $\widehat u\in\Ckin^{2+\beta}(\Qplus_r)$ and
\begin{align}
&r\tnorm[\big]{\nabla_\eta\widehat u}_{L^\infty(\Qplus_r)}+r^2\left(\tnorm[\Big]{D_\eta^2\widehat u}_{L^\infty(\Qplus_r)}+\tnorm[\big]{\Yhat\widehat u}_{L^\infty(\Qplus_r)}\right)+r^{2+\beta}[\widehat u]_{\Ckin^{2+\beta}(\Qplus_r)} \label{lin:order-two-bound}\\
&\qquad\lesssim\left(1+\frac{\abs{\eta_{0,3}}}{r}\right)^N\Bigl(\tnorm[\big]{\widehat u}_{L^\infty(\Qplus_{16r})}+r^2\tnorm[\big]{\widehat S}_{L^\infty(\Qplus_{16r})}+r^{2+\beta}[\widehat S]_{\Ckin^\beta(\Qplus_{16r})}\Bigr).\notag
\end{align}
The exponent $N$ and the implicit constant depend only on $\lambda,\Lambda,\beta,M$.
\end{proposition}

\begin{remark}
Proposition~\ref{lin:boundary-inputs}\ref{lin:gradient} uses incoming data in $\Ckin^{1+\beta}$. For second derivatives, we will construct a lift satisfying
\begin{align}
\trm G=g,\qquad G,\ \nabla_\eta G,\ D_\eta^2G,\ \Yhat G\in\Ckin^\beta(\Qplus_{16r}). \label{lin:lift-required-regularity}
\end{align}
If $\locOp\widehat f=\widehat S$ and $\trm\widehat f=g$, then
\begin{align}
\locOp(\widehat f-G)=\widehat S-\locOp G,\qquad \trm(\widehat f-G)=0. \label{lin:lift-source-reduction}
\end{align}
When $\widehat S\in\Ckin^\beta$, \eqref{lin:coefficients} and \eqref{lin:lift-required-regularity} give $\widehat S-\locOp G\in\Ckin^\beta$, permitting Proposition~\ref{lin:order-two} for the difference. Sections~\ref{cb:inflow-graph-section} and \ref{cb:section-graph} will verify \eqref{lin:lift-required-regularity} for the in-flow and diffuse-reflection lifts.
\end{remark}

\para{Idea of the proof}
We will prove \eqref{lin:order-two-bound} by reducing to unit radius and constructing a polynomial at each center. Fix an arbitrary interior point $\widehat Z$ close to the boundary, and let $\delta$ be the radius of its largest full interior cylinder. We will choose a point $\widehat Z_b$ where this cylinder touches the boundary. Lemma~\ref{lin:boundary-expansion-lemma} will give a polynomial $P_{\widehat Z_b}$ with error of order $\delta^{2+\beta}$ on the neighboring interior cylinder. We will estimate $\widehat u-P_{\widehat Z_b}$ there by the interior Schauder estimate and obtain its polynomial $q_{\widehat Z}$ centered at $\widehat Z$. Then $P_{\widehat Z}:=P_{\widehat Z_b}+q_{\widehat Z}$ will approximate $\widehat u$ at $\widehat Z$ on a small cylinder, with a constant independent of $\delta$.

We will next extend this remainder estimate to every evaluation point in the smaller half-cylinder, keeping $P_{\widehat Z}$ fixed. For points outside the small interior cylinder, we will use the boundary approximation together with the coefficient bounds for $q_{\widehat Z}=P_{\widehat Z}-P_{\widehat Z_b}$. This will establish \eqref{eq:kinetic-holder-polynomial} on the whole smaller half-cylinder. Boundary centers will use Lemma~\ref{lin:boundary-expansion-lemma} directly, and centers farther from the boundary will use the interior estimate.

To implement the interior estimate, we must control the full source of the remainder equation. For $P\in\mathcal P_2$, we have
\begin{align}
\locOp(\widehat u-P)=\widehat S-\Yhat P+\locA:D_\eta^2P+\locb\cdot\nabla_\eta P+\locc P. \label{lin:polynomial-source}
\end{align}
The H\"older seminorm bounds the oscillation of this source. The next lemma will control its constant part using the size of the solution. We will apply it both to the rescaled remainders in Lemma~\ref{lin:boundary-expansion-lemma} and to $\widehat u-P_{\widehat Z_b}$ in the proof of Proposition~\ref{lin:order-two}.

\begin{lemma}[Control of a constant source]\label{lin:constant-source}
Let $\Qc_2=\Qc_2(\widehat Z_*)$ be a full cylinder contained in the domain of the transformed equation. Assume \eqref{lin:coefficients} with $\Qplus_{16r}$ replaced by $\Qc_2$, and use the operator \eqref{lin:local-operator}. Let $C_0\in\R$ and $\widehat S\in\Ckin^\beta(\Qc_2)$. If $\widehat h\in L^\infty(\Qc_2)$ satisfies $\locOp\widehat h=C_0+\widehat S$ in the weak sense of \eqref{lin:normalized-weak}, then this solution necessarily satisfies
\begin{align}
\abs{C_0}\lesssim\tnorm[\big]{\widehat h}_{L^\infty(\Qc_2)}+\tnorm[\big]{\widehat S}_{\Ckin^\beta(\Qc_2)}. \label{lin:source-constant-control}
\end{align}
The implicit constant depends only on $\lambda,\Lambda,\beta$, the $\Ckin^\beta(\Qc_2)$ norms of $\locA,\locb,\locc$ and $\tnorm[\big]{\nabla_\eta\locA}_{L^\infty(\Qc_2)}$.
\end{lemma}

\begin{proof}
Apply \cite[Lemma~2.27, proof of (2.36)]{Ros-Oton.Weidner2025} with polynomial degree $k=0$, polynomial $p=C_0$, and remainder $g=\widehat S$, after a fixed dilation of $\Qc_2$. The ellipticity and $\Ckin^\beta$ coefficient assumptions are \eqref{lin:coefficients}; the weak nondivergence equation has the meaning fixed above. The cited estimate bounds $\tnorm{p}_{L^\infty}$ by $\tnorm[\big]{\widehat h}_{L^\infty}+\tnorm{g}_{\Ckin^\beta}$, which is \eqref{lin:source-constant-control}.
\end{proof}

For a polynomial $P\in\mathcal P_2$ expressed in $\theta-\theta_0$ and $\eta-\eta_0$, write $p_{\widehat Z_0,0}$ for its constant coefficient, $p_{\widehat Z_0,1}$ for its velocity-linear coefficients, and $p_{\widehat Z_0,2}$ for its time and velocity-quadratic coefficients. We use the Euclidean norm of each coefficient list.

\begin{lemma}[Boundary polynomial expansion]\label{lin:boundary-expansion-lemma}
Assume the hypotheses of Proposition~\ref{lin:order-two} with $(\xi_0)_3=0$. There is a polynomial $P_{\widehat Z_0}\in\mathcal P_2$, expressed in $\theta-\theta_0$ and $\eta-\eta_0$, such that for $0<\varrho\leq r$,
\begin{align}
\tnorm[\big]{\widehat u-P_{\widehat Z_0}}_{L^\infty\bigl(\Qplus_\varrho(\widehat Z_0)\bigr)} &\lesssim\left(\frac{\varrho}{r}\right)^{2+\beta}\Bigl(\tnorm[\big]{\widehat u}_{L^\infty\bigl(\Qplus_{16r}(\widehat Z_0)\bigr)} \label{lin:boundary-expansion}\\
&\qquad+r^2\tnorm[\big]{\widehat S}_{L^\infty\bigl(\Qplus_{16r}(\widehat Z_0)\bigr)}+r^{2+\beta}[\widehat S]_{\Ckin^\beta\bigl(\Qplus_{16r}(\widehat Z_0)\bigr)}\Bigr).\notag
\end{align}
For $\nu=0,1,2$, its coefficient lists satisfy
\begin{align}
r^{\nu}\abs[\big]{p_{\widehat Z_0,\nu}} &\lesssim\Bigl(\tnorm[\big]{\widehat u}_{L^\infty\bigl(\Qplus_{16r}(\widehat Z_0)\bigr)}+r^2\tnorm[\big]{\widehat S}_{L^\infty\bigl(\Qplus_{16r}(\widehat Z_0)\bigr)} +r^{2+\beta}[\widehat S]_{\Ckin^\beta\bigl(\Qplus_{16r}(\widehat Z_0)\bigr)}\Bigr). \label{lin:boundary-expansion-coefficients}
\end{align}
The implicit constants have the dependence in Proposition~\ref{lin:order-two} and are independent of $\eta_{0,3}$.
\end{lemma}

\begin{proof}
We will prove the estimate by contradiction. If the polynomial approximation fails, we will subtract the $L^2$ projection onto $\mathcal P_2$ at a shrinking scale and normalize the remainder. A uniform bound for polynomial coefficients will make this construction independent of the boundary center. The remainders will be orthogonal to $\mathcal P_2$, have size at least $\frac12$ on the unit cylinder, and satisfy the growth bound \eqref{lin:growth}. We will first bound their rescaled sources and then obtain a nonzero limit by boundary H\"older estimates. A local energy estimate will allow us to pass the weak equation and, at incoming points, the initial trace to this limit. The Liouville theorems will show that the limit belongs to $\mathcal P_2$, whereas its inherited orthogonality to $\mathcal P_2$ forces it to vanish.

We will use the construction in the proofs of \cite[Lemmas~4.1 and 4.3]{Ros-Oton.Weidner2025} at degree $k=2$. As explained after their Lemma~4.1, the restriction $k\geq3$ supplies the divergence formulation needed for boundary H\"older estimates. Here $\nabla_\eta\locA\in L^\infty$ in \eqref{lin:coefficients} gives $\locB=\locb-\divv_\eta\locA\in L^\infty$. We will check the rescaled source and incoming values before applying Proposition~\ref{lin:boundary-inputs}\ref{lin:rough}, and the local $L^2$ bounds and limiting trace before applying the Liouville theorems.

\proofstep{Step 1: A uniform bound for polynomial coefficients.}
We will bound the coefficients of $P\in\mathcal P_2$ by its $L^2$ norm, uniformly in the cylinder center. This will control the polynomial projections used in Step~2 even when the normal velocity varies. For the cylinders \eqref{lin:cylinders}--\eqref{lin:half-cylinders}, we claim that, whenever $(\xi_0)_3\geq0$,
\begin{align}
\tnorm{P}_{L^\infty(\Qplus_\varrho(\widehat Z_0))}+\sum_{\nu=0}^2\varrho^\nu\abs[\big]{p_{\widehat Z_0,\nu}}&\lesssim\frac{\tnorm{P}_{L^2(\Qplus_\varrho(\widehat Z_0))}}{\abs[\big]{\Qc_\varrho(\widehat Z_0)}^{\frac{1}{2}}}. \label{lin:polynomial-coefficient-bound}
\end{align}
Translation and dilation send the domain to $\Qc_1(0)\cap\{\xi_3+\varrho^{-3}(\xi_0)_3+\varrho^{-1}(\eta_0)_3\theta>0\}$. It contains the fixed product set $\frac{1}{4}<\xi_3<\frac{1}{2}$, $\abs{\xi'}<\frac{1}{4}$, $\abs{\eta}<\frac{1}{2}$ and $\frac{1}{4}<\theta<\frac{1}{2}$ when $(\eta_0)_3\geq0$, or $-\frac{1}{2}<\theta<-\frac{1}{4}$ when $(\eta_0)_3<0$. Norm equivalence for $\mathcal P_2$ on these two fixed sets bounds the rescaled coefficients by the rescaled $L^2$ norm. Those coefficients also bound the polynomial on $\Qc_1(0)$. This proves \eqref{lin:polynomial-coefficient-bound}, with a constant independent of $\widehat Z_0$ and $\varrho$. The same estimate holds on full cylinders.

\proofstep{Step 2: Failure of the estimate and the rescaled remainder.}
We will turn failure of the estimate into remainders satisfying the nonzero bound, growth and orthogonality in \eqref{lin:growth}--\eqref{lin:orthogonality}. Fix $\lambda,\Lambda,\beta,M$. A kinetic dilation preserving $\{\xi_3>0\}$ reduces the radius to $r=1$: the equation is then posed on $\Qplus_{16}(\widehat Z_0)$, and $\Qc_{64}(\widehat Z_0)$ avoids grazing. Dividing the solution and source by their combined norm reduces the desired bound to the class satisfying
\begin{align}
\tnorm[\big]{\widehat u}_{L^\infty(\Qplus_{16}(\widehat Z_0))}+\tnorm[\big]{\widehat S}_{\Ckin^\beta(\Qplus_{16}(\widehat Z_0))}\leq1. \label{lin:boundary-expansion-normalization}
\end{align}
The coefficients and boundary center may vary within this class, while \eqref{lin:coefficients}--\eqref{lin:coefficient-norm} hold with the fixed parameters. We must prove one constant valid for the entire class.

Assume that the conclusion fails uniformly over this class. By \cite[Lemma~2.14]{Ros-Oton.Weidner2025} and \eqref{lin:polynomial-coefficient-bound}, a uniform bound of order $s^{2+\beta}$ for the $L^\infty(\Qplus_s)$ norms of the errors to the $L^2(\Qplus_s)$ projections would give one limiting polynomial and its coefficient bounds. Thus no constant bounds all these projection errors by $Cs^{2+\beta}$. We will use this failure to choose the solutions and scales for the normalized remainders.

We use the selection and dyadic estimates in \cite[proof of Lemma~4.1, Step~1, (4.7)--(4.10), and Lemma~2.15]{Ros-Oton.Weidner2025}. This part of the proof uses boundedness of the normalized functions, uniform norm equivalence for polynomial projections, and $2+\beta-\nu>0$ for $\nu=0,1,2$. These are supplied by \eqref{lin:boundary-expansion-normalization} and Step~1. The construction gives normalized solutions $\widehat u_j$, sources and coefficients in the same class, boundary centers $\widehat Z_j=(\theta_j,\xi_j,\eta_j)$, radii $s_j\to0$, and positive normalizing factors $\Theta_j\to\infty$. Let $P_{j,s_j}$ be the $L^2(\Qplus_{s_j}(\widehat Z_j))$ projection onto $\mathcal P_2$, with coefficient lists $p_{j,s_j,\nu}$ expressed about $\widehat Z_j$, and set
\begin{align}
\widehat h_j:=\frac{(\widehat u_j-P_{j,s_j})\circ T_{\widehat Z_j}\circ S_{s_j}}{s_j^{2+\beta}\Theta_j}. \label{lin:blowup}
\end{align}
The spatial half-space becomes
\begin{align}
\mathcal U_j:=\{\xi_3+k_j\theta>0\},\qquad k_j:=\frac{\eta_{j,3}}{s_j}. \label{lin:tilted-domain}
\end{align}
The quantitative conclusions of the cited construction are
\begin{align}
\frac{1}{2}\leq\tnorm[\big]{\widehat h_j}_{L^\infty(\mathcal U_j\cap\Qc_1(0))}\leq1,\qquad \tnorm[\big]{\widehat h_j}_{L^\infty(\mathcal U_j\cap\Qc_\varrho(0))}\lesssim\varrho^{2+\beta}\quad\left(1\leq\varrho\leq\frac{1}{2s_j}\right), \label{lin:growth}\\
\frac{\abs[\big]{p_{j,s_j,\nu}}}{\Theta_j}\to0\qquad(\nu=0,1,2). \label{lin:coefficient-smallness}
\end{align}
The projection property also gives
\begin{align}
\int_{\mathcal U_j\cap\Qc_1(0)}\widehat h_jP\,\dd\widehat Z=0\qquad(P\in\mathcal P_2). \label{lin:orthogonality}
\end{align}
Since $(\xi_j)_3=0$ and $\Qc_{64}(\widehat Z_j)$ avoids grazing, $\abs{\eta_{j,3}}\geq64$. After taking a subsequence, $k_j\to+\infty$ at incoming points or $k_j\to-\infty$ at outgoing points. Thus $\mathcal U_j$ tends to $\{\theta>0\}$ or $\{\theta<0\}$, respectively; see Figure~\ref{lin:incoming-rescaling}.

In Steps~3--6, we will control the source, pass to a limit $\widehat h\in\mathcal P_2$, and contradict the nonzero bound in \eqref{lin:growth} by passing \eqref{lin:orthogonality} to the limit and taking $P=\widehat h$.
\begin{figure}[tbp]
\centering
\begin{tikzpicture}[x=1.25cm,y=1.25cm,>=Stealth]
\fill[blue!6] (-0.4,1.6)--(2.1,1.6)--(2.1,-1.6)--(0.4,-1.6)--cycle;
\draw[->] (-2.1,0)--(2.35,0) node[right] {$\theta$};
\draw[dashed,->] (0,-1.6)--(0,1.85) node[above] {$\xi_3$};
\node[anchor=east] at (-0.1,-1.45) {$\theta=0$};
\draw[blue!70!black,thick] (-0.4,1.6)--(0.4,-1.6);
\node at (1.15,0.95) {$\mathcal U_j$};
\node[align=center] at (-1.12,-0.92) {$\xi_3+k_j\theta=0$};
\draw[thin] (-0.62,-0.72)--(0.18,-0.72);
\draw[->] (-0.57,1.15)--(-0.06,1.15);
\node[anchor=east] at (-0.62,1.15) {$k_j\to\infty$};

\end{tikzpicture}
\caption{The rescaled domains \eqref{lin:tilted-domain}: locally away from $\theta=0$, $\mathcal U_j$ tends to $\{\theta>0\}$ as $k_j\to+\infty$ and to $\{\theta<0\}$ as $k_j\to-\infty$.}
\label{lin:incoming-rescaling}
\end{figure}
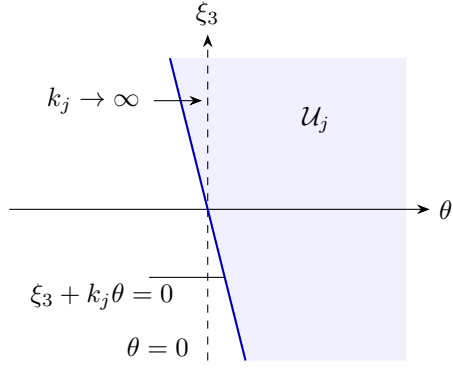

\proofstep{Step 3: The rescaled equation and its source.}
We will show that the equation for $\widehat h_j$ has diffusion and source converging to constants, with both lower-order coefficients tending to zero. In particular, the source must remain bounded so that the boundary H\"older estimate can be used in Step~4. Let $\locOp_j$ be \eqref{lin:local-operator} with the coefficients indexed by $j$. The function $\widehat h_j$ solves
\begin{align}
\Yhat \widehat h_j-\locA_j^{(s)}:D_{\eta}^2\widehat h_j &=\locb_j^{(s)}\cdot\nabla_{\eta}\widehat h_j+\locc_j^{(s)}\widehat h_j+\widehat S_j^{(s)}, \label{lin:blowup-equation}
\end{align}
where $\locA_j^{(s)}=\locA_j\circ(T_{\widehat Z_j}\circ S_{s_j})$, $\locb_j^{(s)}=s_j\locb_j\circ(T_{\widehat Z_j}\circ S_{s_j})$, $\locc_j^{(s)}=s_j^2\locc_j\circ(T_{\widehat Z_j}\circ S_{s_j})$, and
\begin{align}
\widehat S_j^{(s)}&:=\frac{\bigl(\widehat S_j-\locOp_jP_{j,s_j}\bigr)\circ\bigl(T_{\widehat Z_j}\circ S_{s_j}\bigr)}{s_j^\beta\Theta_j}. \label{lin:source}
\end{align}
We verify the source calculation in \cite[proof of Lemma~4.1, Step~2]{Ros-Oton.Weidner2025} at degree two, using \eqref{lin:polynomial-source}. The derivatives $\Yhat P_{j,s_j}$ and $D_\eta^2P_{j,s_j}$ are constant. Coefficient oscillations gain $s_j^\beta$, while velocity-linear and quadratic increments gain $s_j$ and $s_j^2$. Hence \eqref{lin:source} gives, for every fixed $\varrho\geq1$ and all sufficiently large $j$,
\begin{align}
[\widehat S_j^{(s)}]_{\Ckin^\beta(\mathcal U_j\cap\Qc_\varrho(0))} &\lesssim_\varrho\frac{1}{\Theta_j}\left(1+\sum_{\nu=0}^2\abs[\big]{p_{j,s_j,\nu}}\right) \to0 \label{lin:source-seminorm}
\end{align}
by \eqref{lin:coefficient-smallness}. To turn this vanishing oscillation into a bound for the full source, we will control its value at one fixed point. Choose a full cylinder $\Qc_2(\widehat Z_*)$ compactly contained in $\{\theta>0\}$ for incoming points or $\{\theta<0\}$ for outgoing points. It lies in $\mathcal U_j$ for all large $j$. The coefficient assumptions of Lemma~\ref{lin:constant-source} follow by scaling \eqref{lin:coefficients}--\eqref{lin:coefficient-norm}. Apply that lemma to \eqref{lin:blowup-equation}, with constant source $\widehat S_j^{(s)}(\widehat Z_*)$ and remainder $\widehat S_j^{(s)}-\widehat S_j^{(s)}(\widehat Z_*)$. The remainder vanishes at the center, so its full $\Ckin^\beta$ norm is bounded by its seminorm on this fixed cylinder. Equations~\eqref{lin:growth} and \eqref{lin:source-seminorm} give
\begin{align}
\abs[\big]{\widehat S_j^{(s)}(\widehat Z_*)}\lesssim\tnorm[\big]{\widehat h_j}_{L^\infty(\Qc_2(\widehat Z_*))}+[\widehat S_j^{(s)}]_{\Ckin^\beta(\Qc_2(\widehat Z_*))}\lesssim1. \label{lin:rescaled-source-value}
\end{align}
The source values therefore have a convergent subsequence, and \eqref{lin:source-seminorm} makes the limit constant. Similarly, $[\locA_j^{(s)}]_{\Ckin^\beta}\lesssim s_j^\beta$ and ellipticity give a constant diffusion limit, while the factors $s_j$ and $s_j^2$ make the lower-order coefficients vanish. After taking a common subsequence,
\begin{align}
\widehat S_j^{(s)}&\to \widehat S_*\quad\text{locally uniformly}, \qquad \locA_j^{(s)}\to \locA_*, \qquad \locb_j^{(s)},\ \locc_j^{(s)}\to0, \label{lin:constant-limit}
\end{align}
where $\widehat S_*$ and $\locA_*$ are constant and $\locA_*$ is elliptic.

The divergence formulation has first-order coefficient $\locb_j^{(s)}-\divv_\eta\locA_j^{(s)}$. Its second term satisfies
\begin{align}
\tnorm[\Big]{\divv_{\eta}\locA_j^{(s)}}_{L^\infty(\mathcal U_j\cap\Qc_\varrho(0))} &\lesssim s_j\tnorm[\big]{\nabla_{\eta}\locA_j}_{L^\infty\bigl(\Qplus_{16}(\widehat Z_j)\bigr)} \to0, \label{lin:drift-scaling}
\end{align}
Here the derivative on the right is taken before scaling. The bounds \eqref{lin:constant-limit}--\eqref{lin:drift-scaling} supply the coefficients and source for the boundary H\"older estimate in Step~4.

\proofstep{Step 4: Boundary H\"older bounds and convergence.}
We will obtain local uniform convergence up to $\theta=0$, so that the limit retains the nonzero size and growth in \eqref{lin:growth}. We will apply Proposition~\ref{lin:boundary-inputs}\ref{lin:rough} after checking its incoming-data assumption. The ellipticity, source and first-order bounds required by \eqref{lin:local-ellipticity} and \eqref{lin:interior-rough-hypotheses} were proved in Step~3. The local velocity-gradient condition \eqref{lin:boundary-local-hypotheses} holds for each remainder by subtraction and scaling. We apply the proposition on the original spatial half-space and then rescale, using its constant independent of the center velocity.

At outgoing points, for fixed $\varrho$ and all large $j$, every boundary velocity in the original cylinder corresponding to $\mathcal U_j\cap\Qc_{2\varrho}(0)$ is outgoing, so there is no incoming term. Together with \eqref{lin:growth}, \eqref{lin:constant-limit} and \eqref{lin:drift-scaling}, this gives a uniform positive H\"older exponent for $\widehat h_j$ up to the moving boundary.

At incoming points, each original solution $\widehat u_j$ is continuous up to the incoming boundary by Proposition~\ref{lin:boundary-inputs}\ref{lin:rough}, since its incoming value is zero and its original source is bounded. The trace identification in Proposition~\ref{lin:boundary-inputs}\ref{lin:boundary-trace-identification} gives the same zero boundary value for this continuous representative. Hence $\widehat h_j$ is individually continuous up to that boundary, and its boundary value is the negative of the rescaled projection polynomial. On the plane $\theta=-\frac{\xi_3}{k_j}$ it has the form
\begin{align}
\widehat h_j\left(-\frac{\xi_3}{k_j},\xi,\eta\right)&=q_j(\xi,\eta):=q_j(0,0)+\ell_j\cdot \eta+\eta\cdot M_j\eta+e_j\xi_3. \label{lin:incoming-polynomial}
\end{align}
By continuity, \eqref{lin:growth} bounds $q_j$ on fixed boundary patches. The parametrization $(\xi,\eta)\mapsto(-\xi_3/k_j,\xi,\eta)$ covers a fixed open box for large $j$, so polynomial norm equivalence bounds all its coefficients uniformly.

These coefficient bounds give a uniform kinetic H\"older seminorm for the incoming boundary values. Let $\widehat Z,\widehat Z'$ lie on $\{\theta=-\frac{\xi_3}{k_j}\}$ with $\abs{\eta},\abs{\eta'}\leq\varrho$ and $\dkin(\widehat Z,\widehat Z')\leq1$. Put $\Delta\theta:=\theta-\theta'$ and $\Delta\xi:=\xi-\xi'$. Definition~\eqref{eq:kinetic-distance} gives
\begin{align}
\abs{\eta-\eta'}\lesssim\dkin(\widehat Z,\widehat Z'),\qquad \abs{\Delta\theta}\lesssim\bigl(\dkin(\widehat Z,\widehat Z')\bigr)^2, \abs{\Delta\xi-\Delta\theta\,\eta'}\lesssim\bigl(\dkin(\widehat Z,\widehat Z')\bigr)^3. \label{lin:incoming-kinetic-increments}
\end{align}
Since $\Delta\xi_3=-k_j\Delta\theta$, it follows that
\begin{align}
(k_j-\varrho)\abs{\Delta\theta}&\lesssim\bigl(\dkin(\widehat Z,\widehat Z')\bigr)^3,\qquad \abs{\Delta\xi_3}\lesssim\bigl(\dkin(\widehat Z,\widehat Z')\bigr)^3\qquad(k_j\geq2\varrho). \label{eq:I:display:135}
\end{align}
Consequently,
\begin{align}
\abs[\big]{q_j(\xi,\eta)-q_j(\xi',\eta')}&\lesssim_\varrho\Bigl(\dkin(\widehat Z,\widehat Z')+\bigl(\dkin(\widehat Z,\widehat Z')\bigr)^3\Bigr). \label{lin:incoming-polynomial-increment}
\end{align}
This bound is uniform in $j$ and controls the $\Ckin^{\frac{\beta_0}{2}}$ seminorm, where $\beta_0$ is the exponent in Proposition~\ref{lin:boundary-inputs}\ref{lin:rough}. For pairs at distance greater than one, \eqref{lin:growth} bounds the difference of incoming values. Thus \eqref{lin:boundary-rough-hypotheses} holds uniformly on each fixed boundary patch. Formula~\eqref{lin:rough-bound} gives uniform boundary H\"older control of $\widehat h_j$ also at incoming points.

To place the functions on a fixed domain, extend $\widehat h_j$ across its moving plane by evaluation at $\left(\max\{\theta,-\xi_3/k_j\},\xi,\eta\right)$ in the incoming case, and use the minimum in the outgoing case. These maps have uniformly bounded Euclidean Lipschitz constants on compact sets. The kinetic H\"older bounds give a Euclidean H\"older modulus on every bounded set. A diagonal subsequence of the extensions therefore converges uniformly on compact sets to a continuous function whose restriction to $\{\pm\theta\geq0\}$ we denote by $\widehat h$. Here and below, the positive sign is for incoming points and the negative sign is for outgoing points. This convergence reaches $\theta=0$ and preserves the nonzero lower bound and polynomial growth in \eqref{lin:growth}. In particular, for every $\varrho\geq1$,
\begin{align}
\widehat h\in C(\{\pm\theta\geq0\}),\qquad \tnorm[\big]{\widehat h}_{L^\infty(\Qc_1\cap\{\pm\theta>0\})}\geq\frac{1}{2},\qquad \tnorm[\big]{\widehat h}_{L^\infty(\Qc_\varrho\cap\{\pm\theta>0\})}\lesssim\varrho^{2+\beta}. \label{lin:compactness-limit}
\end{align}

\proofstep{Step 5: The local energy bound and the limiting trace.}
We will prove that $\widehat h$ solves the constant-coefficient weak equation \eqref{lin:limit-equation} and, at incoming points, attains the initial value in \eqref{lin:limiting-initial-identity}. The diffusion term requires a uniform local $L^2$ bound for $\nabla_\eta\widehat h_j$ up to the moving plane. Let $\chi$ be a fixed smooth cutoff in a compact cylinder in the rescaled variables. In \eqref{trace:beta-Green}, choose a $C^2$ function with bounded first and second derivatives that agrees with $\frac{s^2}{2}$ on the bounded range of $\widehat h_j$ near $\supp\chi$. Apply the localized identity before translation and dilation, on the fixed spatial half-space, choosing the cutoff whose rescaled form is $\chi^2$. Its hypotheses follow from \eqref{lin:boundary-local-hypotheses}, the bounded coefficients and source, and the boundary continuity proved in Step~4. Changing variables justifies the energy test $\chi^2\widehat h_j$ on $\mathcal U_j$. Ellipticity and Young's inequality give
\begin{align}
\int_{\mathcal U_j}\chi^2\abs[\big]{\nabla_\eta\widehat h_j}^2&\lesssim_\chi\tnorm[\big]{\widehat h_j}_{L^\infty(\mathcal U_j\cap\supp\chi)}^2+\tnorm[\Big]{\widehat S_j^{(s)}}_{L^\infty(\mathcal U_j\cap\supp\chi)}^2 \label{lin:blowup-boundary-energy}\\
&\quad+\int_{\{\xi_3=-k_j\theta\},\,\mathrm{incoming}}\chi^2\abs{q_j}^2\abs{k_j+\eta_3}\,\dd\theta\dd\xi'\dd\eta.\notag
\end{align}
The outgoing energy term is nonnegative on the left and has been dropped. Here the spacetime transport flux on the moving plane has density $-(k_j+\eta_3)\dd\theta\dd\xi'\dd\eta$. In the incoming case, parametrizing the plane by $\xi_3$ gives
\begin{align}
(k_j+\eta_3)\abs{\dd \theta}&=\left(1+\frac{\eta_3}{k_j}\right)\abs{\dd \xi_3}. \label{eq:I:display:136}
\end{align}
The last density is uniformly bounded on fixed velocity patches. The incoming polynomial bounds therefore control the last term in \eqref{lin:blowup-boundary-energy}. In the outgoing case the boundary is outgoing on each fixed patch for large $j$. Together with \eqref{lin:growth} and \eqref{lin:constant-limit}, this proves the required uniform local energy bound, including sets touching the limiting time plane.

Step~4 and the almost-everywhere convergence of the domain indicators give
\begin{align}
\mathbf1_{\mathcal U_j}\widehat h_j&\to\mathbf1_{\{\pm\theta>0\}}\widehat h\quad\text{strongly in }L^2_{\mathrm{loc}}(\R^7). \label{lin:zero-extension-convergence}
\end{align}
The indicators are independent of $\eta$, and hence
\begin{align}
\nabla_\eta\bigl(\mathbf1_{\mathcal U_j}\widehat h_j\bigr)=\mathbf1_{\mathcal U_j}\nabla_\eta\widehat h_j. \label{lin:extended-gradient}
\end{align}
The bound \eqref{lin:blowup-boundary-energy}, together with \eqref{lin:zero-extension-convergence} and \eqref{lin:extended-gradient}, identifies the weak $L^2$ limit of these gradients as the gradient of the zero-extended limit. Thus $\nabla_\eta\widehat h$ belongs to local $L^2$ up to $\theta=0$. For smooth tests supported in $\{\pm\theta>0\}$, we can now pass to the limit in the weak equation. Strong convergence of the diffusion matrices and weak convergence of the velocity gradients pass the diffusion pairing. The first-order and zeroth-order terms vanish by \eqref{lin:constant-limit} and \eqref{lin:drift-scaling}, and \eqref{lin:zero-extension-convergence} passes the transport term. The source converges by \eqref{lin:constant-limit}. Together with the growth retained in Step~4, this gives
\begin{align}
\Yhat \widehat h-\locA_*:D_{\eta}^2\widehat h&=\widehat S_*, \qquad \tnorm[\big]{\widehat h}_{L^\infty(\Qc_\varrho\cap\{\pm \theta>0\})}\lesssim \varrho^{2+\beta}, \label{lin:limit-equation}
\end{align}
for every $\varrho\geq1$. We next identify the incoming initial value by allowing tests to meet $\theta=0$. At incoming points, pass to a subsequence so that the coefficients of \eqref{lin:incoming-polynomial} converge to those of a polynomial $q_*$. For a compactly supported smooth test function $\varphi$, the incoming term of the scaled identity \eqref{lin:boundary-weak-setting} satisfies
\begin{align}
\int q_j(\xi,\eta)\varphi\left(-\frac{\xi_3}{k_j},\xi,\eta\right)\left(1+\frac{\eta_3}{k_j}\right)\dd\xi\dd\eta\to\int q_*(\xi,\eta)\varphi(0,\xi,\eta)\dd\xi\dd\eta. \label{lin:incoming-trace-limit}
\end{align}
The volume terms converge by the bounds just established, which hold up to $\theta=0$. Combining their limits with \eqref{lin:incoming-trace-limit} in the scaled weak identity \eqref{lin:boundary-weak-setting} gives
\begin{align}
-\int_{\theta>0}\widehat h\Yhat\varphi+\int_{\theta>0}\locA_*\nabla_\eta\widehat h\cdot\nabla_\eta\varphi=\int_{\theta>0}\widehat S_*\varphi+\int q_*\varphi(0,\cdot,\cdot). \label{lin:limiting-initial-identity}
\end{align}
Uniform convergence up to the boundary also gives $\widehat h(\theta)\to q_*$ in local $L^2$ as $\theta\downarrow0$. Together with \eqref{lin:blowup-boundary-energy} and \eqref{lin:limit-equation}, this verifies the equation, growth and trace hypotheses used in Step~6.

\proofstep{Step 6: Classification and contradiction.}
We will apply the Liouville theorems to prove $\widehat h\in\mathcal P_2$, then use the projection orthogonality to obtain $\widehat h=0$, contradicting its nonzero size. The outgoing and incoming limits require different time domains. With $k=2$, \cite[Proposition~3.1]{Ros-Oton.Weidner2025} states that a weak solution of \eqref{lin:limit-equation} on $(-\infty,0)\times\R^6$, with constant elliptic $\locA_*$, constant $\widehat S_*$ and the growth bound there for every $\varrho\geq1$, belongs to $\mathcal P_2$. In particular, for constants $a,b\in\R$, $\ell\in\R^3$ and a symmetric matrix $H\in\R^{3\times3}$, the conclusion is
\begin{align}
\widehat h(\theta,\xi,\eta)=a+b\theta+\ell\cdot\eta+\eta\cdot H\eta\in\mathcal P_2. \label{lin:liouville-polynomial}
\end{align}
At outgoing points, the constant-coefficient equation and the growth of order $2+\beta<3$ required by this result are exactly \eqref{lin:constant-limit} and \eqref{lin:limit-equation}, so $\widehat h\in\mathcal P_2$.

With $k=2$, \cite[Proposition~3.3]{Ros-Oton.Weidner2025} gives the same conclusion on $(0,\infty)\times\R^6$ under the same equation and growth assumptions, provided $\widehat h,\nabla_\eta\widehat h\in L^2_{\mathrm{loc}}([0,\infty)\times\R^6)$ and its initial value is the restriction of a polynomial in $\mathcal P_2$, attained in local $L^2$ as $\theta\downarrow0$. At incoming points, the gradient bound and initial convergence were proved in \eqref{lin:blowup-boundary-energy}--\eqref{lin:limiting-initial-identity}. It remains to check that $q_*\in\mathcal P_2$.

The initial value is $q_*(\xi,\eta)=q_*(0,0)+\ell_*\cdot \eta+\eta\cdot M_*\eta+e_*\xi_3$. The growth bound in \eqref{lin:limit-equation} has a constant independent of $\varrho$ and extends to $\theta=0$ by continuity. Evaluating at $\eta=0$ and $\xi=\pm\frac{\varrho^3}{2}e_3$ gives $\abs{e_*}\varrho^3\lesssim\varrho^{2+\beta}$ for every $\varrho\geq1$. Since $\beta<1$, letting $\varrho\to\infty$ proves $e_*=0$. Thus $q_*\in\mathcal P_2$, so all hypotheses of the positive-time classification hold and $\widehat h\in\mathcal P_2$.

In both cases, Step~4 preserves the lower bound in \eqref{lin:growth}, and bounded convergence in \eqref{lin:orthogonality} gives
\begin{align}
\tnorm[\big]{\widehat h}_{L^\infty(\Qc_1(0)\cap\{\pm\theta>0\})}\geq\frac{1}{2},\qquad \int_{\Qc_1(0)\cap\{\pm\theta>0\}}\widehat hP\,\dd\widehat Z=0\quad(P\in\mathcal P_2). \label{lin:limit-contradiction}
\end{align}
Taking $P=\widehat h$ contradicts the first inequality. Consequently, the projection errors have a uniform bound of order $s^{2+\beta}$. The dyadic conclusion of \cite[Lemma~2.14]{Ros-Oton.Weidner2025}, using \eqref{lin:polynomial-coefficient-bound}, supplies one polynomial and its coefficient bounds as described in Step~2. Restoring the data norm and the radius $r$ proves \eqref{lin:boundary-expansion}--\eqref{lin:boundary-expansion-coefficients}.

\end{proof}

\begin{proof}[Proof of Proposition~\ref{lin:order-two}]
For each center $\widehat Z$, we will construct one polynomial $P_{\widehat Z}$ satisfying \eqref{eq:kinetic-holder-polynomial} for every evaluation point in the smaller half-cylinder. Steps~1--2 choose a nearby boundary point $\widehat Z_b$ and its polynomial $P_{\widehat Z_b}$. Step~3 estimates the interior remainder. Step~4 constructs $P_{\widehat Z}$ and controls $P_{\widehat Z}-P_{\widehat Z_b}$ to prove the approximation near $\widehat Z$. Step~5 proves the estimate for all pairs of points in the smaller half-cylinder and restores the scale. We follow the boundary-to-interior argument of \cite[proof of Lemma~4.4]{Ros-Oton.Weidner2025}, with the source estimates and the dependence on normal velocity given below.

\proofstep{Step 1: Normalization and the choice of the boundary point.}
We will reduce to unit radius and choose $\widehat Z_b$ so that the boundary expansion applies on a cylinder containing a smaller full cylinder about $\widehat Z$.
To reduce to unit radius while preserving $\{\xi_3>0\}$, use the map $T_{(\theta_0,(\xi_0',0),(\eta_0',0))}\circ S_r$ from \eqref{eq:I:display:108}--\eqref{eq:intro-kinetic-scaling}. It removes the tangential center velocity and changes the normal center velocity to $\frac{\eta_{0,3}}{r}$. More generally, composition with $T_{\widehat Z}\circ S_s$ transforms the coefficient and source tuple as
\begin{align}
(\locA,\locb,\locc,\widehat S)\mapsto(\locA,s\locb,s^2\locc,s^2\widehat S)\circ T_{\widehat Z}\circ S_s. \label{lin:cover:scaled-equation}
\end{align}
The scaled assumptions remain bounded by \eqref{lin:coefficient-norm}. We may therefore take $r=1$ and, by homogeneity, assume
\begin{align}
\tnorm[\big]{\widehat u}_{L^\infty(\Qplus_{16})}+\tnorm[\big]{\widehat S}_{\Ckin^\beta(\Qplus_{16})}\leq1. \label{lin:cover:normalization}
\end{align}
We will apply Lemma~\ref{lin:boundary-expansion-lemma} near the boundary with the fixed radius
\begin{align}
r_*:=\frac{1}{128\left(1+\abs{\eta_{0,3}}\right)}. \label{lin:cover:scale}
\end{align}
For an arbitrary interior point $\widehat Z=(\theta,\xi,\eta)\in\Qplus_1$, let $\delta=\delta(\widehat Z)$ be the largest radius for which $\Qc_\delta(\widehat Z)$ lies in $\{\xi_3>0\}$. The normal position in this cylinder is $\xi_3+s\eta_3+y_3$, with $\abs{s}<\delta^2$ and $\abs{y_3}<\delta^3$. Hence
\begin{align}
\xi_3=\delta^2\abs{\eta_3}+\delta^3,\qquad \delta>0. \label{eq:I:display:139}
\end{align}
The largest cylinder touches $\{\xi_3=0\}$ at the point
\begin{align}
\widehat Z_b:=\widehat Z\circ\left(-\operatorname{sgn}(\eta_3)\delta^2,-\delta^3e_3,0\right)\in\{\xi_3=0\}. \label{lin:nearest-boundary-point}
\end{align}
It satisfies $\widehat Z_b\in\Qc_{2\delta}(\widehat Z)$ and $\widehat Z\in\Qc_{2\delta}(\widehat Z_b)$. Figure~\ref{lin:interior-boundary-geometry} illustrates the interior cylinder and its nearby boundary center.

\begin{figure}[tbp]
\centering
\begin{tikzpicture}[x=1.8cm,y=1.4cm,>=Stealth]
\fill[blue!6] (-0.3,0) rectangle (2.4,3.1);
\node at (1.15,2.98) {$\Qplus_{\frac{5\delta}{2}}(\widehat Z_b)$};
\draw[->] (-0.3,0)--(2.6,0) node[right] {$\theta$};
\draw[->] (0,0)--(0,3.3) node[above] {$\xi_3$};
\draw[thick] (0,0)--(0,2)--(2,2.8)--(2,0.8)--cycle;
\node at (0.78,2.5) {$\Qc_\delta(\widehat Z)$};
\fill[blue!15] (0.75,1.175)--(0.75,1.425)--(1.25,1.625)--(1.25,1.375)--cycle;
\draw[blue!70!black,thick] (0.75,1.175)--(0.75,1.425)--(1.25,1.625)--(1.25,1.375)--cycle;
\fill (1,1.4) circle (1.3pt);
\node[anchor=west] at (1.29,1.25) {$\widehat Z$};
\fill (0,0) circle (1.3pt);
\node[anchor=north east] at (-0.02,-0.04) {$\widehat Z_b$};
\node at (1.65,1.96) {$\Qc_{\frac{\delta}{2}}(\widehat Z)$};
\draw[thin] (1.46,1.76)--(1.17,1.6);
\node[anchor=west] at (1.75,-0.22) {$\xi_3=0$};
\end{tikzpicture}
\caption{A fixed-velocity section with $\eta_3>0$: $\Qc_\delta(\widehat Z)$ touches $\{\xi_3=0\}$ at $\widehat Z_b$, and $\Qc_{\frac{\delta}{2}}(\widehat Z)\subset\Qplus_{\frac{5\delta}{2}}(\widehat Z_b)$. Light shading shows only the visible portion of the larger boundary cylinder.}
\label{lin:interior-boundary-geometry}
\end{figure}
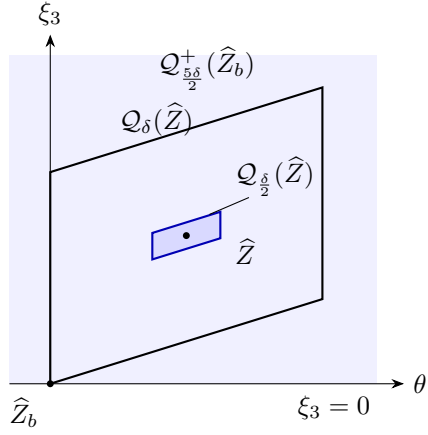

The product law \eqref{eq:I:display:108} gives
\begin{align}
\Qc_{s_2}(\widehat Z)\subset\Qc_{s_1+s_2}(\widehat Z_0)\qquad\text{if }\widehat Z\in\Qc_{s_1}(\widehat Z_0). \label{lin:cylinder-composition}
\end{align}
The time, velocity and spatial increments are bounded by $s_1^2+s_2^2$, $s_1+s_2$ and $s_1^3+s_2^3+s_1s_2^2\leq(s_1+s_2)^3$, respectively. If $\delta\leq\frac{r_*}{8}$, this yields
\begin{align}
\Qc_{16r_*}(\widehat Z_b)\subset\Qc_{1+17r_*}(\widehat Z_0)\subset\Qc_2(\widehat Z_0),\qquad \Qc_{64r_*}(\widehat Z_b)\subset\Qc_{1+65r_*}(\widehat Z_0)\subset\Qc_2(\widehat Z_0). \label{lin:boundary-radius-inclusions}
\end{align}
The same inclusions hold for any boundary center $\widehat Z_b\in\Qc_1(\widehat Z_0)$. The first inclusion in \eqref{lin:boundary-radius-inclusions} supplies the equation domain, and the second verifies \eqref{lin:second-order-nongrazing}. When $\delta\leq\frac{r_*}{8}$, Lemma~\ref{lin:boundary-expansion-lemma} therefore applies at $\widehat Z_b$ with radius $r_*$. We will use its polynomial to construct the interior polynomial at $\widehat Z$. For $\delta>\frac{r_*}{8}$, a full cylinder of fixed radius comparable to $r_*$ is available, and Step~4 will use the interior estimate directly.

\proofstep{Step 2: The boundary polynomial and its remainder.}
For $\delta\leq\frac{r_*}{8}$, we will construct $P_{\widehat Z_b}$ and prove \eqref{lin:interior-remainder-amplitude}, which bounds $\widehat u-P_{\widehat Z_b}$ on the full interior cylinder by $Cr_*^{-2-\beta}\delta^{2+\beta}$.

At each boundary center just described, Lemma~\ref{lin:boundary-expansion-lemma} gives one polynomial $P_{\widehat Z_b}\in\mathcal P_2$, independent of $\varrho$, with
\begin{align}
\tnorm[\big]{\widehat u-P_{\widehat Z_b}}_{L^\infty(\Qplus_\varrho(\widehat Z_b))}\lesssim r_*^{-2-\beta}\varrho^{2+\beta},\qquad \abs[\big]{p_{\widehat Z_b,\nu}}\lesssim r_*^{-\nu}\quad(\nu=0,1,2). \label{lin:cover:boundary}
\end{align}
Here $0<\varrho\leq r_*$. The scaled coefficient bounds follow from \eqref{lin:coefficient-norm} and $r_*\leq1$. The velocity-linear, time and velocity-quadratic terms of this polynomial give
\begin{align}
[\nabla_\eta P_{\widehat Z_b}]_{\Ckin^\beta(\Qplus_{r_*}(\widehat Z_b))}+[P_{\widehat Z_b}]_{\Ckin^\beta(\Qplus_{r_*}(\widehat Z_b))}\lesssim r_*^{-2}. \label{lin:cover:polynomial-holder}
\end{align}
The respective bounds are $r_*^{-1-\beta}$ and $r_*^{-\beta}$, both at most $r_*^{-2}$, by \eqref{lin:cover:boundary}, $\abs{\Delta\eta}\lesssim\dkin$ and $\abs{\Delta\theta}\lesssim\dkin^2$. For example, the quadratic velocity term contributes at most $Cr_*^{-2}(r_*)^{1-\beta}$ to the gradient seminorm. The same coefficient bounds give $\tnorm[\big]{P_{\widehat Z_b}}_{L^\infty}\lesssim1$ and $\tnorm[\big]{\nabla_\eta P_{\widehat Z_b}}_{L^\infty}\lesssim r_*^{-1}$ on $\Qplus_{r_*}(\widehat Z_b)$. The derivatives $D_\eta^2P_{\widehat Z_b}$ and $\Yhat P_{\widehat Z_b}$ are constant.

For $\delta\leq\frac{r_*}{8}$, \eqref{lin:cylinder-composition} gives $\Qc_{\frac{\delta}{2}}(\widehat Z)\subset\Qplus_{\frac{5\delta}{2}}(\widehat Z_b)$, with $\frac{5\delta}{2}<r_*$. Therefore
\begin{align}
\tnorm[\big]{\widehat u-P_{\widehat Z_b}}_{L^\infty(\Qc_{\frac{\delta}{2}}(\widehat Z))}\lesssim r_*^{-2-\beta}\delta^{2+\beta}. \label{lin:interior-remainder-amplitude}
\end{align}

\proofstep{Step 3: Interior Schauder estimates for $\widehat u-P_{\widehat Z_b}$.}
We will prove \eqref{lin:cover:scaled-interior-estimate}. Its seminorm bound will give a polynomial approximation to the remainder centered at $\widehat Z$, and its derivative bounds will control the coefficients of that polynomial. To apply the interior estimate, we must first bound both the oscillation and the value of the source in \eqref{lin:polynomial-source}.

The equation for the remainder is \eqref{lin:polynomial-source}. Separate the value of its source at $\widehat Z$ from its oscillation:
\begin{align}
\widehat S-\locOp P_{\widehat Z_b}=C_\delta+\widehat S_\delta,\qquad C_\delta:=\bigl(\widehat S-\locOp P_{\widehat Z_b}\bigr)(\widehat Z),\qquad \widehat S_\delta(\widehat Z)=0. \label{eq:I:display:141}
\end{align}
The product inequality, the coefficient assumptions \eqref{lin:coefficients}--\eqref{lin:coefficient-norm}, and \eqref{lin:cover:boundary}--\eqref{lin:cover:polynomial-holder} bound the source seminorm. Since $\widehat S_\delta(\widehat Z)=0$, they imply
\begin{align}
\tnorm[\big]{\widehat S_\delta}_{L^\infty(\Qc_{\frac{\delta}{2}}(\widehat Z))}+\delta^\beta[\widehat S_\delta]_{\Ckin^\beta(\Qc_{\frac{\delta}{2}}(\widehat Z))}\lesssim r_*^{-2-\beta}\delta^\beta. \label{lin:cover:source-oscillation}
\end{align}
To control $C_\delta$, apply Lemma~\ref{lin:constant-source} after composition with $T_{\widehat Z}\circ S_{\frac{\delta}{4}}$. This map sends $\Qc_2(0)$ onto $\Qc_{\frac{\delta}{2}}(\widehat Z)$; by \eqref{lin:cover:scaled-equation}, the source is multiplied by $\frac{\delta^2}{16}$ and the scaled coefficients remain uniformly bounded. Thus \eqref{lin:cover:source-oscillation} and \eqref{lin:interior-remainder-amplitude} give
\begin{align}
\delta^2\abs{C_\delta}&\lesssim\tnorm[\big]{\widehat u-P_{\widehat Z_b}}_{L^\infty(\Qc_{\frac{\delta}{2}}(\widehat Z))}+r_*^{-2-\beta}\delta^{2+\beta}\lesssim r_*^{-2-\beta}\delta^{2+\beta}. \label{lin:source-constant-scaled}
\end{align}
Combining \eqref{lin:interior-remainder-amplitude}, \eqref{lin:cover:source-oscillation} and \eqref{lin:source-constant-scaled} now bounds the amplitude and full scaled source. We apply Proposition~\ref{lin:interior-inputs}\ref{lin:interior} to the same rescaled remainder and return to $\Qc_{\frac{\delta}{4}}(\widehat Z)$. The remainder satisfies the weak equation \eqref{lin:polynomial-source}, and \eqref{lin:coefficient-norm} verifies the scaled coefficient assumptions. Estimate~\eqref{lin:interior-bound} therefore gives, with all norms now on $\Qc_{\frac{\delta}{4}}(\widehat Z)$,
\begin{align}
&\delta\tnorm[\big]{\nabla_\eta(\widehat u-P_{\widehat Z_b})}_{L^\infty}+\delta^2\left(\tnorm[\Big]{D_\eta^2(\widehat u-P_{\widehat Z_b})}_{L^\infty}+\tnorm[\big]{\Yhat(\widehat u-P_{\widehat Z_b})}_{L^\infty}\right) \label{lin:cover:scaled-interior-estimate}\\
&\qquad+\delta^{2+\beta}[\widehat u-P_{\widehat Z_b}]_{\Ckin^{2+\beta}}\lesssim r_*^{-2-\beta}\delta^{2+\beta}.\notag
\end{align}

\proofstep{Step 4: Construction of the interior polynomial and its local remainder.}
We will construct $P_{\widehat Z}$ centered at the arbitrary interior point $\widehat Z$ and prove the remainder estimate \eqref{lin:cover:local-centered-remainder}. We will also bound the coefficients of $P_{\widehat Z}-P_{\widehat Z_b}$, which will be needed to extend that estimate beyond the small interior cylinder.

Suppose first that $\delta\leq\frac{r_*}{8}$. By \eqref{lin:cover:scaled-interior-estimate}, the function $\widehat u-P_{\widehat Z_b}$ belongs to $\Ckin^{2+\beta}(\Qc_{\frac{\delta}{4}}(\widehat Z))$. The characterization \eqref{eq:kinetic-holder-polynomial} at center $\widehat Z$ supplies $q_{\widehat Z}\in\mathcal P_2$. Define
\begin{align}
P_{\widehat Z}:=P_{\widehat Z_b}+q_{\widehat Z},\qquad \widehat u-P_{\widehat Z}=\bigl(\widehat u-P_{\widehat Z_b}\bigr)-q_{\widehat Z}. \label{lin:cover:polynomial-construction}
\end{align}
By \eqref{eq:kinetic-holder-polynomial} and \eqref{lin:cover:scaled-interior-estimate}, this fixed polynomial satisfies the required estimate centered at $\widehat Z$:
\begin{align}
\abs[\big]{\widehat u(\widehat Z')-P_{\widehat Z}(\widehat Z')}\lesssim r_*^{-2-\beta}\bigl(\dkin(\widehat Z',\widehat Z)\bigr)^{2+\beta}\qquad\bigl(\widehat Z'\in\Qc_{\frac{\delta}{4}}(\widehat Z)\bigr). \label{lin:cover:local-centered-remainder}
\end{align}

We next identify and bound the coefficients of this polynomial. At an interior center, the approximation in \eqref{eq:kinetic-holder-polynomial} determines the coefficients uniquely: restricting to velocity increments gives the value, velocity gradient and velocity Hessian, and restricting to the transport line gives the transport derivative. The derivative regularity is also supplied by \cite[Lemma~2.12]{Ros-Oton.Weidner2025}. Since $P_{\widehat Z_b}\in\mathcal P_2$, its expansion about $\widehat Z$ is exact. Thus \eqref{lin:cover:polynomial-construction}, expressed about $\widehat Z=(\theta,\xi,\eta)$, is
\begin{align}
P_{\widehat Z}(\widehat Z')&=\widehat u(\widehat Z)+\nabla_\eta\widehat u(\widehat Z)\cdot(\eta'-\eta)+\Yhat\widehat u(\widehat Z)(\theta'-\theta)+\frac{1}{2}(\eta'-\eta)\cdot D_\eta^2\widehat u(\widehat Z)(\eta'-\eta). \label{lin:second-order-polynomial}
\end{align}
The constant coefficient of $q_{\widehat Z}$ is $(\widehat u-P_{\widehat Z_b})(\widehat Z)$, bounded by \eqref{lin:interior-remainder-amplitude}. Its remaining coefficients are the corresponding derivatives of $\widehat u-P_{\widehat Z_b}$ at $\widehat Z$, bounded by \eqref{lin:cover:scaled-interior-estimate}. Consequently,
\begin{align}
\abs[\big]{(P_{\widehat Z}-P_{\widehat Z_b})_{\widehat Z,\nu}}=\abs[\big]{(q_{\widehat Z})_{\widehat Z,\nu}}\lesssim r_*^{-2-\beta}\delta^{2+\beta-\nu}\qquad(\nu=0,1,2). \label{lin:cover:interior-jets}
\end{align}
Both polynomials in the difference are expressed about $\widehat Z$. Translating the boundary polynomial from $\widehat Z_b$ to $\widehat Z$ uses $\abs{\theta-\theta_b}\leq\delta^2$ and $\eta_b=\eta$. Combining \eqref{lin:cover:interior-jets} with \eqref{lin:cover:boundary} therefore gives $\abs[\big]{p_{\widehat Z,\nu}}\lesssim r_*^{-\nu}$.

If $\delta>\frac{r_*}{8}$, use $T_{\widehat Z}\circ S_{\frac{r_*}{32}}$. It maps $\Qc_2(0)$ into the full cylinder $\Qc_{\frac{r_*}{16}}(\widehat Z)\subset\Qplus_2(\widehat Z_0)$, and \eqref{lin:cover:scaled-equation} applies with $s=\frac{r_*}{32}$. The direct interior estimate gives $[\widehat u]_{\Ckin^{2+\beta}(\Qc_{\frac{r_*}{32}}(\widehat Z))}\lesssim r_*^{-2-\beta}$. Its polynomial at $\widehat Z$ is again \eqref{lin:second-order-polynomial}, with $\abs[\big]{p_{\widehat Z,\nu}}\lesssim r_*^{-\nu}$, and
\begin{align}
\abs[\big]{\widehat u(\widehat Z')-P_{\widehat Z}(\widehat Z')}\lesssim r_*^{-2-\beta}\bigl(\dkin(\widehat Z',\widehat Z)\bigr)^{2+\beta}\qquad\bigl(\widehat Z'\in\Qc_{\frac{r_*}{32}}(\widehat Z)\bigr). \label{lin:cover:fixed-interior-remainder}
\end{align}
At a boundary center $\widehat Z$, take $P_{\widehat Z}$ from \eqref{lin:cover:boundary}.

\proofstep{Step 5: Extension of the remainder estimate to the whole half-cylinder.}
We will prove \eqref{lin:cover:all-pairs} for every $\widehat Z,\widehat Z'\in\Qplus_1$, using the polynomial $P_{\widehat Z}$ already fixed in Step~4. The local estimates \eqref{lin:cover:local-centered-remainder} and \eqref{lin:cover:fixed-interior-remainder} cover evaluation points near the center. For the other points, we will use the boundary remainder and the coefficient difference \eqref{lin:cover:interior-jets}.

Let $\widehat Z,\widehat Z'\in\Qplus_1$ and put $\varrho=\dkin(\widehat Z',\widehat Z)$, using \eqref{eq:kinetic-distance}. The estimate is immediate if $\widehat Z'=\widehat Z$, so assume $\varrho>0$. First suppose that $\widehat Z$ is interior and $\varrho\leq\frac{r_*}{128}$. We have $\widehat Z'\in\Qc_{2\varrho}(\widehat Z)$. If $\delta>\frac{r_*}{8}$, then $2\varrho\leq\frac{r_*}{64}<\frac{r_*}{32}$; if $\delta\leq\frac{r_*}{8}$ and $\varrho\leq\frac{\delta}{8}$, then $2\varrho\leq\frac{\delta}{4}$. In either case, \eqref{lin:cover:local-centered-remainder} or \eqref{lin:cover:fixed-interior-remainder} gives
\begin{align}
\abs[\big]{\widehat u(\widehat Z')-P_{\widehat Z}(\widehat Z')}\lesssim r_*^{-2-\beta}\varrho^{2+\beta}. \label{lin:cover:interior-pair}
\end{align}
It remains to consider nearby evaluation points outside the range just covered. Then $\delta\leq\frac{r_*}{8}$ and $\delta<8\varrho$. The inclusion \eqref{lin:cylinder-composition} places both points within the range of the boundary approximation:
\begin{align}
\widehat Z'\in\Qc_{2\varrho}(\widehat Z)\subset\Qc_{2\delta+2\varrho}(\widehat Z_b)\subset\Qc_{18\varrho}(\widehat Z_b),\qquad18\varrho<r_*. \label{lin:cover:near-pair}
\end{align}
Using \eqref{lin:cover:boundary} at radius $18\varrho$ and \eqref{lin:cover:interior-jets}, we obtain
\begin{align}
\abs[\big]{\widehat u(\widehat Z')-P_{\widehat Z}(\widehat Z')}&\leq\abs[\big]{\widehat u(\widehat Z')-P_{\widehat Z_b}(\widehat Z')}+\abs[\big]{(P_{\widehat Z_b}-P_{\widehat Z})(\widehat Z')} \label{lin:cover:far-pair}\\
&\lesssim r_*^{-2-\beta}\left(\varrho^{2+\beta}+\sum_{\nu=0}^2\delta^{2+\beta-\nu}\varrho^\nu\right)\lesssim r_*^{-2-\beta}\varrho^{2+\beta}.\notag
\end{align}
The factors $1$, $\varrho$ and $\varrho^2$ correspond to the constant, velocity-linear, and time or velocity-quadratic terms. The last inequality uses $\delta<8\varrho$ and $2+\beta-\nu>0$. Thus the local polynomial at the interior center also gives the required remainder beyond its interior cylinder.

If $\widehat Z$ is a boundary point and $\varrho\leq\frac{r_*}{128}$, apply \eqref{lin:cover:boundary} on $\Qplus_{2\varrho}(\widehat Z)$. For every center, the remaining case $\varrho>\frac{r_*}{128}$ follows from the amplitude and polynomial coefficient bounds:
\begin{align}
\abs[\big]{\widehat u(\widehat Z')-P_{\widehat Z}(\widehat Z')}\lesssim1+\sum_{\nu=1}^2r_*^{-\nu}\varrho^\nu\lesssim r_*^{-2-\beta}\varrho^{2+\beta}. \label{lin:cover:separated-pair}
\end{align}
We have therefore proved, with one polynomial at each center and one constant for all pairs,
\begin{align}
\abs[\big]{\widehat u(\widehat Z')-P_{\widehat Z}(\widehat Z')}\lesssim r_*^{-2-\beta}\bigl(\dkin(\widehat Z',\widehat Z)\bigr)^{2+\beta}\qquad(\widehat Z,\widehat Z'\in\Qplus_1). \label{lin:cover:all-pairs}
\end{align}
This is \eqref{eq:kinetic-holder-polynomial} on the whole half-cylinder $\Qplus_1$. Together with the derivative and coefficient bounds from Steps~3--4, it gives
\begin{align}
\tnorm[\big]{\nabla_\eta\widehat u}_{L^\infty(\Qplus_1)}+\tnorm[\Big]{D_\eta^2\widehat u}_{L^\infty(\Qplus_1)}+\tnorm[\big]{\Yhat\widehat u}_{L^\infty(\Qplus_1)}+[\widehat u]_{\Ckin^{2+\beta}(\Qplus_1)}\lesssim\left(1+\abs{\eta_{0,3}}\right)^N. \label{lin:cover:normal-center}
\end{align}
Here $r_*^{-2-\beta}\lesssim(1+\abs{\eta_{0,3}})^N$ for an integer $N\geq2+\beta$. Restore the amplitude in \eqref{lin:cover:normalization} and undo the map in Step~1. Formula \eqref{lin:cover:scaled-equation} gives the source factors $r^2$ and $r^{2+\beta}$, while the velocity derivative, second velocity derivative and transport derivative acquire the factors $r$, $r^2$ and $r^2$. The normal velocity is $\frac{\eta_{0,3}}{r}$. This proves \eqref{lin:order-two-bound}.

\end{proof}

\para{Application to the Landau equation}\label{lin:subtract-lift}
Near an incoming boundary point, let $\glift$ be an extension of the incoming value $g$, and set
\begin{align}
u:=f-\glift,\qquad \trm\glift=g,\qquad \trm u=0. \label{lin:lift-subtraction}
\end{align}
The coefficients are still determined by $F=\m+\mh f$. Subtraction in \eqref{cb:scalar} gives
\begin{align}
\Y u-\Af{}:D_v^2u&=-\bigl(\Af{}v\bigr)\cdot\nabla_vu+\Cf{}u+S, \label{lin:physical-lifted-equation}\\
S&:=\Kop f+\sourceR-\Y\glift+\Af{}:D_v^2\glift-\bigl(\Af{}v\bigr)\cdot\nabla_v\glift+\Cf{}\glift. \label{lin:Landau-source}
\end{align}
Set $\widehat u:=u\circ\Psi_0$ and $\widehat S:=\vradius^2S\circ\Psi_0$. The change-of-variables calculation in the proof of Lemma~\ref{cb:normalized-equation}, applied to the linear equation \eqref{lin:physical-lifted-equation}, gives \eqref{lin:Landau-equation} and \eqref{lin:Landau-nondivergence} with $\widehat f$ replaced by $\widehat u$, the weak formulation \eqref{lin:normalized-weak} with source $\widehat S$, and zero incoming trace. Near an outgoing boundary point, we take $\glift=0$, so $u=f$.

For prescribed in-flow, Lemma~\ref{cb:inflow-lift} supplies the extension \eqref{eq:I:display:293}. For diffuse-reflection, Section~\ref{cb:section-graph} constructs a local extension in transformed variables through \eqref{lift:construction} and Lemma~\ref{lift:characteristic}. In that case, we will first estimate $f$ near outgoing boundary points, obtain regularity of its boundary flux, construct the extension, and apply Proposition~\ref{lin:order-two} to $\widehat u$.


\section{Regularity and Strong Solutions for the Landau Equation}\label{cb:section}

We will establish regularity for the weak equation \eqref{cb:equation-residual}, retaining the additional source $\sourceR$. We will first obtain H\"older bounds for $f$ and its coefficients, then Gaussian-weighted $L^\infty$ estimates for its first and second derivatives with grazing factors $1+\graz^{-\frac{1}{2}}$ and $1+\graz^{-\frac{3}{2}}$, where $\graz$ is defined in \eqref{eq:rho}. Sections~\ref{cb:section-initial} and \ref{lin:section-integration} will apply these estimates to the Landau solutions, extend them through the initial time under the compatibility conditions, and verify the $L^2$ integrability in \eqref{eq:I:display:016} required for strong solutions.

Section~\ref{cb:inflow-section} will obtain the incoming H\"older bound directly from the prescribed data; Section~\ref{cb:section-flux} will obtain it for diffuse-reflection by estimating the unknown outgoing flux. Both will then improve the regularity of $f$ and its coefficients. For second derivatives, Sections~\ref{cb:inflow-graph-section} and \ref{cb:section-graph} will construct boundary lifts and apply Proposition~\ref{lin:order-two} to the differences with zero incoming trace.

Throughout this section, $I\subset\R$ is an open time interval. For each compact interval $J\subset\subset I$, fix a compact interval $J'$ with $J\subset\operatorname{int}J'$ and $J'\subset\subset I$. We assume that $f$ satisfies \eqref{cb:equation-residual}, the Gaussian smallness condition \eqref{cb:setting}, and
\begin{align}
f\in C(I;L^2)\cap L^2_{\mathrm{loc}}(I;\Vspace). \label{cb:regularity-class}
\end{align}
The source satisfies \eqref{lin:source-Gaussian}; we will impose \eqref{lin:source-holder} only where needed. The kinetic trace is that of Lemma~\ref{trace:CM}, and each result will specify the boundary condition. The equation, \eqref{cvel:spatial-L}--\eqref{cvel:spatial-dual}, and the Gaussian source bound give $\Y f\in L^2(J';\Vspace')$, hence \eqref{lin:bulk-hypotheses}. The local cylinders used for estimates on $J$ will lie in $J'$. Constants denoted by $C_J$ or implicit in $\lesssim_J$ may depend on the stated bounds on $J'$ and on $\operatorname{dist}(J,\partial J')^{-1}$, and are independent of the velocity center $v_0$. For lifted differences, the coefficients will always be computed from the original $f$.

\subsection{In-flow: First-Order Estimates}\label{cb:inflow-section}

Starting from the Gaussian bound \eqref{cb:setting}, supplied for Landau solutions by Section~\ref{sec:energy}, we will verify the incoming-data hypothesis in Lemma~\ref{cb:inflow-incoming-regularity}. Lemma~\ref{cb:inflow-coefficients} will obtain $f\in\Ckin^{\beta_0}$ for a small $\beta_0>0$ and transfer this regularity to the coefficients by convolution. Proposition~\ref{cb:inflow-gradient} will then give $f\in\Ckin^{\frac12}$ and the weighted first-order estimate.

\begin{lemma}[Regularity of the incoming data]\label{cb:inflow-incoming-regularity}
Let $g\in C^5(I\times\overline\Om\times\R^3)$ satisfy the Gaussian derivative bound \eqref{in:eq:data-derivatives-g} on every compact subinterval of $I$. Then, for every $0<\varepsilon<1$ and $0<\alpha_g<\alpha$,
\begin{align}
\tnorm{g}_{\Ckin^{1+\varepsilon}\left((J\times\Sigma_-)\cap\{\abs{v-v_0}<\bracket{v_0}^{-1}\}\right)}\lesssim_J\ue^{-\alpha_g\abs{v_0}^2}\qquad(v_0\in\R^3). \label{cb:inflow-incoming-holder}
\end{align}
The implicit constant is independent of the velocity center $v_0$.
\end{lemma}

\begin{proof}
Fix an incoming point $Z=(t,x,v)$ in the indicated velocity ball and define the degree-one velocity polynomial
\begin{align}
P_Z(t',x',v'):=g(Z)+\nabla_vg(Z)\cdot(v'-v). \label{cb:inflow-incoming-polynomial}
\end{align}
If another incoming point $Z'=(t',x',v')$ has kinetic distance $r:=\dkin(Z',Z)\leq1$, then
\begin{align}
\abs{t'-t}\leq r^2,\qquad \abs{v'-v}\leq r,\qquad\abs{x'-x}\lesssim r^3+\bracket{v_0}r^2. \label{cb:inflow-incoming-increments}
\end{align}
First assume $r\leq c\bracket{v_0}^{-2}$, with $c>0$ fixed by a finite boundary covering. Use a local boundary parametrization and an orthonormal velocity frame whose third vector is the normal. Its incoming normal component lies in a half-line, so Taylor expansion in this half-space uses only incoming values. Along these paths, the physical velocities remain within $C\bracket{v_0}^{-1}$ of $v_0$, preserving the Gaussian bound with exponent $\alpha$. Returning to the physical velocity changes the expansion by a term bounded by a polynomial in $\bracket{v_0}$ times $\abs{x'-x}$. The derivative assumption and \eqref{cb:inflow-incoming-increments} therefore give
\begin{align}
\abs{g(Z')-P_Z(Z')}\lesssim_J\bracket{v_0}^N\ue^{-\alpha\abs{v_0}^2}r^2\lesssim_J\ue^{-\alpha_g\abs{v_0}^2}r^{1+\varepsilon}. \label{cb:inflow-incoming-remainder}
\end{align}
The polynomial coefficients have the same Gaussian bound. For the remaining pairs, $r\geq c\bracket{v_0}^{-2}$, use the endpoint amplitudes and these coefficients. Dividing by $r^{1+\varepsilon}$ costs at most $C\bracket{v_0}^4$, which is absorbed by decreasing the Gaussian exponent to $\alpha_g$. The polynomial characterization \eqref{eq:kinetic-holder-polynomial} proves \eqref{cb:inflow-incoming-holder}. Only two derivatives of the incoming data enter this argument.
\end{proof}

\begin{lemma}[H\"older estimates for the solution and its coefficients]\label{cb:inflow-coefficients}
Let $f$ be a weak solution of \eqref{cb:equation-residual} with incoming trace $\trm f=g$. Assume the solution class \eqref{cb:regularity-class}, the Gaussian bound \eqref{cb:setting}, and the source bound \eqref{lin:source-Gaussian}. Assume that \eqref{cb:inflow-incoming-holder} holds on every compact subinterval of $I$, for some $\varepsilon\in(0,1)$ and $\alpha_g>0$ on that interval. Then there are $0<\beta_0<\frac12$, $\beta\in(0,1)$, $\alpha_*>0$ and $N\geq0$ such that the following bounds hold uniformly for every $v_0\in\R^3$.
\begin{enumerate}
\item The solution has a continuous representative up to the boundary, agreeing with its kinetic trace, and
\begingroup\usetagform{enumerated}
\begin{align}
\tnorm{f}_{\Ckin^{\beta_0}\left(J\times\overline\Om\times B_{\bracket{v_0}^{-1}}(v_0)\right)}\lesssim_J\ue^{-\alpha_*\abs{v_0}^2}. \label{cb:inflow-rough-solution-holder}
\end{align}
\endgroup
\item With all norms taken on the set in part (i), the coefficients satisfy
\begingroup\usetagform{enumerated}
\begin{align}
\tnorm{\Af{}}_{\Ckin^\beta}+\tnorm{\Bf{}}_{\Ckin^\beta}+\tnorm{\Cf{}}_{\Ckin^\beta}+\tnorm{\nabla_v\Af{}}_{L^\infty}\lesssim_J\bracket{v_0}^{N}. \label{cb:inflow-coefficient-holder}
\end{align}
\endgroup
\item On the same set,
\begingroup\usetagform{enumerated}
\begin{align}
\tnorm{\Kop f}_{\Ckin^\beta}\lesssim_J\ue^{-\alpha_*\abs{v_0}^2}. \label{cb:inflow-K-holder}
\end{align}
\endgroup
If \eqref{lin:source-holder} also holds, then, after decreasing $\beta$ and $\alpha_*$,
\begingroup\usetagform{enumerated}
\begin{align}
\tnorm{\Kop f+\sourceR}_{\Ckin^\beta}\lesssim_J\ue^{-\alpha_*\abs{v_0}^2}. \label{cb:inflow-source-holder}
\end{align}
\endgroup
\end{enumerate}
\end{lemma}

\begin{proof}
We will first verify the hypotheses of Proposition~\ref{lin:boundary-inputs}\ref{lin:rough} and obtain a local H\"older estimate. We will then pass to the full set in part (i) by estimating nearby and separated pairs. The resulting increments of $f$ will give parts (ii)--(iii) by convolution.

\proofstep{The weak equation and the local H\"older estimate} Fix $Z_0=(t_0,x_0,v_0)$ with $t_0\in J$ and $d(x_0)<\frac{d_0}{2}$. Use the flattening map $\Psi_0$ of Lemma~\ref{cb:normalization}, and set $\widehat f=f\circ\Psi_0$ and $\widehat g=g\circ\Psi_0$ on the incoming boundary. Lemma~\ref{cb:normalized-equation} gives the equation \eqref{lin:Landau-equation}, with source \eqref{lin:normalized-f-source}, in the weak sense \eqref{lin:boundary-weak-setting}, and gives $\trm\widehat f=\widehat g$. The solution class \eqref{cb:regularity-class} supplies the local $L^2$ velocity derivative in \eqref{lin:boundary-local-hypotheses}. Lemma~\ref{cb:normalized-bounds} supplies the ellipticity and bounded coefficients in \eqref{lin:normalized-coefficient-bounds}, together with the Gaussian amplitude and source bound \eqref{cb:Gaussian-bound}. Thus the common assumptions \eqref{lin:local-ellipticity}, \eqref{lin:interior-rough-hypotheses} and \eqref{lin:boundary-local-hypotheses} are verified.

Choose $0<\beta_0<\frac12$ below the exponent allowed by Proposition~\ref{lin:boundary-inputs}\ref{lin:rough}. Only the incoming $\Ckin^{\beta_0}$ norm is needed here. The lower-order norm from \eqref{cb:inflow-incoming-holder} and the distance estimate \eqref{lin:kinetic-change-bound} give
\begin{align}
\tnorm[\big]{\widehat g}_{\Ckin^{\beta_0}(\partial_-\Qplus_2)}\lesssim_J\bracket{v_0}^N\ue^{-\alpha_g\abs{v_0}^2}. \label{cb:inflow-normalized-rough-datum}
\end{align}
Here all physical velocities belong to $B_{\bracket{v_0}^{-1}}(v_0)$ by \eqref{eq:I:display:265}; the constants in passing to the lower incoming exponent cost at most a power of $\bracket{v_0}$. We will repeatedly absorb such factors using
\begin{align}
\la^N\ue^{-\alpha_1\abs{v}^2}\lesssim_{N,\alpha_1,\alpha_2}\ue^{-\alpha_2\abs{v}^2},\qquad N\geq0,\quad0<\alpha_2<\alpha_1. \label{lin:Gaussian-absorption}
\end{align}
Proposition~\ref{lin:boundary-inputs}\ref{lin:rough}, with incoming condition \eqref{lin:boundary-rough-hypotheses} verified by \eqref{cb:inflow-normalized-rough-datum}, now gives
\begin{align}
\tnorm[\big]{\widehat f}_{\Ckin^{\beta_0}(\Qplus_1)}&\lesssim_J\tnorm[\big]{\widehat f}_{L^\infty(\Qplus_2)}+\tnorm[\big]{\widehat S_f}_{L^\infty(\Qplus_2)}+\tnorm[\big]{\widehat g}_{\Ckin^{\beta_0}(\partial_-\Qplus_2)}\lesssim_J\ue^{-\alpha_*\abs{v_0}^2}. \label{cb:inflow-normalized-holder}
\end{align}
The coefficient bounds are independent of $v_0$, so $\beta_0$ is also independent of $v_0$. Proposition~\ref{lin:boundary-inputs}\ref{lin:boundary-trace-identification} identifies the continuous boundary values with the kinetic trace. For centers with $d(x_0)\geq\frac{d_0}{2}$, apply the interior linear change of velocity from Lemma~\ref{cb:normalization} and Proposition~\ref{lin:interior-inputs}\ref{lin:rough-full}. The same coefficient and source bounds verify \eqref{lin:common-local-hypotheses}, \eqref{lin:local-ellipticity} and \eqref{lin:interior-rough-hypotheses}, giving the corresponding interior estimate. All time intervals lie in $J'$.

\proofstep{The H\"older estimate in physical variables} The local continuous representatives agree on overlapping interiors because they agree almost everywhere with $f$; continuity then identifies their boundary limits. They therefore define the representative asserted in part (i), with trace already identified in the first step. Fix $v_0$ and consider any two points of the set in part (i). At a center $Z=(t,x,w)$ in this set, we have $\bracket{w}\simeq\bracket{v_0}$ and $\abs{\abs{w}^2-\abs{v_0}^2}\leq3$. We first specify a radius on which the preceding local estimate applies uniformly at every such center.

Choose finitely many boundary neighborhoods with smaller neighborhoods covering $\partial\Om$, so that the flattening maps have uniform bounds and a fixed margin inside their domains. Write $\Psi_Z$ for the map centered at $Z$, with scale $c_J\bracket{w}^{-N}$, and $\widehat Z=\Psi_Z^{-1}(Z)$. For two physical points at kinetic distance $s\leq1$, the definition \eqref{eq:kinetic-distance} gives
\begin{align}
\abs{t'-t}\leq s^2,\qquad\abs{x'-x}\leq s^3+C\bracket{v_0}s^2,\qquad\abs{w'-w}\leq s. \label{cb:inflow-physical-increments}
\end{align}
Consequently, for $s$ smaller than a fixed negative power of $\bracket{v_0}$, both spatial points lie in the larger boundary neighborhood and the time interval lies in $J'$. In the explicit inverse map \eqref{cb:normalizing-map}, the time and velocity differences are bounded by the corresponding terms in \eqref{cb:inflow-physical-increments}. For the spatial characteristic difference, \eqref{lin:inverse-position-cancellation} bounds its numerator by $C\bigl(s^3+(s^3+\bracket{v_0}s^2)^2\bigr)$. Taking the required square and cube roots, using $s\leq1$ and using $\vradius^{-1}\lesssim_J\bracket{v_0}^N$ to absorb the inverse scale into a power of $\bracket{v_0}$, we obtain
\begin{align}
\dkin\bigl(\Psi_Z^{-1}(Z'),\widehat Z\bigr)\leq C_J\bracket{v_0}^{M}\dkin(Z',Z). \label{cb:inflow-physical-inverse-distance}
\end{align}
This calculation uses the explicit inverse on the fixed spatial neighborhood, before verifying that the transformed points lie in $\Qplus_1(\widehat Z)$. Keep the map $\Psi_Z$ fixed for this center. Choose $M$ and $C_J$ also to respect the spatial and time margins above. For the radius below, $s<3r_{v_0}$ makes the transformed distance less than $\frac12$, so
\begin{align}
r_{v_0}:=\frac{1}{6C_J}\bracket{v_0}^{-M},\qquad \Qc_{r_{v_0}}(Z)\cap\bigl(J\times\overline\Om\times\R^3\bigr)\subset\Psi_Z\bigl(\Qc_{\frac12}(\widehat Z)\cap\{\xi_3\geq0\}\bigr). \label{cb:inflow-uniform-local-radius}
\end{align}
On the interior region the same inclusion holds with the interior linear map and the full cylinder. Compactness is used only to fix the finitely many spatial neighborhoods and their uniform margins; no count of velocity-dependent neighborhoods enters the estimate.

Let $s=\dkin(Z',Z)$. If $s<r_{v_0}$, inclusion \eqref{cb:inflow-uniform-local-radius} places both points in the domain of the local estimate centered at $Z$. Apply \eqref{cb:inflow-normalized-holder} and the inverse distance bound \eqref{lin:inverse-kinetic-change-bound}. If $s\geq r_{v_0}$, use the Gaussian amplitude \eqref{cb:setting}. The two bounds are
\begin{align}
\abs{f(Z')-f(Z)}&\lesssim_J\ue^{-a\abs{v_0}^2}\bracket{v_0}^{M\beta_0}s^{\beta_0}\qquad(s<r_{v_0}), \label{cb:inflow-near-pair}\\
\frac{\abs{f(Z')-f(Z)}}{s^{\beta_0}}&\lesssim_J\ue^{-a\abs{v_0}^2}r_{v_0}^{-\beta_0}\qquad(s\geq r_{v_0}), \label{cb:inflow-separated-pair}
\end{align}
for some $a>0$. Absorb the polynomial factors by \eqref{lin:Gaussian-absorption} and take the supremum over both points. This proves \eqref{cb:inflow-rough-solution-holder}. 

\proofstep{Increments for the convolution estimates} At a fixed velocity $w$ and for $\Delta:=\abs{t-s}+\abs{x-y}\leq1$, the kinetic distance satisfies
\begin{align}
\dkin\bigl((t,x,w),(s,y,w)\bigr)\lesssim\bracket{w}^{\frac13}\Delta^{\frac13}. \label{cb:inflow-fixed-velocity-distance}
\end{align}
Use \eqref{cb:inflow-rough-solution-holder} with velocity center $w$ and absorb the resulting factor $\bracket{w}^{\frac{\beta_0}{3}}$ into the Gaussian. After decreasing $\alpha_*$, we have
\begin{align}
\abs[\big]{f(t,x,w)-f(s,y,w)}\lesssim_J\ue^{-\alpha_*\abs{w}^2}\Delta^{\frac{\beta_0}{3}}. \label{eq:I:display:296}
\end{align}
For a velocity increment $\abs{h}\leq1$, the same estimate applies directly when $\abs{h}<\bracket{w}^{-1}$. For larger increments, use the Gaussian amplitude and $\abs{h}^{-\beta_0}\leq\bracket{w}^{\beta_0}$. Since $\abs{w+h}^2\geq\frac12\abs{w}^2-1$, another decrease of the Gaussian exponent gives
\begin{align}
\abs[\big]{f(t,x,w+h)-f(t,x,w)}\lesssim_J\ue^{-\alpha_*\abs{w}^2}\abs{h}^{\beta_0}\qquad(\abs{h}\leq1). \label{cb:inflow-velocity-increments}
\end{align}
The smooth Maxwellian factors show that \eqref{eq:I:display:296} and \eqref{cb:inflow-velocity-increments} also hold for $F$ and $\mh f$, with possibly smaller Gaussian exponents.

\proofstep{The coefficient and source estimates} The kernels $a,b,\nabla a$ are locally integrable for every $-3\leq\gamma\leq1$. Splitting the integral into $\abs{v-w}<1$ and its complement gives
\begin{align}
\int_{\R^3}\bigl(\abs{a(v-w)}+\abs{b(v-w)}+\abs{\nabla a(v-w)}\bigr)\ue^{-\alpha_*\abs{w}^2}\dd w\lesssim\la^N. \label{cb:inflow-convolution-kernel-bound}
\end{align}
For $\gamma>-3$, the same bound holds with $\abs{c(v-w)}$ in the integrand. Apply \eqref{eq:I:display:296} to the difference of the data inside these convolutions. We obtain
\begin{align}
\abs[\big]{\Af{}(t,x,v)-\Af{}(s,y,v)}+\abs[\big]{b_F(t,x,v)-b_F(s,y,v)}+\abs[\big]{c_F(t,x,v)-c_F(s,y,v)}\lesssim_J\la^N\Delta^{\frac{\beta_0}{3}}. \label{cb:inflow-coefficient-time-space-increments}
\end{align}
For $\Delta\geq1$, the same bound follows from the coefficient amplitudes. At $\gamma=-3$, the $c_F$ term follows directly from $c_F=-8\pi F$ in \eqref{cvel:coefficient-Coulomb}.

For velocity increments, place the translation on the data:
\begin{align}
\Af{}(t,x,v+h)-\Af{}(t,x,v)=\int a(v-w)\bigl(F(t,x,w+h)-F(t,x,w)\bigr)\dd w. \label{eq:I:display:297}
\end{align}
Equations \eqref{cb:inflow-velocity-increments} and \eqref{cb:inflow-convolution-kernel-bound} give the bound $C_J\la^N\abs{h}^{\beta_0}$ for this difference and for the corresponding differences of $b_F,c_F$. Also, $\nabla_v\Af{}=(\nabla a)*F$ and \eqref{cb:inflow-convolution-kernel-bound} give its polynomial $L^\infty$ bound. For points in the velocity ball of part (i), a kinetic separation $s\leq1$ gives ordinary time and position separation at most $C\bracket{v_0}s^2+s^3$ by \eqref{cb:inflow-physical-increments}. Combining the two types of increments therefore gives kinetic H\"older bounds for $\Af{},b_F,c_F$. The formulas \eqref{cvel:scalar-coefficients} then prove \eqref{cb:inflow-coefficient-holder}, for example with any $0<\beta<\min\{\varepsilon,\frac{\beta_0}{6}\}$. Pairs with kinetic separation at least one are controlled by the amplitude bounds in Lemma~\ref{cvel:coefficients}.

Apply the same convolution estimates with data $\mh f$. In formula \eqref{cvel:K}, the convolutions have polynomial H\"older bounds, whereas the factors $\mh$ and $\mh(v\otimes v-I)$ have Gaussian H\"older bounds on the indicated velocity ball. The product inequality $[pq]_{\Ckin^\beta}\leq\tnorm{p}_{L^\infty}[q]_{\Ckin^\beta}+\tnorm{q}_{L^\infty}[p]_{\Ckin^\beta}$ and \eqref{lin:Gaussian-absorption} prove \eqref{cb:inflow-K-holder}. At $\gamma=-3$, its local term is $8\pi\m f$, with the same estimate. This implication from the rough H\"older bound of $f$ to the coefficient and $\Kop f$ bounds uses only \eqref{cb:inflow-rough-solution-holder} and \eqref{cb:setting}.

If \eqref{lin:source-holder} also holds, its seminorm bound and \eqref{lin:source-Gaussian} control the full H\"older norm of $\sourceR$. Adding \eqref{cb:inflow-K-holder} proves \eqref{cb:inflow-source-holder}, after decreasing $\beta,\alpha_*$ below the source exponents. The coefficient estimates require only \eqref{lin:source-Gaussian}.
\end{proof}

\begin{proposition}\label{cb:inflow-gradient}
Let $f$ be a weak solution of \eqref{cb:equation-residual} on $I$, with incoming trace $\trm f=g$, satisfying \eqref{cb:regularity-class}, \eqref{cb:setting} and \eqref{lin:source-Gaussian}. Assume that \eqref{cb:inflow-incoming-holder} holds on every compact subinterval of $I$, for some $\varepsilon\in(0,1)$ and $\alpha_g>0$ on that interval. Then, for some $\alpha_*\in(0,\alpha)$, the following holds.
\begin{enumerate}
\item The solution has a continuous representative on $J\times\overline\Om\times\R^3$, agreeing with its kinetic trace, and
\begingroup\usetagform{enumerated}
\begin{align}
\tnorm{f}_{\Ckin^{\frac12}\left(J\times\overline\Om\times B_{\bracket{v_0}^{-1}}(v_0)\right)}\lesssim_J\ue^{-\alpha_*\abs{v_0}^2}\qquad(v_0\in\R^3). \label{cb:inflow-holder-norm}
\end{align}
\endgroup
\item The distributional velocity derivative has a locally continuous representative away from grazing and satisfies
\begingroup\usetagform{enumerated}
\begin{align}
\tnorm{\frac{\ue^{\alpha_*\abs{v}^2}\nabla_vf}{1+\graz^{-\frac12}}}_{L^\infty\left(J\times\{d<\frac{d_0}{2}\}\times\R^3\right)}\lesssim_J1, \label{cb:pointwise-grazing}\\
\tnorm[\Big]{\ue^{\alpha_*\abs{v}^2}\nabla_vf}_{L^\infty\left(J\times\{d\geq\frac{d_0}{2}\}\times\R^3\right)}\lesssim_J1. \label{cb:inflow-interior-gradient}
\end{align}
\endgroup
\end{enumerate}
Here $\graz$ is the grazing weight defined in \eqref{eq:rho}.
\end{proposition}

\begin{proof}
Lemma~\ref{cb:inflow-coefficients} has supplied the preliminary H\"older bound of $f$ and the improved coefficient bound \eqref{cb:inflow-coefficient-holder}. We will use the latter in Proposition~\ref{lin:boundary-inputs}\ref{lin:sharp} to prove part (i). Apply the lemma on a compact enlargement of $J$ inside $J'$, so that its coefficient bounds hold on all cylinders used below. The resulting $\Ckin^{\frac12}$ estimate will then control the oscillation on small cylinders and yield part (ii) by Proposition~\ref{lin:boundary-inputs}\ref{lin:gradient}.

\proofstep{The incoming norm and the sharp H\"older estimate} Use the maps and weak formulation already verified in Lemma~\ref{cb:inflow-coefficients}. The coefficient estimate \eqref{cb:inflow-coefficient-holder} verifies \eqref{lin:physical-coefficient-holder}, so Lemma~\ref{cb:normalized-bounds} gives \eqref{lin:normalized-coefficient-holder} after decreasing the radius $\vradius$ in \eqref{lin:inverse-map} by a polynomial factor in $\bracket{v_0}$.

The first-order estimate will require the stronger incoming norm $\Ckin^{1+\beta}$, so we now transfer \eqref{cb:inflow-incoming-holder} at that order. At each incoming physical point $Z$, let $P_Z$ be the degree-one velocity polynomial in \eqref{eq:kinetic-holder-polynomial}. Its constant coefficient is $g(Z)$. Evaluating it at incoming velocities in a half-ball of radius comparable to $\bracket{v_0}^{-1}$ bounds its linear coefficients by $C_J\bracket{v_0}^N\ue^{-\alpha_g\abs{v_0}^2}$, using the amplitude and polynomial remainder in \eqref{cb:inflow-incoming-holder}.

For $Z=\Psi_0(\zeta)$, compose $P_Z$ with $\Psi_0$ and retain its constant and linear $\eta$ terms at $\zeta$. The physical velocity is linear in $\eta$ for fixed $\xi$. Its additional variation through $\xi$ is $O(\bracket{v_0}^Ns^2)$ when $s=\dkin(\zeta',\zeta)\leq1$, by \eqref{lin:inverse-map} and the cylinder increments. Thus the resulting polynomial $\widehat P_\zeta\in\mathcal P_1$ satisfies
\begin{align}
\abs[\big]{\widehat g(\zeta')-\widehat P_\zeta(\zeta')}\lesssim_J\bracket{v_0}^N\ue^{-\alpha_g\abs{v_0}^2}\left(s^{1+\varepsilon}+s^2\right)\lesssim_J\ue^{-\alpha_*\abs{v_0}^2}s^{1+\varepsilon}. \label{cb:inflow-incoming-pullback-remainder}
\end{align}
For each center use \eqref{cb:inflow-incoming-holder} on the physical velocity ball centered there. If a pair falls outside that ball, its separation is bounded below by a negative power of $\bracket{v_0}$, and the amplitude and polynomial-coefficient bounds give the same estimate. This verifies the polynomial characterization for $\widehat g$. Combining it with \eqref{cb:Gaussian-bound}, we have
\begin{align}
\tnorm[\big]{\widehat f}_{L^\infty(\Qplus_2)}+\tnorm[\big]{\widehat S_f}_{L^\infty(\Qplus_2)}+\tnorm[\big]{\widehat g}_{\Ckin^{1+\varepsilon}(\partial_-\Qplus_2)}\lesssim_J\ue^{-\alpha_*\abs{v_0}^2}. \label{cb:inflow-normalized-amplitudes}
\end{align}
Choose $0<\beta<\min\{\varepsilon,\frac12\}$ below the coefficient exponent. Equations \eqref{lin:normalized-coefficient-holder} and \eqref{cb:inflow-normalized-amplitudes} verify the diffusion and incoming assumptions \eqref{lin:boundary-sharp-hypotheses}. Proposition~\ref{lin:boundary-inputs}\ref{lin:sharp} gives
\begin{align}
\tnorm[\big]{\widehat f}_{\Ckin^{\frac12}(\Qplus_{\frac12}(\widehat Z_0))}\lesssim_J\ue^{-\alpha_*\abs{v_0}^2}. \label{cb:inflow-normalized-sharp}
\end{align}
For interior cylinders the same estimate has an empty incoming set. Apply the nearby- and separated-pair estimates \eqref{cb:inflow-near-pair}--\eqref{cb:inflow-separated-pair}, with exponent $\frac12$ and a smaller uniform local radius if necessary. They convert \eqref{cb:inflow-normalized-sharp} to \eqref{cb:inflow-holder-norm}. The continuous representative and trace were established in Lemma~\ref{cb:inflow-coefficients}.

\proofstep{The first-order estimate on cylinders avoiding grazing} Take $\zeta=(\theta,\xi,\eta)\in\Qplus_{\frac14}(\widehat Z_0)$ and define
\begin{align}
\widehat\graz:=\xi_3^{\frac13}+\abs{\eta_3},\qquad \zeta_*:=\bigl(\theta,(\xi',0),(\eta',0)\bigr). \label{eq:I:display:298}
\end{align}
The set $\gamma_0$ in \eqref{lin:boundary-sets} is the flat grazing set. With the orientation in \eqref{eq:kinetic-distance},
\begin{align}
\inf_{\zeta_b\in\gamma_0}\dkin(\zeta,\zeta_b)=\dkin(\zeta,\zeta_*)=\widehat\graz. \label{cb:inflow-distance-to-grazing}
\end{align}
At a grazing base point the normal position and velocity increments are $\xi_3$ and $\eta_3$, so their contributions to the distance are at least $\xi_3^{\frac13}$ and $\abs{\eta_3}$. The point $\zeta_*$ attains both lower bounds and has zero time and tangential increments. Thus the power $\frac13$ comes from the spatial scale $r^3$ in the kinetic dilation.

Fix a sufficiently small geometric constant $\delta_0>0$. For $0<\widehat\graz<\delta_0$, choose $r=c\widehat\graz$ with a fixed small $c>0$. This radius is comparable to the distance from grazing. The cylinder composition \eqref{lin:cylinder-composition} and \eqref{cb:inflow-distance-to-grazing} give, after decreasing $\delta_0$ and $c$,
\begin{align}
\zeta_*\in\overline{\Qplus_{\frac12}(\widehat Z_0)},\qquad\Qplus_{64r}(\zeta)\subset\Qplus_{\frac12}(\widehat Z_0),\qquad\Qc_{64r}(\zeta)\cap\{\xi_3=\eta_3=0\}=\varnothing. \label{cb:inflow-grazing-cylinder-geometry}
\end{align}
To verify that the enlarged cylinder avoids grazing, if $\abs{\eta_3}\geq\xi_3^{\frac13}$, then $64r\leq\frac12\abs{\eta_3}$ for small $c$, so the enlarged cylinder has nonzero normal velocity. Otherwise its normal spatial displacement is bounded by $C(r^2\abs{\eta_3}+r^3)\leq Cc^2\xi_3<\frac12\xi_3$, so it is interior. Moreover, for $\zeta'\in\Qplus_{16r}(\zeta)$ the same cylinder increments give
\begin{align}
\dkin(\zeta',\zeta_*)\lesssim\widehat\graz,\qquad \abs{\eta_3}/r\leq c^{-1}. \label{cb:inflow-grazing-distance}
\end{align}
These inclusions specify the domain on which we will use \eqref{cb:inflow-normalized-sharp}.

Set $h:=\widehat f-\widehat f(\zeta_*)$. The equation and incoming trace are
\begin{align}
\Yhat h-\divv_\eta(\locA\nabla_\eta h)&=\locB\cdot\nabla_\eta h+\locc h+\widehat S_f+\locc\widehat f(\zeta_*), \label{cb:inflow-constant-subtraction}\\
\trm h&=\widehat g-\widehat f(\zeta_*).\notag
\end{align}
The sharp estimate, \eqref{cb:inflow-grazing-distance}, and $r=c\widehat\graz$ give the small amplitude. The bounded coefficients and \eqref{cb:inflow-normalized-amplitudes} control the source and incoming seminorm:
\begin{align}
\tnorm{h}_{L^\infty(\Qplus_{16r}(\zeta))}&\lesssim_J\ue^{-\alpha_*\abs{v_0}^2}r^{\frac12}, \label{cb:inflow-first-order-oscillation}\\
\tnorm[\big]{\widehat S_f+\locc\widehat f(\zeta_*)}_{L^\infty(\Qplus_{16r}(\zeta))}+[\trm h]_{\Ckin^{1+\beta}(\partial_-\Qplus_{16r}(\zeta))}&\lesssim_J\ue^{-\alpha_*\abs{v_0}^2}. \label{cb:inflow-first-order-data}
\end{align}
Subtracting the constant does not change the incoming seminorm. Thus \eqref{lin:boundary-gradient-hypotheses} and the common weak-solution assumptions of Proposition~\ref{lin:boundary-inputs}\ref{lin:gradient} hold for $h$. The estimate gives
\begin{align}
r\tnorm[\big]{\nabla_\eta\widehat f}_{L^\infty(\Qplus_r(\zeta))}\lesssim_J\ue^{-\alpha_*\abs{v_0}^2}\bigl(r^{\frac12}+r^2+r^{1+\beta}\bigr)\lesssim_J\ue^{-\alpha_*\abs{v_0}^2}r^{\frac12}. \label{cb:inflow-scaled-gradient}
\end{align}
The three powers on the right come from the amplitude, the source and the incoming seminorm, respectively. Dividing by $r$ shows that the $\Ckin^{\frac12}$ oscillation in \eqref{cb:inflow-first-order-oscillation} gives $r^{-\frac12}\simeq\widehat\graz^{-\frac12}$. The source and incoming contributions are smaller for $r\leq1$. The estimate also gives local $\Ckin^{1+\beta}$ regularity of $\widehat f$ and $\Ckin^\beta$ regularity of $\nabla_\eta\widehat f$. On an interior cylinder the incoming seminorm is zero.

For $\widehat\graz\geq\delta_0$, choose the fixed radius $r=c\delta_0$. If $\abs{\eta_3}\geq\frac{\delta_0}{2}$, then $64r<\frac12\abs{\eta_3}$ keeps the enlarged cylinder away from zero normal velocity. Otherwise $\xi_3\geq\bigl(\frac{\delta_0}{2}\bigr)^3$, and its normal spatial displacement is at most $C(r^2\delta_0+r^3)<\frac12\xi_3$. In both cases \eqref{lin:cylinder-composition} also gives $\Qplus_{64r}(\zeta)\subset\Qplus_{\frac12}(\widehat Z_0)$ for small $c$. Proposition~\ref{lin:boundary-inputs}\ref{lin:gradient}, applied directly to $\widehat f$ with \eqref{cb:inflow-normalized-amplitudes}, now gives a uniform Gaussian derivative bound. Together with \eqref{cb:inflow-scaled-gradient}, this proves
\begin{align}
\abs[\big]{\nabla_\eta\widehat f}\lesssim_J\ue^{-\alpha_*\abs{v_0}^2}\bigl(1+\widehat\graz^{-\frac12}\bigr)\quad\text{on }\Qplus_{\frac14}(\widehat Z_0). \label{eq:I:display:300}
\end{align}

\proofstep{The physical first-order estimates} We first compare $\widehat\graz$ with the physical weight $\graz$ in \eqref{eq:rho}. The normal-position and normal-velocity identities \eqref{eq:I:display:279} give
\begin{align}
\widehat\graz=\frac{1}{\vradius}\left(s_n^{-\frac13}\bigl(d(x)\bigr)^{\frac13}+s_n^{-1}\abs[\Big]{n\bigl(\bdproj(x)\bigr)\cdot v}\right). \label{cb:inflow-grazing-comparison}
\end{align}
Consequently
\begin{align}
\min\{s_n^{-\frac13},s_n^{-1}\}\graz\leq\vradius\widehat\graz\leq\max\{s_n^{-\frac13},s_n^{-1}\}\graz. \label{cb:inflow-grazing-weight-comparison}
\end{align}
Since $\vradius=c_J\bracket{v_0}^{-N}\leq1$, \eqref{cb:matrix-bounds} implies, for $p\in\{\frac12,\frac32\}$,
\begin{align}
\frac{1+\widehat\graz^{-p}}{1+\graz^{-p}}\leq\max\{1,\vradius^p s_n^{\frac{p}{3}},\vradius^p s_n^p\}\lesssim_J\bracket{v_0}^N. \label{cb:inflow-grazing-weight-ratio}
\end{align}
The derivative identity \eqref{lin:physical-gradient-conversion} contributes another polynomial factor by \eqref{cb:matrix-bounds}. Combining it with \eqref{eq:I:display:300} and \eqref{cb:inflow-grazing-weight-ratio} for $p=\frac12$, we obtain on $\Psi_0(\Qplus_{\frac14}(\widehat Z_0))$
\begin{align}
\frac{\abs{\nabla_v f}}{1+\graz^{-\frac12}}\lesssim_J\bracket{v_0}^N\ue^{-\alpha_1\abs{v_0}^2}\lesssim_J\ue^{-\alpha_*\abs{v}^2},\qquad0<\alpha_*<\alpha_1. \label{cb:inflow-physical-local-first-order}
\end{align}
The last inequality uses \eqref{lin:Gaussian-absorption} and $\abs{\abs{v}^2-\abs{v_0}^2}\leq3$ from \eqref{eq:I:display:265}. Each point of the boundary collar can be the center $Z_0$, with constants independent of that center, so this proves \eqref{cb:pointwise-grazing}. On $d(x)\geq\frac{d_0}{2}$, the interior change of velocity and Proposition~\ref{lin:boundary-inputs}\ref{lin:gradient}, applied with empty incoming boundary, give the corresponding derivative bound without a grazing factor and proves \eqref{cb:inflow-interior-gradient}. Equivalently, take a countable subcover of these local domains and then their essential supremum; no sum of local bounds is involved. All cylinders have their time intervals in $J'$. The local continuity of the velocity derivative away from grazing is preserved by the smooth change of variables.
\end{proof}

\subsection{In-flow: Boundary Lift and Second-Order Estimates}\label{cb:inflow-graph-section}

We will subtract a lift of the prescribed incoming data and apply Proposition~\ref{lin:order-two} to the resulting zero incoming equation. Lemma~\ref{cb:inflow-coefficients} supplies the coefficient H\"older estimates, and Proposition~\ref{cb:inflow-gradient} supplies the $\Ckin^{\frac12}$ bound. We will estimate the source created by the lift, then use the vanishing incoming value to obtain a remainder of size $r^{\frac12}$ on cylinders of radius $r$ comparable to the distance from grazing. The second-order estimate will consequently have size $r^{-\frac32}$.

\para{The incoming lift} The Gaussian derivative bounds in \eqref{in:eq:data-derivatives-g} are imposed only on $\Sigma_-$. We will extend these controlled incoming values to the interior and to all velocities, obtaining $\trm\glift=g$ and the weighted derivative estimate \eqref{eq:I:display:294} needed after subtracting $\glift$.

\begin{lemma}[Lift of the incoming data]\label{cb:inflow-lift}
Let $g\in C^5(I\times\overline\Om\times\R^3)$ satisfy \eqref{in:eq:data-derivatives-g} on every compact subinterval of $I$. There is a linear lift $\glift$, supported in the boundary collar, such that
\begin{align}
\trm\glift=g. \label{cb:inflow-lift-trace}
\end{align}
Its mixed derivatives with at most two position derivatives and total order at most five are continuous up to the boundary. For every $J\subset\subset I$ and $0<\alpha_1<\alpha$,
\begin{align}
\sum_{j+\abs{\ell}+\abs{k}\leq5,\,\abs{\ell}\leq2}\tnorm[\big]{\partial_t^j\partial_x^\ell\partial_v^k\glift}_{L^\infty(J;L^{\infty,\alpha_1})}\lesssim_{J,\alpha_1}1. \label{eq:I:display:294}
\end{align}
Here $j\in\mathbb N_0$, $\ell,k\in\mathbb N_0^3$, and the implicit constant may depend on the finite bounds in \eqref{in:eq:data-derivatives-g}. The same construction applies at a finite time endpoint when the assumed derivatives of $g$ extend continuously to that endpoint.
\end{lemma}

\begin{proof}
Choose $\mathfrak c=(21,-70,105,-84,35,-6)$ and a smooth collar cutoff $\chi$ equal to one near zero. With $y=\bdproj(x)$ and $s=n(y)\cdot v$, define
\begin{align}
\glift(t,x,v):=\chi\bigl(d(x)\bigr)\sum_{j=1}^6\mathfrak c_jg\bigl(t,y,v-(j+1)s_+n(y)\bigr), \label{eq:I:display:293}
\end{align}
and set $\glift=0$ outside the collar. The identities
\begin{align}
\sum_{j=1}^6\mathfrak c_j(-j)^r=1\qquad(0\leq r\leq5) \label{cb:inflow-lift-coefficients}
\end{align}
give the incoming trace and match the normal-velocity derivatives through order five at $s=0$. In a local velocity frame with third vector $n(y)$, differentiating these identities in the other variables gives the stated mixed-derivative continuity. Only two derivatives of $n$ and $\bdproj$ are needed.

The reflected velocity satisfies
\begin{align}
\abs{v-(j+1)s_+n(y)}^2=\abs{v-sn(y)}^2+\bigl(s-(j+1)s_+\bigr)^2\geq\abs{v}^2. \label{cb:inflow-reflection-Gaussian}
\end{align}
The chain rule therefore bounds each derivative in \eqref{eq:I:display:294} by a finite sum of Gaussian data derivatives of total order at most five, multiplied by powers of $\la$. The geometric factors are bounded because $\partial\Om$ is $C^3$. Absorbing those powers by \eqref{lin:Gaussian-absorption} proves \eqref{eq:I:display:294}.
\end{proof}

\para{Weighted second-order estimates}

\begin{proposition}\label{cb:strong-upgrade}
Let $f$ be a weak solution of \eqref{cb:equation-residual} on $I$, satisfying \eqref{cb:regularity-class}, \eqref{cb:setting} and \eqref{lin:source-Gaussian}, with $\trm f=g$, where $g$ satisfies the hypotheses of Lemma~\ref{cb:inflow-lift}. Assume also \eqref{lin:source-holder}. Then, for some $\alpha_*>0$,
\begin{align}
\tnorm{\frac{\ue^{\alpha_*\abs{v}^2}\bigl(\abs{D_v^2f}+\abs{\Y f}\bigr)}{1+\graz^{-\frac32}}}_{L^\infty\left(J\times\{d<\frac{d_0}{2}\}\times\R^3\right)}\lesssim_J1, \label{cb:pointwise-grazing-hessian}\\
\tnorm{\ue^{\alpha_*\abs{v}^2}\bigl(\abs{D_v^2f}+\abs{\Y f}\bigr)}_{L^\infty\left(J\times\{d\geq\frac{d_0}{2}\}\times\R^3\right)}\lesssim_J1. \label{cb:inflow-interior-hessian}
\end{align}
The derivatives are distributional derivatives represented by functions on the indicated domains.
\end{proposition}

\begin{proof}
We will verify the weak equation and H\"older source for $u=f-\glift$, estimate its amplitude on the nongrazing cylinders, and apply Proposition~\ref{lin:order-two}. The last step will convert the estimates to physical derivatives, including the first-order term in the transport transformation.

\proofstep{The zero incoming equation and its source} Lemma~\ref{cb:inflow-incoming-regularity} verifies \eqref{cb:inflow-incoming-holder}, so Lemma~\ref{cb:inflow-coefficients} and Proposition~\ref{cb:inflow-gradient} apply. Choose the lift from Lemma~\ref{cb:inflow-lift} and use its derivative bound \eqref{eq:I:display:294}. The following argument uses only its weak derivatives through total order three, with at most two position derivatives; Lemma~\ref{cb:inflow-initial} will verify these bounds for the negative-time extension.

Set $u=f-\glift$. The equation and source are \eqref{lin:physical-lifted-equation}--\eqref{lin:Landau-source}. Choose the map $\Psi_0$ of Lemma~\ref{cb:normalization} so that both \eqref{lin:normalized-coefficient-holder} and \eqref{cb:inflow-normalized-sharp} hold, and define
\begin{align}
\widehat u:=(f-\glift)\circ\Psi_0,\qquad \widehat S:=\vradius^2S\circ\Psi_0. \label{cb:inflow-transformed-source}
\end{align}
The change of variables following \eqref{lin:Landau-source} gives
\begin{align}
\Yhat\widehat u-\divv_\eta(\locA\nabla_\eta\widehat u)=\locB\cdot\nabla_\eta\widehat u+\locc\widehat u+\widehat S,\qquad\trm\widehat u=0. \label{cb:inflow-transformed-equation}
\end{align}
This is the weak identity \eqref{lin:boundary-weak-setting}. The local $L^2$ velocity gradient follows from \eqref{cb:regularity-class} and the bound on $\nabla_v\glift$ in \eqref{eq:I:display:294}. The coefficient estimate \eqref{cb:inflow-coefficient-holder} and Lemma~\ref{cb:normalized-bounds} give \eqref{lin:normalized-coefficient-bounds} and \eqref{lin:normalized-coefficient-holder}. They verify the coefficient conditions \eqref{lin:coefficients}--\eqref{lin:coefficient-norm}, uniformly in $v_0$.

We estimate separately the lift terms in the source. For each $H\in\{\glift,\Y\glift,\nabla_v\glift,D_v^2\glift\}$, \eqref{eq:I:display:294} gives, on the physical neighborhood of $Z_0$ with $\abs{v-v_0}\leq\bracket{v_0}^{-1}$,
\begin{align}
\tnorm{H}_{L^\infty}+\tnorm{\nabla_{t,x,v}H}_{L^\infty}\lesssim_J\bracket{v_0}^N\ue^{-a\abs{v_0}^2}. \label{cb:inflow-lift-first-derivatives}
\end{align}
Differentiating $D_v^2\glift$ once uses three derivatives of the lift, while differentiating $\Y\glift=\partial_t\glift+v\cdot\nabla_x\glift$ once uses at most two position derivatives. The factor $v$ contributes only a polynomial velocity factor. These weak derivative bounds give locally Lipschitz representatives on the fixed spatial neighborhoods. For kinetic separation $q\leq1$, ordinary time, position and velocity increments are bounded by $q^2$, $q^3+C\bracket{v_0}q^2$ and $q$, respectively. Thus the same functions have local $\Ckin^\beta$ bounds for every $0<\beta<1$, with a polynomial velocity factor. Formula~\eqref{lin:kinetic-change-bound} transfers these bounds to $H\circ\Psi_0$.

The H\"older product inequality used in the proof of \eqref{cb:inflow-K-holder} gives, with all norms on $\Qplus_1$,
\begin{align}
\tnorm[\Big]{(\Af{}:D_v^2\glift)\circ\Psi_0}_{\Ckin^\beta}\lesssim\tnorm{\Af{}\circ\Psi_0}_{\Ckin^\beta}\tnorm[\Big]{D_v^2\glift\circ\Psi_0}_{\Ckin^\beta}\lesssim_J\bracket{v_0}^N\ue^{-a\abs{v_0}^2}. \label{cb:inflow-source-product}
\end{align}
The coefficient factor has polynomial growth by \eqref{cb:inflow-coefficient-holder}, and the lift factor has Gaussian decay by \eqref{cb:inflow-lift-first-derivatives}. The product inequality also controls $(\Af{}v)\cdot\nabla_v\glift$ and $\Cf{}\glift$, since $v\circ\Psi_0$ has a polynomial H\"older norm; \eqref{cb:inflow-lift-first-derivatives} directly controls $\Y\glift$. Consequently
\begin{align}
\tnorm{\vradius^2(S-\Kop f-\sourceR)\circ\Psi_0}_{\Ckin^\beta(\Qplus_1)}\lesssim_J\bracket{v_0}^N\ue^{-a\abs{v_0}^2}\lesssim_J\ue^{-\alpha_*\abs{v_0}^2}. \label{cb:inflow-normalized-source}
\end{align}
For $\Kop f+\sourceR$, use \eqref{cb:inflow-source-holder} and the same change of variables. Choosing $\beta$ below the coefficient and source exponents and reducing $\alpha_*$ gives
\begin{align}
\tnorm[\big]{\widehat S}_{\Ckin^\beta(\Qplus_1)}\lesssim_J\ue^{-\alpha_*\abs{v_0}^2}. \label{cb:inflow-full-source-holder}
\end{align}
This verifies the source hypothesis in \eqref{lin:coefficients} and also \eqref{lin:interior-holder-hypotheses}. The radius $\vradius$ in \eqref{lin:inverse-map} remains of the form $c_J\bracket{v_0}^{-N}$.

\proofstep{The small remainder on nongrazing cylinders} By \eqref{cb:inflow-normalized-sharp} and the lift bound \eqref{eq:I:display:294}, $\widehat u$ has a Gaussian $\Ckin^{\frac12}$ bound on $\Qplus_{\frac12}(\widehat Z_0)$. Its incoming value is zero by \eqref{cb:inflow-transformed-equation}. Continuity extends this value to the grazing set, which is in the closure of the incoming set.

Take $\zeta\in\Qplus_{\frac14}(\widehat Z_0)$ and use $\widehat\graz$, $\zeta_*$ and $r=c\widehat\graz$ from \eqref{eq:I:display:298}--\eqref{cb:inflow-distance-to-grazing}. For $0<\widehat\graz<\delta_0$, with the same fixed threshold as in the first-order proof, \eqref{cb:inflow-grazing-cylinder-geometry} places $\zeta_*$ and $\Qplus_{64r}(\zeta)$ in the domain of the preceding H\"older bound and verifies \eqref{lin:second-order-nongrazing}. Formula~\eqref{cb:inflow-grazing-distance} and $\widehat u(\zeta_*)=0$ give
\begin{align}
\tnorm[\big]{\widehat u}_{L^\infty(\Qplus_{16r}(\zeta))}\leq[\widehat u]_{\Ckin^{\frac12}}\sup_{\zeta'\in\Qplus_{16r}(\zeta)}\bigl(\dkin(\zeta',\zeta_*)\bigr)^{\frac12}\lesssim_J\ue^{-\alpha_*\abs{v_0}^2}\widehat\graz^{\frac12}. \label{eq:I:display:299}
\end{align}
Thus the remainder amplitude is of order $r^{\frac12}$ on the cylinders where we will apply the second-order estimate.

\proofstep{Application of the second-order estimates} On a cylinder meeting the spatial boundary, apply Proposition~\ref{lin:order-two} with zero incoming trace, using \eqref{cb:inflow-transformed-equation}, \eqref{cb:inflow-full-source-holder} and \eqref{eq:I:display:299}. Its geometric factor is bounded because $\abs{\eta_3}/r\leq c^{-1}$. On an interior cylinder, Proposition~\ref{lin:interior-inputs}\ref{lin:interior} applies to the same equation and source. In both cases the amplitude, source supremum and source seminorm contribute, respectively,
\begin{align}
r^2\left(\tnorm[\Big]{D_\eta^2\widehat u}_{L^\infty(\Qplus_r)}+\tnorm[\big]{\Yhat\widehat u}_{L^\infty(\Qplus_r)}\right)\lesssim_J\ue^{-\alpha_*\abs{v_0}^2}\bigl(\widehat\graz^{\frac12}+r^2+r^{2+\beta}\bigr). \label{cb:inflow-scaled-hessian}
\end{align}
Since $r=c\widehat\graz$, the leading term after division by $r^2$ is $r^{-2}\widehat\graz^{\frac12}\simeq\widehat\graz^{-\frac32}$. The source terms give $1+r^\beta\lesssim1$, so
\begin{align}
\abs[\big]{D_\eta^2\widehat u}+\abs[\big]{\Yhat\widehat u}\lesssim_J\ue^{-\alpha_*\abs{v_0}^2}\bigl(1+\widehat\graz^{-\frac32}\bigr). \label{eq:I:display:300b}
\end{align}
For $\widehat\graz\geq\delta_0$, use the nongrazing cylinders of fixed radius $r=c\delta_0$ constructed in the first-order proof. The factor $(1+\abs{\eta_3}/r)^N$ in Proposition~\ref{lin:order-two} grows at most polynomially in $\bracket{v_0}$ by \eqref{lin:normalized-center} and \eqref{cb:matrix-bounds}. The Gaussian amplitude of $\widehat u$ and \eqref{cb:inflow-full-source-holder} therefore give \eqref{eq:I:display:300b} there as well, after \eqref{lin:Gaussian-absorption} reduces the exponent. We have obtained this estimate throughout $\Qplus_{\frac14}(\widehat Z_0)$.

\proofstep{The physical second-order estimate} Apply the derivative identities \eqref{lin:physical-gradient-conversion}--\eqref{lin:physical-transport-conversion} to $u=f-\glift$. Their matrix and radius bounds give
\begin{align}
\bigl(\abs{D_v^2u}+\abs{\Y u}\bigr)\circ\Psi_0\lesssim\bracket{v_0}^N\left(\abs[\big]{D_\eta^2\widehat u}+\abs[\big]{\Yhat\widehat u}+\abs[\big]{\nabla_\eta\widehat u}\right). \label{cb:pullback-estimate}
\end{align}
The last term is necessary: the transport transformation contains $\vradius^{-1}\bigl(S_0^{-1}D^2\phi(x)[v,v]\bigr)\circ\Psi_0\cdot\nabla_\eta\widehat u$. By \eqref{eq:I:display:300} and \eqref{eq:I:display:294},
\begin{align}
\abs[\big]{\nabla_\eta\widehat u}\leq\abs[\big]{\nabla_\eta\widehat f}+\abs[\big]{\nabla_\eta(\glift\circ\Psi_0)}\lesssim_J\ue^{-\alpha_*\abs{v_0}^2}\bigl(1+\widehat\graz^{-\frac12}\bigr). \label{cb:inflow-remainder-gradient}
\end{align}
Combine this with \eqref{eq:I:display:300b} in \eqref{cb:pullback-estimate}, and add the physical lift derivatives bounded by \eqref{eq:I:display:294}. Use $1+\widehat\graz^{-\frac12}\lesssim1+\widehat\graz^{-\frac32}$ and \eqref{cb:inflow-grazing-weight-ratio} with $p=\frac32$ to divide by the physical grazing factor. Thus, on each physical image of the smaller half-cylinder,
\begin{align}
\frac{\abs{D_v^2f}+\abs{\Y f}}{1+\graz^{-\frac32}}\lesssim_J\bracket{v_0}^N\ue^{-a\abs{v_0}^2}\lesssim_J\ue^{-\alpha_*\abs{v}^2}. \label{cb:inflow-physical-local-hessian}
\end{align}
The constants are uniform in the center, so these local estimates give \eqref{cb:pointwise-grazing-hessian} by taking the essential supremum over the boundary collar. On $\{d\geq\frac{d_0}{2}\}$, the interior change of velocity and Proposition~\ref{lin:interior-inputs}\ref{lin:interior} give the same estimate without the grazing denominator. Their radii are bounded below by a negative power of $\bracket{v_0}$; \eqref{lin:Gaussian-absorption} absorbs this loss and proves \eqref{cb:inflow-interior-hessian}.
\end{proof}

\subsection{Diffuse-Reflection: Flux and First-Order Estimates}\label{cb:section-flux}

For prescribed in-flow, the incoming regularity in \eqref{cb:inflow-incoming-holder} is an assumption on $g$. For diffuse-reflection, the incoming value depends on the unknown solution through the flux:
\begin{align}
\trm f&=\sqrt{2\pi}\,\mh\,\flux{f},\qquad \flux{f}(t,x):=\int_{n(x)\cdot v>0}\mh(v)\,\trp f(t,x,v)\bigl(n(x)\cdot v\bigr)\dd v. \label{cb:diffuse-law}
\end{align}
We will obtain the incoming regularity required by Proposition~\ref{cb:inflow-gradient} through the flux argument of \cite[Proposition~7.2]{Kim.Weidner2026}. Starting from the Gaussian bound \eqref{cb:setting} and the source bound \eqref{lin:source-Gaussian}, Lemma~\ref{cb:rough-flux-coefficients} will first estimate the outgoing trace and integrate in velocity to obtain $\flux{f}\in C^{\kappa_0}$. The diffuse law and the convolution formulas will then give
\begin{align}
\flux{f}\in C^{\kappa_0}\quad\to\quad f\in\Ckinloc^{\beta_0}\quad\to\quad(\Af{},\Bf{},\Cf{},\Kop f)\in\Ckinloc^\beta. \label{cb:diffuse-bootstrap}
\end{align}
Here $\kappa_0$ is an ordinary H\"older exponent in $(t,x)$, while $\beta_0$ and $\beta$ are the kinetic exponents for the solution and its coefficients, respectively, as in Lemma~\ref{cb:inflow-coefficients}.

With the additional source assumption \eqref{lin:source-holder}, Proposition~\ref{cb:flux} will apply the outgoing second-order estimate and improve the flux exponent to $\kappa>\frac{2}{3}$. The cylinders used in that estimate have empty incoming boundary, so Proposition~\ref{lin:order-two} applies before any incoming second-order estimate is established. The improved flux will give \eqref{cb:diffuse-incoming-holder}, and Proposition~\ref{cb:gradient} will then apply Proposition~\ref{cb:inflow-gradient} to obtain sharp $\Ckin^{\frac12}$ regularity and the weighted first-order estimates. Section~\ref{cb:section-graph} will use the same $C^\kappa$ flux bound to construct the lift for the incoming second-order estimate.

\para{Preliminary H\"older estimates}

\begin{lemma}[Initial H\"older estimates for diffuse-reflection]\label{cb:rough-flux-coefficients}
Let $f$ be a weak solution of \eqref{cb:equation-residual}, with diffuse trace \eqref{cb:diffuse-law}. Assume the solution class \eqref{cb:regularity-class}, the Gaussian bound \eqref{cb:setting} and the source bound \eqref{lin:source-Gaussian}. For every $J\subset\subset I$, there are $\alpha_*>0$, $N\geq0$ and $\kappa_0,\beta_0,\beta>0$, with $\beta_0<\frac{1}{2}$, such that the following estimates hold uniformly in $v_0\in\R^3$.
\begin{enumerate}
\item The flux satisfies the ordinary H\"older estimate
\begingroup\usetagform{enumerated}
\begin{align}
\norm{\flux{f}}_{C^{\kappa_0}(J\times\partial\Om)}&\lesssim_J 1. \label{cb:rough-flux-norm}
\end{align}
\endgroup
\item The solution has a continuous representative on $J\times\overline\Om\times\R^3$, agreeing with its kinetic trace at the boundary, and
\begingroup\usetagform{enumerated}
\begin{align}
\tnorm{f}_{\Ckin^{\beta_0}(J\times\overline\Om\times B_{\bracket{v_0}^{-1}}(v_0))}&\lesssim_J\ue^{-\alpha_*\abs{v_0}^2}\qquad(v_0\in\R^3). \label{cb:rough-solution-holder}
\end{align}
\endgroup
\item The coefficients defined by \eqref{cvel:coefficient-convolutions}--\eqref{cvel:scalar-coefficients} satisfy, with all norms on the set in \eqref{cb:rough-solution-holder},
\begingroup\usetagform{enumerated}
\begin{align}
\tnorm{\Af{}}_{\Ckin^{\beta}}+\tnorm{\Bf{}}_{\Ckin^{\beta}}+\tnorm{\Cf{}}_{\Ckin^{\beta}}+\tnorm{\nabla_v\Af{}}_{L^\infty}&\lesssim_J\bracket{v_0}^{N}. \label{cb:diffuse-coefficient-holder}
\end{align}
\endgroup
\item The nonlocal source satisfies
\begingroup\usetagform{enumerated}
\begin{align}
\tnorm{\Kop f}_{\Ckin^{\beta}(J\times\overline\Om\times B_{\bracket{v_0}^{-1}}(v_0))}&\lesssim_J\ue^{-\alpha_*\abs{v_0}^2}. \label{cb:diffuse-K-holder}
\end{align}
\endgroup
If \eqref{lin:source-holder} also holds, then, after decreasing $\beta$ and $\alpha_*$,
\begingroup\usetagform{enumerated}
\begin{align}
\tnorm{\Kop f+\sourceR}_{\Ckin^{\beta}(J\times\overline\Om\times B_{\bracket{v_0}^{-1}}(v_0))}&\lesssim_J\ue^{-\alpha_*\abs{v_0}^2}. \label{cb:diffuse-source-holder}
\end{align}
\endgroup
\end{enumerate}
\end{lemma}

\begin{proof}
We will first estimate the outgoing trace on cylinders with empty incoming boundary, then integrate in velocity to prove (i). The diffuse law will supply the incoming H\"older norm needed for (ii). Finally, the convolution estimates from Lemma~\ref{cb:inflow-coefficients} will give (iii)--(iv).

The local weak equation is \eqref{lin:boundary-weak-setting}, with $\widehat f=f\circ\Psi_0$ and source \eqref{lin:normalized-f-source}, by Lemma~\ref{cb:normalized-equation}. The solution class \eqref{cb:regularity-class} gives the local $L^2$ velocity gradient. Lemma~\ref{cb:normalized-bounds} gives \eqref{lin:normalized-coefficient-bounds} and \eqref{cb:Gaussian-bound}. These verify \eqref{lin:local-ellipticity}, \eqref{lin:interior-rough-hypotheses} and \eqref{lin:boundary-local-hypotheses}, with constants independent of the center velocity. Thus Proposition~\ref{lin:boundary-inputs}\ref{lin:rough} applies as soon as the incoming norm is controlled; that norm is zero on the outgoing cylinders used first.

\proofstep{Proof of (i): the outgoing estimate and the flux integral} Fix a boundary neighborhood on which a $C^2$ orthogonal matrix $\mathsf O(x)$ has third column $n(x)$. Put $v=\mathsf O(x)\nu$, with $\nu=(\nu',s)$ and $s>0$. Orthogonality preserves both the Maxwellian and the velocity measure, so
\begin{align}
\flux{f}(t,x)=\int_{\nu_3>0}\mh(\nu)\trp f\bigl(t,x,\mathsf O(x)\nu\bigr)\nu_3\dd\nu. \label{cb:flat-flux}
\end{align}
The Gaussian trace bound in Lemma~\ref{trace:gaussian} immediately bounds this integral in $L^\infty$. We will estimate its increments for two nearby boundary points, writing $\Delta:=\abs{t-t'}+\abs{x-x'}\leq1$.

For $0<s\leq1$, center $\Psi_0$ at $\bigl(t,x,\mathsf O(x)\nu\bigr)$. Formula~\eqref{eq:I:display:279} gives $\eta_{0,3}=-\frac{s}{\vradius s_n}$. A radius $r=c\min\{1,\abs{\eta_{0,3}}\}$, with fixed $c>0$ sufficiently small, satisfies $\Qc_{2r}(\widehat Z_0)\subset\Qc_1(\widehat Z_0)$ and $\eta_3<0$ on that cylinder. Hence its incoming set is empty. The radius and its inverse differ from $s$ and $s^{-1}$ by at most polynomial factors in $\bracket{\nu}$, by \eqref{cb:matrix-bounds}. Choose $0<\kappa_0<\frac13$ so that $3\kappa_0$ is below the uniform exponent in Proposition~\ref{lin:boundary-inputs}\ref{lin:rough}. That proposition gives, for some $a>0$,
\begin{align}
[\widehat f]_{\Ckin^{3\kappa_0}(\Qplus_r(\widehat Z_0))}\lesssim_J r^{-3\kappa_0}\left(\tnorm[\big]{\widehat f}_{L^\infty(\Qplus_{2r})}+r^2\tnorm[\big]{\widehat S_f}_{L^\infty(\Qplus_{2r})}\right)\lesssim_J r^{-3\kappa_0}\ue^{-a\abs{\nu}^2}. \label{cb:diffuse-rough-outgoing-local}
\end{align}
Its continuous outgoing values agree with the kinetic trace by Proposition~\ref{lin:boundary-inputs}\ref{lin:boundary-trace-identification}.

The velocity-frame increment is bounded by $\abs{\mathsf O(x)\nu-\mathsf O(x')\nu}\lesssim\bracket{\nu}\Delta$. Thus the two physical points have kinetic distance at most $C\bracket{\nu}^{N}\Delta^{\frac13}$. The explicit inverse map and \eqref{lin:inverse-kinetic-change-bound}, as in \eqref{cb:inflow-uniform-local-radius}, place their images in the smaller cylinder whenever $\Delta\leq c\bracket{\nu}^{-N}s^3$. Applying \eqref{cb:diffuse-rough-outgoing-local} gives
\begin{align}
\abs[\Big]{\trp f\bigl(t,x,\mathsf O(x)\nu\bigr)-\trp f\bigl(t',x',\mathsf O(x')\nu\bigr)}\lesssim_J\ue^{-\alpha_*\abs{\nu}^2}\Delta^{\kappa_0}s^{-3\kappa_0}. \label{cb:rough-increment}
\end{align}
Here and below, polynomial velocity factors are absorbed by \eqref{lin:Gaussian-absorption}.

For fixed $\nu'$, set $h:=\min\{1,(C\bracket{\nu'}^N\Delta)^{\frac13}\}$, with $C,N$ large enough that \eqref{cb:rough-increment} applies for $h\leq s\leq1$. On $0<s<h$, use the Gaussian trace amplitude. Since $\bracket{\nu}\simeq\bracket{\nu'}$ on $0<s<1$, the contribution to the flux increment is bounded by
\begin{align}
C_J\ue^{-a\abs{\nu'}^2}\left(h^2+\Delta^{\kappa_0}\int_h^1s^{1-3\kappa_0}\dd s\right)\lesssim_J\ue^{-\alpha_*\abs{\nu'}^2}\bigl(\Delta^{\frac23}+\Delta^{\kappa_0}\bigr)\lesssim_J\ue^{-\alpha_*\abs{\nu'}^2}\Delta^{\kappa_0}. \label{eq:I:display:269}
\end{align}
We have $h^2\leq C\bracket{\nu'}^N\Delta^{\frac23}$, and the integral is uniformly bounded because $1-3\kappa_0>-1$. This also covers $h=1$, when the second integral is zero. For $s\geq1$, choose an outgoing radius bounded below by a negative power of $\bracket{\nu}$. The same local estimate gives a Gaussian multiple of $\Delta^{\kappa_0}$ below the corresponding separation threshold; above it, the Gaussian amplitude gives that bound with a polynomial velocity factor. Integration in $\nu'$ and $s$ proves the local $C^{\kappa_0}$ estimate. Finitely many fixed boundary neighborhoods have a common positive margin, so they cover all sufficiently close pairs. The $L^\infty$ flux bound handles the other pairs, proving \eqref{cb:rough-flux-norm}.

\proofstep{Proof of (ii): the incoming norm and solution regularity} Set $g(t,x,v)=\sqrt{2\pi}\mh(v)\flux{f}(t,x)$. For incoming points $Z,Z'$ in the velocity ball of \eqref{cb:rough-solution-holder}, put $q=\dkin(Z',Z)\leq1$. The ordinary time and position separation is at most $C\bracket{v_0}q^2$, while $\abs{v'-v}\leq q$. The flux estimate and the Maxwellian derivative therefore give
\begin{align}
\abs{g(Z')-g(Z)}\lesssim_J\bracket{v_0}^N\ue^{-a\abs{v_0}^2}\left(q^{2\kappa_0}+q\right). \label{cb:diffuse-rough-incoming-increment}
\end{align}
Choose $0<\beta_0<\min\{\frac12,2\kappa_0\}$ below the exponent allowed by Proposition~\ref{lin:boundary-inputs}\ref{lin:rough}. The amplitude handles $q\geq1$. The distance estimate \eqref{lin:kinetic-change-bound} and Gaussian absorption yield
\begin{align}
\tnorm{g\circ\Psi_0}_{\Ckin^{\beta_0}(\partial_-\Qplus_2)}\lesssim_J\ue^{-\alpha_*\abs{v_0}^2}. \label{cb:diffuse-rough-incoming-norm}
\end{align}
The weak equation and coefficient bounds were verified at the beginning of the proof. Thus \eqref{cb:diffuse-rough-incoming-norm} verifies the remaining assumption \eqref{lin:boundary-rough-hypotheses}, and Proposition~\ref{lin:boundary-inputs}\ref{lin:rough} gives the Gaussian $\Ckin^{\beta_0}$ estimate for $\widehat f$ on $\Qplus_1$. Proposition~\ref{lin:interior-inputs}\ref{lin:rough-full} gives the corresponding interior estimate. Apply the uniform-radius construction \eqref{cb:inflow-uniform-local-radius} and the two-pair estimates \eqref{cb:inflow-near-pair}--\eqref{cb:inflow-separated-pair} from Lemma~\ref{cb:inflow-coefficients}. They use only these local estimates and the Gaussian amplitude, and therefore prove \eqref{cb:rough-solution-holder} on the full physical set, with constants independent of $v_0$.

The local continuous representatives agree on overlaps as in that lemma. Proposition~\ref{lin:boundary-inputs}\ref{lin:boundary-trace-identification} identifies their boundary values with the kinetic trace. The Gaussian trace bound permits passage under the integral in \eqref{cb:flat-flux}; hence its continuous value is the flux in \eqref{cb:diffuse-law}, and this law extends to grazing by continuity.

\proofstep{Proof of (iii)--(iv): coefficient and source estimates} The rough norm \eqref{cb:rough-solution-holder} and Gaussian bound \eqref{cb:setting} are the only inputs used to derive \eqref{eq:I:display:296} and \eqref{cb:inflow-velocity-increments}. Thus both increment estimates hold here. The convolution calculation \eqref{cb:inflow-convolution-kernel-bound}--\eqref{eq:I:display:297}, which uses no boundary condition, proves \eqref{cb:diffuse-coefficient-holder} for $0<\beta<\frac{\beta_0}{6}$. This includes $c_F=-8\pi F$ at $\gamma=-3$ and the bound on $\nabla_v\Af{}$.

Applying that calculation to $\mh f$ in \eqref{cvel:K}, with the Gaussian Maxwellian factors as in \eqref{cb:inflow-K-holder}, proves \eqref{cb:diffuse-K-holder}. If \eqref{lin:source-holder} holds, combine it with \eqref{lin:source-Gaussian} and decrease $\beta,\alpha_*$ below its exponents. The triangle inequality then gives \eqref{cb:diffuse-source-holder}.
\end{proof}

\para{Improved flux regularity}

\begin{proposition}\label{cb:flux}
Let $f$ be a weak solution of \eqref{cb:equation-residual}, with diffuse trace \eqref{cb:diffuse-law}. Assume \eqref{cb:regularity-class}, \eqref{cb:setting} and the source bounds \eqref{lin:source-Gaussian} and \eqref{lin:source-holder}. For every $J\subset\subset I$, there are $\frac{2}{3}<\kappa<1$, $\varepsilon>0$ and $\alpha_*>0$ such that
\begin{enumerate}
\item The flux satisfies the ordinary H\"older estimate
\begingroup\usetagform{enumerated}
\begin{align}
\norm{\flux{f}}_{C^\kappa(J\times\partial\Om)}&\lesssim_J 1. \label{rev14:diffuse-flux-regularity}
\end{align}
\endgroup
\item The incoming value in \eqref{cb:diffuse-law} satisfies
\begingroup\usetagform{enumerated}
\begin{align}
\tnorm{\sqrt{2\pi}\mh\,\flux{f}}_{\Ckin^{1+\varepsilon}\left(\{(t,x,v)\in J\times\Sigma_-:\abs{v-v_0}<\bracket{v_0}^{-1}\}\right)}&\lesssim_J\ue^{-\alpha_*\abs{v_0}^2}\qquad(v_0\in\R^3). \label{cb:diffuse-incoming-holder}
\end{align}
\endgroup
\end{enumerate}
The exponents and the implicit constants are independent of $v_0$.
\end{proposition}

\begin{proof}
\proofstep{Proof of (i): the outgoing second-order estimate} We will apply Proposition~\ref{lin:order-two} on cylinders with empty incoming boundary and integrate its estimate in velocity. The local $L^2$ velocity gradient and weak identity \eqref{lin:boundary-weak-setting} follow from \eqref{cb:regularity-class} and Lemma~\ref{cb:normalized-equation}.

Apply Lemma~\ref{cb:rough-flux-coefficients} on a compact enlargement of $J$ inside $J'$. Its bounds \eqref{cb:diffuse-coefficient-holder} and \eqref{cb:diffuse-source-holder} verify \eqref{lin:physical-coefficient-holder} and \eqref{lin:physical-source-holder}. Lemma~\ref{cb:normalized-bounds}, with radius $\vradius=c_J\bracket{v_0}^{-N}$, gives \eqref{lin:normalized-coefficient-bounds}, \eqref{lin:normalized-coefficient-holder} and \eqref{lin:normalized-source-holder}. These verify \eqref{lin:coefficients}--\eqref{lin:coefficient-norm} uniformly in $v_0$, with Gaussian source decay. We next choose cylinders satisfying \eqref{lin:second-order-nongrazing}.

Fix $0<\beta'<\min\{\beta,\beta_0\}$. In the velocity frame \eqref{cb:flat-flux}, let $\nu=(\nu',s)$ with $0<s\leq1$ and center $\Psi_0$ at $\bigl(t,x,\mathsf O(x)\nu\bigr)$. Write $\widehat Z_0=(0,\xi_0,\eta_0)$ for its image and $\widehat f=f\circ\Psi_0$. By \eqref{eq:I:display:279}, $\eta_{0,3}=-\frac{s}{\vradius s_n}$. Choose
\begin{align}
r_s&:=c\min\{1,\abs{\eta_{0,3}}\},\qquad c\bracket{\nu}^{-N}s\leq r_s\leq C\bracket{\nu}^{N}s, \label{cb:outgoing-radius}
\end{align}
with $c>0$ sufficiently small that $\Qc_{64r_s}(\widehat Z_0)\subset\Qc_1(\widehat Z_0)$ and $\eta_3<0$ throughout that cylinder. These inequalities imply \eqref{lin:second-order-nongrazing} and $\partial_-\Qplus_{16r_s}(\widehat Z_0)=\varnothing$. The second comparison follows from \eqref{cb:matrix-bounds} and the formula for $\vradius$.

If $\abs{\eta_{0,3}}$ is smaller than a fixed constant, subtract the value of $\widehat f$ at the nearby grazing point $\zeta_*$ with the same time, boundary position and tangential velocity. This point lies in $\Qc_1(\widehat Z_0)\cap\{\xi_3=0\}$ after fixing the threshold sufficiently small. The tangential center velocity cancels in the kinetic spatial residual relative to $\zeta_*$, so \eqref{cb:rough-solution-holder} bounds the remainder on $\Qplus_{16r_s}(\widehat Z_0)$ by $C_J\bracket{\nu}^{N}\ue^{-\alpha_*\abs{\nu}^2}s^{\beta_0}$. Otherwise $r_s$ is bounded below, and $s=\vradius s_n\abs{\eta_{0,3}}\geq c\bracket{\nu}^{-N}$ by \eqref{cb:matrix-bounds}. The Gaussian amplitude bound then gives the same estimate without subtracting a value, after increasing $N$. The source for the remainder is $\widehat S_f+\locc\widehat f(\zeta_*)$, with $\widehat S_f$ from \eqref{lin:normalized-f-source}. When no value is subtracted, the source remains $\widehat S_f$. Their $\Ckin^{\beta'}$ norms follow from \eqref{cb:diffuse-coefficient-holder}--\eqref{cb:diffuse-source-holder}. Here we use \eqref{lin:source-holder}.

The geometric factor in Proposition~\ref{lin:order-two} is bounded when $\abs{\eta_{0,3}}\leq1$ and by $C\bracket{\nu}^N$ when the radius saturates. Its right-hand side consists of the remainder amplitude, $r_s^2$ times the source supremum, and $r_s^{2+\beta'}$ times its seminorm. The source terms are bounded by a polynomial velocity factor times $\ue^{-a\abs{\nu}^2}s^2$, using \eqref{cb:outgoing-radius}. Since $0<s\leq1$ and $\beta_0<1$, they are bounded by the same factor times $\ue^{-a\abs{\nu}^2}s^{\beta_0}$.

The outgoing part of Proposition~\ref{lin:order-two}, the radius comparisons \eqref{cb:outgoing-radius}, and a decrease of $\alpha_*$ give
\begin{align}
&s\tnorm[\big]{\nabla_\eta\widehat f}_{L^\infty(\Qplus_{r_s}(\widehat Z_0))}+s^2\Bigl(\tnorm[\Big]{D_\eta^2\widehat f}_{L^\infty(\Qplus_{r_s}(\widehat Z_0))}+\tnorm[\big]{\Yhat\widehat f}_{L^\infty(\Qplus_{r_s}(\widehat Z_0))}\Bigr) \label{cb:outgoing-order-two}\\
&\qquad+s^{2+\beta'}[\widehat f]_{\Ckin^{2+\beta'}(\Qplus_{r_s}(\widehat Z_0))}\lesssim_J\ue^{-\alpha_*\abs{\nu}^2}s^{\beta_0}.\notag
\end{align}
The polynomial estimate also holds with a center on the outgoing boundary. Translate both the center and its comparison point by $(0,he_3,0)$, $h>0$, within a smaller cylinder. Their kinetic increments are unchanged. The polynomial coefficients are bounded by \eqref{cb:outgoing-order-two}. Passing to a convergent subsequence of these finitely many coefficients and using continuity of $f$ gives the same remainder bound at the boundary.

\proofstep{Proof of (i): boundary increments} Choose $\frac{2}{3}<\kappa<\frac{2+\beta'}{3}$. Compare the two points in \eqref{cb:rough-increment} when $\Delta\leq c\bracket{\nu}^{-N}s^3$, using the same map $\Psi_0$ for both. For these two boundary points, \eqref{eq:I:display:279} and $n(x)\cdot\mathsf O(x)\nu=s$ give
\begin{align}
\xi_{0,3}=\xi_{1,3}=0,\qquad \eta_{0,3}=\eta_{1,3}=-\frac{s}{\vradius s_n}. \label{cb:flux-normal-comparison}
\end{align}
Here $\widehat Z_i=(\theta_i,\xi_i,\eta_i)$, $i=0,1$, denote their images under the same $\Psi_0^{-1}$. Their transformed time and velocity increments are bounded by $C\bracket{\nu}^{N}\Delta$, and \eqref{lin:inverse-kinetic-change-bound} bounds their transformed kinetic distance by $C\bracket{\nu}^{N}\Delta^{\frac{1}{3}}$. The polynomial \eqref{lin:second-order-polynomial} and \eqref{cb:outgoing-order-two} bound the linear velocity term by $C_J\ue^{-\alpha_*\abs{\nu}^2}\Delta s^{\beta_0-1}$, the quadratic velocity term by $C_J\ue^{-\alpha_*\abs{\nu}^2}\Delta^2s^{\beta_0-2}$, and the time term by $C_J\ue^{-\alpha_*\abs{\nu}^2}\Delta s^{\beta_0-2}$. The first two are bounded by the third for $s,\Delta\leq1$. Including the polynomial remainder gives
\begin{align}
C_J\ue^{-\alpha_*\abs{\nu}^2} \bigl(\Delta s^{\beta_0-2} +\Delta^{\frac{2+\beta'}{3}}s^{\beta_0-2-\beta'}\bigr) \lesssim_J \ue^{-\alpha_*\abs{\nu}^2}\Delta^\kappa s^{\beta_0-3\kappa}. \label{cb:improved-increment}
\end{align}
In the first term we used $\Delta^{1-\kappa}\leq C\bracket{\nu}^{N}s^{3(1-\kappa)}$ and $s\leq1$, and in the second $\Delta^{\frac{2+\beta'}{3}-\kappa}\leq C\bracket{\nu}^{N}s^{2+\beta'-3\kappa}$, which uses $3\kappa<2+\beta'$.

\proofstep{Proof of (i): integration in velocity} For each $\nu'$, put $h:=\min\{1,(C\bracket{\nu'}^N\Delta)^{\frac13}\}$, with $C,N$ ensuring that the preceding local comparison applies when $h\leq s\leq1$. On $0<s<h$, \eqref{cb:rough-solution-holder} and the velocity-frame increment give the bound $C_J\ue^{-a\abs{\nu}^2}\Delta^{\frac{\beta_0}{3}}$. For separations below a fixed negative power of $\bracket{\nu}$, both velocities lie in a common ball from \eqref{cb:rough-solution-holder}; above that threshold the Gaussian amplitude gives the same bound with a polynomial velocity factor. On $h\leq s\leq1$, use \eqref{cb:improved-increment}. The factor $s\dd s$ in the flux integral therefore gives
\begin{align}
C_J\ue^{-a\abs{\nu'}^2}\left(h^2\Delta^{\frac{\beta_0}{3}}+\Delta^\kappa\int_h^1s^{\beta_0+1-3\kappa}\dd s\right)\lesssim_J\ue^{-\alpha_*\abs{\nu'}^2}\bigl(\Delta^{\frac{2+\beta_0}{3}}+\Delta^\kappa\bigr)\lesssim_J\ue^{-\alpha_*\abs{\nu'}^2}\Delta^\kappa. \label{eq:I:display:274}
\end{align}
Here $h^2\leq C\bracket{\nu'}^N\Delta^{\frac23}$, and the integral is uniformly bounded because $\beta_0+2-3\kappa>\beta_0-\beta'>0$. The same inequalities cover $h=1$, when the second integral is absent. The last inequality uses $3\kappa<2+\beta'<2+\beta_0$ and $\Delta\leq1$.

For $s\geq1$, the outgoing radius is bounded below by a negative power of $\bracket{\nu}$. The same polynomial estimate, using the Gaussian amplitude instead of $s^{\beta_0}$, gives $C_J\ue^{-\alpha_*\abs{\nu}^2}\Delta^\kappa$ for separations below the resulting polynomial threshold. For larger separations, the Gaussian amplitude gives this bound after multiplication by the inverse threshold to the power $\kappa$. Formula~\eqref{lin:Gaussian-absorption} absorbs that polynomial factor. Integrating in $\nu'$ and $s$ proves \eqref{rev14:diffuse-flux-regularity} for points in one smaller boundary neighborhood. The fixed finite spatial covering and the $L^\infty$ flux bound treat all other pairs, as in the proof of \eqref{cb:rough-flux-norm}.

\proofstep{Proof of (ii): the incoming polynomial estimate} At an incoming point $Z_*=(t_*,x_*,v_*)$ in the indicated velocity ball, use the degree-one kinetic polynomial
\begin{align}
\sqrt{2\pi}\,\flux{f}(t_*,x_*)\left(\mh(v_*)+\nabla_v\mh(v_*)\cdot(v-v_*)\right). \label{cb:incoming-polynomial}
\end{align}
For another incoming point at kinetic distance $r\leq1$, the ordinary time and boundary-position separation is at most $C\bracket{v_0}r^2$. The difference from this polynomial is a flux increment times $\sqrt{2\pi}\mh(v)$ plus the velocity Taylor remainder multiplied by $\flux{f}(t_*,x_*)$, hence is bounded by
\begin{align}
C_J\bracket{v_0}^{N}\ue^{-\alpha_*\abs{v_0}^2}\left(r^{2\kappa}+r^2\right). \label{cb:incoming-polynomial-remainder}
\end{align}
The polynomial coefficients have the same Gaussian bound. The definition of $\Ckin^{1+\varepsilon}$ therefore gives \eqref{cb:diffuse-incoming-holder} for $0<\varepsilon<2\kappa-1$. Larger separations are controlled by the amplitude and polynomial coefficients.

\end{proof}

\para{Sharp H\"older and first-order estimates} We will apply Proposition~\ref{cb:inflow-gradient} with the incoming regularity supplied by Proposition~\ref{cb:flux}.

\begin{proposition}\label{cb:gradient}
Let $f$ be a weak solution of \eqref{cb:equation-residual}, with diffuse trace \eqref{cb:diffuse-law}. Assume \eqref{cb:regularity-class}, \eqref{cb:setting} and the source bounds \eqref{lin:source-Gaussian} and \eqref{lin:source-holder}. For every $J\subset\subset I$, there is $\alpha_*>0$ such that
\begin{enumerate}
\item The solution has a continuous representative on $J\times\overline\Om\times\R^3$, agreeing with its kinetic trace at the boundary, and
\begingroup\usetagform{enumerated}
\begin{align}
\tnorm{f}_{\Ckin^{\frac{1}{2}}(J\times\overline\Om\times B_{\bracket{v_0}^{-1}}(v_0))}&\lesssim_J\ue^{-\alpha_*\abs{v_0}^2}\qquad(v_0\in\R^3). \label{cb:diffuse-sharp-holder}
\end{align}
\endgroup
\item The distributional velocity derivative satisfies
\begingroup\usetagform{enumerated}
\begin{align}
\tnorm{\frac{\ue^{\alpha_*\abs{v}^2}\nabla_v f}{1+\graz^{-\frac12}}}_{L^\infty\left(J\times\{d<\frac{d_0}{2}\}\times\R^3\right)}\lesssim_J1, \label{cb:pointwise-estimates}\\
\tnorm[\Big]{\ue^{\alpha_*\abs{v}^2}\nabla_v f}_{L^\infty\left(J\times\{d\geq\frac{d_0}{2}\}\times\R^3\right)}\lesssim_J1. \label{cb:diffuse-interior-gradient}
\end{align}
\endgroup
\end{enumerate}
The distributional velocity derivative has a continuous representative on compact subsets avoiding the grazing set.
\end{proposition}

\begin{proof}
Set $g=\sqrt{2\pi}\mh\,\flux{f}$. Equation~\eqref{cb:diffuse-incoming-holder} is precisely the incoming bound \eqref{cb:inflow-incoming-holder}. The weak equation, solution class, Gaussian bound and source bound required by Proposition~\ref{cb:inflow-gradient} are \eqref{cb:equation-residual}, \eqref{cb:regularity-class}, \eqref{cb:setting} and \eqref{lin:source-Gaussian}. That proposition gives \eqref{cb:diffuse-sharp-holder}, \eqref{cb:pointwise-estimates}--\eqref{cb:diffuse-interior-gradient}, the continuous boundary values agreeing with the trace, and continuity of $\nabla_vf$ away from grazing. The additional source assumption \eqref{lin:source-holder} enters through the outgoing second-order estimate in Proposition~\ref{cb:flux}, which supplies \eqref{cb:diffuse-incoming-holder}.

For the second-order argument below, we will use the sharp estimate for $\widehat f=f\circ\Psi_0$ on the half-cylinder $\Qplus_1(\widehat Z_0)$. Lemma~\ref{cb:normalized-bounds} transfers \eqref{cb:diffuse-coefficient-holder} to the diffusion coefficient $\locA$ in \eqref{lin:Landau-equation}, and \eqref{cb:diffuse-incoming-holder}, followed by \eqref{cb:inflow-incoming-pullback-remainder}--\eqref{cb:inflow-normalized-amplitudes}, gives the required norm of $g\circ\Psi_0$ on $\partial_-\Qplus_2(\widehat Z_0)$. These verify \eqref{lin:boundary-sharp-hypotheses}. Proposition~\ref{lin:boundary-inputs}\ref{lin:sharp}, applied on $\Qplus_2(\widehat Z_0)$ with $r=1$, gives
\begin{align}
\tnorm[\big]{\widehat f}_{\Ckin^{\frac{1}{2}}(\Qplus_1(\widehat Z_0))}&\lesssim_J\ue^{-\alpha_*\abs{v_0}^2}. \label{cb:diffuse-normalized-sharp}
\end{align}
\end{proof}

\subsection{Diffuse-Reflection: Boundary Lift and Second-Order Estimates}\label{cb:section-graph}

We will construct an extension $\widehat{\glift}$ of the diffuse incoming value whose first and second velocity derivatives and transport derivative satisfy the kinetic H\"older estimate \eqref{lift:quantitative-norm}. Its boundary identity \eqref{lift:trace} will make $\widehat u=\widehat f-\widehat{\glift}$ vanish on the incoming boundary, and its derivative bounds will allow us to apply Proposition~\ref{lin:order-two} to this difference. Proposition~\ref{cb:flux} supplies ordinary $C^\kappa$ regularity of the boundary flux, with $\frac23<\kappa<1$. We will first extend this flux to a regularized function $\fluxlift$, then multiply it by the Maxwellian factor to obtain $\widehat{\glift}$.

For the map $\Psi_0$ in \eqref{lin:inverse-map}, write $\Psi_0(\theta,\xi,\eta)=(t,x,v)$ and define
\begin{align}
\widehat M(\theta,\xi,\eta):=\sqrt{2\pi}\mh(v),\qquad \fluxhat(\theta,\xi'):=\flux{f}(t,x)\quad\text{when }\xi_3=0. \label{lift:transformed-boundary-data}
\end{align}
The pair $(t,x)$ on this boundary is independent of $\eta$, so $\fluxhat$ is well defined and the incoming value of $\widehat f=f\circ\Psi_0$ is $\widehat M\fluxhat$. With the coefficients and source from \eqref{cb:diffuse-coefficient-holder}--\eqref{cb:diffuse-source-holder}, subtraction of the lift gives the source
\begin{align}
\widehat S_f-\Yhat\widehat{\glift}+\locA:D_\eta^2\widehat{\glift}+\locb\cdot\nabla_\eta\widehat{\glift}+\locc\widehat{\glift}\in\Ckin^{\beta'} \label{lift:modified-source-requirement}
\end{align}
for a sufficiently small $\beta'>0$. Lemma~\ref{lift:scalar-jets} will estimate the derivatives of $\fluxlift$, and Lemma~\ref{lift:characteristic} will give \eqref{lift:quantitative-norm}, which verifies \eqref{lift:modified-source-requirement}. Proposition~\ref{cb:graph} will then estimate the difference and its source at the scales needed to obtain the grazing weight in \eqref{eq:rho}.

The condition $3\kappa>2$ reflects the incoming boundary geometry. If two boundary points have normal velocity at least $s>0$ and kinetic distance at most $r\leq cs$, then \eqref{eq:kinetic-distance} and their zero normal positions give
\begin{align}
\abs{\Delta\theta}\lesssim\frac{r^3}{s},\qquad \abs{\Delta\xi'}\lesssim\left(1+\frac{\abs{\eta_0'}}{s}\right)r^3. \label{lift:boundary-increments}
\end{align}
Thus the boundary flux varies by order $r^{3\kappa}$. The extension below will retain enough of this regularity to control two velocity derivatives and one transport derivative.

\para{Construction of the regularized flux} Let $B'_r$ be the ball of radius $r$ in $\R^2$. We will first work on
\begin{align}
\liftpatch&:=(-1,1)\times B'_1\times\left(0,\frac{1}{16}\right)\times\left\{\eta:\abs{\eta'}<1,\ \frac{1}{2}<\eta_3<2\right\}, \label{lift:normalized-patch}
\end{align}
with boundary flux $\fluxhat$ given on $(-2,2)\times B'_2$. The incoming boundary of this set is its face $\xi_3=0$. The proof of Proposition~\ref{cb:graph} will use a kinetic dilation and a tangential Galilean translation to bring small positive normal velocities into this range.

For $\zeta=(\theta,\xi,\eta)\in\overline{\liftpatch}$, the backward travel time to $\xi_3=0$ and the time and tangential position reached there are
\begin{align}
\theta_b:=\frac{\xi_3}{\eta_3},\qquad \zeta_b:=\bigl(\theta-\theta_b,\xi'-\eta'\theta_b\bigr),\qquad0\leq\theta_b\leq\frac18. \label{rev14:lift-parameters}
\end{align}
Along the straight transport field, these functions satisfy
\begin{align}
\Yhat\theta_b=1,\qquad \Yhat\zeta_b=0. \label{lift:characteristic-identities}
\end{align}
Fix a nonnegative $\varphi\in C_c^\infty(B_1\subset\R^3)$ with integral one. Define the regularized extension of the boundary flux by
\begin{align}
\fluxlift(\zeta):=\int_{B_1}\varphi(q)\fluxhat\bigl(\zeta_b(\zeta)-\theta_b(\zeta)^{\frac32}q\bigr)\dd q,\qquad \zeta\in\overline{\liftpatch}. \label{lift:construction}
\end{align}
Every averaging ball lies in $(-\frac32,\frac32)\times B'_{\frac32}$: the time and tangential components of $\zeta_b$ have absolute value at most $\frac98$, and $\theta_b^{\frac32}\leq(\frac18)^{\frac32}$. Thus only the given values of $\fluxhat$ are used. At $\xi_3=0$, the radius is zero and \eqref{lift:construction} equals $\fluxhat(\theta,\xi')$.

The averaging radius is chosen according to the kinetic increments. If $\dkin(\zeta_1,\zeta_0)\leq r$, then \eqref{eq:kinetic-distance} gives
\begin{align}
\abs{\theta_1-\theta_0}\lesssim r^2,\qquad \abs{\eta_1-\eta_0}\lesssim r,\qquad \abs[\Big]{\xi_1-\xi_0-\eta_0\bigl(\theta_1-\theta_0\bigr)}\lesssim r^3. \label{lift:kinetic-increments}
\end{align}
The calculation \eqref{lift:xi-increment} below gives $\abs{\Delta\zeta_b}\lesssim\theta_{b,0}r+r^3$. At the velocity scale $r=\sqrt{\theta_{b,0}}$, both terms have size $\theta_{b,0}^{\frac32}$. This explains the radius in \eqref{lift:construction}. We will verify its derivative and H\"older bounds in the following two lemmas.

\begin{lemma}[Extension of the diffuse flux]\label{lift:scalar-jets}
Let $\frac23<\kappa<1$ and $\fluxhat\in C^\kappa((-2,2)\times B'_2)$. The function in \eqref{lift:construction} satisfies
\begin{align}
\fluxlift\in C^\infty(\liftpatch)\cap C(\overline{\liftpatch}),\qquad \fluxlift(\theta,\xi',0,\eta)=\fluxhat(\theta,\xi'). \label{lift:scalar-space}
\end{align}
With the seminorm taken on $(-2,2)\times B'_2$, we have
\begin{align}
\abs[\big]{\fluxlift-\fluxhat(\zeta_b)}&\lesssim_\kappa[\,\fluxhat\,]_{C^\kappa}\theta_b^{\frac{3\kappa}{2}}, \label{lift:vanishing-jets}\\
\abs[\big]{\nabla_\eta\fluxlift}&\lesssim_\kappa[\,\fluxhat\,]_{C^\kappa}\theta_b^{\frac{3\kappa-1}{2}}, \label{lift:vanishing-jets-1}\\
\abs[\big]{D_\eta^2\fluxlift}+\abs[\big]{\Yhat\fluxlift}&\lesssim_\kappa[\,\fluxhat\,]_{C^\kappa}\theta_b^{\frac{3\kappa-2}{2}}. \label{lift:vanishing-jets-2}
\end{align}
These derivatives extend continuously to $\xi_3=0$ with value zero.
\end{lemma}

\begin{proof}
In this proof and the next, the norm and seminorm of $\fluxhat$ are on $(-2,2)\times B'_2$.

\proofstep{Continuity and the incoming value} We will prove \eqref{lift:scalar-space} and \eqref{lift:vanishing-jets}. For $\theta_b>0$, changing variables in \eqref{lift:construction} gives
\begin{align}
\fluxlift(\zeta)=\theta_b^{-\frac92}\int_{(-2,2)\times B'_2}\varphi\left(\frac{\zeta_b-z}{\theta_b^{\frac32}}\right)\fluxhat(z)\dd z. \label{lift:kernel-representation}
\end{align}
The kernel is smooth in $\zeta$ and supported strictly inside the integration domain. Differentiation of this kernel proves interior smoothness without differentiating $\fluxhat$. Also,
\begin{align}
\abs[\big]{\fluxlift(\zeta)-\fluxhat(\zeta_b)}\leq[\,\fluxhat\,]_{C^\kappa}\theta_b^{\frac{3\kappa}{2}}\int_{B_1}\varphi(q)\abs{q}^\kappa\dd q. \label{lift:mollifier-error}
\end{align}
As $\xi_3\to0$, $\theta_b\to0$ and $\zeta_b\to(\theta,\xi')$. This proves the continuity and incoming value in \eqref{lift:scalar-space}, and \eqref{lift:mollifier-error} proves \eqref{lift:vanishing-jets}.

\proofstep{Velocity derivatives} We will prove \eqref{lift:vanishing-jets-1} and the Hessian bound in \eqref{lift:vanishing-jets-2}. From \eqref{rev14:lift-parameters},
\begin{align}
\abs{D_\eta\theta_b}+\abs{D_\eta^2\theta_b}\lesssim\theta_b,\qquad \abs{D_\eta\zeta_b}+\abs{D_\eta^2\zeta_b}\lesssim\theta_b. \label{lift:foot-derivatives}
\end{align}
More precisely, each derivative on the left is $\theta_b$ times a smooth function of $\eta$ with bounded derivatives on the velocity set in \eqref{lift:normalized-patch}. Differentiating the kernel in \eqref{lift:kernel-representation} and then changing back to $q$ gives
\begin{align}
\partial_{\eta_i}\fluxlift&=\frac{\partial_{\eta_i}\zeta_b}{\theta_b^{\frac32}}\cdot\int_{B_1}\nabla\varphi(q)\fluxhat(\zeta_b-\theta_b^{\frac32}q)\dd q \label{lift:first-kernel-derivative}\\
&\quad-\frac{3\partial_{\eta_i}\theta_b}{2\theta_b}\int_{B_1}\operatorname{div}_q(q\varphi(q))\fluxhat(\zeta_b-\theta_b^{\frac32}q)\dd q.\notag
\end{align}
Both kernels on the right have integral zero. Subtracting $\fluxhat(\zeta_b)$ in each integral and using \eqref{lift:foot-derivatives} therefore yields
\begin{align}
\abs[\big]{\nabla_\eta\fluxlift}\lesssim_\kappa[\,\fluxhat\,]_{C^\kappa}\bigl(\theta_b^{\frac{3\kappa-1}{2}}+\theta_b^{\frac{3\kappa}{2}}\bigr)\lesssim_\kappa[\,\fluxhat\,]_{C^\kappa}\theta_b^{\frac{3\kappa-1}{2}}. \label{lift:first-derivative}
\end{align}

For the second derivative, each term in \eqref{lift:first-kernel-derivative} has the form
\begin{align}
T(\zeta)=\mathfrak c(\eta)\theta_b^p\int_{B_1}\psi(q)\fluxhat(\zeta_b-\theta_b^{\frac32}q)\dd q,\qquad\int_{B_1}\psi=0, \label{lift:term-structure}
\end{align}
where $p\in\{-\frac12,0\}$, $\psi\in C_c^\infty(B_1)$ and $\mathfrak c$ is smooth with bounded derivatives. Differentiating $\mathfrak c(\eta)\theta_b^p$ leaves the same power $p$, by \eqref{lift:foot-derivatives}. Differentiating the integral as in \eqref{lift:first-kernel-derivative} either leaves $p$ unchanged or replaces it by $p-\frac12$; its new kernels again have zero integral. Hence every second velocity derivative is a finite sum of \eqref{lift:term-structure} with $p\in\{-1,-\frac12,0\}$. Kernel cancellation gives
\begin{align}
\abs{T(\zeta)}\lesssim_\kappa[\,\fluxhat\,]_{C^\kappa}\theta_b^{p+\frac{3\kappa}{2}},\qquad p\in\left\{-1,-\frac12,0\right\}. \label{lift:hessian-types}
\end{align}
Since $\theta_b\leq\frac18$, summing these terms proves the Hessian bound in \eqref{lift:vanishing-jets-2}.

\proofstep{Transport derivative and zero boundary values} We will prove the remaining bound in \eqref{lift:vanishing-jets-2}. The identities \eqref{lift:characteristic-identities} remove the derivatives of the boundary point from the transport derivative of the kernel. They give
\begin{align}
\Yhat\fluxlift=-\frac{3}{2\theta_b}\int_{B_1}\operatorname{div}_q(q\varphi(q))\left[\fluxhat(\zeta_b-\theta_b^{\frac32}q)-\fluxhat(\zeta_b)\right]\dd q. \label{lift:transport-kernel}
\end{align}
Consequently
\begin{align}
\abs[\big]{\Yhat\fluxlift}\lesssim_\kappa[\,\fluxhat\,]_{C^\kappa}\theta_b^{\frac{3\kappa-2}{2}}. \label{lift:transport-power}
\end{align}
Together with \eqref{lift:first-derivative} and \eqref{lift:hessian-types}, this proves all derivative bounds. Their powers are positive because $3\kappa>2$, so the derivatives extend continuously to the incoming boundary with value zero.
\end{proof}

\para{H\"older estimates for the lift} We will now prove the derivative H\"older bound needed in \eqref{lift:modified-source-requirement}. The velocity component of $\Psi_0$ in \eqref{lin:inverse-map} is affine in $\eta$ with $C^2$ coefficients in $\xi$. Thus $\widehat M$, $\nabla_\eta\widehat M$ and $D_\eta^2\widehat M$ have $C^2$ position coefficients, and $\Yhat\widehat M$ has $C^1$ position coefficients. The same regularity holds after the affine transformations used to place the target velocity in \eqref{lift:normalized-patch}. In particular, for every $0<\beta<1$,
\begin{align}
\widehat M,\quad\nabla_\eta\widehat M,\quad D_\eta^2\widehat M,\quad\Yhat\widehat M\in\Ckin^\beta(\overline{\liftpatch}). \label{lift:multiplier-regularity}
\end{align}
The map bounds \eqref{cb:matrix-bounds} and Gaussian absorption \eqref{lin:Gaussian-absorption} preserve Gaussian decay in the original center velocity, with a smaller exponent.

\begin{lemma}\label{lift:characteristic}
Let $\frac23<\kappa<1$, $\fluxhat\in C^\kappa((-2,2)\times B'_2)$, and
\begin{align}
0<\beta<3\kappa-2. \label{lift:beta}
\end{align}
For $\fluxlift$ defined by \eqref{lift:construction} and $\widehat M$ satisfying \eqref{lift:multiplier-regularity}, the extension
\begin{align}
\widehat{\glift}:=\widehat M\fluxlift,\qquad \widehat{\glift}(\theta,\xi',0,\eta)=\widehat M(\theta,\xi',0,\eta)\fluxhat(\theta,\xi') \label{lift:trace}
\end{align}
satisfies, with every kinetic norm on $\overline{\liftpatch}$,
\begin{align}
&\tnorm[\Big]{\widehat{\glift}}_{\Ckin^\beta}+\tnorm[\Big]{\nabla_\eta\widehat{\glift}}_{\Ckin^\beta}+\tnorm[\Big]{D_\eta^2\widehat{\glift}}_{\Ckin^\beta}+\tnorm[\Big]{\Yhat\widehat{\glift}}_{\Ckin^\beta} \label{lift:quantitative-norm}\\
&\quad\lesssim_{\kappa,\beta}\norm[\big]{\fluxhat}_{C^\kappa}\left(\tnorm[\big]{\widehat M}_{\Ckin^\beta}+\tnorm[\big]{\nabla_\eta\widehat M}_{\Ckin^\beta}+\tnorm[\Big]{D_\eta^2\widehat M}_{\Ckin^\beta}+\tnorm[\big]{\Yhat\widehat M}_{\Ckin^\beta}\right).\notag
\end{align}
The norm of $\fluxhat$ is on $(-2,2)\times B'_2$.
\end{lemma}

\begin{proof}
Lemma~\ref{lift:scalar-jets} gives the amplitudes and boundary values of $\fluxlift$ and its derivatives. We will first prove, for points at kinetic distance at most $r$,
\begin{align}
\abs[\big]{\Delta D_\eta^2\fluxlift}+\abs[\big]{\Delta\Yhat\fluxlift}\lesssim_\kappa[\,\fluxhat\,]_{C^\kappa}r^{3\kappa-2}. \label{lift:holder-increment-goal}
\end{align}
For close points we will compare the integral kernels; for separated points we will use the vanishing derivative bounds. We will then estimate the lower derivatives and multiply by $\widehat M$ to obtain \eqref{lift:quantitative-norm}. Write $\zeta_i=(\theta_i,\xi_i,\eta_i)$, $i=0,1$, and $\Delta H=H(\zeta_1)-H(\zeta_0)$. It suffices to treat \eqref{lift:kinetic-increments} with $0<r<\frac1{100}$, because the amplitude bounds control larger distances.

\proofstep{Increments of the travel time and boundary point} We will derive the changes in $\theta_b$ and $\zeta_b$ needed for the kernel comparison. Write $\Delta\xi_3=\eta_{0,3}\Delta\theta+\omega_3$, where $\abs{\omega_3}\lesssim r^3$. The quotient formula gives
\begin{align}
\Delta\theta_b=\frac{\eta_{0,3}\Delta\theta+\omega_3-\theta_{b,0}\Delta\eta_3}{\eta_{0,3}+\Delta\eta_3}. \label{eq:I:display:159}
\end{align}
In particular, when $r\leq c\sqrt{\min\{\theta_{b,0},\theta_{b,1}\}}$ and $c$ is sufficiently small,
\begin{align}
\abs{\Delta\theta_b}\lesssim r^2+\theta_{b,0}r,\qquad \frac12\theta_{b,0}\leq\theta_{b,1}\leq2\theta_{b,0}. \label{lift:tau-increment}
\end{align}
For the time component of $\zeta_b$, the leading transport increments cancel:
\begin{align}
\Delta\theta-\Delta\theta_b=\frac{(\theta_{b,0}+\Delta\theta)\Delta\eta_3-\omega_3}{\eta_{0,3}+\Delta\eta_3}. \label{eq:I:display:161}
\end{align}
For the tangential components, use $\Delta\xi'=\eta'_0\Delta\theta+O(r^3)$ and $\Delta(\eta'\theta_b)=\eta'_0\Delta\theta_b+\theta_{b,1}\Delta\eta'$. These identities give, for every pair under consideration,
\begin{align}
\abs{\Delta\zeta_b}\lesssim\theta_{b,0}r+r^3. \label{lift:xi-increment}
\end{align}
In the close range, \eqref{lift:tau-increment} and \eqref{lift:xi-increment} yield
\begin{align}
\abs{\Delta\eta}+\frac{\abs{\Delta\theta_b}}{\theta_{b,0}}+\frac{\abs{\Delta\zeta_b}}{\theta_{b,0}^{\frac32}}\lesssim\frac{r}{\sqrt{\theta_{b,0}}}. \label{lift:relative-increments}
\end{align}
We have $r^2/\theta_{b,0}\lesssim r/\sqrt{\theta_{b,0}}$ and $r^3/\theta_{b,0}^{\frac32}\lesssim r/\sqrt{\theta_{b,0}}$ in this range. This is the relative change of the kernel parameters that we need next.

\proofstep{Close points: comparison of the derivative kernels} Assume $r\leq c\sqrt{\min\{\theta_{b,0},\theta_{b,1}\}}$. We will prove \eqref{lift:holder-increment-goal} in this range by estimating each term \eqref{lift:term-structure}, with $p\in\{-1,-\frac12,0\}$. To justify the comparison without differentiating $\fluxhat$, differentiate the smooth kernel in the equivalent integral over $z$. For a zero-mean $\psi\in C_c^\infty(B_1)$ and positive $\tau$, its derivatives satisfy
\begin{align}
\int\left|\partial_b^\ell\partial_\tau^j\left[\tau^{-\frac92}\psi\left(\frac{b-z}{\tau^{\frac32}}\right)\right]\right|\dd z\lesssim_{\ell,j}\tau^{-\frac{3\abs{\ell}}{2}-j}. \label{lift:mollifier-derivatives}
\end{align}
Here $\ell$ is a multi-index and $j\geq0$. Each kernel has integral zero, including at $\ell=j=0$ because $\int\psi=0$. We will subtract the same value $\fluxhat(\zeta_{b,0})$ in every integral along the parameter segment. By \eqref{lift:relative-increments}, all kernel supports lie within distance $C\theta_{b,0}^{\frac32}$ of $\zeta_{b,0}$. The flux factor is therefore bounded by $C_\kappa[\,\fluxhat\,]_{C^\kappa}\theta_{b,0}^{\frac{3\kappa}{2}}$ throughout the comparison.

Along the segment joining $(\eta_0,\theta_{b,0},\zeta_{b,0})$ to $(\eta_1,\theta_{b,1},\zeta_{b,1})$, the positive travel times are comparable to $\theta_{b,0}$. All averaging balls stay inside $(-\frac32,\frac32)\times B'_{\frac32}$ by the bounds following \eqref{lift:construction}. The mean-value formula can therefore be applied throughout this segment. Differentiating the coefficient, the power of $\theta_b$, or the kernel gives, respectively, the three relative changes
\begin{align}
\abs{\Delta T}\lesssim_\kappa[\,\fluxhat\,]_{C^\kappa}\theta_{b,0}^{p+\frac{3\kappa}{2}}\left(\abs{\Delta\eta}+\frac{\abs{\Delta\theta_b}}{\theta_{b,0}}+\frac{\abs{\Delta\zeta_b}}{\theta_{b,0}^{\frac32}}\right). \label{lift:term-increment-estimate}
\end{align}
Using \eqref{lift:relative-increments}, $p\geq-1$ and $\theta_{b,0}\leq\frac18$, we obtain
\begin{align}
\abs{\Delta T}\lesssim_\kappa[\,\fluxhat\,]_{C^\kappa}r\theta_{b,0}^{\frac{3\kappa-3}{2}}\lesssim_\kappa[\,\fluxhat\,]_{C^\kappa}r^{3\kappa-2}. \label{lift:close-source-increment}
\end{align}
The last inequality follows from $\kappa<1$ and $r\leq c\sqrt{\theta_{b,0}}$. Summing the finitely many Hessian terms and the transport term \eqref{lift:transport-kernel} proves \eqref{lift:holder-increment-goal} for close points.

\proofstep{Separated points, including the boundary} We will prove the same increment bound when
\begin{align}
r\geq c\sqrt{\min\{\theta_{b,0},\theta_{b,1}\}}. \label{lift:far-regime}
\end{align}
Orient the pair so that $\theta_{b,0}$ is the smaller value. Interchanging the points preserves \eqref{lift:kinetic-increments} up to fixed constants, since the change in the position increment is $\Delta\eta\Delta\theta=O(r^3)$. Formula~\eqref{eq:I:display:159} now gives $\theta_{b,1}\lesssim\theta_{b,0}+r^2+\theta_{b,0}r+r^3\lesssim r^2$. Hence the amplitude bounds \eqref{lift:vanishing-jets-2} give
\begin{align}
\abs[\big]{\Delta D_\eta^2\fluxlift}+\abs[\big]{\Delta\Yhat\fluxlift}\lesssim_\kappa[\,\fluxhat\,]_{C^\kappa}\bigl(\theta_{b,0}^{\frac{3\kappa-2}{2}}+\theta_{b,1}^{\frac{3\kappa-2}{2}}\bigr)\lesssim_\kappa[\,\fluxhat\,]_{C^\kappa}r^{3\kappa-2}. \label{lift:far-source-increment}
\end{align}
The derivative values at a boundary point are zero by Lemma~\ref{lift:scalar-jets}, so the estimate includes such pairs. Together with the close-point estimate, this proves \eqref{lift:holder-increment-goal} and the required $\Ckin^\beta$ bounds for the two highest derivatives.

\proofstep{Lower derivatives and the estimate for $\widehat{\glift}$} We will complete \eqref{lift:quantitative-norm}. In the close range, the parameter derivatives of the kernel with $\psi=\varphi$ have zero integral, although $\int\varphi=1$. The preceding mean-value argument thus gives
\begin{align}
\abs[\big]{\Delta\fluxlift}\lesssim_\kappa[\,\fluxhat\,]_{C^\kappa}\theta_{b,0}^{\frac{3\kappa}{2}}\frac{r}{\sqrt{\theta_{b,0}}}\lesssim_\kappa[\,\fluxhat\,]_{C^\kappa}r. \label{eq:I:display:167}
\end{align}
For $\nabla_\eta\fluxlift$, the terms \eqref{lift:term-structure} have $p\geq-\frac12$. The same argument bounds their increments by $C_\kappa[\,\fluxhat\,]_{C^\kappa}r\theta_{b,0}^{\frac{3\kappa-2}{2}}\leq C_\kappa[\,\fluxhat\,]_{C^\kappa}r$.

In the separated range, both $\theta_{b,i}\lesssim r^2$ and \eqref{lift:xi-increment} gives $\abs{\Delta\zeta_b}\lesssim r^3$. Subtracting the two boundary-flux values and using \eqref{lift:vanishing-jets} and \eqref{lift:vanishing-jets-1}, we have
\begin{align}
\abs[\big]{\Delta\fluxlift}&\lesssim_\kappa[\,\fluxhat\,]_{C^\kappa}\bigl(\theta_{b,0}^{\frac{3\kappa}{2}}+\theta_{b,1}^{\frac{3\kappa}{2}}+\abs{\Delta\zeta_b}^\kappa\bigr)\lesssim_\kappa[\,\fluxhat\,]_{C^\kappa}r^{3\kappa}, \label{lift:lower-far-increments}\\
\abs[\big]{\Delta\nabla_\eta\fluxlift}&\lesssim_\kappa[\,\fluxhat\,]_{C^\kappa}\bigl(\theta_{b,0}^{\frac{3\kappa-1}{2}}+\theta_{b,1}^{\frac{3\kappa-1}{2}}\bigr)\lesssim_\kappa[\,\fluxhat\,]_{C^\kappa}r^{3\kappa-1}. \label{lift:gradient-far-increments}
\end{align}
Since $3\kappa-1>1$, the close and separated estimates prove the lower-derivative $\Ckin^\beta$ bounds throughout the range \eqref{lift:beta}. Including the amplitude estimates, we have
\begin{align}
\tnorm[\big]{\fluxlift}_{\Ckin^\beta}+\tnorm[\big]{\nabla_\eta\fluxlift}_{\Ckin^\beta}+\tnorm[\Big]{D_\eta^2\fluxlift}_{\Ckin^\beta}+\tnorm[\big]{\Yhat\fluxlift}_{\Ckin^\beta}\lesssim_{\kappa,\beta}\norm[\big]{\fluxhat}_{C^\kappa}. \label{lift:flux-holder-norm}
\end{align}
The product rule gives $\Yhat\widehat{\glift}=(\Yhat\widehat M)\fluxlift+\widehat M\Yhat\fluxlift$ and the corresponding first and second velocity derivative formulas. Apply the product inequality $\norm{ab}_{\Ckin^\beta}\lesssim\norm{a}_{\Ckin^\beta}\norm{b}_{\Ckin^\beta}$ to these formulas, using \eqref{lift:flux-holder-norm} and \eqref{lift:multiplier-regularity}. This proves \eqref{lift:quantitative-norm}; the trace identity in \eqref{lift:trace} follows from \eqref{lift:scalar-space}.
\end{proof}

\para{Weighted second-order estimate} We will use the lift estimates to control the second velocity and transport derivatives near the incoming boundary, and then convert these bounds to the physical grazing weight \eqref{eq:rho}.

\begin{proposition}\label{cb:graph}
Let $f$ be a weak solution of \eqref{cb:equation-residual} on $I$, satisfying \eqref{cb:regularity-class}, \eqref{cb:setting} and the diffuse boundary condition \eqref{cb:diffuse-law}. Assume the source bounds \eqref{lin:source-Gaussian} and \eqref{lin:source-holder}. For every $J\subset\subset I$, there is $\alpha_*>0$ such that
\begin{align}
\tnorm{\frac{\ue^{\alpha_*\abs{v}^2}\bigl(\abs{D_v^2f}+\abs{\Y f}\bigr)}{1+\graz^{-\frac{3}{2}}}}_{L^\infty(J\times\{x\in\Om:d(x)<\frac{d_0}{2}\}\times\R^3)}&\lesssim_J 1, \label{cb:pointwise-hessian}\\
\tnorm{\ue^{\alpha_*\abs{v}^2}\bigl(\abs{D_v^2f}+\abs{\Y f}\bigr)}_{L^\infty(J\times\{x\in\Om:d(x)\geq\frac{d_0}{2}\}\times\R^3)}&\lesssim_J 1. \label{cb:diffuse-interior-hessian}
\end{align}
Here $\graz$ is defined in \eqref{eq:rho}. The distributional derivatives are represented by functions on the indicated domains and are locally H\"older continuous away from grazing.
\end{proposition}

\begin{proof}
Fix the change of variables $\Psi_0$ in \eqref{lin:inverse-map} at a physical center $Z_0=(t_0,x_0,v_0)$, and set $\widehat f=f\circ\Psi_0$. Proposition~\ref{cb:gradient} has given the sharp H\"older bound \eqref{cb:diffuse-normalized-sharp}, Proposition~\ref{cb:flux} has given $\fluxhat\in C^\kappa$ with $\kappa>\frac23$, and Lemmas~\ref{lift:scalar-jets}--\ref{lift:characteristic} control the lift of this flux. We will use these estimates to prove, for $\zeta=(\theta,\xi,\eta)\in\Qplus_{\frac14}(\widehat Z_0)$ away from grazing,
\begin{align}
\abs[\big]{D_\eta^2\widehat f(\zeta)}+\abs[\big]{\Yhat\widehat f(\zeta)}&\lesssim_J\ue^{-a\abs{v_0}^2}\bigl(1+\widehat\graz(\zeta)^{-\frac32}\bigr),\qquad \widehat\graz=\xi_3^{\frac13}+\abs{\eta_3}. \label{cb:diffuse-normalized-hessian}
\end{align}
The derivative identities \eqref{lin:physical-hessian-conversion}--\eqref{lin:physical-transport-conversion} and weight comparison \eqref{cb:inflow-grazing-weight-ratio} will then give \eqref{cb:pointwise-hessian}; the interior estimate will give \eqref{cb:diffuse-interior-hessian}.

Fix small $c_0,c_1>0$ for the geometric choices below. We divide the points of $\Qplus_{\frac14}(\widehat Z_0)$ away from grazing into four cases:
\begin{enumerate}[label=\textup{Case \Roman*:}]
\item $\eta_3\leq0$.
\item $0<\eta_3\leq c_1$ and $\xi_3^{\frac13}>c_0\eta_3$.
\item $0<\eta_3\leq c_1$ and $\xi_3^{\frac13}\leq c_0\eta_3$.
\item $\eta_3>c_1$.
\end{enumerate}
In Cases I and II, we will use interior cylinders or cylinders with empty incoming boundary, subtract a constant, and apply the second-order estimate to the $r^{\frac12}$ oscillation of $\widehat f$. In Case III, we will dilate the small positive normal velocity to one, subtract the regularized lift, and estimate the difference and its source before undoing the dilation. In Case IV, we will use the lift near the boundary and interior estimates at larger heights to obtain a Gaussian bound without a grazing singularity.

\proofstep{The equation for $\widehat f$ and the available bounds} We first verify the coefficient, source and weak-solution hypotheses needed for the local second-order estimates. Lemma~\ref{cb:rough-flux-coefficients} gives \eqref{cb:diffuse-coefficient-holder}--\eqref{cb:diffuse-source-holder}. Choose the radius in $\Psi_0$ as in Lemma~\ref{cb:normalized-bounds}. Equations~\eqref{lin:normalized-coefficient-bounds}--\eqref{lin:normalized-source-holder} then give uniform ellipticity, bounded coefficient H\"older norms, and a Gaussian H\"older bound for $\widehat S_f$ on the equation domain. Lemma~\ref{cb:normalized-equation} gives the weak formulation \eqref{lin:boundary-weak-setting}, and \eqref{cb:regularity-class} gives the local $L^2$ velocity derivative. Thus the equation for $\widehat f$ satisfies \eqref{lin:coefficients}--\eqref{lin:coefficient-norm}.

Take $\frac23<\kappa<1$ from Proposition~\ref{cb:flux}, choose $0<\beta<1$ below the coefficient and source exponents, and fix $0<\beta'<\min\{\beta,3\kappa-2\}$. We will use the sharp bound \eqref{cb:diffuse-normalized-sharp} on the larger half-cylinder.

\proofstep{Cases I and II: interior cylinders and cylinders with empty incoming boundary} We will prove \eqref{cb:diffuse-normalized-hessian} in both cases. In Case I, suppose first that $\widehat\graz$ is small. Take $r=c\widehat\graz$ and the grazing point $\zeta_*$ from \eqref{eq:I:display:298}. The geometry \eqref{cb:inflow-grazing-cylinder-geometry}--\eqref{cb:inflow-grazing-distance} places $\zeta_*$ and $\Qplus_{64r}(\zeta)$ in the domain of the sharp estimate. If $\abs{\eta_3}\geq\xi_3^{\frac13}$, the cylinder has no incoming boundary. Otherwise,
\begin{align}
(64r)^2\abs{\eta_3}+(64r)^3\lesssim(c^2+c^3)\xi_3<\frac12\xi_3, \label{cb:diffuse-cylinder-containment}
\end{align}
so the full cylinder is interior.

In Case II, set $s:=\eta_3\in(0,c_1]$ and take $r=c\xi_3^{\frac13}$, with $c$ sufficiently small depending on $c_0$. If $\xi_3$ is small, then
\begin{align}
(64r)^2s+(64r)^3\lesssim\left(\frac{c^2}{c_0}+c^3\right)\xi_3<\frac12\xi_3. \label{cb:diffuse-interior-containment}
\end{align}
This cylinder is interior, and $r\simeq\widehat\graz$ because $s<c_0^{-1}\xi_3^{\frac13}$. The grazing point and enlarged cylinder again lie in the domain of \eqref{cb:diffuse-normalized-sharp} after decreasing the fixed height threshold.

On each of these small cylinders, subtract $\widehat f(\zeta_*)$. Its difference from $\widehat f$ has $L^\infty$ norm at most $C_J\ue^{-a\abs{v_0}^2}r^{\frac12}$ on $\Qplus_{16r}(\zeta)$, by \eqref{cb:diffuse-normalized-sharp} and \eqref{cb:inflow-grazing-distance}. The source becomes $\widehat S_f+\locc\widehat f(\zeta_*)$, with $\Ckin^\beta$ norm at most $C_J\ue^{-a\abs{v_0}^2}$. Proposition~\ref{lin:order-two} applies in the case of empty incoming boundary; Proposition~\ref{lin:interior-inputs}\ref{lin:interior} applies on full interior cylinders. The geometric factor in the former is bounded because $\abs{\eta_3}/r$ is bounded. Their conclusions give
\begin{align}
r^2\Bigl(\abs[\big]{D_\eta^2\widehat f(\zeta)}+\abs[\big]{\Yhat\widehat f(\zeta)}\Bigr)\lesssim_J\ue^{-a\abs{v_0}^2}\bigl(r^{\frac12}+r^2+r^{2+\beta}\bigr). \label{cb:diffuse-outgoing-hessian}
\end{align}
Since $r\simeq\widehat\graz$, division by $r^2$ proves \eqref{cb:diffuse-normalized-hessian} in these small-cylinder cases. When $\widehat\graz$ is bounded below and $\eta_3\leq0$, the fixed-radius cylinders from the proof of Proposition~\ref{cb:inflow-gradient} are interior or have empty incoming boundary. When $0<s\leq c_1$ and the height is bounded below, a fixed smaller radius gives an interior cylinder by the displacement estimate in \eqref{cb:diffuse-interior-containment}. Apply the same second-order estimates directly with the Gaussian amplitude bound. The radius is fixed, and the geometric factor has at most polynomial growth in $\bracket{v_0}$, which \eqref{lin:Gaussian-absorption} absorbs. This proves \eqref{cb:diffuse-normalized-hessian} throughout Cases I and II.

\proofstep{Case III: construct a difference with zero incoming trace} Let $\zeta_0=(\theta_0,\xi_0,\eta_0)$ satisfy $s:=\eta_{0,3}\in(0,c_1]$ and $\xi_{0,3}^{\frac13}\leq c_0s$. In this region, $s\leq\widehat\graz(\zeta_0)\leq(1+c_0)s$. We will construct the difference on the fixed domain \eqref{lift:normalized-patch} and obtain the lift bounds \eqref{eq:I:display:282}--\eqref{cb:diffuse-lift-lower-derivatives}, which will control its source.

Put $\zeta_*=(\theta_0,(\xi_0',0),(\eta_0',0))$ and use the translation and dilation $T_{\zeta_*}\circ S_s$. They preserve the half-space and send the target to $\bigl(0,(0,\frac{\xi_{0,3}}{s^3}),(0,1)\bigr)$. Define
\begin{align}
\widehat f_s:=\widehat f\circ T_{\zeta_*}\circ S_s,\qquad\widehat M_s:=\widehat M\circ T_{\zeta_*}\circ S_s,\qquad\fluxhat_s(\theta,\xi'):=\fluxhat\bigl(\theta_0+s^2\theta,\xi_0'+s^3\xi'+s^2\theta\eta_0'\bigr). \label{cb:diffuse-rescaled-functions}
\end{align}
Choose a fixed $r_*>0$ so that the target's half-cylinder of radius $64r_*$ lies in \eqref{lift:normalized-patch} and avoids grazing, uniformly for $0\leq\frac{\xi_{0,3}}{s^3}\leq c_0^3$. This is possible because its normal velocity equals one and $c_0^3$ is smaller than the upper height of that domain. The differences from $\zeta_0$, before this dilation, satisfy
\begin{align}
\Delta\xi-\eta_0\Delta\theta=s^3(\xi',\xi_3-\theta)-\xi_{0,3}e_3,\qquad\Delta\theta=s^2\theta,\qquad\Delta\eta=s(\eta',\eta_3-1). \label{cb:rescaled-increments}
\end{align}
Consequently the image of $\liftpatch$ and the boundary rectangle for its data lie within kinetic distance $Cs$ of $\zeta_0$, independently of $\eta_0'$. Taking $c_1$ sufficiently small places these sets in the domain of the original equation and the sharp H\"older estimate.

On $(-2,2)\times B'_2$, the ordinary flux estimate \eqref{rev14:diffuse-flux-regularity} gives
\begin{align}
\tnorm[\big]{\fluxhat_s-\fluxhat_s(0,0)}_{C^\kappa}\lesssim_J\left(s^2(1+\abs{\eta_0'})+s^3\right)^\kappa\tnorm[\big]{\fluxhat}_{C^\kappa}\lesssim_J\bracket{v_0}^N s^{2\kappa}. \label{cb:diffuse-rescaled-flux-bound}
\end{align}
The affine map in \eqref{cb:diffuse-rescaled-functions} has ordinary Lipschitz constant at most $C\bigl(s^2(1+\abs{\eta_0'})+s^3\bigr)$, and subtraction of the value at zero bounds the supremum by the same seminorm. The formula for $\Psi_0$ bounds $\abs{\eta_0'}$ by a fixed power of $\bracket{v_0}$.

Use \eqref{lift:construction} for the data $\fluxhat_s-\fluxhat_s(0,0)$, denote the resulting regularized extension by $\fluxlift_s$, and set
\begin{align}
\widehat{\glift}_s:=\widehat M_s\bigl(\fluxhat_s(0,0)+\fluxlift_s\bigr),\qquad\widehat u_s:=\widehat f_s-\widehat{\glift}_s. \label{cb:diffuse-rescaled-lift}
\end{align}
The trace identity \eqref{lift:trace} gives $\trm\widehat u_s=0$. Lemma~\ref{lift:characteristic} and \eqref{cb:diffuse-rescaled-flux-bound} give the factor $s^{2\kappa}$ for the varying flux. The constant part differentiates only through $\widehat M_s$, for which
\begin{align}
\nabla_\eta\widehat M_s=s(\nabla_\eta\widehat M)\circ T_{\zeta_*}\circ S_s,\qquad D_\eta^2\widehat M_s=s^2(D_\eta^2\widehat M)\circ T_{\zeta_*}\circ S_s,\qquad\Yhat\widehat M_s=s^2(\Yhat\widehat M)\circ T_{\zeta_*}\circ S_s. \label{cb:diffuse-Maxwellian-scaling}
\end{align}
The three derivatives on the left have full $\Ckin^{\beta'}$ norms bounded by a common Gaussian factor times $s,s^2,s^2$, respectively. Combining these identities with \eqref{lift:quantitative-norm} and absorbing polynomial velocity factors gives, on $\overline{\liftpatch}$,
\begin{align}
\tnorm[\Big]{D_\eta^2\widehat{\glift}_s}_{\Ckin^{\beta'}}+\tnorm[\Big]{\Yhat\widehat{\glift}_s}_{\Ckin^{\beta'}}\lesssim_J\ue^{-a\abs{v_0}^2}(s^{2\kappa}+s^2), \label{eq:I:display:282}\\
\tnorm[\Big]{\nabla_\eta\widehat{\glift}_s}_{\Ckin^{\beta'}}\lesssim_J\ue^{-a\abs{v_0}^2}(s^{2\kappa}+s),\qquad\tnorm[\Big]{\widehat{\glift}_s}_{\Ckin^{\beta'}}\lesssim_J\ue^{-a\abs{v_0}^2}. \label{cb:diffuse-lift-lower-derivatives}
\end{align}
\proofstep{Case III: estimate the source and the difference} We will derive the weak equation for $\widehat u_s$ and prove the two bounds \eqref{cb:diffuse-lifted-bounds}. Dilation of \eqref{lin:Landau-nondivergence} and subtraction of the lift give
\begin{align}
\Yhat\widehat u_s-\locA_s:D_\eta^2\widehat u_s=\locb_s\cdot\nabla_\eta\widehat u_s+\locc_s\widehat u_s+\widehat S_s, \label{cb:diffuse-lifted-equation}\\
\widehat S_s:=s^2\widehat S_f\circ T_{\zeta_*}\circ S_s-\Yhat\widehat{\glift}_s+\locA_s:D_\eta^2\widehat{\glift}_s+\locb_s\cdot\nabla_\eta\widehat{\glift}_s+\locc_s\widehat{\glift}_s, \label{cb:diffuse-lifted-source}
\end{align}
where $\locA_s=\locA\circ T_{\zeta_*}\circ S_s$, $\locb_s=s\locb\circ T_{\zeta_*}\circ S_s$ and $\locc_s=s^2\locc\circ T_{\zeta_*}\circ S_s$. These coefficients are determined by the original $f$.

To justify subtraction in the weak formulation, the lift is $C^2$ for $\xi_3>0$, continuous up to $\xi_3=0$, and its transport, first velocity and second velocity derivatives are bounded by \eqref{eq:I:display:282}--\eqref{cb:diffuse-lift-lower-derivatives}. Integrate its equation on $\{\xi_3>\delta\}$ and let $\delta\to0$. The boundary term converges to the trace in \eqref{lift:trace}, and the bulk terms converge by these derivative bounds. Subtraction from the weak identity for $\widehat f_s$ proves \eqref{lin:boundary-weak-setting} for $\widehat u_s$ with zero incoming trace. Its first-order divergence coefficient is $s\locB\circ T_{\zeta_*}\circ S_s$. The weak-solution class supplies the local $L^2$ velocity gradient of $\widehat f_s$, and the bounded lift gradient gives this property for $\widehat u_s$, also up to the incoming boundary.

We next estimate each source term. The full $\Ckin^{\beta'}$ norms of $\locA_s,\locb_s,\locc_s$ are bounded by $C,Cs,Cs^2$, and $\tnorm[\big]{\nabla_\eta\locA_s}_{L^\infty}\leq Cs$. The first term in \eqref{cb:diffuse-lifted-source} has norm at most $C_Js^2\ue^{-a\abs{v_0}^2}$. Equation~\eqref{eq:I:display:282} and the product estimate \eqref{cb:inflow-source-product} bound the transport and diffusion terms by $C_J\ue^{-a\abs{v_0}^2}(s^{2\kappa}+s^2)$. The remaining two terms satisfy
\begin{align}
\tnorm[\Big]{\locb_s\cdot\nabla_\eta\widehat{\glift}_s}_{\Ckin^{\beta'}}\lesssim_J s\ue^{-a\abs{v_0}^2}(s^{2\kappa}+s),\qquad\tnorm[\Big]{\locc_s\widehat{\glift}_s}_{\Ckin^{\beta'}}\lesssim_J s^2\ue^{-a\abs{v_0}^2}. \label{cb:lifted-source-orders}
\end{align}
In particular, the constant flux gives an $s$ term in $\nabla_\eta\widehat{\glift}_s$; multiplication by $\locb_s=O(s)$ makes this contribution $O(s^2)$.

For the difference, continuity of the diffuse law gives $\widehat f(\zeta_*)=\widehat M(\zeta_*)\fluxhat_s(0,0)$. Consequently
\begin{align}
\widehat u_s=\widehat f_s-\widehat f(\zeta_*)-\bigl(\widehat M_s-\widehat M(\zeta_*)\bigr)\fluxhat_s(0,0)-\widehat M_s\fluxlift_s. \label{cb:diffuse-remainder-decomposition}
\end{align}
The first difference has supremum at most $C_J\ue^{-a\abs{v_0}^2}s^{\frac12}$ by \eqref{cb:diffuse-normalized-sharp} and \eqref{cb:rescaled-increments}. The smooth Maxwellian factor gives a bound $C_J\ue^{-a\abs{v_0}^2}s$ for the next term. Positivity and unit integral of the fixed averaging kernel, together with \eqref{cb:diffuse-rescaled-flux-bound}, bound the last term by $C_J\ue^{-a\abs{v_0}^2}s^{2\kappa}$. Since $s\leq1$ and $2\kappa>\frac12$, we have proved
\begin{align}
\tnorm[\big]{\widehat u_s}_{L^\infty(\liftpatch)}\lesssim_J\ue^{-a\abs{v_0}^2}s^{\frac12},\qquad\tnorm[\big]{\widehat S_s}_{\Ckin^{\beta'}(\overline{\liftpatch})}\lesssim_J\ue^{-a\abs{v_0}^2}(s^{2\kappa}+s^2). \label{cb:diffuse-lifted-bounds}
\end{align}
\proofstep{Case III: obtain the second-order bound} We will now prove \eqref{cb:diffuse-normalized-hessian} for $0<s\leq c_1$, $\xi_{0,3}^{\frac13}\leq c_0s$. The coefficient bounds above and \eqref{cb:diffuse-lifted-bounds} verify \eqref{lin:coefficients}--\eqref{lin:coefficient-norm}; the weak formulation and zero incoming trace were verified in the preceding step. The radius-$64r_*$ half-cylinder lies in the lift domain and avoids grazing. Proposition~\ref{lin:order-two} applies at radius $r_*$; if the full cylinder is interior, Proposition~\ref{lin:interior-inputs}\ref{lin:interior} gives the same estimate. The target normal velocity is one and $r_*$ is fixed, so the geometric and inverse-radius factors are bounded. Add the lift derivatives from \eqref{eq:I:display:282} to obtain, at the dilated target,
\begin{align}
\abs[\big]{D_\eta^2\widehat f_s}+\abs[\big]{\Yhat\widehat f_s}\lesssim_J\ue^{-a\abs{v_0}^2}(s^{\frac12}+s^{2\kappa}+s^2). \label{cb:diffuse-rescaled-hessian}
\end{align}
Both derivatives acquire a factor $s^{-2}$ when returning to $\widehat f$. Therefore
\begin{align}
\abs[\big]{D_\eta^2\widehat f(\zeta_0)}+\abs[\big]{\Yhat\widehat f(\zeta_0)}&\lesssim_J\ue^{-a\abs{v_0}^2}(s^{-\frac32}+s^{2\kappa-2}+1)\lesssim_J\ue^{-a\abs{v_0}^2}\widehat\graz(\zeta_0)^{-\frac32}. \label{eq:I:display:283}
\end{align}
The last inequality uses $\kappa>\frac23$, $s\leq1$ and $s\simeq\widehat\graz(\zeta_0)$. This proves \eqref{cb:diffuse-normalized-hessian} in Case III.

\proofstep{Case IV: construct the lift and bound its source} Let $s:=\eta_{0,3}>c_1$. Our goal in this and the next step is
\begin{align}
\abs[\big]{D_\eta^2\widehat f(\zeta_0)}+\abs[\big]{\Yhat\widehat f(\zeta_0)}\lesssim_J\ue^{-a'\abs{v_0}^2} \label{cb:diffuse-large-incoming-goal}
\end{align}
for some $a'>0$. This implies \eqref{cb:diffuse-normalized-hessian} on the whole region $s>c_1$. We first construct a lift where the height is small and prove the source bound \eqref{cb:diffuse-large-source-bound}.

Translate time and tangential position, and subtract $\eta_0'$ from the tangential velocity by a Galilean transformation. The target becomes $(0,(0,h),(0,s))$. This transformation preserves the half-space, $D_\eta$ and $\Yhat$; in this step we keep the same symbols for the composed functions and coefficients. Choose fixed $0<c_2<\min\{c_1/4,1/64\}$ and sufficiently small $H>0$. For $h\leq H/4$, work on
\begin{align}
\abs{\theta}<\frac{H}{s},\qquad\abs{\xi'}<\frac1{64},\qquad0<\xi_3<H,\qquad\abs{\eta'}<c_2,\qquad\abs{\eta_3-s}<c_2. \label{cb:large-incoming-slab}
\end{align}
The normal displacement from transport relative to the target is bounded by $3H$, and its tangential center velocity is zero. Taking $H$ small places this set inside $\Qc_{\frac12}(\zeta_0)$ in the variables before the Galilean transformation, and hence inside $\Qc_1(\widehat Z_0)$ by \eqref{lin:cylinder-composition}. The backward travel time and averaging radius satisfy
\begin{align}
0\leq\theta_b\lesssim\frac{H}{s},\qquad\theta_b^{\frac32}\leq\frac{H}{4s}, \label{cb:large-incoming-exit}
\end{align}
after choosing $H$ small relative to $c_1$. Thus the boundary points $\zeta_b$ and their averaging balls lie in $(-4H/s,4H/s)\times B'_{1/32}$, inside the domain of the boundary data.

Use the same integral \eqref{lift:construction} on this set and put $\widehat{\glift}=\widehat M\fluxlift$. The derivative and H\"older estimates in Lemmas~\ref{lift:scalar-jets}--\ref{lift:characteristic} remain applicable: their coefficient formulas involve reciprocal powers of $\eta_3$ and bounded tangential velocities, while here $\eta_3\geq c_1/2$ and $\abs{\eta'}<c_2$. The cancellation \eqref{eq:I:display:161}, the close-point estimate \eqref{lift:term-increment-estimate}, and the separated-point estimate \eqref{lift:far-source-increment} therefore give the same derivative orders. The flux norm after the tangential transformation has at most polynomial growth in $\bracket{v_0}$, and the Maxwellian factors have Gaussian decay. Hence the lift and its first two velocity and transport derivatives have full $\Ckin^{\beta'}$ norms bounded by $C_J\bracket{v_0}^N\ue^{-a\abs{v_0}^2}$ on \eqref{cb:large-incoming-slab}.

Subtracting the lift gives the weak formulation \eqref{lin:boundary-weak-setting} with zero incoming trace, by the same integration on $\{\xi_3>\delta\}$ used above. The difference has a local $L^2$ velocity gradient. Its source consists of the five terms in \eqref{cb:diffuse-lifted-source}, without the dilation by $s$. The product estimate \eqref{cb:inflow-source-product} and the lift bounds give
\begin{align}
\tnorm[\Big]{\widehat f-\widehat{\glift}}_{L^\infty}+\tnorm[\Big]{\widehat S_f-\Yhat\widehat{\glift}+\locA:D_\eta^2\widehat{\glift}+\locb\cdot\nabla_\eta\widehat{\glift}+\locc\widehat{\glift}}_{\Ckin^{\beta'}}\lesssim_J\bracket{v_0}^N\ue^{-a\abs{v_0}^2}. \label{cb:diffuse-large-source-bound}
\end{align}
All norms are on \eqref{cb:large-incoming-slab}. We have therefore obtained the zero incoming equation and the amplitude and source bounds needed for the second-order estimate when $h\leq H/4$.

\proofstep{Case IV: prove the derivative estimate} We will prove \eqref{cb:diffuse-large-incoming-goal} both near the boundary and at larger heights. Set
\begin{align}
r_s=\varepsilon_0(1+s)^{-\frac12},\qquad(64r_s)^2s+(64r_s)^3<\frac{H}{8},\qquad(64r_s)^2<\frac{H}{4s}, \label{cb:large-incoming-radius}
\end{align}
where $\varepsilon_0>0$ is fixed sufficiently small and $64r_s<c_2/2$. The first inequality controls normal transport displacement; the second places the time interval inside \eqref{cb:large-incoming-slab}. The tangential displacement and velocity radius are also contained in that set. Thus for $h\leq H/4$, the enlarged half-cylinder lies in the lift domain and has normal velocity bounded below by $c_1/2$. The coefficient bounds, weak formulation and source estimate \eqref{cb:diffuse-large-source-bound} verify the hypotheses of Proposition~\ref{lin:order-two}. Its inverse-radius and geometric factors satisfy
\begin{align}
r_s^{-2}\lesssim1+s,\qquad\left(1+\frac{s}{r_s}\right)^N\lesssim(1+s)^{\frac{3N}{2}}. \label{cb:diffuse-large-radius-loss}
\end{align}
Apply \eqref{lin:order-two-bound} and add the lift derivatives. If $h>H/4$, \eqref{cb:large-incoming-radius} instead places the enlarged full cylinder in $\{\xi_3>H/8\}$. Proposition~\ref{lin:interior-inputs}\ref{lin:interior} then applies directly to $\widehat f$, using its Gaussian amplitude and source bounds. In both cases we obtain
\begin{align}
\abs[\big]{D_\eta^2\widehat f(\zeta_0)}+\abs[\big]{\Yhat\widehat f(\zeta_0)}&\lesssim_J\bracket{v_0}^N\ue^{-a\abs{v_0}^2}r_s^{-2}\left(1+\frac{s}{r_s}\right)^N\lesssim_J\ue^{-a'\abs{v_0}^2}. \label{cb:diffuse-large-incoming-hessian}
\end{align}
The last inequality uses \eqref{cb:diffuse-large-radius-loss}, the polynomial bound on $s$ from \eqref{lin:normalized-center}, and \eqref{lin:Gaussian-absorption}. This proves \eqref{cb:diffuse-large-incoming-goal}. Since $\widehat\graz\geq s>c_1$, this completes Case IV.

\proofstep{Obtain the physical estimates} Cases I--IV prove \eqref{cb:diffuse-normalized-hessian} on $\Qplus_{\frac14}(\widehat Z_0)$ away from grazing. We will now deduce the two estimates in the proposition. The derivative identities \eqref{lin:physical-hessian-conversion}--\eqref{lin:physical-transport-conversion} give \eqref{cb:pullback-estimate} with $f$ and $\widehat f$ in place of $u$ and $\widehat u$. The transport identity also contains $\nabla_\eta\widehat f$. Proposition~\ref{cb:gradient}, applied in the variables of $\Psi_0$ as in \eqref{eq:I:display:300}, bounds this term by $C_J\ue^{-a\abs{v_0}^2}(1+\widehat\graz^{-\frac12})$, which is controlled by the right-hand side of \eqref{cb:diffuse-normalized-hessian}. The weight ratio \eqref{cb:inflow-grazing-weight-ratio} with exponent $\frac32$ and Gaussian absorption therefore give
\begin{align}
\frac{\abs{D_v^2f}+\abs{\Y f}}{1+\graz^{-\frac32}}\lesssim_J\bracket{v_0}^N\ue^{-a\abs{v_0}^2}\lesssim_J\ue^{-\alpha_*\abs{v}^2}. \label{cb:diffuse-physical-local-hessian}
\end{align}
The constants are uniform in the physical center. Taking the essential supremum over the boundary collar proves \eqref{cb:pointwise-hessian}. On $\{d\geq\frac{d_0}{2}\}$, use the interior change of variables and Proposition~\ref{lin:interior-inputs}\ref{lin:interior}. The inverse radii grow at most polynomially in $\bracket{v_0}$, so \eqref{lin:Gaussian-absorption} gives \eqref{cb:diffuse-interior-hessian} without the grazing factor.

Finally, the local $\Ckin^{2+\beta'}$ estimates give H\"older continuous second velocity and transport derivatives away from grazing. The lift derivatives have this regularity by Lemma~\ref{lift:characteristic}, and the first-order bound controls the additional term in the transport transformation. These local representatives agree on overlaps with the distributional derivatives of $f$, completing the proof.
\end{proof}

\subsection{The Initial Time and the Compatibility Conditions}\label{cb:section-initial}

From now on, $f$ will solve the Landau equation \eqref{eq:normalized-problem}, so $\sourceR=0$ in \eqref{cb:equation-residual} on positive times. Under the corresponding data assumptions, the preceding subsections give the following estimates on every $J\subset\subset(0,T_1)$. For in-flow, we have the H\"older and first-order bounds \eqref{cb:inflow-holder-norm}--\eqref{cb:inflow-interior-gradient} and the second-order bounds \eqref{cb:pointwise-grazing-hessian}--\eqref{cb:inflow-interior-hessian}. For diffuse-reflection, we have the flux bound \eqref{rev14:diffuse-flux-regularity}, the H\"older and first-order bounds \eqref{cb:diffuse-sharp-holder}--\eqref{cb:diffuse-interior-gradient}, and the second-order bounds \eqref{cb:pointwise-hessian}--\eqref{cb:diffuse-interior-hessian}.

We will extend these estimates to $[0,T]$, for every $T<T_1$, under the compatibility conditions at $t=0$. We will construct $\bar f$ on $(-t_*,T_1)$ with the prescribed boundary law and $\bar f=f$ for $t\geq0$. Its additional source $\sourceR$ may be nonzero for $t<0$. We will verify \eqref{lin:source-Gaussian} for the first-order estimate and \eqref{lin:source-holder} for the second-order estimate, then apply the preceding results on an interval containing $[0,T]$.

Assume in addition to \eqref{cb:setting} that
\begin{align}
\delta\leq\frac{1}{2}\delta_A. \label{cb:initial-smallness}
\end{align}
We will choose $t_*$ so that $\gnorm[\big]{\infty}{\frac{\alpha}{2}}{\bar f(t)}\leq2\delta\leq\delta_A$ on $(-t_*,0)$, which supplies the coefficient bounds required by the regularity estimates.

\begin{lemma}[In-flow estimates at the initial time]\label{cb:inflow-initial}
Let $f$ be a weak solution of \eqref{eq:normalized-problem} on $[0,T_1]$ in the sense of Definition~\ref{def:energy}, with incoming data $g$ and initial data $f_0$. Assume \eqref{cb:setting} and \eqref{cb:initial-smallness}. Suppose that $f_0\in C^5(\overline\Om\times\R^3)$ and $g\in C^5([0,T_1]\times\overline\Om\times\R^3)$ satisfy the Gaussian derivative bounds \eqref{in:eq:data-derivatives} and \eqref{in:eq:data-derivatives-g}. Let $a_0$ be the expression defined in \eqref{in:eq:compat1}.

For every $0<T<T_1$, the following holds, with constants depending on $T$ and the finite derivative bounds of the data.
\begin{enumerate}
\item If $\trm f_0=g(0)$, then $f\in\Ckinloc^{\frac{1}{2}}([0,T]\times\overline\Om\times\R^3)$. There are $\alpha_T>0$ and $0<C_T<\infty$ such that
\begingroup\usetagform{enumerated}
\begin{align}
\tnorm{\frac{\ue^{\alpha_T\abs{v}^2}\nabla_vf}{1+\graz^{-\frac{1}{2}}}}_{L^\infty([0,T]\times\left\{d<\frac{d_0}{2}\right\}\times\R^3)}\leq C_T. \label{cb:initial-inflow-gradient}
\end{align}
\endgroup
\item If also $\trm a_0=\partial_tg(0)$, then, after decreasing $\alpha_T$,
\begingroup\usetagform{enumerated}
\begin{align}
\tnorm{\frac{\ue^{\alpha_T\abs{v}^2}\bigl(\abs{D_v^2f}+\abs{\Y f}\bigr)}{1+\graz^{-\frac{3}{2}}}}_{L^\infty([0,T]\times\left\{d<\frac{d_0}{2}\right\}\times\R^3)}\leq C_T. \label{cb:initial-inflow-hessian}
\end{align}
\endgroup
\end{enumerate}
On $\left\{d\geq\frac{d_0}{2}\right\}$, the same bounds hold with the denominators replaced by $1$. Away from grazing at $t=0$, the gradient in part (i) is $\nabla_vf_0$. In part (ii), the second velocity derivatives are $D_v^2f_0$ and $\Y f(0)=\Ga(f_0,f_0)-\Lop f_0$. On compact positive-time intervals, Propositions~\ref{cb:inflow-gradient} and \ref{cb:strong-upgrade} give these conclusions without initial compatibility conditions.
\end{lemma}

\begin{proof}
We will construct an extension satisfying the weak equation \eqref{cb:equation-residual} on $(-t_*,T_1)$. Proposition~\ref{cb:inflow-gradient} gives its first-order estimate once we verify \eqref{cb:regularity-class}, Gaussian smallness, the incoming-data bound \eqref{cb:inflow-incoming-holder}, and \eqref{lin:source-Gaussian}. For the second-order estimate, we will use the proof of Proposition~\ref{cb:strong-upgrade}. It uses only Gaussian bounds on the weak derivatives of the lift through total order three, with at most two position derivatives, and the source bound \eqref{lin:source-holder}. We will verify these requirements across $t=0$, with the lift estimate stated in \eqref{cb:extended-lift-bounds}.

\proofstep{Extend the incoming data with controlled derivatives} We will obtain the lift bound \eqref{cb:extended-lift-bounds} and the incoming estimate \eqref{cb:inflow-incoming-holder} across $t=0$. The Gaussian derivative assumption \eqref{in:eq:data-derivatives-g} controls $g$ only on $\Sigma_-$. We therefore use the extension $\glift$ supplied by Lemma~\ref{cb:inflow-lift}, whose derivatives satisfy the bulk bound \eqref{eq:I:display:294}. For negative times define
\begin{align}
\glift(t):=\glift(0)+t\partial_t\glift(0)+\frac{t^2}{2}\partial_t^2\glift(0). \label{cb:negative-lift}
\end{align}
Set $\alpha_1:=\frac{\alpha}{2}$. The extension matches the value and first two time derivatives of the positive-time lift. Its third weak time derivative therefore has no measure at $t=0$. On negative times, a mixed derivative $\partial_t^j\partial_x^\ell\partial_v^k$ of total order at most three is a sum of derivatives $\partial_t^r\partial_x^\ell\partial_v^k\glift(0)$ with $j\leq r\leq2$. These require at most $r+\abs{\ell}+\abs{k}\leq5$ derivatives of the original lift. Consequently \eqref{eq:I:display:294} gives, for $t_*\leq1$,
\begin{align}
\sum_{j+\abs{\ell}+\abs{k}\leq3,\,\abs{\ell}\leq2}\tnorm[\Big]{\partial_t^j\partial_x^\ell\partial_v^k\glift}_{L^\infty(-t_*,T;L^{\infty,\alpha_1})}\lesssim_{T,\alpha_1} 1. \label{cb:extended-lift-bounds}
\end{align}
Thus the extended lift satisfies the derivative bounds used in the proof of Proposition~\ref{cb:strong-upgrade}, including across $t=0$. We will use it to construct the negative-time solution and, in part (ii), to subtract its boundary value. To verify the incoming-data hypothesis of the first-order estimate, set $g=\trm\glift$ for $t<0$. The value and first two time derivatives match those of the original $g$ at zero. For $Z=(t,x,v)$ and $Z'=(t',x',v')$ with $\abs{v-v_0},\abs{v'-v_0}\leq\bracket{v_0}^{-1}$ and $r:=\dkin(Z',Z)\leq1$, we have
\begin{align}
\abs{t-t'}\leq r^2,\qquad \abs{v-v'}\leq r,\qquad \abs{x-x'}\leq r^3+C\bracket{v_0}r^2. \label{cb:initial-kinetic-increments}
\end{align}
Taylor's formula and \eqref{cb:extended-lift-bounds} therefore give, for incoming boundary points $Z,Z'$ in one of the fixed boundary neighborhoods and some $0<\alpha_2<\alpha_1$,
\begin{align}
\abs[\big]{g(Z')-g(Z)-\nabla_v\glift(Z)\cdot(v'-v)}&\lesssim_T\ue^{-\alpha_2\abs{v_0}^2}r^2. \label{cb:extended-incoming-polynomial}
\end{align}
Here and below, polynomial velocity factors are absorbed by decreasing the Gaussian exponent as in \eqref{lin:Gaussian-absorption}. The Gaussian bound controls separated points. The polynomial characterization \eqref{eq:kinetic-holder-polynomial} thus verifies \eqref{cb:inflow-incoming-holder} for every fixed $0<\epsilon<1$.

\proofstep{Estimate the initial time derivative} We will prove \eqref{cb:initial-derivative-Gaussian}, which controls the source of the negative-time extension in part (ii). Recall from \eqref{in:eq:compat1} that $a_0=-v\cdot\nabla_xf_0-\Lop f_0+\Ga(f_0,f_0)$. Equivalently,
\begin{align}
\mh a_0=-v\cdot\nabla_x(\mh f_0)+Q(\m,\mh f_0)+Q(\mh f_0,\m)+Q(\mh f_0,\mh f_0). \label{eq:I:display:309}
\end{align}
The identity $Q(G,H)=(a*G):D_v^2H-(c*G)H$ shows that three further derivatives use at most five derivatives of $f_0$. Velocity derivatives of a convolution are placed on its Gaussian factor. The kernels $a$ and $c$ are locally integrable when $\gamma>-3$, and at $\gamma=-3$ the second term is $8\pi GH$. After division by $\mh$, the resulting terms have Gaussian decay with finite polynomial factors. Thus
\begin{align}
\sum_{\abs{k}\leq3}\gnorm{\infty}{\alpha_1}{\partial_{x,v}^ka_0}<\infty. \label{cb:initial-derivative-Gaussian}
\end{align}
\proofstep{Construct the extension for the first-order estimate} Assume $\trm f_0=g(0)$. We will preserve this boundary condition and prove the bounded-source estimate \eqref{lin:source-Gaussian}. Define
\begin{align}
\bar f(t):=\begin{cases}f_0+\glift(t)-\glift(0),&-t_*<t<0,\\ f(t),&0\leq t<T_1.\end{cases} \label{eq:I:display:306}
\end{align}
For $t<0$, zeroth-order compatibility gives $\trm\bar f(t)=g(0)+\trm\glift(t)-\trm\glift(0)=\trm\glift(t)$. Initial attainment in $L^2$ and \eqref{cb:setting} give $\gnorm{\infty}{\alpha_1}{f_0}\leq\delta$. For $t_*\leq1$, \eqref{cb:extended-lift-bounds} therefore gives
\begin{align}
\gnorm[\big]{\infty}{\alpha_1}{\bar f(t)-f_0}\leq C\abs{t},\qquad \gnorm[\big]{\infty}{\alpha_1}{\bar f(t)}\leq\delta+C t_*,\qquad -t_*<t<0. \label{cb:negative-amplitude}
\end{align}
Choose $t_*$ so that $Ct_*\leq\delta$. By \eqref{cb:initial-smallness}, $\bar f$ satisfies the Gaussian smallness required for the coefficient comparison at exponent $\alpha_1$.

We next compute the source. Since $\partial_t\bar f=\partial_t\glift$ on negative times and $a_0$ is given by \eqref{in:eq:compat1},
\begin{align}
\sourceR_-(t)&:=\partial_t\glift(t)+v\cdot\nabla_x\bar f(t)+\Lop\bar f(t)-\Ga(\bar f(t),\bar f(t)),\qquad t<0, \label{cb:inflow-residual-zero}\\
\sourceR_-(0^-)&=\partial_t\glift(0)-a_0. \label{cb:inflow-residual-zero-limit}
\end{align}
Set $\sourceR=\sourceR_-$ for $t<0$ and $\sourceR=0$ for $t\geq0$. Applying \eqref{cvel:scalar-graph-pointwise} to $\bar f$, the bounds on $\partial_t\glift$, $\nabla_x\bar f$ and $D_v^2\bar f$ give, after decreasing $\alpha_2$ if necessary,
\begin{align}
\tnorm{\sourceR}_{L^\infty(-t_*,T;L^{\infty,\alpha_2})}\lesssim_T 1. \label{cb:inflow-residual-bounded}
\end{align}
This verifies \eqref{lin:source-Gaussian}; Proposition~\ref{cb:inflow-gradient} allows the jump of $\sourceR$ at $t=0$.

\proofstep{Verify the weak equation and conclude the first-order estimate} We will verify \eqref{cb:regularity-class} and the weak equation on $(-t_*,T_1)$, then obtain part (i). The extension has matching $L^2$ values at $t=0$. When the equations on negative and positive times are integrated against a test function, the time-boundary terms cancel because
\begin{align}
\bigl(\bar f(0^-)-f_0,\varphi(0)\bigr)_{L^2}=0. \label{cb:initial-time-cancellation}
\end{align}
Consequently $\Y\bar f$ has no measure supported on $\{t=0\}$, and $\bar f$ satisfies the weak equation \eqref{cb:equation-residual} on $(-t_*,T_1)$ with the source just defined. Its negative-time derivatives have Gaussian bounds, so it also satisfies \eqref{cb:regularity-class}. The incoming bound follows from \eqref{cb:extended-incoming-polynomial}, and the source bound is \eqref{cb:inflow-residual-bounded}. Proposition~\ref{cb:inflow-gradient} therefore applies on an interval containing $[0,T]$ and gives \eqref{cb:initial-inflow-gradient} and the H\"older regularity in part (i).

\proofstep{Construct an extension whose source vanishes at zero} Assume also $\trm a_0=\partial_tg(0)$. To obtain the second-order estimate, we will modify the extension so that $\sourceR_-(0^-)=0$, then prove the H\"older bound \eqref{lin:source-holder}. Replace \eqref{eq:I:display:306} on negative times by
\begin{align}
\bar f(t):=f_0+ta_0+\glift(t)-\glift(0)-t\partial_t\glift(0)=f_0+ta_0+\frac{t^2}{2}\partial_t^2\glift(0). \label{eq:I:display:308}
\end{align}
The identities \eqref{in:eq:compat0}--\eqref{in:eq:compat1} give $\trm\bar f(t)=\trm\glift(t)$. The bounds \eqref{cb:initial-derivative-Gaussian} and \eqref{eq:I:display:294} again give \eqref{cb:negative-amplitude}, so the same smallness argument applies after reducing $t_*$. The identity $a_0=-v\cdot\nabla_xf_0-\Lop f_0+\Ga(f_0,f_0)$ in \eqref{in:eq:compat1} gives
\begin{align}
\sourceR_-(t)&:=a_0+t\partial_t^2\glift(0)+v\cdot\nabla_x\bar f(t)+\Lop\bar f(t)-\Ga(\bar f(t),\bar f(t)), \label{cb:inflow-residual-first}\\
\sourceR_-(0^-)&=a_0+v\cdot\nabla_xf_0+\Lop f_0-\Ga(f_0,f_0)=0. \label{cb:inflow-residual-first-limit}
\end{align}
\proofstep{Prove the H\"older source bound across zero} We will estimate the source and its first derivatives, and then its kinetic increments. The source in \eqref{cb:inflow-residual-first} and one derivative in $(t,x,v)$ use at most three derivatives of $a_0$. For the term $\partial_t^2\glift(0)$, they use at most three further derivatives, with at most two in position. All these derivatives are included in \eqref{eq:I:display:294} and follow from the assumed $C^5$ bound on $g$. After a further Gaussian exponent loss, the product and convolution formulas give
\begin{align}
\tnorm{\ue^{\alpha_2\abs{v}^2}\bigl(\abs{\sourceR_-}+\abs{\partial_t\sourceR_-}+\abs{\nabla_x\sourceR_-}+\abs{\nabla_v\sourceR_-}\bigr)}_{L^\infty((-t_*,0)\times\phase)}\lesssim_T 1. \label{cb:inflow-residual-Lipschitz}
\end{align}
Extend this source by zero to $t\geq0$. By \eqref{cb:inflow-residual-first-limit} and the time-derivative bound in \eqref{cb:inflow-residual-Lipschitz},
\begin{align}
\gnorm{\infty}{\alpha_2}{\sourceR_-(t)}\lesssim_T\abs{t},\qquad -t_*<t<0. \label{cb:inflow-residual-vanishing}
\end{align}
For points $Z,Z'$ at kinetic distance $r\leq1$ with $t,t'<0$, the derivative bound and \eqref{cb:initial-kinetic-increments} give an increment bounded by $C_T\ue^{-\alpha_3\abs{v_0}^2}(r+r^2+r^3)$, after absorbing the polynomial velocity factor in a smaller exponent $\alpha_3>0$. Local boundary neighborhoods and a finite covering give the same bound up to the spatial boundary. For $t<0\leq t'$, \eqref{cb:inflow-residual-vanishing} gives an increment at most $C_T\ue^{-\alpha_3\abs{v_0}^2}\abs{t}\leq C_T\ue^{-\alpha_3\abs{v_0}^2}r^2$. For $t,t'\geq0$, the increment is zero. Larger separations are controlled by \eqref{cb:inflow-residual-bounded}. Hence, for every fixed $0<\beta<1$,
\begin{align}
[\sourceR]_{\Ckin^\beta\left((-t_*,T)\times\overline\Om\times B_{\bracket{v_0}^{-1}}(v_0)\right)}\lesssim_{T,\beta}\ue^{-\alpha_3\abs{v_0}^2}. \label{cb:initial-residual-holder}
\end{align}
This proves \eqref{lin:source-holder} on an interval containing $[0,T]$.

\proofstep{Conclude the second-order estimate} We will apply the second-order argument of Proposition~\ref{cb:strong-upgrade} to $\bar f-\glift$. The cancellation \eqref{cb:initial-time-cancellation} still applies, so $\bar f$ satisfies the weak equation with precisely this source. The incoming bound \eqref{cb:extended-incoming-polynomial} permits application of Lemma~\ref{cb:inflow-coefficients} and Proposition~\ref{cb:inflow-gradient}; in particular, \eqref{cb:inflow-source-holder} follows from \eqref{cb:initial-residual-holder}. The lift bound \eqref{cb:extended-lift-bounds} is the order-three weak derivative condition used in the proof of Proposition~\ref{cb:strong-upgrade}. Applying that argument to $\bar f-\glift$, with \eqref{cb:inflow-residual-bounded} and \eqref{cb:initial-residual-holder}, gives \eqref{cb:initial-inflow-hessian} and the corresponding interior bound. The local derivative representatives given by that argument are continuous away from grazing on an interval containing zero. On negative times, their values tend to the classical velocity derivatives of $f_0$ and to $a_0+v\cdot\nabla_xf_0=\Ga(f_0,f_0)-\Lop f_0$. This identifies their values at $t=0$ and completes part (ii).
\end{proof}

\begin{lemma}[Diffuse-reflection estimates at the initial time]\label{cb:initial}
Let $f$ be a weak solution of \eqref{eq:normalized-problem} on $[0,T_1]$ in the sense of Definition~\ref{def:energy}, with diffuse-reflection and initial data $f_0$. Assume \eqref{cb:setting} and \eqref{cb:initial-smallness}. Suppose that $f_0\in C^5(\overline\Om\times\R^3)$ satisfies \eqref{eq:data-derivatives}, and let $a_0$ be defined by \eqref{eq:a0}. Assume both identities in \eqref{eq:compatibility}.

For every $0<T<T_1$, the following holds.
\begin{enumerate}
\item For some $\kappa>\frac{2}{3}$,
\begingroup\usetagform{enumerated}
\begin{align}
f\in\Ckinloc^{\frac{1}{2}}([0,T]\times\overline\Om\times\R^3),\qquad \flux{f}\in C^\kappa([0,T]\times\partial\Om). \label{cb:initial-diffuse-holder}
\end{align}
\endgroup
\item There are $\alpha_T>0$ and $0<C_T<\infty$ such that both \eqref{cb:initial-inflow-gradient} and \eqref{cb:initial-inflow-hessian} hold. On $\left\{d\geq\frac{d_0}{2}\right\}$, their denominators are replaced by $1$.
\end{enumerate}
The constants may depend on $T$ and the finite derivative bounds of $f_0$. Away from grazing, the velocity derivatives at $t=0$ are those of $f_0$, and $\Y f(0)=\Ga(f_0,f_0)-\Lop f_0$.
\end{lemma}

\begin{proof}
We will extend $f$ to negative times while preserving diffuse reflection. We will then prove that the additional source vanishes at $t=0$ and satisfies \eqref{lin:source-Gaussian} and \eqref{lin:source-holder} across zero. After verifying the weak equation on the extended interval, Propositions~\ref{cb:flux}, \ref{cb:gradient} and \ref{cb:graph} will give the flux, first-order and second-order conclusions on $[0,T]$.

\proofstep{Construct the extension and preserve the boundary law} We will construct $\bar f$ satisfying the Gaussian smallness and the diffuse law on $(-t_*,0)$. Recall from \eqref{eq:a0} that $a_0=-v\cdot\nabla_xf_0-\Lop f_0+\Ga(f_0,f_0)$. Formula~\eqref{eq:I:display:309} and its proof give \eqref{cb:initial-derivative-Gaussian}, with $\alpha_1=\frac{\alpha}{2}$. Set
\begin{align}
\bar f(t):=\begin{cases}f_0+ta_0,&-t_*<t<0,\\ f(t),&0\leq t<T_1.\end{cases} \label{cb:diffuse-negative-extension}
\end{align}
Initial attainment in $L^2$ and \eqref{cb:setting} give $\gnorm{\infty}{\alpha_1}{f_0}\leq\delta$. Choose $t_*\leq1$ so that $t_*\gnorm{\infty}{\alpha_1}{a_0}\leq\delta$. Then
\begin{align}
\gnorm[\big]{\infty}{\alpha_1}{\bar f(t)}\leq\delta+\abs{t}\gnorm{\infty}{\alpha_1}{a_0}\leq2\delta\leq\delta_A,\qquad -t_*<t<0. \label{cb:diffuse-negative-amplitude}
\end{align}
The linearity of $\Pm$ and the two compatibility identities \eqref{eq:compatibility} also give
\begin{align}
\trm\bar f(t)=\trm f_0+t\trm a_0=\Pm\bigl(\trp f_0+t\trp a_0\bigr)=\Pm\trp\bar f(t),\qquad -t_*<t<0. \label{cb:diffuse-negative-boundary}
\end{align}
Thus the extension has the required small amplitude and exactly preserves diffuse reflection.

\proofstep{Compute the source and its first derivatives} We will prove that the negative-time source is $O(\abs{t})$ in a Gaussian-weighted $L^\infty$ norm and has bounded first derivatives. Substituting \eqref{cb:diffuse-negative-extension} into \eqref{cb:equation-residual}, the constant term is $a_0+v\cdot\nabla_xf_0+\Lop f_0-\Ga(f_0,f_0)=0$ by \eqref{eq:a0}. The remaining terms are
\begin{align}
\sourceR(t)=t\Bigl[v\cdot\nabla_xa_0+\Lop a_0-\Ga(a_0,f_0)-\Ga(f_0,a_0)\Bigr]-t^2\Ga(a_0,a_0),\qquad t<0. \label{cb:initial-residual}
\end{align}
The bounds \eqref{eq:data-derivatives} and \eqref{cb:initial-derivative-Gaussian} control all derivatives of $f_0$ and $a_0$ needed here. In \eqref{cvel:Q}, place derivatives of a convolution on its Gaussian factor and use \eqref{cvel:Gaussian-convolution-bound}; at $\gamma=-3$, use $c*G=-8\pi G$. After division by $\mh$, formula~\eqref{eq:I:display:071} expresses the remaining derivatives as Gaussian factors times polynomials in $v$. Absorbing these polynomials in a smaller exponent $0<\alpha_2<\alpha_1$ gives
\begin{align}
&\sum_{\abs{k}\leq1}\gnorm{\infty}{\alpha_2}{\partial_{x,v}^k(v\cdot\nabla_xa_0+\Lop a_0)}+\sum_{G,H\in\{f_0,a_0\}}\sum_{\abs{k}\leq1}\gnorm{\infty}{\alpha_2}{\partial_{x,v}^k\Ga(G,H)}\lesssim_T1. \label{cb:diffuse-residual-coefficients}
\end{align}
The largest derivative order is three on $a_0$, which uses at most five derivatives of $f_0$. Applying \eqref{cb:diffuse-residual-coefficients} to the explicit powers of $t$ in \eqref{cb:initial-residual}, with $t_*\leq1$, yields
\begin{align}
\gnorm{\infty}{\alpha_2}{\sourceR(t)}\lesssim_T\abs{t},\qquad \gnorm{\infty}{\alpha_2}{\abs{\partial_t\sourceR(t)}+\abs{\nabla_x\sourceR(t)}+\abs{\nabla_v\sourceR(t)}}\lesssim_T1,\qquad -t_*<t<0. \label{cb:diffuse-residual-derivatives}
\end{align}
In particular, the source tends to zero in this weighted norm as $t\uparrow0$.

\proofstep{Verify the source bounds across zero} Set $\sourceR=0$ for $t\geq0$. Equation~\eqref{cb:diffuse-residual-derivatives} gives \eqref{lin:source-Gaussian}. It also gives exactly the two bounds used in the proof of \eqref{cb:initial-residual-holder}: Gaussian bounds for the first derivatives on $t<0$ and $\gnorm{\infty}{\alpha_2}{\sourceR(t)}\lesssim_T\abs{t}$. For kinetic separation $r\leq1$, \eqref{cb:initial-kinetic-increments} consequently bounds source increments by $C_T\ue^{-\alpha_3\abs{v_0}^2}r$ when both times are negative and by $C_T\ue^{-\alpha_3\abs{v_0}^2}r^2$ when they straddle zero, for some $0<\alpha_3<\alpha_2$. The amplitude controls larger separations. Thus, for every fixed $0<\beta<1$,
\begin{align}
[\sourceR]_{\Ckin^\beta((-t_*,T)\times\overline\Om\times B_{\bracket{v_0}^{-1}}(v_0))}\lesssim_{T,\beta}\ue^{-\alpha_3\abs{v_0}^2}. \label{cb:diffuse-initial-source-holder}
\end{align}
This verifies \eqref{lin:source-holder} across zero.

\proofstep{Verify the weak equation through zero} The Gaussian negative-time derivative bounds give $\bar f\in L^2_{\mathrm{loc}}((-t_*,T_1);\Vspace)$. Initial attainment and \eqref{cb:diffuse-negative-extension} give $\bar f(0^-)=f(0^+)=f_0$ in $L^2$, so \eqref{cb:regularity-class} holds across zero. The time-boundary terms cancel by \eqref{cb:initial-time-cancellation}. Adding the two weak equations therefore gives \eqref{cb:equation-residual} with the source \eqref{cb:initial-residual}, extended by zero, and no additional term at $t=0$. Together with \eqref{cb:diffuse-negative-boundary}, this verifies the extended equation and diffuse law.

\proofstep{Obtain the estimates on $[0,T]$} The bounds \eqref{cb:diffuse-negative-amplitude}, \eqref{cb:diffuse-residual-derivatives} and \eqref{cb:diffuse-initial-source-holder}, and the preceding weak equation, verify the hypotheses of Propositions~\ref{cb:flux}, \ref{cb:gradient} and \ref{cb:graph} on an interval containing $[0,T]$. Proposition~\ref{cb:flux} gives the flux bound \eqref{rev14:diffuse-flux-regularity}. Proposition~\ref{cb:gradient} gives \eqref{cb:diffuse-sharp-holder} and \eqref{cb:pointwise-estimates}--\eqref{cb:diffuse-interior-gradient}. Restricting these estimates to $[0,T]$ proves \eqref{cb:initial-diffuse-holder} and \eqref{cb:initial-inflow-gradient}. Proposition~\ref{cb:graph} gives \eqref{cb:pointwise-hessian}--\eqref{cb:diffuse-interior-hessian}, hence \eqref{cb:initial-inflow-hessian} and its interior version. We decrease the Gaussian exponents to one common $\alpha_T>0$.

The derivative representatives are continuous away from grazing on the enlarged time interval. From the negative-time formula, $\nabla_v\bar f(0)=\nabla_vf_0$, $D_v^2\bar f(0)=D_v^2f_0$, and $\Y\bar f(0)=a_0+v\cdot\nabla_xf_0=\Ga(f_0,f_0)-\Lop f_0$. This proves the stated values at $t=0$.
\end{proof}

\subsection{Regularity and Strong Solution}\label{lin:section-integration}

We will collect the regularity obtained in Sections~\ref{cb:inflow-section}--\ref{cb:section-initial} and prove that the global weak Landau solution is strong. The first- and second-order estimates hold through $t=0$ by the preceding lemmas. We will integrate their grazing factors to obtain weighted $L^2$ derivative bounds, then apply Lemma~\ref{cvel:weighted-collisions} to bound $\Lop f$ and $\Ga(f,f)$ separately. These estimates will verify the integrability condition \eqref{eq:I:display:016}.

\begin{proposition}[Regularity and strong solution]\label{cb:integration}
Let $f$ be a global weak solution of \eqref{eq:normalized-problem} in the sense of Definition~\ref{def:energy}, furnished by Proposition~\ref{ci:global-energy} in the in-flow case or Proposition~\ref{cex:CM-input} in the diffuse-reflection case. Suppose that $f$ satisfies \eqref{cb:setting} and \eqref{cb:initial-smallness}. Assume the data satisfy Assumption~\ref{in:data:regularity} in the in-flow case, or Assumption~\ref{def:data} in the diffuse-reflection case. The following holds.
\begin{enumerate}
\item We have
\begingroup\usetagform{enumerated}
\begin{align}
f\in\Ckinloc^{\frac12}([0,\infty)\times\overline\Om\times\R^3). \label{cb:integration-holder}
\end{align}
\endgroup
For every $T>0$, there are $\alpha_T\in(0,\alpha)$ and $C_T<\infty$ such that \eqref{cb:initial-inflow-gradient} and \eqref{cb:initial-inflow-hessian} hold, and
\begingroup\usetagform{enumerated}
\begin{align}
\tnorm{\ue^{\alpha_T\abs{v}^2}\bigl(\abs{\nabla_vf}+\abs{D_v^2f}+\abs{\Y f}\bigr)}_{L^\infty([0,T]\times\{d\geq\frac{d_0}{2}\}\times\R^3)}\leq C_T. \label{cb:integration-interior}
\end{align}
\endgroup
For diffuse-reflection, there is also $\kappa_T\in(\frac23,1)$ such that
\begingroup\usetagform{enumerated}
\begin{align}
\flux{f}\in C^{\kappa_T}([0,T]\times\partial\Om),\qquad\norm{\flux{f}}_{C^{\kappa_T}([0,T]\times\partial\Om)}\leq C_T. \label{cb:integration-flux}
\end{align}
\endgroup
\item For every $T>0$ and $m\geq0$,
\begingroup\usetagform{enumerated}
\begin{align}
\la^mf,\la^m\nabla_vf,\la^mD_v^2f,\la^m\Y f,\la^m\Lop f,\la^m\Ga(f,f)\in L^2((0,T)\times\Om\times\R^3), \label{cb:strong-graph}
\end{align}
\endgroup
and
\begingroup\usetagform{enumerated}
\begin{align}
&\sum_{\abs{k}\leq2}\tnorm[\Big]{\la^m\partial_v^kf}_{L^2((0,T)\times\phase)}^2+\tnorm[\Big]{\la^m\Y f}_{L^2((0,T)\times\phase)}^2 \label{cb:graph-L2}\\
&\qquad+\tnorm[\Big]{\la^m\Lop f}_{L^2((0,T)\times\phase)}^2+\tnorm[\big]{\la^m\Ga(f,f)}_{L^2((0,T)\times\phase)}^2\lesssim_{T,m}T.\notag
\end{align}
\endgroup
The equation holds as
\begingroup\usetagform{enumerated}
\begin{align}
\Y f+\Lop f=\Ga(f,f)\qquad\text{in }L^2((0,T)\times\Om\times\R^3). \label{cb:strong-equation}
\end{align}
\endgroup
In particular, $f$ is a strong solution on every $[0,T]$, and hence on $[0,\infty)$.
\end{enumerate}
The constants and exponents may depend on the finite derivative bounds of the data on $[0,T+1]$.
\end{proposition}

\begin{proof}
Fix $T>0$ and $m\geq0$. The first step will prove part~(i). For part~(ii), we will integrate the derivative bounds in $x$ and $v$, estimate the collision terms separately, and obtain the equation in $L^2$.

\proofstep{Obtain the regularity estimates on $[0,T]$} Take $T_1=T+1$. The assumed smallness and the data assumptions verify all hypotheses of Lemma~\ref{cb:inflow-initial}, including both compatibility identities, in the in-flow case. That lemma gives the H\"older bound and \eqref{cb:initial-inflow-gradient}--\eqref{cb:initial-inflow-hessian}, together with the corresponding interior bounds. For diffuse-reflection, Lemma~\ref{cb:initial} gives these conclusions and the flux bound. Taking a common smaller Gaussian exponent gives \eqref{cb:integration-interior}, and decreasing the flux exponent if needed gives \eqref{cb:integration-flux}. Since $T$ is arbitrary, the H\"older conclusions give \eqref{cb:integration-holder}. This proves part~(i).

\proofstep{Prove the weighted $L^2$ derivative bounds} For part~(ii), we will show \eqref{cb:integrated-second-derivatives}--\eqref{cb:integrated-first-derivatives}. For $d(x)<\frac{d_0}{2}$, rotate velocity so that $s=n(\bdproj(x))\cdot v$. This preserves the Gaussian and velocity measure. By \eqref{eq:rho}, the singular normal integral is
\begin{align}
\int_\R\bigl(d^{\frac13}+\abs{s}\bigr)^{-3}\dd s=d^{-\frac23},\qquad d>0. \label{eq:I:display:304}
\end{align}
The bounded Jacobian of the boundary collar and
\begin{align}
\int_0^{\frac{d_0}{2}}\int_\R\bigl(d^{\frac13}+\abs{s}\bigr)^{-3}\dd s\dd d=\int_0^{\frac{d_0}{2}}d^{-\frac23}\dd d=3\left(\frac{d_0}{2}\right)^{\frac13}<\infty \label{eq:I:display:285}
\end{align}
therefore control the singularity in the square of \eqref{cb:initial-inflow-hessian}. More precisely, first use $\la^{2m}\ue^{-2\alpha_T\abs{v}^2}\lesssim_{m,\alpha_T}\ue^{-\alpha_T\abs{v}^2}$ to absorb the polynomial weight, then integrate the tangential Gaussian. For the factor $\graz^{-3}$, discard the normal Gaussian and apply \eqref{eq:I:display:304}--\eqref{eq:I:display:285}; the remaining term is a Gaussian integral. On $\{d\geq\frac{d_0}{2}\}$, use \eqref{cb:integration-interior}. Integrating also in time gives
\begin{align}
\tnorm[\Big]{\la^mD_v^2f}_{L^2((0,T)\times\phase)}^2+\tnorm[\Big]{\la^m\Y f}_{L^2((0,T)\times\phase)}^2\lesssim_{T,m}T. \label{cb:integrated-second-derivatives}
\end{align}
The first-order estimate \eqref{cb:initial-inflow-gradient} has squared grazing factor bounded by $C(1+\graz^{-1})\leq C(2+\graz^{-3})$. The same integration, together with the Gaussian amplitude \eqref{cb:setting}, consequently gives
\begin{align}
\tnorm[\Big]{\la^mf}_{L^2((0,T)\times\phase)}^2+\tnorm[\Big]{\la^m\nabla_vf}_{L^2((0,T)\times\phase)}^2\lesssim_{T,m}T. \label{cb:integrated-first-derivatives}
\end{align}
\proofstep{Prove the separate collision bounds} We will estimate $\la^m\Lop f$ and $\la^m\Ga(f,f)$ in $L^2((0,T)\times\phase)$. The Gaussian smallness \eqref{cb:setting} verifies \eqref{cvel:weighted-Gaussian}, and \eqref{cb:integrated-second-derivatives}--\eqref{cb:integrated-first-derivatives}, for every weight, verify \eqref{cvel:scalar-graph-assumption}. Lemma~\ref{cvel:weighted-collisions} therefore identifies the collision distributions with $L^2$ functions. Its estimate \eqref{cvel:scalar-graph-pointwise} gives
\begin{align}
\la^m\bigl(\abs{\Lop f}+\abs{\Ga(f,f)}\bigr)\lesssim\la^{m+5}\bigl(\abs{D_v^2f}+\abs{\nabla_vf}+\abs{f}\bigr)+\delta\la^{m+5}\ue^{-\frac14\abs{v}^2}. \label{cb:integration-collision-bound}
\end{align}
The first three terms have the preceding $L^2$ bounds with weight $m+5$, and the last term has $L^2$ norm bounded by $C_m\sqrt{T}$. Hence
\begin{align}
\tnorm[\Big]{\la^m\Lop f}_{L^2((0,T)\times\phase)}+\tnorm[\Big]{\la^m\Ga(f,f)}_{L^2((0,T)\times\phase)}\lesssim_{T,m}\sqrt{T}. \label{cb:integrated-collisions}
\end{align}
Together with \eqref{cb:integrated-second-derivatives}--\eqref{cb:integrated-first-derivatives}, this proves \eqref{cb:strong-graph} and \eqref{cb:graph-L2} throughout $-3\leq\gamma\leq1$.

\proofstep{Conclude the strong equation and global result} All terms of the distributional equation \eqref{eq:intro-interior-weak} now belong to $L^2((0,T)\times\phase)$. Their distributional equality is therefore the $L^2$ equality \eqref{cb:strong-equation}. Taking $m=0$ in \eqref{cb:graph-L2} verifies the integrability condition \eqref{eq:I:display:016}, so Definition~\ref{def:strong} applies on $[0,T]$. Since $T>0$ was arbitrary, $f$ is a global strong solution with the stated regularity. This proves part~(ii).
\end{proof}


\section{Continuous Dependence and Uniqueness}\label{sec:uniqueness}

We will first prove quantitative continuous dependence and then deduce uniqueness for weak solutions in the small Gaussian class. We will use the first-order bound \eqref{cb:integration-gradient} for the reference solution, supplied by Proposition~\ref{cb:integration}. This bound is additional to the $L^2$ conditions in Definition~\ref{def:strong}.

\subsection{Logarithmic Estimates}\label{sec:uniq-velocity-log}

We write $\logconst:=2\operatorname{diam}(\Om)$, so $\frac{\logconst}{d(x)}\geq2$ and the logarithmic weight in \eqref{eq:intro-log-bound} is positive.
\begin{lemma}\label{uniq:velocity-log}
Let $f$ be a weak solution of \eqref{eq:normalized-problem} on $[0,T]$ in the sense of Definition~\ref{def:energy}. Suppose that, for some $\alpha_T>0$ and $C_T>0$,
\begin{align}
&\tnorm{f}_{L^\infty(0,T;L^{\infty,\alpha_T})}+\tnorm{\frac{\ue^{\alpha_T\abs{v}^2}\nabla_vf}{1+\graz^{-\frac{1}{2}}}}_{L^\infty\left((0,T)\times\{d<\frac{d_0}{2}\}\times\R^3\right)}+\tnorm[\Big]{\ue^{\alpha_T\abs{v}^2}\nabla_vf}_{L^\infty\left((0,T)\times\{d\geq\frac{d_0}{2}\}\times\R^3\right)}\leq C_T. \label{cb:integration-gradient}
\end{align}
Then we have, for almost every $(t,x)\in (0,T)\times\Om$,
\begin{align}
\vnorm[\big]{f(t,x,\cdot)}_{\Vv}^2\lesssim_{\alpha_T}C_T^2\left(1+\ln\left(\frac{\logconst}{d(x)}\right)\right). \label{cb:strong-log}
\end{align}
\end{lemma}
\begin{proof}
By \eqref{cvel:rough-coefficients} with $F=\m$, we have $\abs{\sig(v)}\lesssim\la^{\gamma+2}$. On the collar $d(x)<\frac{d_0}{2}$, the definition \eqref{eq:V}, the assumed bound \eqref{cb:integration-gradient}, and $\left(1+\graz^{-\frac12}\right)^2\leq2\left(1+\graz^{-1}\right)$ give
\begin{align}
\vnorm[\big]{f(t,x,\cdot)}_{\Vv}^2&\lesssim C_T^2\int_{\R^3}\la^{\gamma+2}\ue^{-2\alpha_T\abs{v}^2}\left(1+\graz^{-1}\right)\dd v \label{uniq:velocity-Gaussian-reduction}\\
&\lesssim_{\alpha_T}C_T^2\int_{\R^3}\ue^{-\frac{\alpha_T}{2}\abs{v}^2}\left(1+\graz^{-1}\right)\dd v.\notag
\end{align}
The term without $\graz^{-1}$ is a Gaussian integral and also controls the zeroth-order part of \eqref{eq:V}. To estimate the remaining term, rotate velocity so that $s=n(\bdproj(x))\cdot v$. By \eqref{eq:rho}, integration in the two tangential velocity variables gives
\begin{align}
\int_{\R^3}\frac{\ue^{-\frac{\alpha_T}{2}\abs{v}^2}}{\graz(x,v)}\dd v=\frac{2\pi}{\alpha_T}\int_\R\frac{\ue^{-\frac{\alpha_T}{2}s^2}}{\bigl(d(x)\bigr)^{\frac13}+\abs{s}}\dd s. \label{uniq:normal-velocity-integral}
\end{align}
Splitting at $\abs{s}=1$ and using $1\leq s$ in the Gaussian tail, we obtain
\begin{align}
\int_\R\frac{\ue^{-\frac{\alpha_T}{2}s^2}}{\bigl(d(x)\bigr)^{\frac13}+\abs{s}}\dd s&\leq2\int_0^1\frac{\dd s}{\bigl(d(x)\bigr)^{\frac13}+s}+2\int_1^\infty\ue^{-\frac{\alpha_T}{2}s^2}\dd s \label{eq:I:display:284}\\
&\leq2\ln\left(1+\bigl(d(x)\bigr)^{-\frac13}\right)+\frac{2}{\alpha_T}\ue^{-\frac{\alpha_T}{2}}\lesssim_{\alpha_T}1+\ln\left(\frac{\logconst}{d(x)}\right).\notag
\end{align}
Combining \eqref{uniq:velocity-Gaussian-reduction}--\eqref{eq:I:display:284} proves \eqref{cb:strong-log} on the collar. On $d(x)\geq\frac{d_0}{2}$, the last term in \eqref{cb:integration-gradient} gives the same bound with a uniform Gaussian integral.
\end{proof}

In the difference argument, \eqref{eq:I:display:319} leaves the integral $\int_\Om\vnorm{f(x,\cdot)}_{\Vv}^2\vnorm{h(x,\cdot)}_{L_v^2}^2\dd x$. Lemma~\ref{uniq:velocity-log} bounds the first factor by the spatial logarithm. The following lemma estimates the resulting integral in terms of the $L^2$ norm of $h$.
\begin{lemma}\label{uniq:logarithm}
For every $p\geq2$, we have
\begin{align}
\norm{1+\ln\left(\frac{\logconst}{d(x)}\right)}_{L^p(\Om)}\lesssim p. \label{uniq:logLp}
\end{align}
For $h\in L_x^\infty L_v^2$ and $p\geq2$, we have
\begin{align}
\int_\Om\left(1+\ln\left(\frac{\logconst}{d(x)}\right)\right)\vnorm[\big]{h(x,\cdot)}_{L_v^2}^2\dd x\lesssim p\norm{h}_{L_x^\infty L_v^2}^{\frac2p}\norm{h}_{L^2}^{2-\frac2p}. \label{uniq:multiplier}
\end{align}
\end{lemma}
\begin{proof}
We first have
\begin{align}
\abs{\{x\in\Om:d(x)<r\}}\lesssim r,\qquad r>0. \label{uniq:collar-volume}
\end{align}
For $0<r<d_0$, this follows by integrating the bounded Jacobian in the boundary collar. For $r\geq d_0$, boundedness of $\Om$ gives the same estimate. Consequently, for $s\geq1$,
\begin{align}
\abs{\left\{1+\ln\left(\frac{\logconst}{d(x)}\right)>s\right\}}=\abs{\left\{d(x)<\logconst\ue^{1-s}\right\}}\lesssim\ue^{-s}. \label{uniq:logarithm-distribution}
\end{align}
Integration over the levels $s$, split at $s=1$, yields
\begin{align}
\int_\Om\left(1+\ln\left(\frac{\logconst}{d(x)}\right)\right)^p\dd x&=p\int_0^\infty s^{p-1}\abs{\left\{1+\ln\left(\frac{\logconst}{d(x)}\right)>s\right\}}\dd s \label{eq:I:display:312}\\
&\leq C+Cp\int_1^\infty s^{p-1}\ue^{-s}\dd s\leq C+C\int_0^\infty s^p\ue^{-s}\dd s\notag\\
&\leq C+2C\left(\frac{2p}{\ue}\right)^p\leq(Cp)^p.\notag
\end{align}
The last two inequalities use integration by parts and $s^p\ue^{-\frac{s}{2}}\leq\left(\frac{2p}{\ue}\right)^p$, with $\int_0^\infty\ue^{-\frac{s}{2}}\dd s=2$. Taking the $p$th root proves \eqref{uniq:logLp} for every real $p\geq2$.

For almost every $x$, we have
\begin{align}
\vnorm[\big]{h(x,\cdot)}_{L_v^2}^{\frac{2p}{p-1}}\leq\norm{h}_{L_x^\infty L_v^2}^{\frac{2}{p-1}}\vnorm[\big]{h(x,\cdot)}_{L_v^2}^{2}. \label{uniq:interpolation-power}
\end{align}
H\"older's inequality with exponents $p$ and $\frac{p}{p-1}$ then gives
\begin{align}
\int_\Om\left(1+\ln\left(\frac{\logconst}{d(x)}\right)\right)\vnorm[\big]{h(x,\cdot)}_{L_v^2}^2\dd x&\leq\norm{1+\ln\left(\frac{\logconst}{d(x)}\right)}_{L^p(\Om)}\left(\int_\Om\vnorm[\big]{h(x,\cdot)}_{L_v^2}^{\frac{2p}{p-1}}\dd x\right)^{\frac{p-1}{p}} \label{uniq:logarithm-interpolation}\\
&\lesssim p\norm{h}_{L_x^\infty L_v^2}^{\frac2p}\norm{h}_{L^2}^{2-\frac2p}.\notag
\end{align}
This proves \eqref{uniq:multiplier}.
\end{proof}

\subsection{In-flow: Continuous Dependence and Uniqueness}\label{sec:uniq-osgood}

We will estimate the difference of two solutions in terms of their initial and incoming data. Logarithmic estimates for differences of solutions also appear in \cite[Section~3.2]{Loeper2006} and \cite[Lemma~8 and Theorem~9]{Fournier2010}.

\begin{proposition}\label{uniq:osgood}
There is $\delta_*>0$ such that the following holds. Let $f$ and $\widetilde f$ be weak solutions of \eqref{eq:normalized-problem} on $[0,T]$, with in-flow data $g$ and $\widetilde g$, respectively. Assume that $f$ satisfies \eqref{cb:integration-gradient} and that
\begin{align}
\tnorm[\big]{\widetilde f}_{L^\infty(0,T;L^{\infty,\alpha})}\leq\delta_*. \label{uniq:smallness}
\end{align}
Then we have, for some $C_T\geq1$ and $\theta_T\in(0,1)$,
\begin{align}
\tnorm[\big]{\widetilde f-f}_{C([0,T];L^2)}^2\leq C_T\left(\norm[\big]{\widetilde f(0)-f(0)}_{L^2}^2+\tnorm[\big]{\widetilde g-g}_{L^2((0,T)\times\Sigma_-)}^2\right)^{\theta_T}. \label{uniq:modulus}
\end{align}
The constants depend additionally on $T$ and the bounds in \eqref{cb:integration-gradient}. In particular, equal initial and incoming data imply $\widetilde f=f$.
\end{proposition}
\begin{proof}
Set $h:=\widetilde f-f$. We will derive its energy inequality, apply the logarithmic estimates and integrate the resulting differential inequality.

\proofstep{The difference energy inequality} We will prove \eqref{uniq:energy}, including the dissipation needed for convergence of derivatives and traces. Subtracting the weak equations gives
\begin{align}
\Y h+\Lop h=\Ga(\widetilde f,h)+\Ga(h,f),\qquad \trm h=\widetilde g-g. \label{trace:difference-equation}
\end{align}
Both solutions belong to $C([0,T];L^2)\cap L^2(0,T;\Vspace)$. Their Gaussian bounds imply $h\in L^\infty(0,T;L_x^\infty L_v^2)$, so \eqref{cvel:spatial-L}, \eqref{cvel:spatial-dual} and \eqref{cvel:spatial-reversed} give all terms of \eqref{trace:difference-equation} in $L^2(0,T;\Vspace')$. In particular,
\begin{align}
\tnorm{\Ga(h,f)}_{L^2(0,T;\Vspace')}\lesssim\tnorm{h}_{L^\infty(0,T;L_x^\infty L_v^2)}\tnorm{f}_{L^2(0,T;\Vspace)}. \label{uniq:reversed-source}
\end{align}
Lemmas~\ref{trace:CM}(ii) and \ref{trace:gaussian} give the $L^2$ traces. These facts verify the hypotheses of Proposition~\ref{trace:green}, which yields, for $0\leq s\leq t\leq T$,
\begin{align}
&\frac12\norm{h(t)}_{L^2}^2+\frac12\int_s^t\norm{\trp h(r)}_{L^2_{\Sigma_+}}^2\dd r+\int_s^t\ip{\Lop h}{h}\dd r \label{trace:difference-inflow-energy}\\
&\qquad=\frac12\norm{h(s)}_{L^2}^2+\frac12\int_s^t\norm[\big]{\widetilde g(r)-g(r)}_{L^2_{\Sigma_-}}^2\dd r+\int_s^t\ip{\Ga(\widetilde f,h)+\Ga(h,f)}{h}\dd r.\notag
\end{align}
Every integrand belongs to $L^1(0,T)$, so $t\mapsto\norm{h(t)}_{L^2}^2$ is absolutely continuous.

Integrating \eqref{cvel:localized-garding} in $x$, applying \eqref{cvel:spatial-dual}, and using $\norm[\big]{\widetilde f}_{L_x^\infty L_v^2}\lesssim\delta_*$ from \eqref{uniq:smallness}, we obtain, for sufficiently small $\delta_*$,
\begin{align}
\ip{\Lop h-\Ga(\widetilde f,h)}{h}\geq(c_0-C\delta_*)\norm{h}_{\Vspace}^2-C\norm{h}_{L^2}^2\geq\frac{c_0}{2}\norm{h}_{\Vspace}^2-C\norm{h}_{L^2}^2. \label{cvel:principal}
\end{align}
For the remaining term, integrate \eqref{cvel:trilinear} in $x$ and apply Young's inequality to obtain
\begin{align}
\abs{\ip{\Ga(h,f)}{h}}\leq\frac{c_0}{4}\norm{h}_{\Vspace}^2+C\int_\Om\vnorm{f(t,x,\cdot)}_{\Vv}^2\vnorm{h(t,x,\cdot)}_{L_v^2}^2\dd x. \label{eq:I:display:319}
\end{align}
Substitution into \eqref{trace:difference-inflow-energy} gives, for almost every $t$,
\begin{align}
&\frac{\dd}{\dd t}\norm{h(t)}_{L^2}^2+\frac{c_0}{2}\norm{h(t)}_{\Vspace}^2+\norm{\trp h(t)}_{L^2_{\Sigma_+}}^2 \label{uniq:energy}\\
&\qquad\leq C\norm{h(t)}_{L^2}^2+C\int_\Om\vnorm{f(t,x,\cdot)}_{\Vv}^2\vnorm{h(t,x,\cdot)}_{L_v^2}^2\dd x+\norm[\big]{\widetilde g(t)-g(t)}_{L^2_{\Sigma_-}}^2.\notag
\end{align}

\proofstep{Estimate the remaining collision integral} We will recover a full factor $\norm{h(t)}_{L^2}^2$ in the right-hand side of \eqref{uniq:energy}, with a logarithmic factor that can be integrated in the next step. Applying \eqref{cb:strong-log} and then \eqref{uniq:multiplier} to the integral in \eqref{eq:I:display:319}, we obtain, for every real $p\geq2$,
\begin{align}
\int_\Om\vnorm{f(t,x,\cdot)}_{\Vv}^2\vnorm{h(t,x,\cdot)}_{L_v^2}^2\dd x\lesssim_T p\norm{h(t)}_{L_x^\infty L_v^2}^{\frac{2}{p}}\norm{h(t)}_{L^2}^{2-\frac{2}{p}}. \label{uniq:optimization}
\end{align}
These inequalities hold for all $p\geq2$ at every time outside one null set. The Gaussian bounds, Lemma~\ref{trace:gaussian} and continuity into $L^2$ allow us to fix $C_T\geq1$ sufficiently large that
\begin{align}
&\tnorm{h}_{L^\infty(0,T;L_x^\infty L_v^2)}^2\leq C_T,\qquad \sup_{0\leq t\leq T}\norm{h(t)}_{L^2}^2+\tnorm[\big]{\widetilde g-g}_{L^2((0,T)\times\Sigma_-)}^2\leq\frac{C_T}{\ue^2}. \label{uniq:uniform-bounds}
\end{align}
Choose this constant also larger than the coefficients needed below, and keep it fixed. For $\norm{h(t)}_{L^2}>0$, factoring out its square gives
\begin{align}
p\norm{h(t)}_{L_x^\infty L_v^2}^{\frac{2}{p}}\norm{h(t)}_{L^2}^{2-\frac{2}{p}}\leq p\norm{h(t)}_{L^2}^2\left(\frac{C_T}{\norm{h(t)}_{L^2}^2}\right)^{\frac1p}. \label{uniq:factored-product}
\end{align}
We choose $p$ so that the last power is constant:
\begin{align}
p=\ln\left(\frac{C_T}{\norm{h(t)}_{L^2}^2}\right)\geq2,\qquad \left(\frac{C_T}{\norm{h(t)}_{L^2}^2}\right)^{\frac1p}=\ue. \label{uniq:exponent-choice}
\end{align}
This value minimizes $p\left(\frac{C_T}{\norm{h(t)}_{L^2}^2}\right)^{\frac1p}$ over $p\geq2$. Substitution into \eqref{uniq:factored-product} yields
\begin{align}
p\norm{h(t)}_{L_x^\infty L_v^2}^{\frac{2}{p}}\norm{h(t)}_{L^2}^{2-\frac{2}{p}}\leq\ue\norm{h(t)}_{L^2}^2\ln\left(\frac{C_T}{\norm{h(t)}_{L^2}^2}\right). \label{uniq:optimized-product}
\end{align}
Combining \eqref{uniq:energy}--\eqref{uniq:optimized-product}, and using the lower bound for the logarithm in \eqref{uniq:exponent-choice} to absorb the term $C\norm{h(t)}_{L^2}^2$, gives
\begin{align}
\frac{\dd}{\dd t}\norm{h(t)}_{L^2}^2\leq C_T\norm{h(t)}_{L^2}^2\ln\left(\frac{C_T}{\norm{h(t)}_{L^2}^2}\right)+\norm[\big]{\widetilde g(t)-g(t)}_{L^2_{\Sigma_-}}^2. \label{uniq:differential}
\end{align}
The logarithmic term has continuous value zero when $\norm{h(t)}_{L^2}=0$. We will also use \eqref{uniq:optimization}--\eqref{uniq:optimized-product} with the dissipative terms in \eqref{uniq:energy} retained.

\proofstep{Integrate the inequality and conclude continuous dependence} Notice that
\begin{align}
\frac{\dd}{\dd t}\tnorm[\big]{\widetilde g-g}_{L^2((t,T)\times\Sigma_-)}^2=\frac{\dd}{\dd t}\int_t^T\norm[\big]{\widetilde g(r)-g(r)}_{L^2_{\Sigma_-}}^2\dd r=-\norm[\big]{\widetilde g(t)-g(t)}_{L^2_{\Sigma_-}}^2. \label{uniq:incoming-cancellation}
\end{align}
Also, for all sufficiently small $\eta$, \eqref{uniq:uniform-bounds} gives
\begin{align}
0<\norm{h(t)}_{L^2}^2+\tnorm[\big]{\widetilde g-g}_{L^2((t,T)\times\Sigma_-)}^2+\eta<\frac{C_T}{\ue},\qquad 0\leq t\leq T. \label{uniq:forced-range}
\end{align}
The function $s\mapsto s\ln\left(\frac{C_T}{s}\right)$ is increasing on this range, because its derivative is $\ln\left(\frac{C_T}{s}\right)-1>0$. Consequently, \eqref{uniq:differential} and \eqref{uniq:incoming-cancellation} imply
\begin{align}
&\frac{\dd}{\dd t}\left(\norm{h(t)}_{L^2}^2+\tnorm[\big]{\widetilde g-g}_{L^2((t,T)\times\Sigma_-)}^2+\eta\right) \label{uniq:forced-differential}\\
&\quad\leq C_T\norm{h(t)}_{L^2}^2\ln\left(\frac{C_T}{\norm{h(t)}_{L^2}^2}\right)\notag\\
&\quad\leq C_T\left(\norm{h(t)}_{L^2}^2+\tnorm[\big]{\widetilde g-g}_{L^2((t,T)\times\Sigma_-)}^2+\eta\right)\ln\left(\frac{C_T}{\norm{h(t)}_{L^2}^2+\tnorm[\big]{\widetilde g-g}_{L^2((t,T)\times\Sigma_-)}^2+\eta}\right).\notag
\end{align}
All quantities differentiated here are absolutely continuous. Dividing \eqref{uniq:forced-differential} by its positive parenthesis and changing sign gives a lower bound for the derivative of the logarithm. Multiplication by the integrating factor $\ue^{C_Tt}$ yields
\begin{align}
\frac{\dd}{\dd t}\left[\ue^{C_Tt}\ln\left(\frac{C_T}{\norm{h(t)}_{L^2}^2+\tnorm[\big]{\widetilde g-g}_{L^2((t,T)\times\Sigma_-)}^2+\eta}\right)\right]\geq0. \label{uniq:forced-comparison}
\end{align}
Integration from $0$ to $t$ therefore gives
\begin{align}
\ln\left(\frac{C_T}{\norm{h(t)}_{L^2}^2+\tnorm[\big]{\widetilde g-g}_{L^2((t,T)\times\Sigma_-)}^2+\eta}\right)\geq\ue^{-C_Tt}\ln\left(\frac{C_T}{\norm{h(0)}_{L^2}^2+\tnorm[\big]{\widetilde g-g}_{L^2((0,T)\times\Sigma_-)}^2+\eta}\right). \label{uniq:integrated-comparison}
\end{align}
Exponentiating and taking reciprocals bounds the denominator on the left from above. Since that denominator is at least $\norm{h(t)}_{L^2}^2$, we obtain
\begin{align}
\norm{h(t)}_{L^2}^2\leq C_T\left(\frac{\norm{h(0)}_{L^2}^2+\tnorm[\big]{\widetilde g-g}_{L^2((0,T)\times\Sigma_-)}^2+\eta}{C_T}\right)^{\exp(-C_Tt)}. \label{uniq:integrated-logarithm}
\end{align}
Let $\eta\downarrow0$. The fraction in \eqref{uniq:integrated-logarithm} then lies in $[0,\ue^{-2}]$, and $\exp(-C_Tt)\geq\exp(-C_TT)>0$. Its power is therefore bounded above by the same fraction raised to $\theta_T:=\exp(-C_TT)$. Taking the supremum over $0\leq t\leq T$ proves \eqref{uniq:modulus}, after absorbing $C_T^{1-\theta_T}$ into its prefactor. If the initial and incoming data agree, this fraction tends to zero and the same estimate gives $h=0$.
\end{proof}

\subsection{Diffuse Reflection: Continuous Dependence and Uniqueness}\label{sec:uniq-diffuse}

For diffuse reflection, the boundary identity \eqref{trace:projection-identity} replaces the incoming discrepancy by a nonnegative dissipation term.
\begin{proposition}\label{uniq:osgood-diffuse}
Let $f$ and $\widetilde f$ be weak solutions of \eqref{eq:normalized-problem} on $[0,T]$, both satisfying \eqref{eq:diffuse-normalized}. Assume \eqref{cb:integration-gradient} for $f$ and \eqref{uniq:smallness} for $\widetilde f$, with $\delta_*$ as in Proposition~\ref{uniq:osgood}. Then we have, for some $C_T\geq1$ and $\theta_T\in(0,1)$,
\begin{align}
\tnorm[\big]{\widetilde f-f}_{C([0,T];L^2)}^2\leq C_T\norm[\big]{\widetilde f(0)-f(0)}_{L^2}^{2\theta_T}. \label{uniq:modulus-diffuse}
\end{align}
The constants depend additionally on $T$ and the bounds in \eqref{cb:integration-gradient}. In particular, equal initial data imply $\widetilde f=f$.
\end{proposition}
\begin{proof}
\proofstep{The difference energy inequality} Set $h:=\widetilde f-f$. Its equation is the bulk equation in \eqref{trace:difference-equation}, while trace linearity gives
\begin{align}
\trm h=\Pm\trp h. \label{uniq:diffuse-difference-boundary}
\end{align}
The source and trace estimates verified in Proposition~\ref{uniq:osgood}, including \eqref{uniq:reversed-source}, verify the hypotheses of Proposition~\ref{trace:green}. Applying it and then \eqref{trace:projection-identity}, we obtain
\begin{align}
\frac12\left[\norm{h(r)}_{L^2}^2\right]_{r=s}^{r=t}+\int_s^t\left(\frac12\norm{(I-\Pp)\trp h(r)}_{L^2_{\Sigma_+}}^2+\ip{\Lop h}{h}\right)\dd r=\int_s^t\ip{\Ga(\widetilde f,h)+\Ga(h,f)}{h}\dd r. \label{trace:difference-diffuse-energy}
\end{align}
This identity again makes the squared $L^2$ norm absolutely continuous. Inserting \eqref{cvel:principal} and \eqref{eq:I:display:319} gives
\begin{align}
&\frac{\dd}{\dd t}\norm{h(t)}_{L^2}^2+\frac{c_0}{2}\norm{h(t)}_{\Vspace}^2+\norm{(I-\Pp)\trp h(t)}_{L^2_{\Sigma_+}}^2 \label{uniq:energy-diffuse}\\
&\qquad\leq C\norm{h(t)}_{L^2}^2+C\int_\Om\vnorm{f(t,x,\cdot)}_{\Vv}^2\vnorm{h(t,x,\cdot)}_{L_v^2}^2\dd x.\notag
\end{align}

\proofstep{Continuous dependence and uniqueness} Lemmas~\ref{uniq:velocity-log} and \ref{uniq:logarithm} give \eqref{uniq:optimization} for the last integral. Choose $C_T$ to bound the Gaussian norms as in \eqref{uniq:uniform-bounds}, without its incoming term. Taking $p=\ln\left(\frac{C_T}{\norm{h(t)}_{L^2}^2}\right)$ when $\norm{h(t)}_{L^2}>0$ gives
\begin{align}
\frac{\dd}{\dd t}\norm{h(t)}_{L^2}^2\leq C_T\norm{h(t)}_{L^2}^2\ln\left(\frac{C_T}{\norm{h(t)}_{L^2}^2}\right). \label{uniq:differential-diffuse}
\end{align}
The logarithmic term is zero when the norm vanishes. The regularization and integration in \eqref{uniq:forced-range}--\eqref{uniq:integrated-logarithm}, with the incoming integral omitted, give
\begin{align}
\norm{h(t)}_{L^2}^2\leq C_T\left(\frac{\norm{h(0)}_{L^2}^2}{C_T}\right)^{\exp(-C_Tt)},\qquad 0\leq t\leq T. \label{uniq:diffuse-integrated}
\end{align}
Taking the supremum in time proves \eqref{uniq:modulus-diffuse}. If $h(0)=0$, its right-hand side vanishes, proving uniqueness.
\end{proof}

\subsection{Convergence of Derivatives and Traces}\label{sec:uniq-stability}

We will retain the dissipation in the difference inequalities to prove convergence in $L^2(0,T;\Vspace)$, then use the difference equation to control the transport derivatives. For in-flow, the same energy inequality gives convergence of the outgoing trace. For diffuse reflection, it controls only $(I-\Pp)\trp(f_j-f)$. We will estimate the remaining component through the flux, using \eqref{trace:local-quadratic-Green} with a fixed velocity cutoff.
\begin{lemma}[In-flow]\label{uniq:bulk}
Let $f$ and $f_j$ be weak solutions of \eqref{eq:normalized-problem} on $[0,T]$, with incoming data $g$ and $g_j$, respectively. Assume a common bound in $L^\infty(0,T;L^{\infty,\alpha})$, \eqref{uniq:smallness} with $\widetilde f=f_j$ for every $j$, and \eqref{cb:integration-gradient} for $f$. If
\begin{align}
f_j(0)\to f(0)\quad\text{in }L^2,\qquad g_j\to g\quad\text{in }L^2((0,T)\times\Sigma_-), \label{uniq:inflow-data-convergence}
\end{align}
then we have
\begin{align}
f_j&\to f\quad\text{in }C([0,T];L^2)\cap L^2(0,T;\Vspace), \label{uniq:bulk-convergence}\\
\Y(f_j-f)&\to0\quad\text{in }L^2(0,T;\Vspace'), \label{uniq:transport-convergence}\\
\trpm f_j&\to\trpm f\quad\text{in }L^2\bigl((0,T)\times\Sigma_\pm;\dd t\dd\gamma\bigr). \label{uniq:inflow-trace-convergence}
\end{align}
\end{lemma}
\begin{proof}
\proofstep{Convergence in $C([0,T];L^2)$} Set $h_j:=f_j-f$. The reference solution is fixed, and the Gaussian bounds are uniform in $j$. Thus the constants in Proposition~\ref{uniq:osgood} can be chosen independently of $j$, and \eqref{uniq:modulus} gives
\begin{align}
\tnorm{h_j}_{C([0,T];L^2)}^2\leq C_T\left(\norm{h_j(0)}_{L^2}^2+\tnorm{g_j-g}_{L^2((0,T)\times\Sigma_-)}^2\right)^{\theta_T}\to0,\qquad\theta_T>0. \label{uniq:bulk-uniform}
\end{align}

\proofstep{Convergence of the velocity derivatives and outgoing trace} We will retain the dissipation in \eqref{uniq:energy} to prove the $L^2(0,T;\Vspace)$ convergence in \eqref{uniq:bulk-convergence} and the outgoing convergence in \eqref{uniq:inflow-trace-convergence}. Choose one sufficiently large logarithmic constant $C_T>1$ for all $j$, as in the proof of Proposition~\ref{uniq:osgood}. The common Gaussian bound allows us to require $\tnorm{h_j}_{C([0,T];L^2)}^2\leq\frac{C_T}{\ue^2}$. Integrating \eqref{uniq:energy} and applying \eqref{uniq:optimization}--\eqref{uniq:optimized-product} gives
\begin{align}
&\tnorm{h_j}_{L^2(0,T;\Vspace)}^2+\tnorm{\trp h_j}_{L^2((0,T)\times\Sigma_+)}^2 \label{eq:I:display:330}\\
&\qquad\lesssim_T\norm{h_j(0)}_{L^2}^2+\tnorm{g_j-g}_{L^2((0,T)\times\Sigma_-)}^2+\int_0^T\norm{h_j(t)}_{L^2}^2\ln\left(\frac{C_T}{\norm{h_j(t)}_{L^2}^2}\right)\dd t.\notag
\end{align}
We give the logarithmic expression its continuous value zero when $h_j(t)=0$. Since $s\ln\left(\frac{C_T}{s}\right)$ is increasing on $\left[0,\frac{C_T}{\ue}\right]$ and tends to zero as $s\downarrow0$, \eqref{uniq:bulk-uniform} yields
\begin{align}
\int_0^T\norm{h_j(t)}_{L^2}^2\ln\left(\frac{C_T}{\norm{h_j(t)}_{L^2}^2}\right)\dd t\leq T\tnorm{h_j}_{C([0,T];L^2)}^2\ln\left(\frac{C_T}{\tnorm{h_j}_{C([0,T];L^2)}^2}\right)\to0. \label{uniq:bulk-vanishing}
\end{align}
Together with \eqref{uniq:inflow-data-convergence}, this proves \eqref{uniq:bulk-convergence} and the outgoing convergence in \eqref{uniq:inflow-trace-convergence}. The incoming convergence follows from $\trm h_j=g_j-g$.

\proofstep{Convergence of the transport derivatives} We will prove \eqref{uniq:transport-convergence}. Equation~\eqref{trace:difference-equation} reads
\begin{align}
\Y h_j=-\Lop h_j+\Ga(f_j,h_j)+\Ga(h_j,f). \label{eq:I:display:331}
\end{align}
For the first two terms, \eqref{cvel:spatial-L}, \eqref{cvel:spatial-dual} and the common Gaussian bound give
\begin{align}
\tnorm{\Lop h_j}_{L^2(0,T;\Vspace')}+\tnorm{\Ga(f_j,h_j)}_{L^2(0,T;\Vspace')}\lesssim_T\tnorm{h_j}_{L^2(0,T;\Vspace)}\to0. \label{uniq:linear-source-convergence}
\end{align}
For the last term, \eqref{cvel:spatial-reversed} gives
\begin{align}
\norm[\big]{\Ga(h_j,f)}_{\Vspace'}^2\lesssim\int_\Om\vnorm[\big]{f(t,x,\cdot)}_{\Vv}^2\vnorm[\big]{h_j(t,x,\cdot)}_{L_v^2}^2\dd x. \label{eq:I:display:332}
\end{align}
Applying \eqref{uniq:optimization}--\eqref{uniq:optimized-product} and then \eqref{uniq:bulk-vanishing}, we obtain
\begin{align}
\tnorm[\big]{\Ga(h_j,f)}_{L^2(0,T;\Vspace')}^2\lesssim_T\int_0^T\norm{h_j(t)}_{L^2}^2\ln\left(\frac{C_T}{\norm{h_j(t)}_{L^2}^2}\right)\dd t\to0. \label{eq:I:display:333}
\end{align}
This proves \eqref{uniq:transport-convergence}.
\end{proof}
\begin{lemma}[Diffuse-reflection]\label{uniq:bulk-diffuse}
Let $f$ and $f_j$ be weak solutions of \eqref{eq:normalized-problem} on $[0,T]$, all satisfying \eqref{eq:diffuse-normalized}. Assume a common bound in $L^\infty(0,T;L^{\infty,\alpha})$, \eqref{uniq:smallness} with $\widetilde f=f_j$ for every $j$, and \eqref{cb:integration-gradient} for $f$. If $f_j(0)\to f(0)$ in $L^2$, then we have
\begin{align}
f_j&\to f\quad\text{in }C([0,T];L^2)\cap L^2(0,T;\Vspace), \label{uniq:diffuse-bulk-convergence}\\
\Y(f_j-f)&\to0\quad\text{in }L^2(0,T;\Vspace'), \label{uniq:diffuse-transport-convergence}\\
\trpm f_j&\to\trpm f\quad\text{in }L^2\bigl((0,T)\times\Sigma_\pm;\dd t\dd\gamma\bigr), \label{rev14:trace-convergence}\\
\flux{f_j}&\to\flux{f}\quad\text{in }L^2\bigl((0,T)\times\partial\Om\bigr). \label{rev14:flux-convergence}
\end{align}
\end{lemma}
\begin{proof}
\proofstep{Convergence in $C([0,T];L^2)$} Set $h_j:=f_j-f$. The common Gaussian bound and the fixed reference estimate allow uniform constants in \eqref{uniq:modulus-diffuse}, which gives $\tnorm{h_j}_{C([0,T];L^2)}^2\leq C_T\norm{h_j(0)}_{L^2}^{2\theta_T}\to0$, with $\theta_T>0$ independent of $j$.

\proofstep{Convergence of the velocity derivatives and boundary dissipation} We will retain the dissipative terms in \eqref{uniq:energy-diffuse} to prove the $L^2(0,T;\Vspace)$ convergence in \eqref{uniq:diffuse-bulk-convergence} and estimate $(I-\Pp)\trp h_j$. Choose a common logarithmic constant $C_T$ as in Lemma~\ref{uniq:bulk}. Integration and \eqref{uniq:optimization}--\eqref{uniq:optimized-product} give
\begin{align}
&\tnorm{h_j}_{L^2(0,T;\Vspace)}^2+\tnorm[\big]{(I-\Pp)\trp h_j}_{L^2((0,T)\times\Sigma_+)}^2 \label{uniq:diffuse-dissipation-convergence}\\
&\qquad\lesssim_T\norm{h_j(0)}_{L^2}^2+\int_0^T\norm{h_j(t)}_{L^2}^2\ln\left(\frac{C_T}{\norm{h_j(t)}_{L^2}^2}\right)\dd t\to0.\notag
\end{align}
The integral tends to zero by \eqref{uniq:bulk-vanishing}, proving \eqref{uniq:diffuse-bulk-convergence}.

\proofstep{Convergence of the transport derivatives} The bounds \eqref{uniq:linear-source-convergence}--\eqref{eq:I:display:333} use only the common Gaussian bound, the reference first-order estimate and the bulk convergence just proved. Their right-hand sides tend to zero by \eqref{uniq:diffuse-dissipation-convergence} and \eqref{uniq:bulk-vanishing}. All three terms in the difference equation \eqref{eq:I:display:331} therefore converge to zero in $L^2(0,T;\Vspace')$, proving \eqref{uniq:diffuse-transport-convergence}.

\proofstep{Convergence of the flux} We will prove \eqref{rev14:flux-convergence} by controlling the incoming trace on a fixed velocity ball. The Gaussian trace bound in Lemma~\ref{trace:gaussian} gives $\trpm h_j\in L^2((0,T)\times\Sigma_\pm)$. Together with \eqref{uniq:diffuse-bulk-convergence} and \eqref{uniq:diffuse-transport-convergence}, this verifies the hypotheses of Proposition~\ref{trace:green} for $h_j$, with source $\Y h_j+\Lop h_j\in L^2(0,T;\Vspace')$.

Let $\bar n\in C^1(\overline\Om;\R^3)$ be a bounded extension of the unit outward normal. Fix a radial $\chi\in C_c^\infty(B_2)$ with $0\leq\chi\leq1$ and $\chi=1$ on $B_1$, and set $\varphi(x,v):=\bigl(\bar n(x)\cdot v\bigr)\chi(v)^2$. Equation~\eqref{trace:local-quadratic-Green} is linear in its multiplier, so it permits this signed test. Using the cancellation \eqref{trace:cutoff-cancellation} gives
\begin{align}
&\int_0^T\int_{\partial\Om\times\R^3}\chi^2\abs{\trace h_j}^2(n\cdot v)^2\dd S_x\dd v\dd t \label{uniq:normal-multiplier}\\
&\qquad=\int_{\phase}\bigl(h_j(0)^2-h_j(T)^2\bigr)\varphi\dd x\dd v+2\int_0^T\ip*{\Y h_j}{\varphi h_j}\dd t+\int_0^T\int_{\phase}h_j^2\Y\varphi\dd x\dd v\dd t.\notag
\end{align}
Approximation in $C^1$ by smooth multipliers with common compact velocity support justifies \eqref{uniq:normal-multiplier}: the endpoint and trace terms converge by their $L^2$ bounds, while the duality pairing and volume integral converge using the uniform convergence of $\varphi$, $\nabla_v\varphi$ and $\Y\varphi$.

The function $\Y\varphi=v_i(\partial_{x_i}\bar n_k)v_k\chi^2$ is bounded. By \eqref{eq:V} and the local ellipticity \eqref{trace:local-matrix-assumption}, multiplication by $\varphi$ is bounded on $\Vspace$. Hence \eqref{uniq:normal-multiplier} implies
\begin{align}
&\int_0^T\int_{\partial\Om\times\R^3}\chi^2\abs{\trace h_j}^2(n\cdot v)^2\dd S_x\dd v\dd t \label{uniq:normal-multiplier-bound}\\
&\qquad\lesssim_T\tnorm{h_j}_{C([0,T];L^2)}^2+\tnorm{\Y h_j}_{L^2(0,T;\Vspace')}\tnorm{h_j}_{L^2(0,T;\Vspace)}\to0.\notag
\end{align}
Here the bound for the last integral in \eqref{uniq:normal-multiplier} uses $\tnorm{h_j}_{L^2(0,T;L^2)}^2\leq T\tnorm{h_j}_{C([0,T];L^2)}^2$.

The diffuse condition \eqref{eq:diffuse-normalized} gives $\trm h_j=\sqrt{2\pi}\mh\flux{h_j}$. Since $\chi$ and $\m$ are radial, rotation of velocity at each boundary point yields
\begin{align}
&\int_0^T\int_{\Sigma_-}\chi^2\abs{\trm h_j}^2(n\cdot v)^2\dd S_x\dd v\dd t \label{uniq:flux-from-incoming}\\
&\qquad=2\pi\left(\int_{v_3<0}\chi(v)^2\m(v)v_3^2\dd v\right)\tnorm{\flux{h_j}}_{L^2((0,T)\times\partial\Om)}^2.\notag
\end{align}
The velocity integral is a fixed positive number because $\chi=1$ on $B_1$. The left-hand side is bounded by the nonnegative integral in \eqref{uniq:normal-multiplier-bound}. Therefore \eqref{uniq:flux-from-incoming} proves \eqref{rev14:flux-convergence}.

\proofstep{Convergence of both full traces} We will combine the flux convergence with \eqref{uniq:diffuse-dissipation-convergence} to prove \eqref{rev14:trace-convergence}. Integrating \eqref{trace:projection-identity-a} in $(t,x)$ gives
\begin{align}
\tnorm{\Pp\trp h_j}_{L^2((0,T)\times\Sigma_+)}^2=\tnorm{\trm h_j}_{L^2((0,T)\times\Sigma_-)}^2=\sqrt{2\pi}\tnorm{\flux{h_j}}_{L^2((0,T)\times\partial\Om)}^2\to0. \label{uniq:diffuse-projected-convergence}
\end{align}
The orthogonal decomposition \eqref{trace:projection-Pythagoras} then yields
\begin{align}
\tnorm{\trp h_j}_{L^2((0,T)\times\Sigma_+)}^2=\tnorm{\Pp\trp h_j}_{L^2((0,T)\times\Sigma_+)}^2+\tnorm[\big]{(I-\Pp)\trp h_j}_{L^2((0,T)\times\Sigma_+)}^2\to0. \label{uniq:diffuse-full-outgoing-convergence}
\end{align}
This proves \eqref{rev14:trace-convergence} for both signs and completes the proof.
\end{proof}


\section{Proof of the Main Results}\label{sec:assembly}

We will combine the existence, regularity and uniqueness estimates from Sections~\ref{sec:energy}, \ref{cb:section} and \ref{sec:uniqueness} to prove the main results.

\subsection{Proof of Theorem~\ref{main:inflow} and Corollary~\ref{in:cor:decay}}\label{sec:assembly-inflow}

\begin{proof}[Proof of Theorem~\ref{main:inflow}]
We choose the constants for both main theorems. Take $\delta_*>0$ below the thresholds in Propositions~\ref{uniq:osgood} and \ref{uniq:osgood-diffuse}, with $\delta_*\leq\delta_A$. Choose $C_*\geq1$ to dominate the constants in Propositions~\ref{ci:global-energy}, \ref{cex:CM-input} and Corollary~\ref{cex:decay-traces}, using $\alpha'=\frac{\alpha}{2}$ in the latter, and choose $\lambda_*>0$ below their decay rates, including \eqref{cex:density-decay}. Finally, take $\varepsilon_*$ below both weak-existence thresholds with $C_*\varepsilon_*\leq\frac{\delta_*}{2}$. These choices give \eqref{in:eq:thresholds} and \eqref{eq:diffuse-thresholds}; the constructed solutions will satisfy
\begin{align}
\tnorm{f}_{L^\infty(0,\infty;L^{\infty,\alpha})}\leq C_*\varepsilon_*\leq\frac{\delta_*}{2}\leq\frac{\delta_A}{2}. \label{assembly:constructed-smallness}
\end{align}
This verifies \eqref{cb:setting} and \eqref{cb:initial-smallness}. These choices depend only on the fixed parameters; data derivative bounds enter only the regularity estimates.

\proofstep{Part (i): Global weak existence and the energy estimates} Assumption~\ref{in:data:regularity} gives nonnegative initial and incoming densities, and \eqref{in:eq:data} verifies \eqref{eq:I:display:214}. Proposition~\ref{ci:global-energy} therefore supplies a global weak solution with $\m+\mh f\geq0$. Its Gaussian estimate \eqref{ci:uniform-Gaussian} and \eqref{eq:gaussian-energy} give \eqref{in:eq:main-uniform} and \eqref{assembly:constructed-smallness}. The forcing estimate \eqref{ci:forced-energy-decay}, with its rate decreased to $\lambda_*$, gives \eqref{in:eq:main-energy-forcing}, and \eqref{ci:integrated-energy} gives \eqref{in:eq:main-dissipation}. The constants are uniform in the final time.

\proofstep{Part (ii): Regularity and the strong equation} The bound \eqref{assembly:constructed-smallness}, the data derivatives \eqref{in:eq:data-derivatives}--\eqref{in:eq:data-derivatives-g}, and compatibility \eqref{in:eq:compat0}--\eqref{in:eq:compat1} verify the in-flow hypotheses of Proposition~\ref{cb:integration}. For every $T>0$, its conclusion \eqref{cb:integration-holder} gives \eqref{in:eq:main-holder}, while \eqref{cb:initial-inflow-gradient}--\eqref{cb:initial-inflow-hessian} and \eqref{cb:integration-interior} give \eqref{in:eq:main-gradient}--\eqref{in:eq:main-hessian}, with a common smaller Gaussian exponent $\alpha_T\in(0,\alpha)$. The weighted bound \eqref{cb:strong-graph} gives \eqref{in:eq:main-strong} for every $m\geq0$, and \eqref{cb:strong-equation} gives the equation in $L^2((0,T)\times\phase)$.

\proofstep{Part (iii): Continuous dependence and uniqueness} Fix $T>0$. The Gaussian bound \eqref{in:eq:main-uniform}, the first-order estimate \eqref{in:eq:main-gradient} and its interior version verify \eqref{cb:integration-gradient} for $f$. The comparison bound \eqref{in:eq:comparison-class} is \eqref{uniq:smallness}. For the sequence in part~(iii), the same Gaussian bound and \eqref{uniq:inflow-data-convergence} verify the hypotheses of Lemma~\ref{uniq:bulk}. Its conclusions \eqref{uniq:bulk-convergence} and \eqref{uniq:inflow-trace-convergence} give \eqref{in:eq:main-stability} and \eqref{in:eq:main-stability-trace}.

For equal initial and incoming data, both differences in \eqref{uniq:modulus} vanish. Proposition~\ref{uniq:osgood} therefore gives $\widetilde f=f$ on every $[0,T]$, proving global uniqueness in the stated class.
\end{proof}

\begin{proof}[Proof of Corollary~\ref{in:cor:decay}]
\proofstep{Decay of the solution} The global weak solution satisfies the Gaussian smallness required in Proposition~\ref{ci:gaussian-bootstrap}. Under \eqref{in:eq:inflow-decay}, its estimate \eqref{ci:decaying-Gaussian} gives the Gaussian term of \eqref{in:eq:main-decay}; \eqref{eq:gaussian-energy} gives the $L^2$ term. Choose $\lambda_\omega>0$ strictly below $\omega$ and the available decay rates. This choice also absorbs any polynomial time factor from equal rates in the forcing convolution, and proves \eqref{in:eq:main-decay}.

\proofstep{The time-integrated dissipation} Apply \eqref{ci:integrated-energy} on $[t,T]$ and let $T\to\infty$. The incoming bound \eqref{in:eq:inflow-decay}, Gaussian integration in velocity, and monotone convergence give
\begin{align}
\int_t^\infty\left(\norm[\big]{f(s)}_{\Vspace}^2+\norm[\big]{\trp f(s)}_{L^2_{\Sigma_+}}^2\right)\dd s\lesssim\norm[\big]{f(t)}_{L^2}^2+\frac{A^2}{2\omega}\ue^{-2\omega t}. \label{in:eq:tail-energy-proof}
\end{align}
Insert \eqref{in:eq:main-decay} and use $\lambda_\omega\leq\omega$. This proves \eqref{in:eq:decay-dissipation-flux}, with constants independent of the data derivative bounds.

If instead the incoming Gaussian norm only tends to zero, \eqref{ci:vanishing-inflow-Gaussian} in Proposition~\ref{ci:gaussian-bootstrap}(iii) gives \eqref{eq:I:display:029}, again using \eqref{eq:gaussian-energy} for the $L^2$ norm.
\end{proof}

\subsection{Proof of Theorem~\ref{main:diffuse} and Corollary~\ref{main:diffuse-decay}}\label{sec:assembly-diffuse}

\begin{proof}[Proof of Theorem~\ref{main:diffuse}]
Use the constants fixed at the beginning of the proof of Theorem~\ref{main:inflow}.
\proofstep{Part (i): Global weak existence, mass and decay} Assumption~\ref{def:data} gives nonnegativity and zero perturbation mass, and \eqref{eq:small-data} verifies the Gaussian smallness in Proposition~\ref{cex:CM-input}. That proposition supplies a global weak solution with $\m+\mh f\geq0$ and mass conservation \eqref{eq:main-mass}. Its bound \eqref{cex:CM-bounds} and \eqref{eq:gaussian-energy} give \eqref{eq:main-decay} and \eqref{assembly:constructed-smallness}.

Apply Corollary~\ref{cex:decay-traces} with $\alpha'=\frac{\alpha}{2}$. Squaring the $L^2$ term in \eqref{cex:weaker-Gaussian} and decreasing its rate to $\lambda_*$ gives \eqref{eq:main-energy-decay}; the energy and both trace terms in \eqref{cex:global-integrals} give \eqref{eq:main-decay-dissipation}.

\proofstep{Part (ii): Regularity and the strong equation} The bound \eqref{assembly:constructed-smallness}, the data derivatives \eqref{eq:data-derivatives} and compatibility \eqref{eq:compatibility} verify the diffuse-reflection hypotheses of Proposition~\ref{cb:integration}. For every $T>0$, its conclusions \eqref{cb:integration-holder} and \eqref{cb:integration-flux} give \eqref{eq:main-holder}, with $\kappa_T\in(\frac23,1)$. Its derivative bounds \eqref{cb:initial-inflow-gradient}--\eqref{cb:initial-inflow-hessian} and \eqref{cb:integration-interior} give \eqref{eq:main-gradient}--\eqref{eq:main-hessian}, with a common smaller Gaussian exponent $\alpha_T\in(0,\alpha)$. Finally, \eqref{cb:strong-graph} gives \eqref{eq:main-graph} for every $m\geq0$, and \eqref{cb:strong-equation} gives the equation in $L^2((0,T)\times\phase)$.

\proofstep{Part (iii): Continuous dependence and uniqueness} Fix $T>0$. The amplitude \eqref{eq:main-decay}, the first-order estimate \eqref{eq:main-gradient} and its interior version verify \eqref{cb:integration-gradient} for $f$. The comparison bound \eqref{eq:diffuse-comparison-class} is \eqref{uniq:smallness}. For the sequence in part~(iii), the same Gaussian bound and the initial convergence in $L^2$ verify the hypotheses of Lemma~\ref{uniq:bulk-diffuse}. Its conclusions \eqref{uniq:diffuse-bulk-convergence}, \eqref{rev14:trace-convergence} and \eqref{rev14:flux-convergence} give \eqref{eq:main-stability}, \eqref{eq:main-stability-trace} and \eqref{eq:main-stability-flux}.

For equal initial data and the same diffuse boundary law, the right side of \eqref{uniq:modulus-diffuse} vanishes. Proposition~\ref{uniq:osgood-diffuse} therefore gives $\widetilde f=f$ on every $[0,T]$, proving global uniqueness in the stated class.
\end{proof}

\begin{proof}[Proof of Corollary~\ref{main:diffuse-decay}]
The solution constructed in Theorem~\ref{main:diffuse} is the solution of Proposition~\ref{cex:CM-input}, so Corollary~\ref{cex:decay-traces} applies. For each $0<\alpha'<\alpha$, its estimate \eqref{cex:weaker-Gaussian} is \eqref{eq:main-decay-rate}, with $p_\gamma$ given by \eqref{full:decay-exponent}.

For $-2\leq\gamma\leq1$, \eqref{cex:same-weight-decay} retains the initial Gaussian exponent $\alpha$. Choosing the common rate $\lambda_*\leq\lambda_D$ gives \eqref{full:same-weight-decay}. Both estimates have constants independent of the data derivative bounds.
\end{proof}


\bibliographystyle{siam}
\bibliography{Reference}

@article{Albritton.Armstrong.Mourrat.Novack2024,
  author = {Albritton, Dallas and Armstrong, Scott and Mourrat, Jean-Christophe and Novack, Matthew},
  title = {Variational methods for the kinetic {Fokker--Planck} equation},
  journal = {Anal. PDE},
  volume = {17},
  number = {6},
  pages = {1953--2010},
  year = {2024},
  doi = {10.2140/apde.2024.17.1953},
}

@article{Bardos1970,
    AUTHOR = {Bardos, Claude},
     TITLE = {Probl\`emes aux limites pour les \'equations aux d\'eriv\'ees partielles du premier ordre \`a coefficients r\'eels; th\'eor\`emes d'approximation; application \`a l'\'equation de transport},
   JOURNAL = {Ann. Sci. \'Ecole Norm. Sup. (4)},
  FJOURNAL = {Annales Scientifiques de l'\'Ecole Normale Sup\'erieure. Quatri\`eme S\'erie},
    VOLUME = {3},
      YEAR = {1970},
    NUMBER = {2},
     PAGES = {185--233},
       URL = {http://www.numdam.org/item/ASENS_1970_4_3_2_185_0/},
}

@article{Bernou.Carrapatoso.Mischler.Tristani2023,
  author = {Bernou, Armand and Carrapatoso, Kleber and Mischler, St{\'e}phane and Tristani, Isabelle},
  title = {Hypocoercivity for kinetic linear equations in bounded domains with general {Maxwell} boundary condition},
  journal = {Ann. Inst. H. Poincar{\'e} Anal. Non Lin{\'e}aire},
  volume = {40},
  number = {2},
  pages = {287--338},
  year = {2023},
  doi = {10.4171/AIHPC/44},
}

@article{Boyer2005,
  author = {Boyer, Franck},
  title = {Trace theorems and spatial continuity properties for the solutions of the transport equation},
  journal = {Differential Integral Equations},
  volume = {18},
  number = {8},
  pages = {891--934},
  year = {2005},
}

@article{Cameron.Silvestre.Snelson2018,
    AUTHOR = {Cameron, Stephen and Silvestre, Luis and Snelson, Stanley},
     TITLE = {Global a priori estimates for the inhomogeneous {L}andau
              equation with moderately soft potentials},
   JOURNAL = {Ann. Inst. H. Poincar\'e{} C Anal. Non Lin\'eaire},
  FJOURNAL = {Annales de l'Institut Henri Poincar\'e{} C. Analyse Non
              Lin\'eaire},
    VOLUME = {35},
      YEAR = {2018},
    NUMBER = {3},
     PAGES = {625--642},
      ISSN = {0294-1449,1873-1430},
   MRCLASS = {35Q82 (35B45)},
  MRNUMBER = {3778645},
MRREVIEWER = {Laurent\ B\'etermin},
       DOI = {10.1016/j.anihpc.2017.07.001},
       URL = {https://doi.org/10.1016/j.anihpc.2017.07.001},
}

@article{Carrapatoso.Desvillettes.He2017,
    AUTHOR = {Carrapatoso, Kleber and Desvillettes, Laurent and He, Lingbing},
     TITLE = {Estimates for the large time behavior of the {L}andau equation in the {C}oulomb case},
   JOURNAL = {Arch. Ration. Mech. Anal.},
    VOLUME = {224},
      YEAR = {2017},
     PAGES = {381--420},
       DOI = {10.1007/s00205-017-1078-3},
       URL = {https://doi.org/10.1007/s00205-017-1078-3},
}

@article{Carrapatoso.Mischler2017,
    AUTHOR = {Carrapatoso, Kleber and Mischler, St\'ephane},
     TITLE = {{L}andau equation for very soft and {C}oulomb potentials near {M}axwellians},
   JOURNAL = {Ann. PDE},
  FJOURNAL = {Annals of PDE},
    VOLUME = {3},
      YEAR = {2017},
    NUMBER = {1},
     PAGES = {Paper No. 1, 65},
       DOI = {10.1007/s40818-017-0021-0},
       URL = {https://doi.org/10.1007/s40818-017-0021-0},
}

@article{Carrapatoso.Mischler2024,
    AUTHOR = {Carrapatoso, Kleber and Mischler, St\'ephane},
     TITLE = {The {Landau} equation in a domain},
      YEAR = {2024},
   JOURNAL = {Preprint, arXiv:2407.09031v1},
       URL = {https://arxiv.org/abs/2407.09031v1},
}

@article{Carrapatoso.Tristani.Wu2016,
    AUTHOR = {Carrapatoso, Kleber and Tristani, Isabelle and Wu, Kung-Chien},
     TITLE = {Cauchy problem and exponential stability for the inhomogeneous {L}andau equation},
   JOURNAL = {Arch. Ration. Mech. Anal.},
  FJOURNAL = {Archive for Rational Mechanics and Analysis},
    VOLUME = {221},
      YEAR = {2016},
    NUMBER = {1},
     PAGES = {363--418},
       DOI = {10.1007/s00205-015-0963-x},
       URL = {https://doi.org/10.1007/s00205-015-0963-x},
}

@article{Carrapatoso.Tristani.Wu2017Erratum,
    AUTHOR = {Carrapatoso, Kleber and Tristani, Isabelle and Wu, Kung-Chien},
     TITLE = {Erratum to: Cauchy problem and exponential stability for the inhomogeneous {Landau} equation},
   JOURNAL = {Arch. Ration. Mech. Anal.},
    VOLUME = {223},
      YEAR = {2017},
     PAGES = {1035--1037},
       DOI = {10.1007/s00205-016-1064-1},
       URL = {https://doi.org/10.1007/s00205-016-1064-1},
}

@book{Cercignani1990,
    AUTHOR = {Cercignani, Carlo},
     TITLE = {Mathematical methods in kinetic theory},
   EDITION = {Second},
 PUBLISHER = {Plenum Press, New York},
      YEAR = {1990},
     PAGES = {xii+252},
      ISBN = {0-306-43460-1},
   MRCLASS = {82C40 (35Q99 76P05 82-02)},
  MRNUMBER = {1069558},
MRREVIEWER = {Gert\ Due\ Billing},
       DOI = {10.1007/978-1-4899-7291-0},
       URL = {https://doi.org/10.1007/978-1-4899-7291-0},
}

@book{Cercignani.Illner.Pulvirenti1994,
    AUTHOR = {Cercignani, Carlo and Illner, Reinhard and Pulvirenti, Mario},
     TITLE = {The mathematical theory of dilute gases},
    SERIES = {Applied Mathematical Sciences},
    VOLUME = {106},
 PUBLISHER = {Springer-Verlag, New York},
      YEAR = {1994},
     PAGES = {viii+347},
      ISBN = {0-387-94294-7},
   MRCLASS = {82C40 (76-02 76P05 82-02 82B40)},
  MRNUMBER = {1307620},
MRREVIEWER = {Giuseppe\ Toscani},
       DOI = {10.1007/978-1-4419-8524-8},
       URL = {https://doi.org/10.1007/978-1-4419-8524-8},
}

@article{Chen.Nguyen.Yang2025,
    AUTHOR = {Chen, Ke and Nguyen, Quoc-Hung and Yang, Tong},
     TITLE = {Scaling-Critical Theory for the {Boltzmann} and {Landau} Equations},
   JOURNAL = {Preprint, arXiv:2509.14845v2},
      YEAR = {2025},
       DOI = {10.48550/arXiv.2509.14845},
       URL = {https://arxiv.org/abs/2509.14845v2},
}

@article{DengDuan2021,
    AUTHOR = {Deng, Dingqun and Duan, Renjun},
     TITLE = {Spectral gap formation to kinetic equations with soft potentials in bounded domain},
   JOURNAL = {Comm. Math. Phys.},
    VOLUME = {397},
      YEAR = {2023},
     PAGES = {1441--1489},
       DOI = {10.1007/s00220-022-04519-2},
       URL = {https://doi.org/10.1007/s00220-022-04519-2},
}

@article{Desvillettes1992,
    AUTHOR = {Desvillettes, L.},
     TITLE = {On asymptotics of the {B}oltzmann equation when the collisions become grazing},
   JOURNAL = {Transport Theory Statist. Phys.},
  FJOURNAL = {Transport Theory and Statistical Physics},
    VOLUME = {21},
      YEAR = {1992},
    NUMBER = {3},
     PAGES = {259--276},
  MRNUMBER = {1165528},
       DOI = {10.1080/00411459208203923},
       URL = {https://doi.org/10.1080/00411459208203923},
}

@article{Desvillettes2015,
    AUTHOR = {Desvillettes, Laurent},
     TITLE = {Entropy dissipation estimates for the {L}andau equation in the {C}oulomb case and applications},
   JOURNAL = {J. Funct. Anal.},
    VOLUME = {269},
      YEAR = {2015},
    NUMBER = {5},
     PAGES = {1359--1403},
       DOI = {10.1016/j.jfa.2015.05.009},
       URL = {https://doi.org/10.1016/j.jfa.2015.05.009},
}

@article{Desvillettes.Villani2000,
    AUTHOR = {Desvillettes, Laurent and Villani, C{\'e}dric},
     TITLE = {On the spatially homogeneous {L}andau equation for hard potentials. {P}art {I}: Existence, uniqueness and smoothness},
   JOURNAL = {Comm. Partial Differential Equations},
    VOLUME = {25},
      YEAR = {2000},
    NUMBER = {1--2},
     PAGES = {179--259},
       DOI = {10.1080/03605300008821512},
       URL = {https://doi.org/10.1080/03605300008821512},
}

@article{Dong.Guo.Yastrzhembskiy2024,
    AUTHOR = {Dong, Hongjie and Guo, Yan and Yastrzhembskiy, Timur},
     TITLE = {Asymptotic stability for relativistic {Vlasov--Maxwell--Landau} system in bounded domain},
      YEAR = {2024},
   JOURNAL = {Preprint, arXiv:2401.00554v2; Revised December 2024},
       URL = {https://arxiv.org/abs/2401.00554v2},
}

@article{Duan.Liu.Sakamoto.Strain2020,
    AUTHOR = {Duan, Renjun and Liu, Shuangqian and Sakamoto, Shota and
              Strain, Robert M.},
     TITLE = {Global mild solutions of the {L}andau and non-cutoff
              {B}oltzmann equations},
   JOURNAL = {Comm. Pure Appl. Math.},
  FJOURNAL = {Communications on Pure and Applied Mathematics},
    VOLUME = {74},
      YEAR = {2021},
    NUMBER = {5},
     PAGES = {932--1020},
      ISSN = {0010-3640,1097-0312},
   MRCLASS = {35Q35 (82C40)},
  MRNUMBER = {4230064},
       DOI = {10.1002/cpa.21920},
       URL = {https://doi.org/10.1002/cpa.21920},
}

@article{Esposito.Guo.Kim.Marra2013,
    AUTHOR = {Esposito, R. and Guo, Y. and Kim, C. and Marra, R.},
     TITLE = {Non-isothermal boundary in the {B}oltzmann theory and
              {F}ourier law},
   JOURNAL = {Comm. Math. Phys.},
  FJOURNAL = {Communications in Mathematical Physics},
    VOLUME = {323},
      YEAR = {2013},
    NUMBER = {1},
     PAGES = {177--239},
      ISSN = {0010-3616,1432-0916},
   MRCLASS = {35Q20 (35B30 76P05)},
  MRNUMBER = {3085665},
MRREVIEWER = {Calvin\ Tadmon},
       DOI = {10.1007/s00220-013-1766-2},
       URL = {https://doi.org/10.1007/s00220-013-1766-2},
}

@article{Fournier2010,
    AUTHOR = {Fournier, Nicolas},
     TITLE = {Uniqueness of bounded solutions for the homogeneous {L}andau equation with a {C}oulomb potential},
   JOURNAL = {Comm. Math. Phys.},
  FJOURNAL = {Communications in Mathematical Physics},
    VOLUME = {299},
      YEAR = {2010},
    NUMBER = {3},
     PAGES = {765--782},
       DOI = {10.1007/s00220-010-1113-9},
       URL = {https://doi.org/10.1007/s00220-010-1113-9},
}

@article{Fournier.Guerin2009,
    AUTHOR = {Fournier, Nicolas and Gu{\'e}rin, H{\'e}l{\`e}ne},
     TITLE = {Well-posedness of the spatially homogeneous {L}andau equation for soft potentials},
   JOURNAL = {J. Funct. Anal.},
    VOLUME = {256},
      YEAR = {2009},
    NUMBER = {8},
     PAGES = {2542--2560},
       DOI = {10.1016/j.jfa.2008.11.008},
       URL = {https://doi.org/10.1016/j.jfa.2008.11.008},
}

@article{Golding.Gualdani.Loher2025,
    AUTHOR = {Golding, William and Gualdani, Maria and Loher, Am{\'e}lie},
     TITLE = {Global smooth solutions to the {L}andau--{C}oulomb equation in {$L^{3/2}$}},
   JOURNAL = {Arch. Ration. Mech. Anal.},
    VOLUME = {249},
      YEAR = {2025},
     PAGES = {Paper No. 34},
       DOI = {10.1007/s00205-025-02107-x},
       URL = {https://doi.org/10.1007/s00205-025-02107-x},
}

@article{Golse.Imbert.Mouhot.Vasseur2019,
    AUTHOR = {Golse, Fran{\c{c}}ois and Imbert, Cyril and Mouhot, Cl\'ement and Vasseur, Alexis F.},
     TITLE = {Harnack inequality for kinetic {F}okker-{P}lanck equations with rough coefficients and application to the {L}andau
              equation},
   JOURNAL = {Ann. Sc. Norm. Super. Pisa Cl. Sci. (5)},
  FJOURNAL = {Annali della Scuola Normale Superiore di Pisa. Classe di
              Scienze. Serie V},
    VOLUME = {19},
      YEAR = {2019},
    NUMBER = {1},
     PAGES = {253--295},
      ISSN = {0391-173X,2036-2145},
   MRCLASS = {35Q84 (35B45 35B65)},
  MRNUMBER = {3923847},
MRREVIEWER = {Barbara\ Brandolini},
}

@article{Guerand.Mouhot2022,
    AUTHOR = {Guerand, Jessica and Mouhot, Cl\'ement},
     TITLE = {Quantitative {D}e {G}iorgi methods in kinetic theory},
   JOURNAL = {J. \'Ec. polytech. Math.},
  FJOURNAL = {Journal de l'\'Ecole polytechnique. Math\'ematiques},
    VOLUME = {9},
      YEAR = {2022},
     PAGES = {1159--1181},
      ISSN = {2429-7100,2270-518X},
   MRCLASS = {35K70 (35Q84 35R09)},
  MRNUMBER = {4453413},
       DOI = {10.5802/jep.203},
       URL = {https://doi.org/10.5802/jep.203},
}

@article{Guillen.Silvestre2025,
    AUTHOR = {Guillen, Nestor and Silvestre, Luis},
     TITLE = {The {L}andau equation does not blow up},
   JOURNAL = {Acta Math.},
  FJOURNAL = {Acta Mathematica},
    VOLUME = {234},
      YEAR = {2025},
    NUMBER = {2},
     PAGES = {315--375},
       DOI = {10.4310/ACTA.2025.v234.n2.a2},
       URL = {https://doi.org/10.4310/ACTA.2025.v234.n2.a2},
}

@article{Guo.Hwang.Jang.Ouyang2020,
    AUTHOR = {Guo, Yan and Hwang, Hyung Ju and Jang, Jin Woo and Ouyang, Zhimeng},
     TITLE = {The {L}andau equation with the specular reflection boundary condition},
   JOURNAL = {Arch. Ration. Mech. Anal.},
  FJOURNAL = {Archive for Rational Mechanics and Analysis},
    VOLUME = {236},
      YEAR = {2020},
    NUMBER = {3},
     PAGES = {1389--1454},
      ISSN = {0003-9527,1432-0673},
   MRCLASS = {35Q82 (35Q20)},
  MRNUMBER = {4076068},
       DOI = {10.1007/s00205-020-01496-5},
       URL = {https://doi.org/10.1007/s00205-020-01496-5},
}

@article{Guo.Hwang.Jang.Ouyang2021Correction,
    AUTHOR = {Guo, Yan and Hwang, Hyung Ju and Jang, Jin Woo and Ouyang, Zhimeng},
     TITLE = {Correction to: The {Landau} equation with the specular reflection boundary condition},
   JOURNAL = {Arch. Ration. Mech. Anal.},
    VOLUME = {240},
      YEAR = {2021},
    NUMBER = {1},
     PAGES = {605--626},
       DOI = {10.1007/s00205-021-01622-x},
       URL = {https://doi.org/10.1007/s00205-021-01622-x},
}

@article{Henderson.Snelson2020,
    AUTHOR = {Henderson, Christopher and Snelson, Stanley},
     TITLE = {{$C^\infty$} smoothing for weak solutions of the inhomogeneous {L}andau equation},
   JOURNAL = {Arch. Ration. Mech. Anal.},
  FJOURNAL = {Archive for Rational Mechanics and Analysis},
    VOLUME = {236},
      YEAR = {2020},
    NUMBER = {1},
     PAGES = {113--143},
      ISSN = {0003-9527,1432-0673},
   MRCLASS = {35Q20 (35B65 35D30 82D10)},
  MRNUMBER = {4072211},
       DOI = {10.1007/s00205-019-01465-7},
       URL = {https://doi.org/10.1007/s00205-019-01465-7},
}

@article{Henderson.Snelson.Tarfulea2019,
    AUTHOR = {Henderson, Christopher and Snelson, Stanley and Tarfulea, Andrei},
     TITLE = {Local existence, lower mass bounds, and a new continuation criterion for the {L}andau equation},
   JOURNAL = {J. Differential Equations},
  FJOURNAL = {Journal of Differential Equations},
    VOLUME = {266},
      YEAR = {2019},
    NUMBER = {2-3},
     PAGES = {1536--1577},
      ISSN = {0022-0396,1090-2732},
   MRCLASS = {35Q20 (35A01 35B44)},
  MRNUMBER = {3906224},
       DOI = {10.1016/j.jde.2018.08.005},
       URL = {https://doi.org/10.1016/j.jde.2018.08.005},
}

@article{Imbert.Mouhot2021,
    AUTHOR = {Imbert, Cyril and Mouhot, Cl\'ement},
     TITLE = {The {S}chauder estimate in kinetic theory with application to a toy nonlinear model},
   JOURNAL = {Ann. H. Lebesgue},
    VOLUME = {4},
      YEAR = {2021},
     PAGES = {369--405},
       DOI = {10.5802/ahl.75},
}

@article{Imbert.Silvestre2021,
    AUTHOR = {Imbert, Cyril and Silvestre, Luis},
     TITLE = {The {S}chauder estimate for kinetic integral equations},
   JOURNAL = {Anal. PDE},
  FJOURNAL = {Analysis \& PDE},
    VOLUME = {14},
      YEAR = {2021},
    NUMBER = {1},
     PAGES = {171--204},
      ISSN = {2157-5045,1948-206X},
   MRCLASS = {35K70 (35R09)},
  MRNUMBER = {4229202},
       DOI = {10.2140/apde.2021.14.171},
       URL = {https://doi.org/10.2140/apde.2021.14.171},
}

@article{Kim.Guo.Hwang2020,
    AUTHOR = {Kim, Jinoh and Guo, Yan and Hwang, Hyung Ju},
     TITLE = {An {$L^2$} to {$L^\infty$} framework for the {L}andau
              equation},
   JOURNAL = {Peking Math. J.},
  FJOURNAL = {Peking Mathematical Journal},
    VOLUME = {3},
      YEAR = {2020},
    NUMBER = {2},
     PAGES = {131--202},
      ISSN = {2096-6075,2524-7182},
   MRCLASS = {35Qxx},
  MRNUMBER = {4171912},
       DOI = {10.1007/s42543-019-00018-x},
       URL = {https://doi.org/10.1007/s42543-019-00018-x},
}

@article{Kim.Seo2026,
    AUTHOR = {Kim, Chanwoo and Seo, Jinsol},
     TITLE = {Boundary-Driven Nonequilibrium for the {Coulomb} {Landau} Equation},
      YEAR = {2026},
   JOURNAL = {Preprint, arXiv:2610.30111v1},
       URL = {https://arxiv.org/abs/2610.30111v1},
}

@article{Kim.Weidner2026,
    AUTHOR = {Kim, Kyeongbae and Weidner, Marvin},
     TITLE = {Sharp regularity near the grazing set for kinetic {Fokker--Planck} equations},
      YEAR = {2026},
   JOURNAL = {Preprint, arXiv:2603.04121v1},
       URL = {https://arxiv.org/abs/2603.04121v1},
}

@article{Landau1936,
    AUTHOR = {Landau, Lev D.},
     TITLE = {Die kinetische {G}leichung f\"ur den {F}all {C}oulombscher {W}echselwirkung},
   JOURNAL = {Phys. Z. Sowjetunion},
    VOLUME = {10},
      YEAR = {1936},
     PAGES = {154--164},
}

@article{Liu.Yang.Ma2020,
  author = {Liu, Liping and Yang, Hang and Ma, Xuan},
  title = {The {Landau} equation with inflow boundary condition in a finite channel},
  journal = {Acta Math. Sci. Ser. A},
  volume = {40},
  number = {4},
  pages = {904--917},
  year = {2020},
  url = {https://actams.apm.ac.cn/sxwlxbA/EN/abstract/abstract16170.shtml},
}

@article{Loeper2006,
    AUTHOR = {Loeper, Gr{\'e}goire},
     TITLE = {Uniqueness of the solution to the {V}lasov--{P}oisson system with bounded density},
   JOURNAL = {J. Math. Pures Appl. (9)},
  FJOURNAL = {Journal de Math{\'e}matiques Pures et Appliqu{\'e}es},
    VOLUME = {86},
      YEAR = {2006},
    NUMBER = {1},
     PAGES = {68--79},
       DOI = {10.1016/j.matpur.2006.01.005},
       URL = {https://doi.org/10.1016/j.matpur.2006.01.005},
}

@article{Mischler2000,
  author = {Mischler, St{\'e}phane},
  title = {On the trace problem for solutions of the {Vlasov} equation},
  journal = {Comm. Partial Differential Equations},
  volume = {25},
  number = {7--8},
  pages = {1415--1443},
  year = {2000},
  doi = {10.1080/03605300008821554},
}

@article{Mischler2010,
    AUTHOR = {Mischler, St\'ephane},
     TITLE = {Kinetic equations with {M}axwell boundary conditions},
   JOURNAL = {Ann. Sci. \'Ec. Norm. Sup\'er. (4)},
  FJOURNAL = {Annales Scientifiques de l'\'Ecole Normale Sup\'erieure.
              Quatri\`eme S\'erie},
    VOLUME = {43},
      YEAR = {2010},
    NUMBER = {5},
     PAGES = {719--760},
      ISSN = {0012-9593},
       DOI = {10.24033/asens.2132},
       URL = {https://doi.org/10.24033/asens.2132},
}

@article{Mouhot.Strain2007,
    AUTHOR = {Mouhot, Cl\'ement and Strain, Robert M.},
     TITLE = {Spectral gap and coercivity estimates for linearized
              {B}oltzmann collision operators without angular cutoff},
   JOURNAL = {J. Math. Pures Appl. (9)},
  FJOURNAL = {Journal de Math\'ematiques Pures et Appliqu\'ees. Neuvi\`eme S\'erie},
    VOLUME = {87},
      YEAR = {2007},
    NUMBER = {5},
     PAGES = {515--535},
      ISSN = {0021-7824},
       DOI = {10.1016/j.matpur.2007.03.003},
       URL = {https://doi.org/10.1016/j.matpur.2007.03.003},
}

@article{Niebel.Valentini2026,
    AUTHOR = {Niebel, Lukas and Valentini, Lisa},
     TITLE = {Sharp kinetic trace theory},
      YEAR = {2026},
   JOURNAL = {Preprint, arXiv:2607.24708v1},
       URL = {https://arxiv.org/abs/2607.24708v1},
}

@article{Ros-Oton.Weidner2025,
    AUTHOR = {Ros-Oton, Xavier and Weidner, Marvin},
     TITLE = {Optimal regularity for kinetic {Fokker--Planck} equations in domains},
      YEAR = {2025},
   JOURNAL = {Preprint, arXiv:2505.11943v1},
       URL = {https://arxiv.org/abs/2505.11943v1},
}

@article{Silvestre2022,
    AUTHOR = {Silvestre, Luis},
     TITLE = {H\"older estimates for kinetic {F}okker-{P}lanck equations up
              to the boundary},
   JOURNAL = {Ars Inven. Anal.},
  FJOURNAL = {Ars Inveniendi Analytica},
      YEAR = {2022},
     PAGES = {Paper No. 6, 29},
      ISSN = {2769-8505},
   MRCLASS = {35J62 (46E35)},
  MRNUMBER = {4527757},
MRREVIEWER = {Jiabin\ Zuo},
}

@article{Snelson2020,
    AUTHOR = {Snelson, Stanley},
     TITLE = {Gaussian bounds for the inhomogeneous {L}andau equation with hard potentials},
   JOURNAL = {SIAM J. Math. Anal.},
  FJOURNAL = {SIAM Journal on Mathematical Analysis},
    VOLUME = {52},
      YEAR = {2020},
    NUMBER = {2},
     PAGES = {2081--2097},
      ISSN = {0036-1410,1095-7154},
   MRCLASS = {35Q20 (35B45 35B65 35F20 35R09 82D10)},
  MRNUMBER = {4091877},
       DOI = {10.1137/19M1244275},
       URL = {https://doi.org/10.1137/19M1244275},
}

@article{Strain.Guo2006,
    AUTHOR = {Strain, Robert M. and Guo, Yan},
     TITLE = {Almost exponential decay near {M}axwellian},
   JOURNAL = {Comm. Partial Differential Equations},
  FJOURNAL = {Communications in Partial Differential Equations},
    VOLUME = {31},
      YEAR = {2006},
    NUMBER = {1-3},
     PAGES = {417--429},
      ISSN = {0360-5302,1532-4133},
   MRCLASS = {82B40 (82C40 82D05 82D10)},
  MRNUMBER = {2209761},
MRREVIEWER = {Hai-Liang\ Li},
       DOI = {10.1080/03605300500361545},
       URL = {https://doi.org/10.1080/03605300500361545},
}

@article{Strain.Guo2008,
    AUTHOR = {Strain, Robert M. and Guo, Yan},
     TITLE = {Exponential decay for soft potentials near {M}axwellian},
   JOURNAL = {Arch. Ration. Mech. Anal.},
  FJOURNAL = {Archive for Rational Mechanics and Analysis},
    VOLUME = {187},
      YEAR = {2008},
    NUMBER = {2},
     PAGES = {287--339},
      ISSN = {0003-9527,1432-0673},
   MRCLASS = {82B05 (35B45 35F25)},
  MRNUMBER = {2366140},
MRREVIEWER = {Simone\ Calogero},
       DOI = {10.1007/s00205-007-0067-3},
       URL = {https://doi.org/10.1007/s00205-007-0067-3},
}

@incollection{Ukai1986,
    AUTHOR = {Ukai, Seiji},
     TITLE = {Solutions of the {B}oltzmann equation},
 BOOKTITLE = {Patterns and waves},
    SERIES = {Stud. Math. Appl.},
    VOLUME = {18},
     PAGES = {37--96},
 PUBLISHER = {North-Holland, Amsterdam},
      YEAR = {1986},
      ISBN = {0-444-70144-3},
   MRCLASS = {35Q20 (76P05 82A40)},
  MRNUMBER = {882376},
MRREVIEWER = {Nina\ B.\ Maslova},
       DOI = {10.1016/S0168-2024(08)70128-0},
       URL = {https://doi.org/10.1016/S0168-2024(08)70128-0},
}

@article{Villani1996,
    AUTHOR = {Villani, C\'edric},
     TITLE = {On the {C}auchy problem for {L}andau equation: sequential stability, global existence},
   JOURNAL = {Adv. Differential Equations},
  FJOURNAL = {Advances in Differential Equations},
    VOLUME = {1},
      YEAR = {1996},
    NUMBER = {5},
     PAGES = {793--816},
      ISSN = {1079-9389},
   MRCLASS = {82D10 (35Q99)},
  MRNUMBER = {1392006},
MRREVIEWER = {Carlo\ Cercignani},
}

@article{Villani1998,
    AUTHOR = {Villani, C\'edric},
     TITLE = {On a new class of weak solutions to the spatially homogeneous {B}oltzmann and {L}andau equations},
   JOURNAL = {Arch. Rational Mech. Anal.},
  FJOURNAL = {Archive for Rational Mechanics and Analysis},
    VOLUME = {143},
      YEAR = {1998},
    NUMBER = {3},
     PAGES = {273--307},
       DOI = {10.1007/s002050050106},
       URL = {https://doi.org/10.1007/s002050050106},
}

@incollection{Villani2002,
    AUTHOR = {Villani, C\'edric},
     TITLE = {A review of mathematical topics in collisional kinetic theory},
 BOOKTITLE = {Handbook of mathematical fluid dynamics, {V}ol. {I}},
     PAGES = {71--305},
 PUBLISHER = {North-Holland, Amsterdam},
      YEAR = {2002},
      ISBN = {0-444-50330-7},
   MRCLASS = {82C40 (35F20 76P05 82-02)},
  MRNUMBER = {1942465},
MRREVIEWER = {Fran{\c{c}}ois\ Castella},
       DOI = {10.1016/S1874-5792(02)80004-0},
       URL = {https://doi.org/10.1016/S1874-5792(02)80004-0},
}

@article{Yudovich1963,
    AUTHOR = {Yudovich, V. I.},
     TITLE = {Non-stationary flows of an ideal incompressible fluid},
   JOURNAL = {Zh. Vychisl. Mat. Mat. Fiz.},
  FJOURNAL = {Zhurnal Vychislitel'noi Matematiki i Matematicheskoi Fiziki},
    VOLUME = {3},
      YEAR = {1963},
    NUMBER = {6},
     PAGES = {1032--1066},
}

@article{Zhu2024,
    AUTHOR = {Zhu, Yuzhe},
     TITLE = {Regularity of kinetic {F}okker-{P}lanck equations in bounded
              domains},
   JOURNAL = {Ann. H. Lebesgue},
  FJOURNAL = {Annales Henri Lebesgue},
    VOLUME = {7},
      YEAR = {2024},
     PAGES = {1323--1366},
      ISSN = {2644-9463},
   MRCLASS = {35Q84 (35B65 35H10)},
  MRNUMBER = {4883871},
MRREVIEWER = {Felix\ X.-F.\ Ye},
       DOI = {10.5802/ahl.221},
       URL = {https://doi.org/10.5802/ahl.221},
}

@article{Zhu2025,
    AUTHOR = {Zhu, Yuzhe},
     TITLE = {Sharp boundary regularity properties for hypoelliptic kinetic equations},
      YEAR = {2025},
   JOURNAL = {Preprint, arXiv:2509.02536v1},
       URL = {https://arxiv.org/abs/2509.02536v1},
}

@article{Dong.Guo.Ouyang2020,
    AUTHOR = {Dong, Hongjie and Guo, Yan and Ouyang, Zhimeng},
     TITLE = {The {V}lasov-{P}oisson-{L}andau system with the
              specular-reflection boundary condition},
   JOURNAL = {Arch. Ration. Mech. Anal.},
  FJOURNAL = {Archive for Rational Mechanics and Analysis},
    VOLUME = {246},
      YEAR = {2022},
    NUMBER = {2-3},
     PAGES = {333--396},
      ISSN = {0003-9527,1432-0673},
   MRCLASS = {35Q83 (82D10)},
  MRNUMBER = {4514055},
       DOI = {10.1007/s00205-022-01818-9},
       URL = {https://doi.org/10.1007/s00205-022-01818-9},
}

@article{Dong.Guo.Ouyang.Yastrzhembskiy2024,
    AUTHOR = {Dong, Hongjie and Guo, Yan and Ouyang, Zhimeng and
              Yastrzhembskiy, Timur},
     TITLE = {The local well-posedness of the relativistic
              {V}lasov-{M}axwell-{L}andau system with the specular
              reflection boundary condition},
   JOURNAL = {SIAM J. Math. Anal.},
  FJOURNAL = {SIAM Journal on Mathematical Analysis},
    VOLUME = {56},
      YEAR = {2024},
    NUMBER = {5},
     PAGES = {6613--6688},
      ISSN = {0036-1410,1095-7154},
   MRCLASS = {35Q83 (34A12 35H10 35K70 35Q61 35Q84)},
  MRNUMBER = {4802874},
       DOI = {10.1137/23M1608938},
       URL = {https://doi.org/10.1137/23M1608938},
}

\end{document}